\documentclass[12pt]{book}

\usepackage[utf8]{inputenc} 

\usepackage[top=3cm,left=2.5cm,right=2.5cm,bottom=3.5cm]{geometry} 
\usepackage{graphicx} 

 \usepackage[parfill]{parskip} 

\usepackage{booktabs} 
\usepackage{array} 
\usepackage{paralist} 
\usepackage{verbatim} 
\usepackage{subfig} 
\usepackage{hyperref}
\hypersetup{colorlinks,linkcolor=black,urlcolor=blue,pdftex}
\usepackage[english]{babel}

\usepackage{fancyhdr} 
\usepackage{sectsty}

\usepackage[nottoc,notlof,notlot]{tocbibind} 
\usepackage[titles,subfigure]{tocloft} 

\usepackage{fontenc}
\usepackage{amsfonts}
\usepackage{wrapfig}
\usepackage{multirow}
\usepackage{amssymb}
\usepackage{quoting}
\usepackage{rotating}
\quotingsetup{font=normalsize}
\usepackage{amsthm}
\theoremstyle{plain}
\newtheorem{teo}{\\ Theorem}
\usepackage{bm}
\usepackage{amsmath}

\newcommand{\bsig}{\boldsymbol{\sigma}}
\newcommand{\beps}{\boldsymbol{\varepsilon}}
\newcommand{\bepsi}{\boldsymbol{\epsilon}}
\newcommand{\bom}{\boldsymbol{\omega}}
\newcommand{\bnu}{\boldsymbol{\nu}}
\newcommand{\bmu}{\boldsymbol{\mu}}
\newcommand{\btau}{\boldsymbol{\tau}}
\newcommand{\bb}{\boldsymbol{\beta}}
\newcommand{\bg}{\boldsymbol{\gamma}}

\newcommand{\bl}{\boldsymbol{\lambda}}
\newcommand{\bkappa}{\boldsymbol{\kappa}}
\newcommand{\bphi}{\boldsymbol{\varphi}}
\newcommand{\ba}{\boldsymbol{\alpha}}
\newcommand{\bd}{\boldsymbol{\delta}}
\newcommand{\au}{\alpha_1}
\newcommand{\ad}{\alpha_2}
\newcommand{\ai}{\alpha_i}

\renewcommand{\Phi}{\varPhi}
\newcommand{\gr}{\mathbf}
\renewcommand{\L}{\mathbf{L}}
\newcommand{\x}{\mathbf{x}}
\newcommand{\bt}{\mathbf{t}}
\newcommand{\eps}{\varepsilon}
\newcommand{\Ey}{\mathbb{E}}
\newcommand{\R}{\mathbb{R}}
\newcommand{\Q}{\mathbf{Q}}
\newcommand{\tr}{\mathrm{tr}}
\renewcommand{\div}{\mathrm{div}}
\newcommand{\C}{\mathbf{C}}
\newcommand{\X}{\mathbf{X}}
\newcommand{\D}{\mathbf{D}}
\newcommand{\Eu}{\mathcal{E}}
\newcommand{\Ba}{\mathcal{B}}
\newcommand{\Ve}{\mathcal{V}}
\newcommand{\Rep}{\mathcal{R}}
\renewcommand{\S}{\mathcal{S}}
\newcommand{\bu}{\mathbf{u}}
\newcommand{\bv}{\mathbf{v}}
\newcommand{\bw}{\mathbf{w}}
\newcommand{\bo}{\mathbf{o}}
\newcommand{\g}{\mathbf{g}}
\newcommand{\e}{\mathbf{e}}
\newcommand{\p}{\mathbf{p}}
\newcommand{\n}{\mathbf{n}}
\newcommand{\F}{\mathbf{F}}
\newcommand{\f}{\mathbf{f}}
\renewcommand{\O}{\mathbf{O}}
\newcommand{\I}{\mathbf{I}}
\newcommand{\B}{\mathbf{B}}
\newcommand{\U}{\mathbf{U}}
\newcommand{\V}{\mathbf{V}}
\newcommand{\A}{\mathbf{A}}
\renewcommand{\P}{\mathbf{P}}
\renewcommand{\a}{\mathbf{a}}
\renewcommand{\b}{\mathbf{b}}
\renewcommand{\c}{\mathbf{c}}
\renewcommand{\d}{\mathbf{d}}
\newcommand{\M}{\mathbf{M}}
\newcommand{\N}{\mathbf{N}}
\newcommand{\W}{\mathbf{W}}
\newcommand{\Rb}{\mathbf{R}}
\renewcommand{\r}{\mathbf{r}}

\newcommand{\Ro}{\mathbf{R}}
\newcommand{\Sy}{\mathbf{S}}
\newcommand{\Lq}{\mathbb{L}}
\newcommand{\Aq}{\mathbb{A}}
\newcommand{\Bq}{\mathbb{B}}
\newcommand{\Cq}{\mathbb{C}}
\newcommand{\Dq}{\mathbb{D}}
\newcommand{\Iq}{\mathbb{I}}
\newcommand{\Uq}{\mathbb{U}}
\newcommand{\Oq}{\mathbb{O}}
\newcommand{\Sq}{\mathbb{S}}
\newcommand{\Tq}{\mathbb{T}}
\newcommand{\Wq}{\mathbb{W}}
\newcommand{\Nq}{\mathbb{N}}
\newcommand{\Mq}{\mathbb{M}}
\newcommand{\Ax}{\mathcal{A}}
\newcommand{\eu}{\gr{e}_1}
\newcommand{\ed}{\gr{e}_2}
\newcommand{\et}{\gr{e}_3}

\newcommand{\ez}{\gr{e}_\mathit{z}}
\newcommand{\ei}{\gr{e}_\mathit{i}}
\newcommand{\ej}{\gr{e}_\mathit{j}}
\newcommand{\eh}{\gr{e}_\mathit{h}}
\newcommand{\ek}{\gr{e}_\mathit{k}}
\newcommand{\el}{\gr{e}_\mathit{l}}
\newcommand{\ep}{\gr{e}_\mathit{p}}
\newcommand{\eq}{\gr{e}_\mathit{q}}
\newcommand{\er}{\gr{e}_\mathit{r}}
\newcommand{\es}{\gr{e}_\mathit{s}}
\newcommand{\curl}{\mathrm{curl}}
\newcommand{\grad}{\mathrm{grad}}

\renewcommand{\(}{\begin{columns}}
\renewcommand{\)}{\end{columns}}
\newcommand{\<}[1]{\begin{column}{#1}}
\renewcommand{\>}{\end{column}}
\newcommand{\be}{\begin{equation}}
\newcommand{\ee}{\end{equation}}
\newcommand{\bes}{\begin{equation*}}
\newcommand{\ees}{\end{equation*}}
\newcommand{\besp}{\begin{split}}
\newcommand{\esp}{\end{split}}
\renewcommand{\phi}{\varphi}

\begin{document}
    \frontmatter

\title{\begin{figure}[t]
	\begin{center}
         \includegraphics[scale=.15]{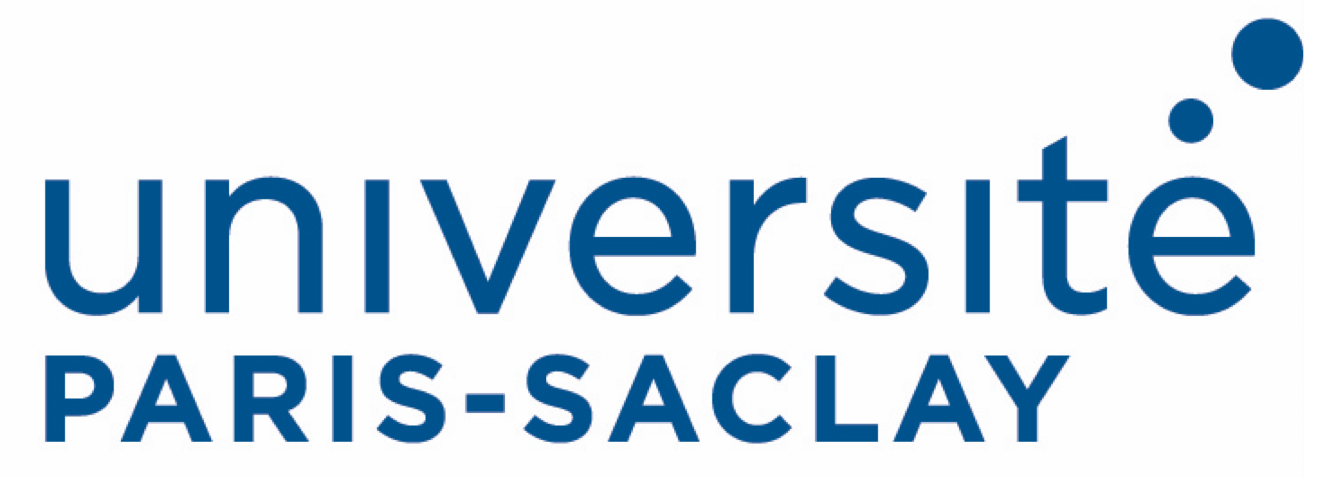}\vspace{-2cm}\\
	\end{center}
	\end{figure}
	{A short course on}\vspace{3mm}\\\huge{\textbf{Tensor algebra and analysis for engineers}}\vspace{2mm}\\\begin{large}with applications to differential geometry\end{large}}

\author{
\\
}

\date{\vspace{.5cm}
\includegraphics[scale=.25]{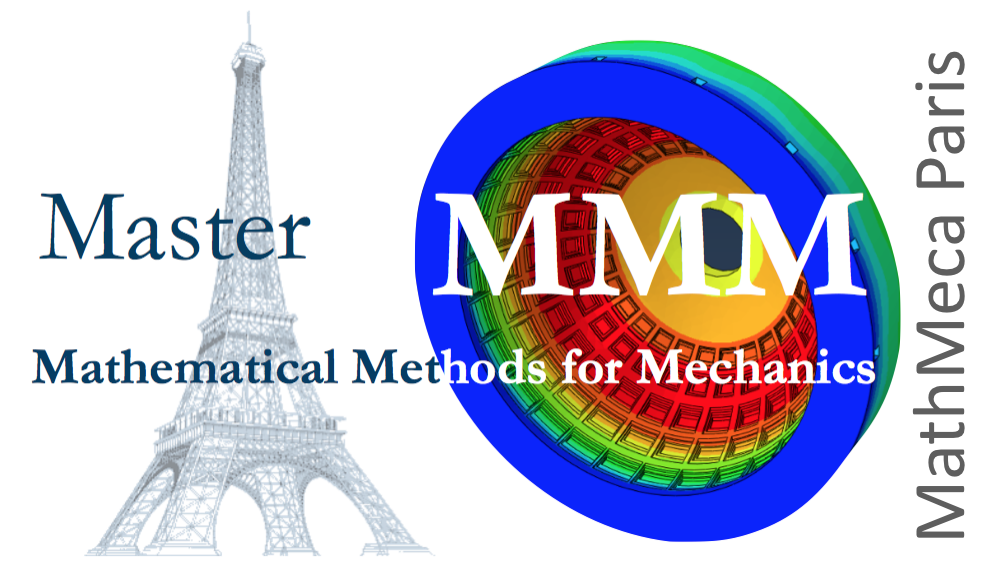}\vspace{1.8cm}\\
Paolo Vannucci\\
\vspace{.2cm}
{\href{mailto:paolo.vannucci@uvsq.fr}{paolo.vannucci@uvsq.fr}}\vspace{1.5cm}}

\maketitle

\vspace{7cm}

\hspace{10cm}{\it To Carla, Bianca and Alessandro}

\newpage

    \chapter{Preface}

      This textbook, addressed to graduate students and young researchers in mechanics, has been  developed from the class notes  of different courses in continuum mechanics that I have been delivering for several years as part of  the Master's program in MMM - Mathematical Methods for Mechanics, at the  University Paris-Saclay.
    
  Far from being exhaustive, as any primer text, the intention of this book is to introduce students  in mechanics and engineering to the  mathematical language and tools that are necessary for a modern approach to continuum theoretical and applied mechanics. The presentation of the matter is hence tailored for this scientific community and necessarily it is different, in terms of  language and  objectives, from that normally proposed to students of other disciplines, such as physics, especially general relativity,  or  pure mathematics.

  What has motivated me to write this textbook is the idea of collecting in a single, introductory book a set of results and tools useful for studies in mechanics and  presenting them in a modern, succinct way. Almost all the results and theorems are proved and the reader is guided along a {\it tour} that starts from vectors and ends at the differential geometry of surfaces, passing through the algebra of tensors of second and fourth orders, the differential geometry of curves, the tensor analysis for fields and deformations and the use of curvilinear coordinates. 
  
  Some topics are specially treated, such as rotations, the algebra of fourth order tensors, fundamental for the mechanics of modern materials, or the properties of differential operators. 
  Some other topics are intentionally omitted because they are less important to continuum mechanics or too advanced for an introductory text. Though, in some modern texts, tensors are directly presented in the most general setting of curvilinear coordinates, I preferred here to choose a more traditional approach, introducing first tensors in Cartesian coordinates, normally used for classical problems.  Then, an entire chapter is devoted to the passage to curvilinear coordinates and to the formalism of co- and contravariant components.

    The tensor theory and results  are specially applied to introduce some subjects concerning differential geometry of curves and surfaces. Also in this case, the presentation is mainly intended for applications to continuum mechanics and, in particular, in view of courses on slender beams or thin shells. All the presentation of the topics of differential geometry is extensively based on tensor algebra and analysis.
    
    More than a hundred  exercises are proposed to the reader, many of them completing the theoretical part through new results and proofs. All the exercises are entirely  developed and solved  at the end of the book, in order to provide the reader with  thorough support for his learning.
    
    In   Chapter 1, vectors and points are introduced and also, with a small anticipation of some results of the second Chapter, also applied vectors are visited.
   Chapter 2 is completely devoted to the algebra of second-rank tensors and the succeeding  Chapter 3 to that of fourth-rank tensors. Intentionally, these are the only two types of tensors introduced in the book: They are the most important types of  tensors in mechanics, by them we can represent deformation, stress and the constitutive laws. 
    I preferred not to introduce tensors in an absolutely general way but to go directly to the most important tensors for applications in mechanics; for the same reason, the algebra of other tensors, namely of third-rank tensors, is not presented in this primer text. 
   
   The analysis of tensors is first done using  differential geometry of curves, in Chapter 4, for differentiation and integration with respect to only one variable, then introducing the differential operators for fields and deformations, in Chapter 5.
   
   Then, a generalization of second-rank tensor algebra and analysis in the sense of the use of curvilinear coordinates is presented in Chapter 6, where the notion of metric tensor, co- and contravariant components and Christoffel's symbols are introduced. 
   
   Finally, Chapter 7 is entirely devoted to an introduction to the differential geometry of surfaces. Classical topics such as the first and second fundamental forms of a surface, the different types of curvatures, the Gauss-Weingarten equations or the concepts of minimal surfaces,  geodesics and the Gauss-Codazzi conditions are presented, with all these topics being of a great interest in mechanics.  
   
   I tried to write a coherent, almost self-contained manual of mathematical tools for graduate students in  mechanics, with the hope of helping young students  progress in their studies. The exposition is as simple as possible, sober, sometimes minimalist.     
   I intentionally avoided burdening the language and the text  with nonessential details and considerations, but I have always tried to grasp the essence of a result and its usefulness. 
   
   It is my most sincere hope that the reader who dares to persevere through the pages of this book will find a benefit for his studies in continuum mechanics. This is, eventually, the goal of this primer text.   
   \bigskip\bigskip\\
    {\it Versailles, June 16, 2022}
\chapter{Acknowledgments}
I am indebted to many persons for the topics of this book. Professor E. G. G. Virga, University of Pavia, introduced me, the first, to tensor algebra, during my PhD at the University of Pisa, many years ago. 

Then, I have had the privilege of collaborating with Professor G. Verchery, at the University of Burgundy, who introduced me to the representation methods based upon tensor invariants. This has been very useful for developing some results in the algebra of fourth rank tensors. 

I wish also to thank  Professor P. M. Mariano, University of Florence, who always pushed me to go forward and to consider problems of modern mechanics; many of the discussions I had with him have been very important to me.

I have also had many interesting and useful discussions  with Professor J. Lerbet, University of Evry, and with Doctor C. Fourcade, of Renault S.A.; I wish to thank them sincerely.

I am also grateful to Prof. A. Frediani, who has been for me an example of honesty in teaching and science, and I cannot forget Prof. P. Villaggio, my PhD director, who passed by some years ago: his teaching and personality leaved an indelebile trace in all my life of scientist.

Finally, I wish to thank my wife, Carla, my daughter, Bianca, and my son, Alessandro. Without them, nothing would have been possible; because of them, many things happen.

\chapter{List of symbols}
\label{ch:list}

$:=$ : definition symbol\smallskip\\
$|$: such that\smallskip\\
$\exists !$: exists and is unique\smallskip\\
$\R$: set of real numbers\smallskip\\
$\Eu$: ordinary 3D Euclidean space\smallskip\\
$\Ve$: vector space of translations, associated with $\Eu$\smallskip\\
$Lin(\Ve)$: vector space of second-rank tensors\smallskip\\
$\Lq$in$(\Ve)$: vector space of fourth-rank tensors\smallskip\\
$x,y,z$ etc.: scalars (elements of $\R$)\smallskip\\
$p,q,r$ etc.: points (elements of $\Eu$)\smallskip\\
$\bu,\bv,\bw$ etc.: vectors (elements of $\Ve$)\smallskip\\
$\bu=p-q$ etc.: vector  difference of two points of $\Eu$\smallskip\\
$\L,\M,\N$ etc.: second-rank tensors (elements of $Lin(\Ve)$)\smallskip\\
$\Lq,\Mq,\Nq$ etc.: fourth-rank tensors (elements of $\Lq$in$(\Ve)$)\smallskip\\
$\bu\cdot\bv$ etc.: scalar product of vectors\smallskip\\
$\L\cdot\M$ etc.: scalar product of second-rank tensors\smallskip\\
$\bu\times\bv$ etc.: cross product of vectors\smallskip\\
$\bu\otimes\bv$ etc.: dyad (second-rank tensor) of vectors\smallskip\\
$u=|\bu|$: norm of a vector\smallskip\\
$L=|\L|$: norm of a tensor\smallskip\\
$|p-q|,|\bu-\bv|,|\L-\M|$: distance of  points, vectors or tensors\smallskip\\
$\det\L$: determinant of a  second-rank tensor $\L$\smallskip\\
$\L^\top$: transpose of a second-rank tensor $\L$\smallskip\\
$\L^{-1}$: inverse of a second-rank tensor $\L$\smallskip\\
$\L^{-\top}=(\L^\top)^{-1}=(\L^{-1})^\top$\smallskip\\
$Sym(\Ve)$: subspace of $Lin(\Ve)$ of symmetric second-rank tensors\smallskip\\
$Skw(\Ve)$: subspace of $Lin(\Ve)$ of skew second-rank tensors\smallskip\\
$Sph(\Ve)$: subspace of $Lin(\Ve)$ of spherical second-rank tensors\smallskip\\
$Dev(\Ve)$: subspace of $Lin(\Ve)$ of deviatoric second-rank tensors\smallskip\\
$Orth(\Ve)$: subspace of $Lin(\Ve)$ of orthogonal second-rank tensors\smallskip\\
$Orth(\Ve)^+$: subspace of $Lin(\Ve)$ of rotation second-rank tensors\smallskip\\
$\S$: unit sphere of $\Ve$: $\S=\{\bu\in\Ve|\ |\bu|=1\}$\smallskip\\
$\delta_{ij}$: Kronecker's delta\smallskip\\
$\Re$: real part of a complex quantity\smallskip\\
$\Im$: imaginary part of a complex quantity\smallskip\\
$\Lq=\A\boxtimes\B$: conjugation product of two second-rank tensors\smallskip\\
$\Sq^{sph}$: spherical projector\smallskip\\
$\Dq^{dev}$: deviatoric projector\smallskip\\
$\Tq^{trp}$: transpose projector\smallskip\\
$\Sq^{sym}$: symmetry projector\smallskip\\
$\Wq^{skw}$: antisymmetry projector\smallskip\\
$\Iq$: identity of $\Lq$in($\Ve$)\smallskip\\
$\Iq^s$: restriction of $\Iq$ to symmetric tensors of $Lin(\Ve)$

\tableofcontents

    \mainmatter

\chapter{Points and vectors}
\label{ch:1}
\section{Points and vectors}

\begin{figure}[ht]
	\begin{center}
         \includegraphics[scale=1]{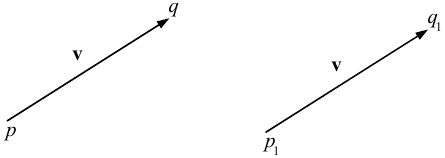}
	\caption{Same translation applied to two different points.}
	\label{fig:1}
	\end{center}
\end{figure}

We consider in the following a {\it point space} $\Eu$, whose elements are points $p$. In classical mechanics, $\Eu$ is to be identified with the {\it Euclidean three-dimensional space}, wherein   events are intended to be set. On $\Eu$, we admit the existence of an operation, the {\it difference} of any couple of its  elements:
\bes
q-p,\ \ p,q\in\Eu.
\ees
We associate to $\Eu$ a {\it vector space} $\Ve$ whose dimension is dim$\Ve=3$ and whose elements are vectors $\bv$ representing {\it translations over} $\Eu$:
\bes
\forall p,q\in\Eu,\ \exists!\ \bv\in\Ve|\ q-p=\bv.
\ees
Any element $\bv\in\Ve$ is hence a transformation over $\Eu$ that, using the previous definition, can be written as :
\bes
\forall\bv\in\Ve,\ \bv:\Eu\rightarrow\Eu|\  q=\bv(p)\ \rightarrow\ q=p+\bv.
\ees
We remark that the result of the application of the translation $\bv$ depends upon the argument $p$: 
\bes
q=p+\bv\neq p_1+\bv=q_1,
\ees
whose geometric meaning is depicted in Fig. \ref{fig:1}.
Unlike  difference, the sum of two points is not defined and  is meaningless.

\begin{figure}[h]
	\begin{center}
         \includegraphics[scale=1]{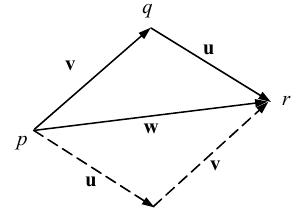}
	\caption{Sum of two vectors: the parallelogram rule.}
	\label{fig:2}
	\end{center}
\end{figure}

We define the {\it sum of two vectors} $\bu$ and $\bv$ as the vector $\bw$ such that
\bes
(\bu+\bv)(p)=\bu(\bv(p))=\bw(p)
\ees
This means that, if
\bes
q=\bv(p)=p+\bv,
\ees
then
\bes
r=\bu(q)=q+\bu=\bw(p),
\ees
see Fig. \ref{fig:2}, which shows that the above definition actually coincides with the {\it parallelogram rule} and that
\bes
\bu+\bv=\bv+\bu,
\ees
as obvious, for the sum over a vector space commutes.
It is evident that the sum of more than two vectors can be defined iteratively, summing up a vector at a time to the sum of the previous vectors.

The {\it null vector} $\bo$ is defined as the difference of any two coincident points:
\bes
\bo:=p-p\ \ \forall p\in\Eu;
\ees
$\bo$ is unique and the only vector such that
\bes
\bv+\bo=\bv\ \ \forall\bv\in\Ve.
\ees

In fact,
\bes
\forall p\in\Eu,\ \ \bv+\bo=\bv+p-p\ \rightarrow\ p+\bv+\bo=p+\bv\ \iff\ \bv+\bo=\bv.
\ees

A {\it linear combination} of $n$ vectors $\bv_i$ is defined as the vector\footnote{\label{note:sumeinstein}We adopt here and in the following the {\it Einstein notation} for summations: All the times when an index is repeated in a monomial, then the summation with respect to that index, called the {\it dummy index}, is understood, e.g., $k_i\bv_i=\sum_i k_i\bv_i$. We then say that the index $i$ is {\it saturated}. If a repeated index is underlined, then it is not a dummy index, i.e. there is no summation.}
\bes
\bw:=k_i\bv_i,\ \ k_i\in\R,\ \ i=1,...,n.
\ees
The $n+1$ vectors $\bw,\ \bv_i,\ i=1,...,n,$ are said to be {\it linear independent} if there does not exist a set of $n$ scalars $k_i$ such that the above equation is satisfied and are said to be {\it linear dependent} in the opposite case.

\section{Scalar product, distance, orthogonality}
A {\it scalar product} on a vector space is a positive definite, symmetric, bilinear form. A {\it form} $\omega$ is a function 
\bes
\omega:\Ve\times\Ve\rightarrow\R,
\ees 
i.e. $\omega$ operates on a couple of vectors to give a real number, a {\it scalar}. We will indicate the scalar product of two vectors $\bu$ and $\bv$ as\footnote{The scalar product $\omega(\bu,\bv)$ is also  indicated as $<\bu,\bv>$.}
\bes
\omega(\bu,\bv)=\bu\cdot\bv.
\ees
The properties of {\it bilinearity} prescribe that, $\forall \bu,\bv\in\Ve$ and $\forall \alpha,\beta\in\R$,
\bes
\besp
&\bu\cdot(\alpha\bv+\beta\bw)=\alpha\bu\cdot\bv+\beta\bu\cdot\bw,\\
&(\alpha\bu+\beta\bv)\cdot\bw=\alpha\bu\cdot\bw+\beta\bv\cdot\bw,
\end{split}
\ees
while {\it symmetry} implies that
\bes
\bu\cdot\bv=\bv\cdot\bu\ \ \forall \bu,\bv\in\Ve.
\ees
Finally, the {\it positive definiteness} means that
\bes
\bv\cdot\bv>0\ \ \forall \bv\in\Ve,\ \  \bv\cdot\bv=0\ \iff\ \bv=\bo.
\ees

Any two vectors $\bu,\bv\in\Ve$ are said to be {\it orthogonal} $\iff$
\bes
 \bu\cdot\bv=0.
\ees
Thanks to the properties of the scalar product, we can define the {\it Euclidean norm} of a vector $\bv$ as the nonnegative scalar, denoted equivalently by $v$ or $|\bv|$, 
\bes
v=|\bv|:=\sqrt{\bv\cdot\bv}
\ees
\begin{teo} The norm of a vector has the following properties: $\forall \bu,\bv\in\Ve,k\in\R$,
\bes
\besp
&|\bu\cdot\bv|\leq u\ v\ \ \mathrm{(Schwarz's\ inequality)};\\
&|\bu+\bv|\leq u+v\ \ \mathrm{(Minkowski's\ triangular\ inequality)};\\
&|k\bv|=|k| v.
\end{split}
\ees
\begin{proof}
{\it Schwarz's inequality}: It is sufficient to prove that 
\bes
(\bu\cdot\bv)^2\leq\bu\cdot\bu\ \bv\cdot\bv.
\ees
Let $x=\bv\cdot\bv$ and $y=-\bu\cdot\bv$. Then, by the positive definiteness of the scalar product, we get
\bes
(x\bu+y\bv)\cdot(x\bu+y\bv)\geq0,
\ees
which implies that
\bes
x^2\bu\cdot\bu+2xy\bu\cdot\bv+y^2\bv\cdot\bv=(\bv\cdot\bv)^2\bu\cdot\bu-2\bv\cdot\bv(\bu\cdot\bv)^2+\bv\cdot\bv(\bu\cdot\bv)^2\geq0;
\ees
supposing $\bv\neq\bo$ (otherwise, the proof is trivial), we get the thesis on dividing by $\bv\cdot\bv$.

{\it Minkowski's inequality}: Because the two members of the inequality to be proved are nonnegative, it is sufficient to prove that
\bes
(\bu+\bv)\cdot(\bu+\bv)\leq(u+v)^2=u^2+2uv+v^2.
\ees
This can be proved easily:
\bes
\besp
(\bu+\bv)\cdot(\bu+\bv)&=\bu\cdot\bu+2\bu\cdot\bv+\bv\cdot\bv=u^2+2\bu\cdot\bv+v^2\\
&\leq u^2+2|\bu\cdot\bv|+v^2\leq u^2+2uv+v^2,
\end{split}
\ees
in which the last operation follows from  the Schwarz's inequality.

The proof of the third property is immediate, it is sufficient to use the same definition of norm.
\end{proof}
\end{teo}
We define {\it distance between any two  points} $p$ and $q\ \in\Eu$ the scalar
\bes
d(p,q):=|p-q|=|q-p|.
\ees
Similarly, the {\it distance between two any vectors} $\bu$ and $\bv\ \in\Ve$ is defined as
\bes
d(\bu,\bv):=|\bu-\bv|=|\bv-\bu|.
\ees
Two points or two vectors are {\it coincident} if and only if their distance is null.

The {\it unit sphere $\S$ of} $\Ve$ is defined as the set of all the vectors whose norm is one:
\bes
\S:=\{\bv\in\Ve|\ v=1\}.
\ees

\section{Basis of $\Ve$, expression of the scalar product}
\begin{figure}[h]
	\begin{center}
         \includegraphics[scale=.7]{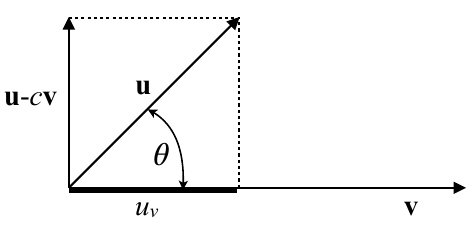}
	\caption{Angle between two vectors.}
	\label{fig:3}
	\end{center}
\end{figure}

There is a general way to define a {\it basis} for a vector space of any kind. Here, we limit the introduction of the concept of basis to the case of $\Ve$ only, of interest in classical mechanics. Generally speaking, a {\it basis $\Ba$} of $\Ve$  is any set of three linearly independent vectors $\e_i,i=1,2,3$ of $\Ve$:
\bes
\Ba=\{\e_1,\e_2,\e_3\}.
\ees
The introduction  of a  basis for $\Ve$ is useful for representing vectors. In fact, once a basis $\Ba$ is fixed, any  vector $\bv\in\Ve$ can be represented as a linear combination of the vectors of the basis, where the coefficients $v_i$ of the linear combination are the {\it Cartesian components} of $\bv$:
\bes
\bv=v_i\e_i=v_1\e_1+v_2\e_2+v_3\e_3.
\ees
Though the choice of the elements of a basis is completely arbitrary, the only condition being their linear independency, we will use in the following only {\it orthonormal bases}, which are  bases composed of mutually orthogonal vectors of $\S$, i.e. satisfying
\bes
\e_i\cdot\e_j=\delta_{ij},
\ees
where the symbol $\delta_{ij}$ is the so-called {\it Kronecker's delta}:
\bes
\delta_{ij}=\left\{
\besp
&1 \ \mathrm{if}\ i=j,\\
&0 \ \mathrm{if}\ i\neq j.
\end{split}
\right.
\ees
The use of orthonormal bases has great advantages; namely, it allows us to give a very simple rule for the calculation of the scalar product:
\bes
\bu\cdot\bv=u_i\e_i\cdot v_j\e_j=u_iv_j\delta_{ij}=u_iv_i=u_1v_1+u_2v_2+u_3v_3.
\ees
In particular, it is
\bes
\bv\cdot\e_i=v_k\e_k\cdot\e_i=v_k\delta_{ik}=v_i,\ \ i=1,2,3.
\ees
So, the  Cartesian components of a vector are the projection of the vector on the three vectors of the basis $\Ba$;  such quantities are the {\it director cosines} of $\bv$ in the basis $\Ba$. In fact, if $\theta$ is the angle formed by two vectors $\bu$ and $\bv$, then 
\bes
\bu\cdot\bv=u\ v\ \cos\theta.
\ees
This relation is used to define the angle between two vectors,
\bes
\theta=\arccos\frac{\bu\cdot\bv}{u\ v},
\ees
which can be proved easily: Given two vectors $\bu$ and $\bv$, we look for $c\in\R$ such that the vector $\bu-c\bv$ is orthogonal to $\bv$:
\bes
(\bu-c\bv)\cdot\bv=0\ \iff\ c=\frac{\bu\cdot\bv}{\bv\cdot\bv}=\frac{\bu\cdot\bv}{v^2}.
\ees
Now, if $\bu$ is inclined of $\theta$ on $\bv$, its projection $u_v$ on the direction of $\bv$ is
\bes
u_v=u\ \cos\theta,
\ees
and, by construction (see Fig. \ref{fig:3}), it is also
\bes
u_v=c\ v.
\ees
So,
\bes
c=\frac{u}{v}\cos\theta\ \rightarrow\ \frac{u}{v}\cos\theta=\frac{\bu\cdot\bv}{v^2}\ \Rightarrow\ \cos\theta=\frac{\bu\cdot\bv}{u\ v}.
\ees

We remark that while the scalar product, being an intrinsic operation, does not change for a change of basis, the components $v_i$ of a vector are {\it not} intrinsic quantities, but they are basis-dependent: A change of the basis makes  the components change. The way this change is done will be introduced in Section \ref{sec:changeofbasis}.

A {\it frame} $\Rep$ for $\Eu$ is composed of a point $o\in\Eu$, the {\it origin}, and a basis $\Ba$ of $\Ve$:
\bes
\Rep:=\{o,\Ba\}=\{o;\e_1,\e_2,\e_3\}.
\ees
The use of a frame for $\Eu$ is useful for determining the position of a point $p$, which can be done through its {\it Cartesian coordinates} $x_i$, defined as the components in $\Ba$ of the vector $p-o$:
\bes
x_i:=(p-o)\cdot\e_i,\ \ i=1,2,3.
\ees
Of course, the coordinates $x_i$ of a point $p\in\Eu$ depend  upon  the choice of $o$ and  $\Ba$.

\section{Applied vectors}
We introduce now a set of definitions, concepts and results that are widely used in physics, especially in mechanics. For that, we need to anticipate some results that are introduced in the next Chapter, namely that of cross product, in Section \ref{sec:crossprod}, and of complementary projector, in Exercise 2, Chapter 2. 
This slight deviation from the good rule of consistent progression in stating the results is justified by the fact that, actually, the matter exposed hereafter is still that of vectors. The readers can, of course,  come back to the topics of this section once they have studied Chapter \ref{ch:2}. 

We call {\it applied vector} $\bv^p$ a vector $\bv$ associated to a point $p\in\Eu$. In physics the concept of applied vector\footnote{In the literature, applied vectors are also called {\it bound vectors}.} is often employed, for example, to represent forces\footnote{The fact that in classical mechanics forces can be represented by vectors is actually a fundamental postulate of physics. Forces are vectors that cannot be considered belonging to the translation space $\Ve$; nevertheless,  the definitions and results found earlier are also  valid  for vectors $\notin\Ve$.}. We define the {\it resultant} of a system of $n$ applied vectors $\bv^p_i$ as the  vector
\bes
\Ro:=\sum_{i=1}^n\bv^p_i.
\ees
 We define {\it moment of an applied vector $\bv^p$} about a point $o$, called the {\it center of the moment},  the vector
\bes
\M_o:=(p-o)\times\bv^p,
\ees
and {\it resultant moment of a system of $n$ applied vectors $\bv^p_i$} about a point $o$ the vector
\bes
\M^r_o:=\sum_{i=1}^n(p_i-o)\times\bv^p_i.
\ees
We remark that $\Ro,\M_o$ and $\M^r_o$ are {\it not}  applied vectors.

\begin{figure}[h]
\begin{center}
\includegraphics[width=.23\textwidth]{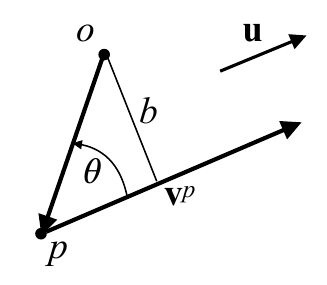}
\caption{Moment arm of an applied vector.}
\label{fig:55}
\end{center}
\end{figure}

If $\bu\in\S| \ \bu\times\bv^p=\bo$, then
\bes
b:=|(\I-\bu\otimes\bu)(p-o)|=|p-o|\sin\theta
\ees
is called the {\it moment arm} of $\bv^p$ with respect to the center $o$. It measures the distance of $o$ from the {\it line of action}, i.e. the line passing  through $p$ and parallel to $\bv^p$, cf. Fig. \ref{fig:55}.

\begin{figure}[h]
\begin{center}
\includegraphics[width=.3\textwidth]{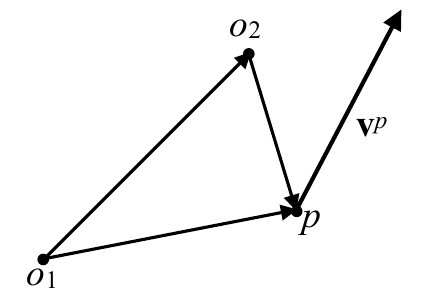}
\caption{Scheme for the transposition of the moment.}
\label{fig:56}
\end{center}
\end{figure}

\begin{teo}{\bf{(Transport of the moment).}} If $\M_{o_1}$ is the moment of an applied vector $\bv^p$ about a center $o_1$, the moment $\M_{o_2}$ of $\bv^p$ about to another center $o_2$ is
\bes
\M_{o_2}=\M_{o_1}+(o_1-o_2)\times\bv^p.
\ees
\begin{proof}
Referring to Fig. \ref{fig:56},
\bes
\besp
\M_{o_2}&=(p-o_2)\times\bv^p=(p-o_1+o_1-o_2)\times\bv^p\\
&=(p-o_1)\times\bv^p+(o_1-o_2)\times\bv^p\\
&=\M_{o_1}+(o_1-o_2)\times\bv^p.
\end{split}
\ees
\end{proof}
\end{teo}
A consequence of this theorem is that $\M_{o_1}=\M_{o_2}\iff\bv^p\times(o_1-o_2)$, i.e. if $\bv^p$ and $o_1-o_2$ are parallel. It follows from this that the moment of an applied vector does not change when calculated about the points of a straight line parallel to the vector itself or, more importantly, if $\bv^p$ is translated along its line of action.

The above theorem can be extended to the resultant moment of a system of applied vectors to give (the proof is quite similar)
\be
\label{eq:traspmomres}
\M^r_{o_2}=\M^r_{o_1}+(o_1-o_2)\times\Ro.
\ee
Also in this case, the resultant moment does not change when $\Ro$ and $o_1-o_2$ are parallel vectors, but not exclusively, as another possibility is that $\Ro=\bo$: For the systems of applied vectors with null resultant, the resultant moment is invariant with respect to the center of the moment. 

An interesting relation can be found if the two members of the last equation are projected onto $\Ro$, which gives
\be
\label{eq:momresproj}
\M^r_{o_1}\cdot\Ro=\M^r_{o_2}\cdot\Ro:
\ee
The projection of the resultant moment onto the direction of $\Ro$ does not depend upon the center of the moment.

\begin{figure}[h]
\begin{center}
\includegraphics[width=.4\textwidth]{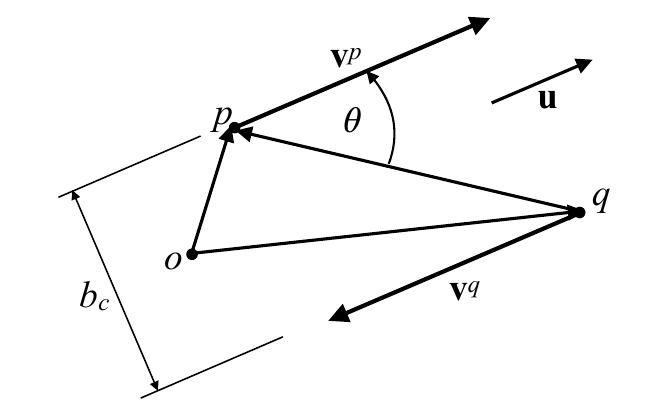}
\caption{Scheme of a couple.}
\label{fig:57}
\end{center}
\end{figure}

A particularly important case of system with null resultant is that of a {\it couple}, which is composed by two opposite vectors $\bv,-\bv$, applied to two points $p$ and $q$:
\bes
\bv^p=-\bv^q.
\ees
Of course, by definition, $\Ro=\bo$ for any couple and, as a consequence, the resultant moment $\M^r$  of a couple, called the {\it moment of the couple} and simply denoted by $\M$, is independent of the center of the moment (that is why the index denoting the center of the moment is omitted): Referring to Fig. \ref{fig:57},
\bes
\besp
\M&=(p-o)\times\bv^p+(q-o)\times\bv^q=(p-o)\times\bv-(q-o)\times\bv\\
&=((p-o)-(q-o))\times\bv=(p-q)\times\bv.
\end{split}
\ees
If $\bu\in\S| \ \bu\times\bv=\bo$, then
\bes
b_c:=|(\I-\bu\otimes\bu)(p-q)|=|p-q|\sin\theta
\ees
is  the {\it couple arm}. We then have
\bes
M=|(p-q)\times\bv|=|p-q|v\sin\theta=b_cv.
\ees

The {\it central axis} $\A$ of a system of $n$ applied vectors with $\Ro\neq\bo$ is the axis such that
\bes
\M^r_a\times\Ro=\bo\ \ \ \forall a\in\Ax.
\ees
\begin{teo} {\bf{(Existence and uniqueness of the central axis).}}  The central axis of a system of $n$ vectors exists and is unique.
\begin{proof}
Existence: We need  at least a point $a\in\Eu|\ \M^r_a=k\Ro,\ \ k\in\R\Rightarrow\M^r_a\times\Ro=\bo$.
From Eq. (\ref{eq:traspmomres}), $\forall o\in\Eu$ we get
\bes
\M^r_a\times\Ro=\M^r_o\times\Ro+((o-a)\times\Ro)\times\Ro=\M^r_o\times\Ro-\Ro^2(o-a)+(\Ro\cdot(o-a))\Ro.
\ees
Then, if we take for $o-a$ the vector
\bes
o-a=\frac{\M^r_o\times\Ro}{\Ro^2},
\ees
 it is evident that we get
\bes
\M^r_a\times\Ro=\bo.
\ees
Hence, the point
\bes
a=o-\frac{\M^r_o\times\Ro}{\Ro^2}\ \in\Ax.
\ees
So, because $\M^r$ does not change when calculated with respect to the points of an axis parallel to $\Ro$, $\Ax$ is the axis passing through $a$ and parallel to $\Ro$. Its equation is
\bes
p=a+t\ \Ro=o-\frac{\M^r_o\times\Ro}{\Ro^2}+t\ \Ro,\ \ \ t\in\R.
\ees
Uniqueness: Suppose another axis $\hat{\Ax}\neq\Ax$ exists, which is necessarily parallel to $\Ax$. If $q\in\hat{\Ax}$, again using Eq. (\ref{eq:traspmomres}) we get
\bes
\M^r_q\times\Ro=\M^r_a\times\Ro+((a-q)\times\Ro)\times\Ro.
\ees
In this equation the left-hand side and the first term on the right-hand side are null by the definition of central axis. Because $(a-q)\times\Ro$ is perpendicular to $\Ro$ and $\Ro\neq\bo$ by hypothesis, the left-hand side is null if and only if $a=q\Rightarrow\hat{\Ax}=\Ax$.
\end{proof}
\end{teo}
The central axis has another remarkable property:
\begin{teo} {\bf{(Property of minimum of the central axis).}}  The points of the central axis minimize the resultant moment.
\begin{proof}
When $\M^r$ is calculated about a point $a\in\Ax$, it is parallel to $\Ro$, which is not the case for any point $q\notin\Ax$. In this last case, hence, $\M^r$ has also a component orthogonal to $\Ro$. Then, by virtue of the invariance of the projection of $\M^r$ onto $\Ro$, Eq. (\ref{eq:momresproj}), $\M^r$ gets its minimum value when calculated about the points of $\Ax$.
\end{proof}
\end{teo}
Let us now consider the case of systems for which
\bes
\M^r_o\cdot\Ro=0\ \ \ \forall o\in\Eu.
\ees
This is namely the case of systems of coplanar or parallel vectors (cf. Exercise 7). 
Because in this case, for the points $a\in\Ax$, it must be at the same time
$\M^r_a\cdot\Ro=0$ and $\M^r_a\times\Ro=\bo$,
the only possibility is that
\bes
\M^r_a=\bo\ \ \ \forall a\in\Ax,
\ees
i.e., in these cases $\Ax$ is the axis of points that make the resultant moment vanish. 

Two systems of applied vectors are {\it equivalent} if they have the same resultant $\Ro$ and the same resultant moment $\M^r_o$ about any center $o\in\Eu$. The equivalence does not depend upon the center $o$. In fact, by Eq. (\ref{eq:traspmomres}), if two systems have the same $\Ro$ and the same $\M^r_{o_1}$, with $o_1$ a given point, then also $\M^r_{o_2}$ will be the same, $\forall o_2\in\Eu$.

\begin{teo} {\bf{(Reduction of a system of applied vectors).}}  A system of applied vectors is always equivalent to the system composed by the resultant $\Ro$ applied at a point $o$ and by a couple with moment $\M=\M^r_o$, with $o$ any point of $\Eu$. 
\begin{proof} By construction, $\Ro$ is the same for the two systems; moreover, for the equivalent system (resultant plus couple) it is
\bes
\M+(o-o)\times\Ro=\M.
\ees
So, if the couple has a moment $\M=\M^r_o$, the two systems are equivalent. 

\end{proof}
\end{teo}
In practice, this theorems affirms that it is always possible to reduce a system of $n$ applied vectors to only an applied vector equal to $\Ro$ and to a couple or, if one of the two vectors composing the couple is applied to the same point of $\Ro$, to two applied vectors. It is worth noting that the equivalence of two systems is preserved if a vector is translated along its line of action, because in such a case $\Ro$ and $\M^r_o$ do not change.

Finally, a system of $n$ applied vectors is said to be {\it equilibrated} if 
\bes
\Ro=\bo,\ \ \ \M^r_o=\bo\ \ \ \forall o\in\Eu.
\ees
We note that, because $\Ro=\bo$, the center $o$ can be any point of $\Eu$.

\section{Exercises}
\begin{enumerate}
\item Prove that the null vector is unique.
\item Prove that the null vector is orthogonal to any vector.
\item Prove that the norm of the null vector is zero.
\item Prove that 
\bes
\bu\cdot\bv=0\ \iff\ |\bu-\bv|=|\bu+\bv|\ \ \forall\bu,\bv\in\Ve.
\ees
\item Prove the {\it linear forms representation theorem}: Let $\psi:\Ve\rightarrow\R$ be a linear function. Then, $\exists!\ \bu\in\Ve$ such that
\bes
\psi(\bv)=\bu\cdot\bv\ \ \forall\bv\in\Ve.
\ees
\item Consider a point $p$ and  two noncollinear vectors $\bu,\bv\in\S$ at $p$. Show that a vector $\bw$ is the bisector of the angle formed by $\bu$ and $\bv$ if and only if $\bw\cdot\bu=\bw\cdot\bv$.

\item \label{ex:momnullch1} Show that in the case of systems composed of coplanar or parallel applied vectors with $\Ro\neq\bo,\M^r_o\cdot\Ro=0\ \forall o\in\Eu.$
\item Prove that any system of applied vectors with $\Ro=\bo$ is equivalent to a couple.
\item \label{ex:9ch1} Prove that a system of applied vectors all passing through a point $p$ is equivalent to $\Ro$ applied to $p$. 
\item Prove that if for a system of applied vectors $\M^r_o=\bo$, then the system is equivalent to $\Ro$ applied to $o$. Then, show that if $o\in\Ax$, this is the case of coplanar  or  parallel vectors.

\item Prove that a system of applied vectors is equilibrated if and only if any equivalent system is equilibrated.
\item Prove that two applied vectors form an equilibrated system if and only if they are two opposite vectors applied to the same point.
\item Prove that a system of  applied vectors is equilibrated if all the vectors pass through the same point and $\Ro=\bo$.

\end{enumerate}

\chapter{Second rank tensors}
\label{ch:2}
\section{Second-rank tensors}

A {\it second-rank tensor} $\L$ is any linear application from $\Ve$ to $\Ve$:
\bes
\L:\Ve\rightarrow\Ve\ |\ \L(\alpha_i\bu_i)=\alpha_i\L\bu_i\ \forall\alpha_i\in\R,\ \bu_i\in\Ve,\ i=1,...,n.
\ees
Though here $\Ve$ indicates the vector space of translations over $\Eu$, the definition of tensor\footnote{We consider, for the time being, only second-rank tensors, that constitute a very important set of operators in classical and continuum mechanics.  In the following, we also introduce  fourth-rank tensors.} 
is more general and, in particular, $\Ve$ can be any vector space.

Defining the {\it sum of two tensors} as
\be
\label{eq:tenssum}
(\L_1+\L_2)\bu=\L_1\bu+\L_2\bu\ \ \forall\bu\in\Ve,
\ee
the {\it product of a scalar by a tensor} as
\bes
(\alpha\L)\bu=\alpha(\L\bu)\ \ \forall\alpha\in\R,\bu\in\Ve
\ees
and the {\it null tensor} $\O$ as the unique tensor such that
\bes
\O\bu=\bo\ \forall\bu\in\Ve,
\ees
then the set of all the tensors $\L$ that operate on $\Ve$ forms a vector space, denoted by $Lin(\Ve)$. We define the {\it identity tensor} $\I$ as the unique tensor such that
\bes
\I\bu=\bu\ \ \forall\bu\in\Ve.
\ees

Different  operations can be defined for the second-rank tensors. We consider all of them in the following sections.

\section{Dyads, tensor components}
For any couple of vectors $\bu$ and $\bv$, the {\it dyad}\footnote{In some texts, the dyad is also called the {\it tensor product}; we prefer to use the term dyad because the expression tensor product can be ambiguous, as   it is used  to denote the product of two tensors, see Section \ref{sec:tensprod}.} $\bu\otimes\bv$ is the tensor defined by
\bes
(\bu\otimes\bv)\bw:=\bv\cdot\bw\ \bu\ \ \forall\bw\in\Ve.
\ees
The application defined above is actually a tensor because of the bilinearity of the scalar product. The introduction of dyads allows us to express any tensor as a linear combination of dyads. In fact, it can be proved that if $\Ba=\{\e_1,\e_2,\e_3\}$ is a basis of $\Ve$, then the set of nine dyads
\bes
\Ba^2=\{\e_i\otimes\e_j,\ i,j=1,2,3\},
\ees
is a basis of $Lin(\Ve)$, so that dim$(Lin(\Ve))=9$. This implies that any tensor $\L\in Lin(\Ve)$ can be expressed as
\bes
\L=L_{ij}\ \e_i\otimes\e_j,\ \ i,j=1,2,3,
\ees
where the $L_{ij}$s are the nine {\it Cartesian components} of $\L$ with respect to $\Ba^2$. The $L_{ij}$s can be calculated easily:
\bes
\e_i\cdot\L\e_j=\e_i\cdot L_{hk}\e_h\otimes\e_k\ \e_j=L_{hk} \e_i\cdot\e_h\ \e_k\cdot\e_j=L_{hk}\delta_{ih}\delta_{jk}=L_{ij}.
\ees
The above expression is sometimes called the {\it canonical decomposition of a tensor}.
The components of a dyad can be computed as follows:
\be
\label{eq:compdyad}
(\bu\otimes\bv)_{ij}=\e_i\cdot(\bu\otimes\bv)\ \e_j=\bu\cdot\e_i\ \bv\cdot\e_j=u_i\ v_j.
\ee

The components of a vector $\bv$ resulting from the application of a tensor $\L$ on a vector $\bu$, can now be  calculated:
\be
\label{eq:Lu}
\bv=\L\bu=L_{ij}(\e_i\otimes\e_j)(u_k\e_k)=L_{ij}u_k\delta_{jk}\e_i=L_{ij}u_j\e_i\ \rightarrow\ v_i= L_{ij}u_j.
\ee

Depending upon two indices, any second-rank tensor $\L$ can be represented by a matrix, whose entries are the Cartesian components of $\L$ in the basis $\Ba$:
\bes
\L=\left[
\begin{array}{ccc}
L_{11} & L_{12} & L_{13}\\
L_{21} & L_{22} & L_{23} \\
L_{31} & L_{32} & L_{33}
\end{array}
\right].
\ees
Because any $\bu\in\Ve$, depending upon only one index, can be represented by a column vector, Eq. (\ref{eq:Lu}) represents actually the classical operation of the multiplication of a $3\times3$ matrix by a $3\times1$ vector.

\section{Tensor product}
\label{sec:tensprod}
The {\it tensor product} of $\L_1$ and $\L_2\in Lin(\Ve)$ is defined by
\bes
(\L_1\L_2)\bv:=\L_1(\L_2\bv)\ \ \forall\bv\in\Ve.
\ees
By linearity and Eq. (\ref{eq:tenssum}), $\forall \L,\L_1,\L_2\in Lin(\Ve),\bu\in\Ve$, we get
\bes
\besp
[\L(\L_1+\L_2)]\bv&=\L[(\L_1+\L_2)\bv]=\L(\L_1\bv+\L_2\bv)\\
&=\L\L_1\bv+\L\L_2\bv=(\L\L_1+\L\L_2)\bv\ \rightarrow\\ 
&\ \ \  \ \L(\L_1+\L_2)=\L\L_1+\L\L_2.
\end{split}
\ees
We remark that the tensor product is not symmetric:
\bes
\L_1\L_2\neq\L_2\L_1;
\ees
however, by the same definition of the identity tensor and of tensor product,
\bes
\I\L=\L\I=\L\ \forall\L\in Lin(\Ve).
\ees
The Cartesian components of a tensor $\L=\A\B$ can be  calculated using Eq. (\ref{eq:Lu}):
\bes
\besp
L_{ij}&=\e_i\cdot(\A\B)\e_j=\e_i\cdot\A(\B\e_j)=\e_i\cdot\A(B_{hk}(\e_j)_k\ \e_h)=B_{hk}\delta_{jk}\e_i\cdot\A\e_h\\
&=B_{hk}\delta_{jk}\e_i\cdot(A_{pq}(\e_h)_q\ \e_p)=A_{pq}B_{hk}\delta_{jk}\delta_{qh}\delta_{ip}=A_{ih}B_{hj}.
\end{split}
\ees
The above result simply corresponds to the row-column multiplication  of two matrices. Using that, the following two identities can be readily shown:
\be
\label{eq:propdyad0}
\besp
&(\a\otimes\b)(\c\otimes\d)=\b\cdot\c(\a\otimes\d)\ \ \forall \a,\b,\c,\d\in\Ve,\\
&\A(\a\otimes\b)=(\A\a)\otimes\b\ \ \forall \a,\b\in\Ve,\ \A\in Lin(\Ve).
\end{split}
\ee

Finally, the symbol $\L^2$ is normally used to denote, in short, the product $\L\L,\ \forall\L\in Lin(\Ve)$.

\section{Transpose, symmetric and skew tensors}
For any tensor $\L\in Lin(\Ve)$, there exists just one tensor $\L^\top$, called the {\it transpose} of $\L$, such that
\be
\label{eq:transposedef}
\bu\cdot\L\bv=\bv\cdot\L^\top\bu\ \ \forall\bu,\bv\in\Ve.
\ee
The transpose of the transpose of $\L$ is $\L$:
\bes
\bu\cdot\L\bv=\bv\cdot\L^\top\bu=\bu\cdot(\L^\top)^\top\bv\ \ \Rightarrow\ \ (\L^\top)^\top=\L.
\ees
The Cartesian components of $\L^\top$ are obtained by swapping the indices of the components of $\L$:
\bes
L^\top_{ij}=\e_i\cdot\L^\top\e_j=\e_j\cdot(\L^\top)^\top\e_i=\e_j\cdot\L\e_i=L_{ji}.
\ees
It is immediate to show that 
\bes
(\A+\B)^\top=\A^\top+\B^\top\ \ \forall\A,\B\in Lin(\Ve),
\ees
while
\bes
\bu\cdot(\A\B)\bv=\B \bv\cdot\A^\top\bu=\bv\cdot\B^\top\A^\top\bu\ \ \Rightarrow\ \ (\A\B)^\top=\B^\top\A^\top.
\ees
Moreover,
\be
\label{eq:transposedyad}
\bu\cdot(\a\otimes\b)\bv=\a\cdot\bu\ \b\cdot\bv=\bv\cdot(\b\otimes\a)\bu\ \ \Rightarrow\ \ (\a\otimes\b)^\top=\b\otimes\a.
\ee

A tensor $\L$ is {\it symmetric} $\iff$
\bes
\L=\L^\top.
\ees
In such a case, because $L_{ij}=L^\top_{ij}$, we have
\bes
L_{ij}=L_{ji}.
\ees
A symmetric tensor is hence represented, in a given basis, by a symmetric matrix and has just six independent Cartesian components. Applying Eq. (\ref{eq:transposedef}) to $\I$, it is immediately recognized that  the identity tensor is symmetric: $\I=\I^\top$.

A tensor $\L$ is {\it antisymmetric} or {\it skew} $\iff$
\bes
\L=-\L^\top.
\ees
In this a case, because $L_{ij}=-L^\top_{ij}$, we have (no summation on the index $i$, see footnote \ref{note:sumeinstein}, Chapter \ref{ch:1})
\bes
L_{ij}=-L_{ji}\  \Rightarrow\ L_{\underline{i}\underline{i}}=0\ \forall i=1,2,3.
\ees
A skew tensor is hence represented, in a given basis, by an antisymmetric matrix whose components on the diagonal are identically null in any basis; finally, a skew tensor  only depends upon  three independent Cartesian components.

If we denote by $Sym(\Ve)$ the set of all the symmetric tensors and by $Skw(\Ve)$ that of all the skew tensors, then it is evident that, $\forall\alpha,\beta,\lambda,\mu\in\R$,
\bes
\besp
Sym(\Ve) \cap Sk&w(\Ve)=\O,\\
\alpha\A+\beta\B\in Sym(\Ve)\ &\forall\A,\B\in Sym(\Ve),\\
\lambda\L+\mu\M\in Skw(\Ve)\ &\forall\L,\M\in Skw(\Ve),\\
\end{split}
\ees 
so $Sym(\Ve)$ and $Skw(\Ve)$ are vector subspaces of $Lin(\Ve)$ with dim$(Sym(\Ve))=6$, while   dim$(Skw(\Ve))=3$.

Any tensor $\L$ can be decomposed into the sum of a symmetric, $\L^s$, and an antisymmetric, $\L^a$, tensor:
\bes
\L=\L^s+\L^a, 
\ees
with
\bes
\L^s=\frac{\L+\L^\top}{2}\ \ \in Sym(\Ve)
\ees
and
\bes
\L^a=\frac{\L-\L^\top}{2}\ \ \in Skw(\Ve),
\ees
so that, finally,
\bes
Lin(\Ve)=Sym(\Ve)\oplus Skw(\Ve).
\ees

\section{Trace, scalar product of tensors}
There exists one and only one linear form
\bes
\tr:Lin(\Ve)\rightarrow\R,
\ees
called the {\it trace}, such that
\bes
\tr(\a\otimes\b)=\a\cdot\b\ \ \forall\a,\b\in\Ve.
\ees
For its same definition, that has been given without making use of any basis of $\Ve$, the trace of a tensor is a {\it tensor invariant}, i.e. a quantity, extracted from a tensor, that does not depend upon the basis. 

Linearity implies that
\bes
\tr(\alpha\A+\beta\B)=\alpha\tr\A+\beta\tr\B\ \ \forall\alpha,\beta\in\R,\ \A,\B\in Lin(\Ve).
\ees
It is just  linearity to give the rule for calculating the trace of a tensor $\L$:
\be
\label{eq:tracciaL}
\tr\L=\tr(L_{ij}\e_i\otimes\e_j)=L_{ij}\tr(\e_i\otimes\e_j)=L_{ij}\ \e_i\cdot\e_j=L_{ij}\delta_{ij}=L_{ii}.
\ee
A tensor is hence an operator whose sum of the components on the diagonal,
\bes
 \tr\L=L_{11}+L_{22}+L_{33},
\ees 
is constant, regardless of  the basis.

Following the same procedure above, it is readily seen that
\bes
\tr\L^\top=\tr\L,
\ees
which implies, by linearity, that 
\be
\label{eq:tracciaskew}
\tr \L=0\ \ \forall \L\in Skw(\Ve).
\ee

The {\it scalar product} of tensors $\A$ and $\B$ is the positive definite, symmetric bilinear form defined by
\bes
\A\cdot\B=\tr(\A^\top\B).
\ees
This definition implies that, $\forall\L,\M,\N\in Lin(\Ve),\ \alpha,\beta\in\R$,
\bes
\besp
&\L\cdot(\alpha\M+\beta\N)=\alpha\L\cdot\M+\beta\L\cdot\N,\\
&(\alpha\L+\beta\M)\cdot\N=\alpha\L\cdot\N+\beta\M\cdot\N,\\
&\L\cdot\M=\M\cdot\L,\\
&\L\cdot\L>0\ \ \forall\L\in Lin(\Ve),\ \ \L\cdot\L=0\ \iff\ \L=\O.
\end{split}
\ees
These properties give the rule for computing the scalar product of two tensors $\A$ and $\B$:
\bes
\besp
\A\cdot\B&=A_{ij}(\e_i\otimes\e_j)\cdot B_{hk}(\e_h\otimes\e_k)=A_{ij}B_{hk}(\e_i\otimes\e_j)\cdot(\e_h\otimes\e_k)\\
&=A_{ij}B_{hk}\ \tr[(\e_i\otimes\e_j)^\top(\e_h\otimes\e_k)]=A_{ij}B_{hk}\ \tr[(\e_j\otimes\e_i)(\e_h\otimes\e_k)]\\
&=A_{ij}B_{hk}\ \tr[\e_i\cdot\e_h(\e_j\otimes\e_k)]=A_{ij}B_{hk}\ \e_i\cdot\e_h\ \e_j\cdot\e_k\\
&=A_{ij}B_{hk}\delta_{ih}\delta_{jk}=A_{ij}B_{ij}.
\end{split}
\ees
As in the case of vectors, the scalar product of two tensors is equal to the sum of the products of the corresponding components. 
In a similar manner, or using Eq. (\ref{eq:propdyad0})$_1$, it is easily shown that, $\forall\a,\b,\c,\d\in\Ve$,
\bes
(\a\otimes\b)\cdot(\c\otimes\d)=\a\cdot\c\ \b\cdot\d=a_ib_jc_id_j,
\ees
while by the same definition of the tensor scalar product,
\bes
\tr\L=\I\cdot\L\ \ \forall\L\in Lin(\Ve).
\ees
Similar to vectors, we define {\it Euclidean norm} of a tensor $\L$ the nonnegative scalar, denoted either by $L$ or $|\L|$,
\bes
L=|\L|=\sqrt{\L\cdot\L}=\sqrt{\tr(\L^\top\L)}=\sqrt{L_{ij}L_{ij}},
\ees
and the {\it distance} $d(\L,\M)$ of two tensors $\L$ and $\M$ the norm of the tensor difference:
\bes
d(\L,\M):=|\L-\M|=|\M-\L|.
\ees

\section{Spherical and deviatoric parts}
Let $\L\in Sym(\Ve)$; the {\it spherical part} of $\L$ is defined by
\bes
\L^{sph}:=\frac{1}{3}\tr\L\ \I,
\ees
and the {\it deviatoric part} by
\bes
\L^{dev}:=\L-\L^{sph},
\ees
so that
\bes
\L=\L^{sph}+\L^{dev}.
\ees
We remark that 
\bes
\tr\L^{sph}=\frac{1}{3}\tr\L\ \tr\I=\tr\L\ \Rightarrow\ \tr\L^{dev}=0,
\ees
i.e. the deviatoric part is a traceless tensor.
Let $\A,\B\in Lin(\Ve)$; then
\be
\label{eq:orthsphdev}
\A^{sph}\cdot\B^{dev}=\frac{1}{3}\tr\A\ \I\cdot\B^{dev}=\frac{1}{3}\tr\A\ \tr\B^{dev}=0,
\ee
i.e. any spherical tensor is orthogonal to any deviatoric tensor.

The sets
\bes
\besp
&Sph(\Ve):=\left\{\A^{sph}\in Lin(\Ve)|\ \A^{sph}=\frac{1}{3}\tr\A\I\ \forall\A\in Lin(\Ve)\right\},\\
&Dev(\Ve):=\left\{\A^{dev}\in Lin(\Ve)|\ \A^{dev}=\A-\A^{sph}\ \forall\A\in Lin(\Ve)\right\}
\end{split}
\ees
form two subspaces of $Lin(\Ve)$; the proof is left to the reader. For what is proved above, $Sph(\Ve)$ and $Dev(\Ve)$ are two {\it mutually orthogonal subspaces} of $Lin(\Ve)$.

\section{Determinant, inverse of a tensor}
The reader is probably familiar with the concept of determinant of a matrix. We show here that the determinant of a second-rank tensor can be defined intrinsically and that it corresponds with the determinant of the matrix that represents it in any basis of $\Ve$. 
For this purpose, we first need  to introduce a  mapping:
\bes
\omega:\Ve\times\Ve\times\Ve\rightarrow\R
\ees
is a {\it skew trilinear form} if $\omega(\bu,\bv,\cdot),\omega(\bu,\cdot,\bv)$ and $\omega(\cdot,\bu,\bv)$ are linear forms on $\Ve$ and if
\be
\label{eq:trilin1}
\omega(\bu,\bv,\bw)=-\omega(\bv,\bu,\bw)=-\omega(\bu,\bw,\bv)=-\omega(\bw,\bv,\bu)\ \forall\bu,\bv,\bw\in\Ve.
\ee

Using this definition, we can state the following
\begin{teo}
\label{teo:trilin1}
Three vectors are linearly independent if and only if every skew trilinear form on them is not null.
\begin{proof}
In fact, let $\bu=\alpha\bv+\beta\bw$, then for any skew trilinear form $\omega$,
\bes
\omega(\bu,\bv,\bw)=\omega(\alpha\bv+\beta\bw,\bv,\bw)=\alpha\omega(\bv,\bv,\bw)+\beta\omega(\bw,\bv,\bw)=0
\ees
because of Eq. (\ref{eq:trilin1}) applied to the permutation of the positions of the two $\bv$ and the two $\bw$.
\end{proof}
\end{teo}

It is evident that the set of all the skew trilinear forms  is a vector space, that we denote by $\Omega$, whose null element is the {\it null form} $\omega_0$:
\bes
\omega_0(\bu,\bv,\bw)=0\ \forall\bu,\bv,\bw\in\Ve.
\ees
For a given $\omega(\bu,\bv,\bw)\in\Omega$, any $\L\in Lin(\Ve)$ induces another form $\omega_L(\bu,\bv,\bw)\in\Omega$, defined as
\bes
\omega_L(\bu,\bv,\bw)=\omega(\L\bu,\L\bv,\L\bw)\ \forall\bu,\bv,\bw\in\Ve.
\ees
A key point\footnote{The proof of this statement is rather involved and outside of our scope; the interested reader is referred to the classical textbook by Halmos on linear algebra, Section 31 (see the bibiography). The theory of the determinants is developed in Section 53.} for the following developments is that $\dim\Omega=1$. 

This means that $\forall\omega_1,\omega_2\neq\omega_0\in\Omega,\exists\lambda\in\R$ such that
\bes
\omega_2(\bu,\bv,\bw)=\lambda\omega_1(\bu,\bv,\bw)\ \forall\bu,\bv,\bw\in\Ve.
\ees
So, $\forall\L\in Lin(\Ve)$, there must exist $\lambda_L\in\R$ such that
\be
\label{eq:trilin3}
\omega(\L\bu,\L\bv,\L\bw)=\omega_L(\bu,\bv,\bw)=\lambda_L\ \omega(\bu,\bv,\bw)\ \forall\bu,\bv,\bw\in\Ve.
\ee
The scalar\footnote{More precisely, $\det\L$ is the function that associates a scalar with each tensor (Halmos, Section 53). We can, however, for the sake of practice, identify $\det\L$ with the scalar associated with $\L$, without  consequences for our purposes.} $\lambda_L$ is the {\it determinant of }$\L$ and in the following it will be denoted as $\det\L$. The determinant of a tensor $\L$ is an intrinsic quantity of $\L$, i.e. it does not depend upon the particular form $\omega$, nor on the basis of $\Ve$. In fact, we have never introduced, so far, a basis for defining $\det\L$, hence it cannot depend upon the choice of a basis for $\Ve$, i.e. $\det\L$ is a {\it tensor invariant}. 

Then, if $\omega^a$ and $\omega^b\in\Omega$, because $\dim\Omega=1$, there exists $k\in\R,\ k\neq0$ such that
\bes
\besp
&\omega^b(\bu,\bv,\bw)=k\ \omega^a(\bu,\bv,\bw) \ \forall\bu,\bv,\bw\in\Ve\Rightarrow\\
&\omega^b(\L\bu,\L\bv,\L\bw)=k\ \omega^a(\L\bu,\L\bv,\L\bw)\rightarrow\\
&\omega^b_L(\bu,\bv,\bw)=k\ \omega^a_L(\bu,\bv,\bw).
\end{split}
\ees
Moreover, by Eq. (\ref{eq:trilin3}) we get
\bes
\besp
&\omega^a(\L\bu,\L\bv,\L\bw)=\omega^a_L(\bu,\bv,\bw)=\lambda^a_L\omega^a(\bu,\bv,\bw),\\
&\omega^b(\L\bu,\L\bv,\L\bw)=\omega^b_L(\bu,\bv,\bw)=\lambda^b_L\omega^b(\bu,\bv,\bw),
\end{split}
\ees
so that
\bes
\besp
&\lambda^b_Lk\ \omega^a(\bu,\bv,\bw)=\lambda^b_L\omega^b(\bu,\bv,\bw)=\omega^b_L(\bu,\bv,\bw)=\\
&k\ \omega^a_L(\bu,\bv,\bw)=\lambda^a_Lk\ \omega^a(\bu,\bv,\bw)\iff\lambda^a_L=\lambda^b_L,
\end{split}
\ees
which proves that $\det\L$ does not depend upon the skew trilinear form, but only upon $\L$.

The definition given for $\det\L$ allows us to prove some important properties. First of all,
\bes
\det\O=0;
\ees
in fact, $\forall\omega\in\Omega$, 
\bes
\det\O\ \omega(\bu,\bv,\bw)=\omega(\O\bu,\O\bv,\O\bw)=\omega(\bo,\bo,\bo)=0\ \forall\bu,\bv,\bw\in\Ve
\ees
because $\omega$ operates on three identical, i.e. linearly dependent, vectors.
Moreover, if $\L=\I$, then
\bes
\det\I\ \omega(\bu,\bv,\bw)=\omega(\I\bu,\I\bv,\I\bw)=\omega(\bu,\bv,\bw)
\ees
if and only if
\be
\label{eq:detident}
\det\I=1.
\ee
A third property is that $\forall\a,\b\in\Ve$,
\be
\label{eq:trilin4}
\det(\a\otimes\b)=0.
\ee
In fact, if $\L=\a\otimes\b$, then
\bes
\det\L\ \omega(\bu,\bv,\bw)=\omega(\L\bu,\L\bv,\L\bw)=\omega((\b\cdot\bu)\a,(\b\cdot\bv)\a,(\b\cdot\bw)\a)=0
\ees
because the three vectors on which $\omega\in\Omega$ operates are linearly dependent; because $\bu,\bv$ and $\bw$ are arbitrarily chosen, this   implies Eq. (\ref{eq:trilin4}).

An important result is the 
\begin{teo}
\label{teo:binet}
{\bf (Theorem of Binet).} $\forall\A,\B\in Lin(\Ve)$
\be
\label{eq:thbinet}
\det(\A\B)=\det\A\det\B.
\ee
\begin{proof}
$\forall\omega\in\Omega$ and $\forall\bu,\bv,\bw\in\Ve$,
\bes
\besp
&\lambda_{AB}\omega(\bu,\bv,\bw)=\omega(\A\B\bu,\A\B\bv,\A\B\bw)=\omega(\A(\B\bu),\A(\B\bv),\A(\B\bw))=\\
&\lambda_A\omega(\B\bu,\B\bv,\B\bw)=\lambda_A\lambda_B\omega(\bu,\bv,\bw)\iff\lambda_{AB}=\lambda_A\lambda_B,
\end{split}
\ees
which proves the theorem.
\end{proof}

\end{teo}

A tensor $\L$ is called {\it singular} if $\det\L=0$, otherwise it is {\it non-singular}.

Considering Eq. (\ref{eq:trilin3}), with some effort but without major difficulties, one can  see that, if in a basis $\Ba$ of $\Ve$ it is $\L=L_{ij}\e_i\otimes\e_j$, then
\bes
\det\L=\sum_{\pi\in\mathcal{P}_3}\epsilon_{\pi(1),\pi(2),\pi(3)}L_{1,\pi(1)}L_{2,\pi(2)}L_{3,\pi(3)},
\ees
where $\mathcal{P}_3$ is the set of all the permutations $\pi$ of $\{1,2,3\}$ and the $\epsilon_{i,j,k}$s are the components of the {\it Ricci's alternator}\footnote{
We recall that a {\it permutation} of an ordered set of $n$ objects is {\it even} if it can be obtained as the product of an even number of {\it transpositions}, i.e. exchange of places, of any couple of its elements, it is {\it odd} if the number of transpositions is odd. For the set $\{1,2,3\}$ the even permutations are $\{1,2,3\}, \{3,1,2\},\{2,3,1\}$, while the odd ones are $\{2,1,3\},\{1,3,2\},\{3,2,1\}$;  any triplet having at least a repeated number is not a permutation.}:
\bes
\epsilon_{i,j,k}:=
\left\{
\begin{array}{rll}  
1& \mathrm{if}\ \{i,j,k\}  &\mathrm{is\ an\ even\ permutation\ of} \ \{1,2,3\},   \\  
0& \mathrm{if}\ \{i,j,k\}  & \mathrm{is\ not\ a\ permutation\ of} \ \{1,2,3\}  \\  
-1& \mathrm{if}\ \{i,j,k\}  & \mathrm{is\ an\ odd\ permutation\ of} \ \{1,2,3\}.
 \end{array}
\right.
\ees

The above rule for  $\det\L$  coincides with that for calculating the determinant of the matrix whose entries are the $L_{ij}$s. This shows that, once chosen a basis $\Ba$ for $\Ve$,  $\det\L$ coincides with the determinant of the matrix representing it in  $\Ba$, and finally that
\be
\label{eq:detformula}
\besp
\det \L&=L_{11}L_{22}L_{33}+L_{12}L_{23}L_{31}+L_{13}L_{32}L_{21}\\
&-L_{11}L_{23}L_{32}-L_{22}L_{13}L_{31}-L_{33}L_{12}L_{21}.
\end{split}
\ee
This result shows immediately that $\forall\L\in Lin(\Ve)$, and regardless of $\Ba$, we have
\be
\label{eq:dettransp}
\det\L^\top=\det\L.
\ee
Using Eq. (\ref{eq:detformula}), it is not difficult to show that, $\forall\alpha\in\R$,
\be
\label{eq:detsum}
\det(\I+\alpha\L)=1+\alpha I_1+\alpha^2 I_2+\alpha^3 I_3,
\ee
where $I_1,I_2$ and $I_3$ are the three {\it principal invariants} of $\L$:
\be
\label{eq:princinv}
I_1=\tr \L,\ \ I_2=\frac{\tr^2\L-\tr\L^2}{2},\ \ I_3=\det\L.
\ee
A tensor $\L\in Lin(\Ve)$ is said to be {\it invertible} if there exists a tensor $\L^{-1}\in Lin(\Ve)$, called the {\it inverse} of $\L$, such that
\be
\label{eq:definvert}
\L\L^{-1}=\L^{-1}\L=\I.
\ee
If $\L$ is invertible, then $\L^{-1}$ is unique. By the above definition, if $\L$ is invertible, then
\bes
\bu_1=\L\bu\Rightarrow\bu=\L^{-1}\bu_1.
\ees
\begin{teo} Any invertible tensor maps triples of linearly independent vectors into triples of still linearly independent vectors.
\begin{proof}
Let $\L$ be an invertible tensor and  $\bu_1=\L\bu,\bv_1=\L\bv,\bw_1=\L\bw$, where $\bu,\bv,\bw$ are three linearly independent vectors.
Let us suppose that there exist $h,k\in\R$ such that
\bes
\bu_1=h\bv_1+k\bw_1.
\ees
Then, because $\L$ is invertible,
\bes
\L^{-1}\bu_1=\L^{-1}(h\bv_1+k\bw_1)=h\L^{-1}\bv_1+k\L^{-1}\bw_1=h\bv+k\bw,
\ees
which goes against the hypothesis. Consequently, $\bu_1,\bv_1$ and $\bw_1$ are linearly independent.
\end{proof}
\end{teo}
This result, along with the definition of determinant, Eq. (\ref{eq:trilin3}), and Theorem \ref{teo:trilin1}, proves the 
\begin{teo}{\bf (Invertibility theorem).} $\L\in Lin(\Ve)$ is invertible $\iff\det\L\neq0$.
\end{teo}
Using the theorem of Binet, \ref{teo:binet}, along with Eqs. (\ref{eq:detident}) and (\ref{eq:definvert}), we get
\bes
\det\L^{-1}=\frac{1}{\det\L}.
\ees 
Equation (\ref{eq:definvert}) applied to $\L^{-1}$, along with the uniqueness of the inverse, gives immediately that
\bes
(\L^{-1})^{-1}=\L,
\ees
while
\bes
\B^{-1}\A^{-1}=\B^{-1}\A^{-1}\A\B(\A\B)^{-1}=(\A\B)^{-1}.
\ees
The operations of transpose and inversion commute:
\bes
\begin{array}{c}
\L^\top(\L^\top)^{-1}=\I=\L^{-1}\L=\I^\top=(\L^{-1}\L)^\top=\L^\top(\L^{-1})^\top\Rightarrow\medskip\\
(\L^{-1})^\top=(\L^\top)^{-1}:=\L^{-\top}.
\end{array}
\ees

\section {Eigenvalues and eigenvectors of a tensor}
\label{sec:eigenvaluesvectors}
If there exists a $\lambda\in\R$ and a $\bv\in\Ve$, except the null vector, such that
\be
\label{eq:eig1}
\L\bv=\lambda\bv,
\ee
then $\lambda$ is an {\it eigenvalue} and $\bv$ an {\it eigenvector}, relatif to $\lambda$, of $\L$. It is immediate to observe that, thanks to linearity, any eigenvector $\bv$ of $\L$ is determined to within a multiplier, i.e. that $k\bv$ is an eigenvector of $\L$ too $\forall k\in\R$. Often, the multiplier $k$ is fixed in such a way that $|\bv|=1$.

To determine the eigenvalues and eigenvectors of a tensor, we rewrite Eq. (\ref{eq:eig1}) as
\be
\label{eq:eig2}
(\L-\lambda\I)\bv=\bo.
\ee
The condition for this homogeneous system having a non null solution is
\bes
\det(\L-\lambda\I)=0;
\ees
this is the so-called {\it characteristic} or {\it Laplace's equation}. In the case of a second-rank tensor over $\Ve$, the Laplace's equation is an algebraic equation of degree three with real coefficients. The roots of the Laplace's equation are the eigenvalues of $\L$; because the components of $\L$, and hence the coefficients of the characteristic equation, are all real, then the eigenvalues of $\L$ are all real or one real and two complex conjugate. 

For any eigenvalue $\lambda_i,\ i=1,2,3$, of $\L$, the corresponding eigenvectors $\bv_i$ can be found solving Eq. (\ref{eq:eig2}), once set $\lambda=\lambda_i$.

The {\it proper space} of $\L$ relatif to $\lambda$ is the subspace of $\Ve$ composed of all the vectors that satisfy Eq. (\ref{eq:eig2}). The {\it multiplicity} of $\lambda$ is the dimension of its proper space, while the {\it spectrum} of $\L$ is the set composed by all of its eigenvalues, each one with its multiplicity.

$\L^\top$ has the same eigenvalues of $\L$, because the Laplace's equation is the same in both the cases:
\bes
\det(\L^\top-\lambda\I)=\det(\L^\top-\lambda\I^\top)=\det(\L-\lambda\I)^\top=\det(\L-\lambda\I).
\ees
However, this is not the case for the eigenvectors, that are generally  different, as a numerical example can show.

Developing the Laplace's equation, it is easy to show that it can be written as
\bes
\det(\L-\lambda\I)=-\lambda^3+I_1\lambda^2-I_2\lambda+I_3=0,
\ees
which is merely an application of Eq. (\ref{eq:detsum}). If we denote $\L^3=\L\L\L$, using Eq. (\ref{eq:princinv}) one can prove the {\it Cayley-Hamilton theorem}:
\begin{teo}
{\bf (Cayley-Hamilton theorem).} $\forall \L\in Lin(\Ve)$,
\bes
\L^3-I_1\L^2+I_2\L-I_3\I=\O.
\ees
\end{teo}
A {\it quadratic form}  defined by $\L$ is any form $\omega:\Ve\times\Ve\rightarrow\R$  of the type
\bes
\omega=\bv\cdot\L\bv;
\ees
if $\omega>0\ \forall\bv\in\Ve,\ \omega=0\iff\bv=\bo$, then $\omega$ and $\L$ are said to be {\it positive definite}. The eigenvalues of a positive definite tensor are positive. In fact, if $\lambda$ is an eigenvalue of $\L$, which is positive definite, and  $\bv$  its eigenvector, then 
\bes
\bv\cdot\L\bv=\bv\cdot\lambda\bv=\lambda\bv^2>0\iff\lambda>0.
\ees
Let $\bv_1$ and $\bv_2$ be two eigenvectors of a symmetric tensor $\L$ relative to the eigenvalues $\lambda_1$ and $\lambda_2$, respectively, with $\lambda_1\neq\lambda_2$. Then
\bes
\lambda_1\bv_1\cdot\bv_2=\L\bv_1\cdot\bv_2=\L\bv_2\cdot\bv_1=\lambda_2\bv_2\cdot\bv_1\iff\bv_1\cdot\bv_2=0.
\ees
Actually, symmetric tensors have a particular importance, specified by the {\it spectral theorem}:
\begin{teo}
{\bf(Spectral theorem).} The eigenvectors of a symmetric tensor form a basis of $\Ve$.
\end{teo}
This theorem\footnote{The proof of the spectral theorem is omitted here; the interested reader can find a proof of it in the classical text by Halmos, p. 155, see the suggested texts.} is of  paramount importance in linear algebra: It proves that the eigenvalues of a symmetric tensor $\L$ are real valued and, remembering the definition of eigenvalues and eigenvectors, Eq. (\ref{eq:eig1}), that there exists a basis $\Ba_N=\{\bu_1,\bu_2,\bu_3\}$ of $\Ve$ composed of eigenvectors of $\L$, i.e. by vectors that are mutually orthogonal and that remain mutually orthogonal once transformed by $\L$. Such a basis is called the {\it normal basis}.

If  $\lambda_i,i=1,2,3,$ are the eigenvalues of $\L$, then the components of $\L$ in $\Ba_N$ are
\bes
L_{ij}=\bu_i\cdot\L\bu_j=\bu_i\cdot\lambda_j\bu_j=\lambda_j\delta_{ij}
\ees
so finally in $\Ba_N$ we have
\bes
\L=\lambda_i\e_i\otimes\e_i,
\ees
i.e. $\L$ is diagonal and  is completely represented by its eigenvalues. In addition, it is easy to check that
\bes
I_1=\lambda_1+\lambda_2+\lambda_3,\ I_2=\lambda_1\lambda_2+\lambda_2\lambda_3+\lambda_3\lambda_1,\ I_3=\lambda_1\lambda_2\lambda_3.
\ees
A tensor with a unique eigenvalue $\lambda$ of multiplicity three is said to be {\it spherical}; in such a case, any basis of $\Ve$ is $\Ba_N$ and 
\bes
\L=\lambda\I.
\ees
Eigenvalues and eigenvectors have also another important property: Let us consider the quadratic form $\omega:=\bv\cdot\L\bv,\ \forall\bv\in\S$, defined by a symmetric tensor $\L$. We look for the directions $\bv\in\S$ whereupon $\omega$ is stationary. Then, we have to solve the constrained problem
\bes
\nabla_{\bv}(\bv\cdot\L\bv)=\bo,\ \ \bv\in\S.
\ees
Using the Lagrange's multiplier technique, we solve the equivalent problem
\bes
\nabla_{(\bv,\lambda)}(\bv\cdot\L\bv-\lambda(\bv^2-1))=0,
\ees
which restitutes the equation
\bes
\L\bv=\lambda\bv
\ees
and the constraint $|\bv|=1$. The above equation is exactly the one defining the eigenvalue problem for $\L$: The stationary values (i.e. the maximum and minimum) of $\omega$ hence corresponds  to two eigenvalues of $\L$ and the directions $\bv$, whereupon stationarity is get, coincide with the respective eigenvectors.

Two tensors $\A$ and $\B$ are said to be {\it coaxial} if they have the same normal basis $\Ba_N$, i.e. if they share the same eigenvectors. Let $\bu$ be an eigenvector of $\A$, relative to the eigenvalue $\lambda_A$, and of $\B$, relatif to $\lambda_B$. Then,
\bes
\A\B\bu=\A\lambda_B\bu=\lambda_B\A\bu=\lambda_A\lambda_B\bu=\lambda_A\B\bu=\B\lambda_A\bu=\B\A\bu,
\ees
which shows, on the one hand, that also $\B\bu$ is an eigenvector of $\A$, relative to the same eigenvalue $\lambda_A$; in the same way, of course, $\A\bu$ is an eigenvector of $\B$ relative to $\lambda_B$. In other words, this shows that $\B$ leaves unchanged any proper space of $\A$ and vice versa. On the other hand, we see that, at least for what concerns the eigenvectors, two tensors commute if and only if they are coaxial. Because any vector can be written as a linear combination of the vectors of $\Ba_N$, and for the linearity of tensors, we finally have proved the {\it commutation theorem}:

\begin{teo}
\label{teo: commutativity}
{\bf(Commutation theorem).} Two tensors commute if and only if they are coaxial.
\end{teo}

\section{Skew tensors and cross product}
\label{sec:crossprod}
Because $\dim(\Ve)=\dim(Skw(\Ve))=3$, an isomorphism can be established between $\Ve$ and $Skw(\Ve)$, i.e. between vectors and skew tensors. We establish hence a way to associate in a unique way a vector to any skew tensor and inversely. For this purpose, we first introduce the following theorem: 
\begin{teo} The spectrum of any tensor $\W\in Skw(\Ve)$ is $\{0\}$ and the dimension of its proper space is 1.
\label{teo:spectreskew}
\begin{proof}
This theorem states that zero is the only real eigenvalue of any skew tensor  and that its multiplicity is 1.  In fact, let $\bw$ be an eigenvector of $\W$ relative to the eigenvalue $\lambda$. Then
\bes
\besp
\lambda^2\bw^2&=\W\bw\cdot\W\bw=\bw\cdot\W^\top\W\bw=-\bw\cdot\W\W\bw\\
&=-\bw\cdot\W(\lambda\bw)=-\lambda\bw\cdot\W\bw=-\lambda^2\bw^2\iff\lambda=0.
\end{split}
\ees
Then, if $\W\neq\O$ its rank is necessarily 2, because $\det\W=0\ \forall\W\in Skw(\Ve)$; hence, the equation 
\be
\label{eq:isomorph}
\W\bw=\bo
\ee
has $\infty^1$ solutions, i.e. the multiplicity of $\lambda$ is 1, which proves the theorem. 
\end{proof}
\end{teo}
The last equation   also shows the way the isomorphism is constructed: In fact, using Eq. (\ref{eq:isomorph}) it is easy to check that if $\bw=(a,b,c)$, then
\be
\label{eq:isomorphskew}
\bw=(a,b,c)\iff\W=\left[\begin{array}{ccc}0 & -c & b \\c & 0 & -a\\-b & a & 0\end{array}\right].
\ee
The proper space of $\W$ is called the {\it axis of }$\W$ and it is indicated by $\mathcal{A}(\W)$:
\bes
\mathcal{A}(\W):=\{\bu\in\Ve|\ \W\bu=\bo\}.
\ees
The consequence of what shown above is that $\dim\mathcal{A}(\W)=1$. With regard to Eq. (\ref{eq:isomorphskew}), one can easily check  that  the equation
\be
\label{eq:axtensvect1}
\bu\cdot\bu=\frac{1}{2}\W\cdot\W
\ee
is satisfied only by $\bw$ and by its opposite $-\bw$.
Because both these vectors belong to $\mathcal{A}(\W)$, choosing one of them corresponds to choose an orientation for $\Eu$, see the next section. We  always make our choice according to Eq. (\ref{eq:isomorphskew}), which fixes once and for all the isomorphism between $\Ve$ and $Skw(\Ve)$ that makes correspond any vector $\bw$ with one and only one {\it axial tensor} $\W$ and vice-versa, any skew tensor $\W$ with a unique {\it axial vector} $\bw$. 

It is worth noting that the above isomorphism between the vector spaces $\Ve$ and $Skw(\Ve)$ implies that to any linear combination of vectors $\a$ and $\b$ corresponds an equal linear combination of the corresponding axial tensors $\W_a$ and $\W_b$ and vice-versa, i.e. $\forall a,b\in\R$
\be
\label{eq:lincombax}
\bw=\alpha\a+\beta\b\iff\W=\alpha\W_a+\beta\W_b,
\ee
where $\W$ is the axial tensor of $\bw$. Such a property is immediately checked using Eq. (\ref{eq:isomorphskew}).

We define {\it cross product} of two vectors $\a$ and $\b$ the vector
\bes
\a\times\b=\W_a\b,
\ees
where $\W_a$ is the axial tensor of $\a$. If $\a=(a_1,a_2,a_3)$ and $\b=(b_1,b_2,b_3)$, then by Eq. (\ref{eq:isomorphskew}) we get
\bes
\a\times\b=(a_2b_3-a_3b_2,a_3b_1-a_1b_3,a_1b_2-a_2b_1).
\ees
It is immediate to check that such a result can also be obtained  using  Ricci's alternator
\be
\label{eq:prodvectricci}
\a\times\b=\epsilon_{ijk}a_jb_k\ei,
\ee
or even computing the symbolic determinant
\bes
\a\times\b=\det\left[\begin{array}{ccc}\eu & \ed & \et \\a_1 & a_2 & a_3 \\b_1 & b_2 & b_3\end{array}\right].
\ees

The cross product is bilinear: $\forall\a,\b,\bu\in\Ve,\ \alpha,\beta\in\R$,
\bes
\besp
&(\alpha\a+\beta\b)\times\bu=\alpha\a\times\bu+\beta\b\times\bu,\\
&\bu\times(\alpha\a+\beta\b)=\alpha\bu\times\a+\beta\bu\times\b.
\end{split}
\ees
In fact, the first equation above is a consequence of Eq. (\ref{eq:lincombax}), while the second one is a simple application to axial tensors of the same definition of  tensor.

Three important results concerning the cross product are stated by the following theorems.
\begin{teo}{\bf(Condition of parallelism).} Two vectors $\a$ and $\b$ are parallel, i.e. $\b=k\a$, $k\in\R,\iff$
\bes
\a\times\b=\bo.
\ees
\begin{proof}
This property is actually a consequence of the fact that any eigenvalue of a tensor is determined to within a multiplier:
\bes
\a\times\b=\W_a\b=\bo\iff\b=k\a,\ k\in\R,
\ees
for Theorem \ref{teo:spectreskew}.
\end{proof}
\end{teo}
\begin{teo}{\bf(Orthogonality property).}
\label{teo:orthprop}
\be
\label{eq:orthprop}
\a\times\b\cdot\a=\a\times\b\cdot\b=0.
\ee
\begin{proof}
\bes
\besp
&\a\times\b\cdot\a=\W_a\b\cdot\a=\b\cdot\W_a^\top\a=-\b\cdot\W_a\a=-\b\cdot\bo=0,\\
\a\times&\b\cdot\b=\W_a\b\cdot\b=\b\cdot\W_a^\top\b=-\b\cdot\W_a\b\iff\a\times\b\cdot\b=0.
\end{split}
\ees
\end{proof}
\end{teo}

\begin{teo}
\label{teo:axialvector}
$\a\times\b$ is the axial vector of the tensor $(\b\otimes\a-\a\otimes\b)$.
\begin{proof}
First of all, by Eq. (\ref{eq:transposedyad}) we see that
\bes
(\b\otimes\a-\a\otimes\b)\in Skew(\Ve).
\ees
Then,
\bes
(\b\otimes\a-\a\otimes\b)(\a\times\b)=\a\cdot\a\times\b\ \b-\b\cdot\a\times\b\ \a=0
\ees
for Theorem \ref{teo:orthprop}.
\end{proof}
\end{teo}

\begin{figure}[h]
\begin{center}
\includegraphics[scale=.5]{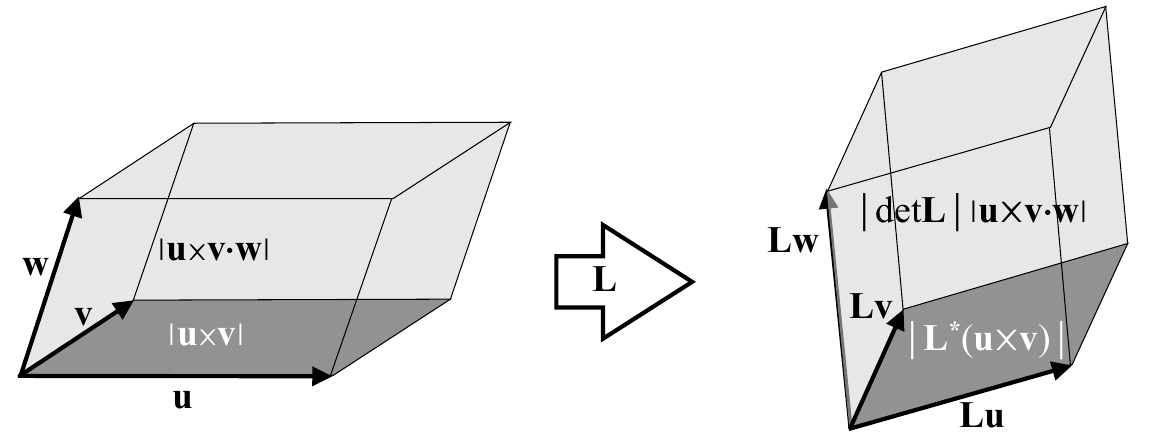}
\caption{Geometrical meaning of the cross and mixed products before (left) and after (right) the application of a tensor $\L$ on the vectors $\bu,\bv,\bw$.}
\label{fig:53}
\end{center}
\end{figure}

Theorem \ref{teo:axialvector} allows us to show another important result about cross product: the {\it antisymmetry of the cross product}:
\begin{teo}{\bf (Antisymmetry of the cross product).} The cross product is antisymmetric:
\be
\label{eq:antisymcrossprod}
\a\times\b=-\b\times\a \ \forall\a,\b\in\Ve.
\ee
\begin{proof}
Let $\W_1=(\b\otimes\a-\a\otimes\b)$ be  the axial tensor of $\a\times\b$ and  $\W_2=(-\a\otimes\b+\b\otimes\a)$  that of $-\b\times\a$. Evidently, $\W_1=\W_2$ which implies Eq. (\ref{eq:antisymcrossprod}) for the isomorphism between $\Ve$ and $Lin(\Ve)$.
\end{proof}
\end{teo}
This property and, again, Theorem \ref{teo:axialvector} lets us derive the formula for the {\it  double cross product}:
\be
\label{eq:doublecrossprod}
\bu\times(\bv\times\bw)=-(\bv\times\bw)\times\bu=-(\bw\otimes\bv-\bv\otimes\bw)\bu=
\bu\cdot\bw\ \bv-\bu\cdot\bv\ \bw.
\ee
Another interesting result concerns the {\it mixed product}:
\be
\label{eq:prodmixte}
\bu\times\bv\cdot\bw=\W_u\bv\cdot\bw=-\bv\cdot\W_u\bw=-\bv\cdot\bu\times\bw=\bw\times\bu\cdot\bv,
\ee
and similarly 
\bes
\bu\times\bv\cdot\bw=\bv\times\bw\cdot\bu.
\ees
Using this last result,  we can obtain a formula for the norm of a cross product; if $\a=a\ \e_a$ and $\b=b\ \e_b$, with $\e_a,\e_b\in\mathcal{S}$, are two vectors forming the angle $\theta$, then
\be
\label{eq:formulevectprod}
\besp
&(\a\times\b)\cdot(\a\times\b)=\a\times\b\cdot(\a\times\b)=(\a\times\b)\times\a\cdot\b=-\a\times(\a\times\b)\cdot\b=\\
&(-\a\cdot\b\ \a+\a^2\ \b)\cdot\b=\b\cdot(\a^2\I-\a\otimes\a)\b=a^2\ \b\cdot(\I-\e_a\otimes\e_a)\b=\\
&a^2b^2\ \e_b\cdot(\I-\e_a\otimes\e_a)\e_b=
a^2b^2(1-\cos^2\theta)=a^2b^2\sin^2\theta\rightarrow
|\a\times\b|=ab\sin\theta.
\end{split}
\ee
So, the norm of a cross product can be interpreted, geometrically, as the area of the parallelogram spanned by the two vectors. As a consequence, the absolute value of the mixed product (\ref{eq:prodmixte}) measures the volume of the prism delimited by three non coplanar vectors, cf. Fig. \ref{fig:53}.

Because the cross product is antisymmetric and the scalar one is symmetric, it is easy to check that the form
\bes
\beta(\bu,\bv,\bw)=\bu\times\bv\cdot\bw
\ees
is a skew trilinear form. Then, Eq. (\ref{eq:trilin3}), we get
\be
\label{eq:skewmixed}
\L\bu\times\L\bv\cdot\L\bw=\det\L\ \bu\times\bv\cdot\bw.
\ee
Following the interpretation given above for the absolute value of the mixed product, we can conclude that $|\det\L|$ can be interpreted as a coefficient of volume expansion\footnote{This result is classical and fundamental for the analysis of deformation in continuum mechanics.}, cf. again Fig. \ref{fig:53}. A geometrical interpretation can then be given to the case of a non invertible tensor, i.e. of $\det\L=0$: It crushes a prism into a flat region (the three original vectors become coplanar, i.e. linearly dependent).

The {\it adjugate} of $\L$ is the tensor 
\bes
\L^*:=(\det\L)\L^{-\top}.
\ees
From Eq. (\ref{eq:skewmixed}) we get hence
\bes
\begin{array}{c}
\det\L\ \bu\times\bv\cdot\bw=\L\bu\times\L\bv\cdot\L\bw=\L^\top(\L\bu\times\L\bv)\cdot\bw\ \ \forall\bw\Rightarrow \medskip\\\L\bu\times\L\bv=\L^*(\bu\times\bv).
\end{array}
\ees

It is useful, for further development, to calculate the powers of $\W$:
\be
\label{eq:w2sym}
{\W^2}=\W\W=-\W^\top(-\W^\top)=(\W\W)^\top{=(\W^2)^\top},
\ee
i.e., $\W^2$ is symmetric. Moreover, if we take $\bw\in\S$, which is always possible  because eigenvectors are determined to within an arbitrary multiplier, 
\be
\label{eq:skew2}
\besp
{\W^2\bu}&=\W\W\bu=\bw\times(\bw\times\bu)=\bw\cdot\bu\bw-\bw\cdot\bw\bu\\
&=-(\I-\bw\otimes\bw)\bu\ \Rightarrow\ {\W^2{=-(\I-\bw\otimes\bw)}};
\end{split}
\ee
We remark that $\W^2\bu$ gives the opposite of the projection of any vector $\bu\in\Ve$ onto the direction orthogonal to $\bw$, see Exercise 2.

Applying recursively the previous results,
\be
\label{eq:skew4}
\besp
&{\W^3}=\W\W^2=-\W(\I-\bw\otimes\bw)=-\W+(\W\bw)\otimes\bw{=-\W}\\
&{\W^4}=\W\W^3{=-\W^2}\\
&{\W^5}=\W\W^4{=-\W^3}\\
& etc.
\end{split}
\ee
An important property  of any couple axial tensor $\W-$axial vector $\bw\in\S$ is
\be
\label{eq:normW}
\W\W=-\frac{1}{2}|\W|^2(\I-\bw\otimes\bw),
\ee
while Eq. (\ref{eq:axtensvect1}) can be generalized to any two axial couples $\bw_1,\W_1$ and $\bw_2,\W_2$ :
\bes
\bw_1\cdot\bw_2=\frac{1}{2}\W_1\cdot\W_2.
\ees
The proof of these two last properties is rather easy and left to the reader.

\section{Orientation of a basis}
It is immediate to observe that a basis $\Ba=\{\e_1,\e_2,\e_3\}$ can be oriented in two opposite ways\footnote{It is evident that this is true also for one- and two-dimensional vector spaces.}: For example, once two unit mutually orthogonal vectors $\e_1$ and $\e_2$ are chosen, there are two opposite unit vectors perpendicular to both $\e_1$ and $\e_2$ that can be chosen to form $\Ba$.  

\begin{figure}[h]
	\begin{center}
         \includegraphics[scale=.25]{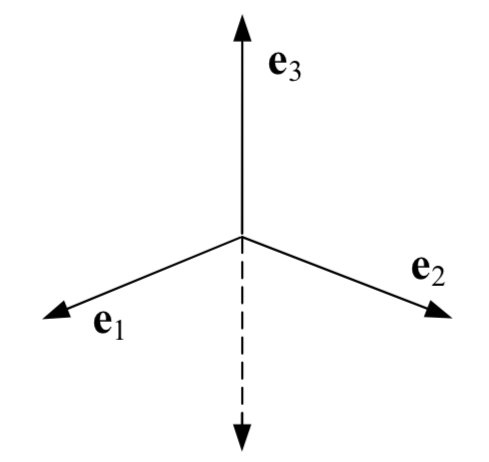}
	\caption{Right- and left-handed bases.}
	\label{fig:4}
	\end{center}
\end{figure}

We say that $\Ba$ is {\it positively oriented} or {\it right-handed} if
\bes
\e_1\times\e_2\cdot\e_3=1,
\ees
while $\Ba$ is {\it negatively oriented} or {\it left-handed} if
\bes
\e_1\times\e_2\cdot\e_3=-1.
\ees
Schematically, a right-handed basis is represented in Fig. \ref{fig:4}, where a left-handed basis is represented too with a dashed $\e_3$.

With a right-handed basis, by definition,  the axial tensors of the three vectors of the basis are
\bes
\besp
&\W_{\e_1}=\e_3\otimes\e_2-\e_2\otimes\e_3,\\
&\W_{\e_2}=\e_1\otimes\e_3-\e_3\otimes\e_1,\\
&\W_{\e_3}=\e_2\otimes\e_1-\e_1\otimes\e_2.\\
\end{split}
\ees

\section{Rotations}
\label{sec:changeofbasis}
In the previous chapter, we have seen that the elements of $\Ve$ represent translations over $\Eu$. A {\it rotation}, i.e. a rigid rotation of the space, is an operation that transforms any two vectors $\bu,\bv\in\Ve$ into two other vectors $\hat{\bu},\hat{\bv}\in\Ve$ in such a way that 
\be
\label{eq:rot1}
u=\hat{u},\ \ v=\hat{v},\ \ \bu\cdot\bv=\hat{\bu}\cdot\hat{\bv},
\ee
i.e., a rotation is a transformation that  preserves norms and angles. Because a rotation is a transformation from $\Ve$ to $\Ve$, rotations are tensors, so we can write
\bes
\hat{\bv}=\Ro\bv,
\ees
with $\Ro$ the {\it rotation tensor} or simply  {\it rotation}.

Conditions (\ref{eq:rot1}) impose some restrictions on $\Ro$:
\bes
\hat{\bu}\cdot\hat{\bv}=\Ro\bu\cdot\Ro\bv=\bu\cdot\Ro^\top\Ro\bv=\bu\cdot\bv\iff\Ro^\top\Ro=\I=\Ro\Ro^\top.
\ees
A tensor that preserves the angles belongs to $Orth(\Ve)$, the subspace of {\it orthogonal tensors}; we leave to the reader the proof that actually $Orth(\Ve)$ is actually a subspace of $Lin(\Ve)$. Replacing in the above equation $\bv$ with $\bu$ shows immediately that an orthogonal tensor also preserves  the norms. By the uniqueness of the inverse, we see that 
\bes
\Ro\in Orth(\Ve)\iff \Ro^{-1}=\Ro^\top.
\ees
The above condition is not sufficient to characterize a rotation; in fact, a rotation must transform a right-handed basis into another right-handed basis, i.e. it must {\it preserve the orientation of the space}. This means that it must be  
\bes
\hat{\e}_1\times\hat{\e}_2\cdot\hat{\e}_3=\Ro\e_1\times\Ro \e_2\cdot\Ro\e_3=\e_1\times\e_2\cdot\e_3.
\ees 
By Eq. (\ref{eq:skewmixed}), we get hence the condition\footnote{From the condition $\Ro^\top\Ro=\I$ and through Eq. (\ref{eq:dettransp}) and the theorem of Binet, we recognize immediately that $\det\Ro=\pm1\ \ \forall\Ro\in Orth(\Ve)$.}
\bes
\det\Ro(\e_1\times\e_2\cdot\e_3)=\e_1\times\e_2\cdot\e_3\iff\det\Ro=1.
\ees
The tensors of $Orth(\Ve)$ that have a determinant equal to 1 form the subspace of  {\it proper rotations} or simply {\it rotations}, indicated by $Orth(\Ve)^+$ or also by $SO(3)$. Only tensors of $Orth(\Ve)^+$ represent rigid rotations of $\Eu$\footnote{A tensor $\Sy\in Orth(\Ve)$ such that $\det\Sy=-1$ represents a transformation that changes the orientation of the space, like mirror symmetries do, see Section  \ref{sec:symmetries}.}.
\begin{teo} Each tensor $\Ro\in Orth(\Ve)$ has the eigenvalue $\pm1$, with $+1$ for rotations.
\begin{proof}
Let $\bu$ be an eigenvector of $\Ro\in Orth(\Ve)$ corresponding to the eigenvalue $\lambda$.
Because  $\Ro$ preserves the norm, we have
\bes
\Ro\bu\cdot\Ro\bu=\lambda^2\bu^2=\bu^2\ \rightarrow\ \lambda^2=1.
\ees
We must now  prove that there exists at least one real eigenvector $\lambda$. To this end, we consider the characteristic equation
\bes
f(\lambda)=\lambda^3+k_1\lambda^2+k_2\lambda+k_3=0, 
\ees
whose coefficients $k_i$ are real-valued, because $\Ro$ has real-valued components. It is immediate to recognize that
\bes
\lim_{\lambda\rightarrow\pm\infty}f(\lambda)=\pm\infty.
\ees
So, because $f(\lambda)$ is a real-valued continuous function, actually a polynomial of $\lambda$, there exists at least one $\lambda_1\in\R$ such that
\bes
f(\lambda_1)=0.
\ees  
In addition, we already know that $\forall\Ro\in Orth(\Ve), \det\Ro=\pm1$ and that, if $\lambda_i, i=1,2,3$ are the eigenvalues of $\Ro$, then $\det\Ro=\lambda_1\lambda_2\lambda_3$. Hence, the following two  are the possible cases:
\begin{enumerate}[i.]
\item $\lambda_1\in\R$ and $\lambda_2,\lambda_3\in\Cq$, with $\lambda_3=\overline{\lambda}_2$, the complex conjugate of $\lambda_2$;
\item $\lambda_i\in\R\ \forall i=1,2,3$.
\end{enumerate}
Let us consider the case of $\Ro\in Orth(\Ve)^+$, i.e. a (proper) rotation $\rightarrow\det\Ro=1$. Then, in the first case above,
\bes
\det\Ro=\lambda_1\lambda_2\overline{\lambda}_2=\lambda_1[\Re^2(\lambda_2)+\Im^2(\lambda_2)].
\ees
But
\bes
\Re^2(\lambda_2)+\Im^2(\lambda_2)=1
\ees 
because it is the square of the modulus of the complex eigenvalue $\lambda_2$. So in this case
\bes
\det\Ro=1\iff\lambda_1=1.
\ees
In the second case, $\lambda_i\in\R\ \forall i=1,2,3$, either $\lambda_1>0, \lambda_2,\lambda_3<0$, or all of them are positive. Because the modulus of each eigenvalue must be equal to 1, either $\lambda_1=1$ or $\lambda_i=1\ \forall i=1,2,3$ (in this case, $\Ro=\I$).

Following the same steps, one can easily  show that $\forall\Sy\in Orth(\Ve)$ with $\det\Sy=-1$, there exists at least one real eigenvalue $\lambda_1=-1$.
\end{proof}
\end{teo}

Generally, a rotation tensor rotates  the basis $\Ba=\{\e_1,\e_2,\e_3\}$ into the basis $\hat{\Ba}=\{\hat{\e}_1,\hat{\e}_2,\hat{\e}_3\}$:
\be
\label{eq:rotcomponents}
\Rb \e_i=\hat{\e}_i\ \ \forall i=1,2,3\ \Rightarrow\ {R_{ij}=\e_i\cdot\Rb\e_j=\e_i\cdot\hat{\e}_j}.
\ee

This result actually means that the $j$th column of $\Rb$ is formed by the components  in the basis $\Ba$ of the vector $\hat{\e}_j$ of $\hat{\Ba}$. Because the two bases are orthonormal, such components are the director cosines of the axes of $\hat{\Ba}$ with respect to $\Ba$.

Geometrically, any rotation is characterized by an {\it axis of rotation}  ${\bw},|\bw|=1$, and by an {\it amplitude} $\phi$, i.e. the angle through which the space is rotated about $\bw$. By definition, $\bw$ is the (only) vector that is left unchanged by $\Ro$, i.e.
\bes
\Ro\bw=\bw,
\ees
or, in other words, it is the eigenvector corresponding to the eigenvalue $+1$.

The question is then: How can a rotation tensor $\Ro$  be expressed by means of its geometrical parameters, $\bw$ and $\phi$? To this end we have a fundamental theorem:
\begin{teo}{\bf(Euler's rotation representation theorem).} $\forall\Ro\in Orth(\Ve)^+$,
\be
\label{eq:rotation1}
{\Rb=\I+\sin\phi\W+(1-\cos\phi)\W^2}
\ee
with $\phi$ the rotation's amplitude and $\W$ the axial tensor of the rotation axis $\bw$.
\begin{proof}
We observe preliminarily that
\be
{\Rb\bw}=\I\bw+\sin\varphi\W\bw+(1-\cos\varphi)\W\W\bw=\I\bw={\bw}
\ee
i.e. that Eq. (\ref{eq:rotation1}) actually defines a transformation that leaves unchanged the axis $\bw$, like a rotation about $\bw$ must do, and that $+1$ is an eigenvalue of $\Ro$.

We need now to prove that Eq. (\ref{eq:rotation1}) actually represents a rotation tensor, i.e. we must prove that 
\bes
\Rb\Rb^\top=\I,\ \ \det\Rb=1.
\ees
Through Eq. (\ref{eq:skew4}) we get
\bes
\besp
{\Rb\Rb^\top}&=(\I+\sin\varphi\W+(1-\cos\varphi)\W^2)(\I+\sin\varphi\W+(1-\cos\varphi)\W^2)^\top\\
&=(\I+\sin\varphi\W+(1-\cos\varphi)\W^2)(\I-\sin\varphi\W+(1-\cos\varphi)\W^2)\\
&=\I+2(1-\cos\varphi)\W^2-\sin^2\varphi\W^2+(1-\cos\varphi)^2\W^4\\
&=\I+2(1-\cos\varphi)\W^2-\sin^2\varphi\W^2-(1-\cos\varphi)^2\W^2{=\I}.
\end{split}
\ees
Then, through Eq. (\ref{eq:skew2}), we obtain
\be
\label{eq:rot2}
\besp
{\Rb}&=\I+\sin\varphi\W+(1-\cos\varphi)\W^2\\
&=\I+\sin\varphi\W-(1-\cos\varphi)(\I-\bw\otimes\bw)\\
&{=\cos\varphi\I+\sin\varphi\W+(1-\cos\varphi)\bw\otimes\bw}.
\end{split}
\ee
To go on, we need to express $\W$ and $\bw\otimes\bw$; if $\bw=(w_1,w_2,w_3)$, then by Eq. (\ref{eq:isomorphskew}) we have
\bes
{\W=\left[\begin{array}{ccc}0 & -w_3 & w_2 \\w_3 & 0 & -w_1 \\-w_2 & w_1 & 0\end{array}\right]}
\ees 
and by Eq. (\ref{eq:compdyad}),
\bes
\bw\otimes\bw=\left[\begin{array}{ccc}w_1^2 &w_1w_2 & w_1w_3 \\  w_1w_2& w_2^2 & w_2w_3 \\ w_1w_3 &w_2w_3  & w_3^2\end{array}\right],
\ees
which on injecting into Eq. (\ref{eq:rot2}) gives
\be
\label{eq:rot3}
\Rb\hspace{-1mm}=\hspace{-1mm}\left[\hspace{-1mm}\begin{array}{ccc}
\cos\varphi\hspace{-1mm}+\hspace{-1mm}(1\hspace{-1mm}-\hspace{-1mm}\cos\varphi)w_1^2 \hspace{-1mm}&-w_3\hspace{-1mm}\sin\varphi\hspace{-1mm}+\hspace{-1mm}w_1w_2(1\hspace{-1mm}-\hspace{-1mm}\cos\varphi) \hspace{-1mm}&w_2\sin\varphi\hspace{-1mm}+\hspace{-1mm} w_1w_3(1\hspace{-1mm}-\hspace{-1mm}\cos\varphi) \\
  w_3\sin\varphi\hspace{-1mm}+\hspace{-1mm}w_1w_2(1\hspace{-1mm}-\hspace{-1mm}\cos\varphi)\hspace{-1mm}& \cos\varphi\hspace{-1mm}+\hspace{-1mm}(1\hspace{-1mm}-\hspace{-1mm}\cos\varphi)w_2^2 \hspace{-1mm}& -\hspace{-.5mm}w_1\sin\varphi\hspace{-1mm}+\hspace{-1mm} w_2w_3(1\hspace{-1mm}-\hspace{-1mm}\cos\varphi) \\
   -\hspace{-.5mm}w_2\sin\varphi\hspace{-1mm}+\hspace{-1mm} w_1w_3(1\hspace{-1mm}-\hspace{-1mm}\cos\varphi) \hspace{-1mm}&w_1\hspace{-.5mm}\sin\varphi\hspace{-1mm}+\hspace{-1mm} w_2w_3(1\hspace{-1mm}-\hspace{-1mm}\cos\varphi)  \hspace{-1mm}& \cos\varphi\hspace{-1mm}+\hspace{-1mm}(1\hspace{-1mm}-\hspace{-1mm}\cos\varphi)w_3^2
   \end{array}\hspace{-1mm}\right]\hspace{-1mm}.
\ee
This formula gives $\Rb$ as a function exclusively of $\bw$ and $\varphi$, the geometrical elements of the rotation. Then
\bes
\det\Rb=(w^2+(1-w^2)\cos\varphi)(\cos^2\varphi+w^2\sin^2\varphi)
\ees
and because $w=1,{\det\Rb=1}$, which proves that Eq. (\ref{eq:rotation1}) actually represents a rotation.

We eventually need to prove that Eq. (\ref{eq:rotation1}) represents the rotation about $\bw$ of amplitude $\varphi$. To this end, we choose an orthonormal basis $\Ba=\{\eu,\ed,\et\}$ of $\Ve$ such that {$\bw=\et$}, i.e. we analyze the particular case of  a rotation of amplitude $\varphi$ about $\et$. This is always possible thanks to the arbitrariness of the basis of $\Ve$. In such a case, Eq. (\ref{eq:rotcomponents}) gives 
\be
\label{eq:rotationZ}
\Rb=\left[\begin{array}{ccc}\cos\varphi & -\sin\varphi &0 \\\sin\varphi & \cos\varphi & 0 \\0 & 0 & 1\end{array}\right].
\ee
Moreover,
\bes
\besp
&{\W=\left[\begin{array}{ccc}0 & -1 & 0 \\1 & 0 &0 \\0 & 0 & 0\end{array}\right],\ \ \
\bw\otimes\bw=\left[\begin{array}{ccc}0 & 0 & 0 \\0 & 0 &0 \\0 & 0 & 1\end{array}\right]},\\
&{\W^2=-(\I-\bw\otimes\bw)=\left[\begin{array}{ccc}-1 & 0 & 0 \\0 & -1 &0 \\0 & 0 & 0\end{array}\right]}.
\end{split}
\ees

Hence
\be
\besp
&{\I+\sin\varphi\W+(1-\cos\varphi)\W^2}=\left[\begin{array}{ccc}1& 0 & 0 \\0 & 1 &0 \\0 & 0 & 1\end{array}\right]+
\sin\varphi\left[\begin{array}{ccc}0 & -1 & 0 \\1 & 0 &0 \\0 & 0 & 0\end{array}\right]+\\
&+(1-\cos\varphi)\left[\begin{array}{ccc}-1 & 0 & 0 \\0 & -1 &0 \\0 & 0 & 0\end{array}\right]=\left[\begin{array}{ccc}\cos\varphi & -\sin\varphi &0 \\\sin\varphi & \cos\varphi & 0 \\0 & 0 & 1\end{array}\right]{=\Rb}.
\end{split}
\ee
\end{proof}
\end{teo}
Equation (\ref{eq:rotation1}) gives another result: To obtain the inverse of $\Ro$ it is sufficient to change the sign of $\phi$. In fact, because $\W\in Skw(\Ve)$ and through Eq. (\ref{eq:w2sym})
\bes
\besp
\Ro^{-1}=\Ro^\top&=(\I+\sin\phi\W+(1-\cos\phi)\W^2)^\top=\I+\sin\phi\W^\top+(1-\cos\phi)(\W^2)^\top\\
&=\I-\sin\phi\W+(1-\cos\phi)\W^2=\I+\sin(-\phi)\W+(1-\cos(-\phi))\W^2.
\end{split}
\ees

The knowledge of the inverse of a rotation also allows us  to perform the operation of {\it change of basis}, i.e. to determine the components of a vector or of a tensor in a basis   $\hat{\Ba}=\{\hat{\e}_1,\hat{\e}_2,\hat{\e}_3\}$ rotated with respect to an original basis $\Ba=\{\eu,\ed,\et\}$ by a rotation $\Ro$ (in the following equations, the symbol $\ \hat{ }\ $ indicates a quantity specified in the basis $\hat{\Ba}$). Considering that
\bes
\e_i=\Ro^{-1}\hat{\e}_i=\Ro^\top\hat{\e}_i=R_{hk}^\top(\hat{\e}_h\otimes\hat{\e}_k)\hat{\e}_i=R^\top_{hk}\delta_{ki}\hat{\e}_h
\ees
we get, for a vector $\bu$,
\bes
\bu=u_i\e_i=R^\top_{ki}u_i\hat{\e}_k
\ees
i.e.
\bes
\hat{u}_k=R^\top_{ki}u_i\ \rightarrow\ \hat{\bu}=\Ro^\top\bu.
\ees
We remark that, because $\Ro^\top=\Ro^{-1}$, the operation of change of basis is just the opposite  of the rotation of the space (and actually, we have seen that it is sufficient to take the opposite of $\phi$ in Eq. (\ref{eq:rotation1}) to get $\Ro^{-1}$).

For a second-rank tensor $\L$ we get
\bes
\L=L_{ij}\e_i\otimes\e_j=L_{ij}R^\top_{mi}\hat{\e}_m\otimes R^\top_{nj}\hat{\e}_n=R^\top_{mi}R^\top_{nj}L_{ij}\hat{\e}_m\otimes\hat{\e}_n,
\ees
i.e.
\bes
\hat{L}_{mn}=R^\top_{mi}R^\top_{nj}L_{ij}\ =R^\top_{mi}L_{ij}R_{jn}\ \rightarrow\ \hat{\L}=\Ro^\top\L\Ro.
\ees

We remark something that is typical of tensors: the components of a $r$-rank tensor in a rotated basis $\hat{\Ba}$ depend upon the $r$th powers of the directors cosines of the axes of $\hat{\Ba}$, i.e. on the $r$th powers of the components $R_{ij}$ of $\Ro$.

\begin{figure}[h]
	\begin{center}
         \includegraphics[scale=.6]{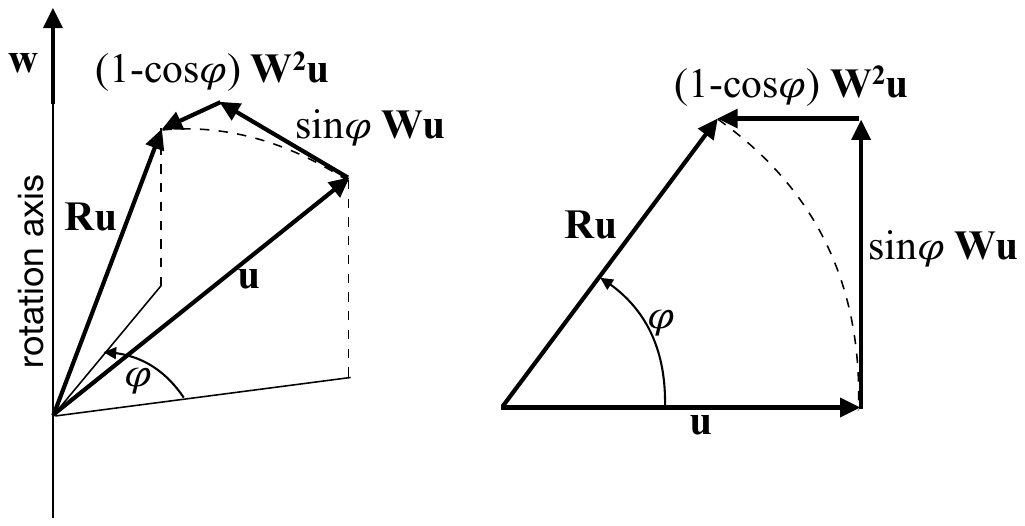}
	\caption{Rotation of a vector.}
	\label{fig:5}
	\end{center}
\end{figure}

If a rotation tensor is known through its Cartesian components in a given basis $\Ba$, it is easy to calculate its geometrical elements: The rotation axis $\bw$ is the eigenvector of $\Ro$ corresponding to the eigenvalue 1, so it is found solving the equation
\bes
\Ro\bw=\bw
\ees
and then normalizing it, while the rotation amplitude $\phi$ can be found  using (\ref{eq:rotation1}) along with (\ref{eq:skew2}): Because the trace of a tensor is an invariant, we get
\bes
\tr\Ro=3+(1-\cos\phi)\tr(-\I+\bw\otimes\bw)=1+2\cos\phi\ \rightarrow\ \phi=\arccos\frac{\tr\Ro-1}{2}.
\ees
It is interesting to consider the   geometrical meaning of Eq. (\ref{eq:rotation1}).  For this purpose, we apply Eq. (\ref{eq:rotation1}) to a vector $\bu$, see Fig. \ref{fig:5},
\bes
\besp
\Rb\bu&=(\I+\sin\varphi\W+(1-\cos\varphi)\W^2)\bu\\
&=\bu+\sin\varphi\bw\times\bu+(1-\cos\varphi)\bw\times(\bw\times\bu)
\end{split}
\ees
The rotated vector $\Rb\bu$ is the sum of three vectors; in particular, $\sin\varphi\W\bu$ is always orthogonal to $\bu$, $\bw$ and  $(1-\cos\varphi)\W^2\bu$.  If $\bu\cdot\bw=0$, then  $(1-\cos\varphi)\W^2\bu$  is also parallel to $\bu$, see the sketch on the right in Fig. \ref{fig:5}.

\begin{figure}[h]
	\begin{center}
         \includegraphics[scale=.45]{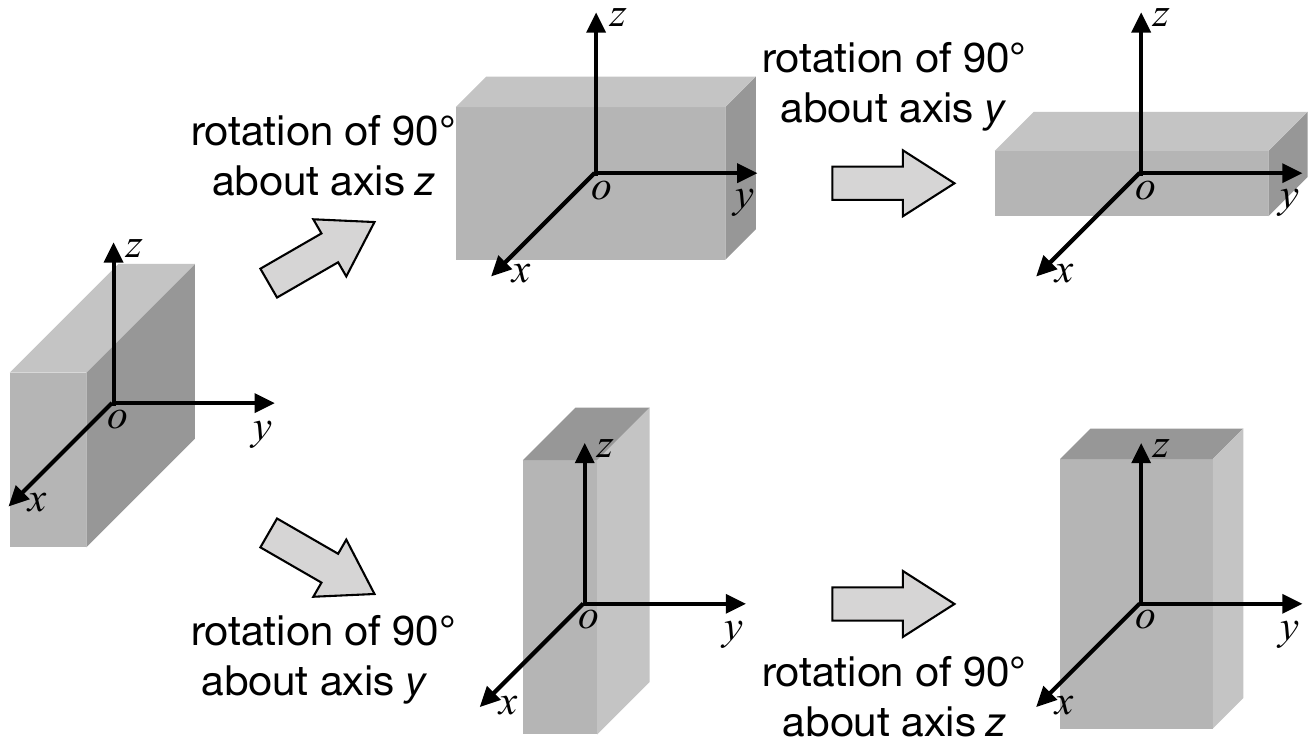}
	\caption{Non-commutativity of the rotations.}
	\label{fig:6}
	\end{center}
\end{figure}

Let us consider now a composition of rotations. 
In particular, let us imagine that a vector $\bu$ is rotated first by $\Rb_1$ around $\bw_1$ through $\varphi_1$, then by $\Rb_2$ around $\bw_2$ through $\varphi_2$.  So, first, the vector $\bu$ becomes  the vector
\bes
\bu_1=\Rb_1\bu.
\ees
Then, the vector $\bu_1$ is  rotated about $\bw_2$ through $\varphi_2$ to become
\bes
\bu_{12}=\Rb_2\bu_1=\Rb_2\Rb_1\bu.
\ees
Let us now suppose  that we change the order of the rotations: $\Rb_2$ first and then $\Rb_1$.
The final result will be the vector
\be
\bu_{21}=\Rb_1\Rb_2\bu.
\ee
Because the tensor product is not symmetric (i.e., it has not the commutativity property), generally\footnote{We have seen in Theorem \ref{teo: commutativity}, that two tensors commute $\iff$ they are coaxial, i.e. if they have the same eigenvectors. Because the rotation axis is always a real eigenvector of a rotation tensor, if two tensors operate a rotation about different axes, they are not coaxial. Hence, the rotation tensors about different axes never commute.}
\bes
{\bu_{12}\neq\bu_{21}}.
\ees
In other words, the order of the rotations matters: Changing the order of the rotations leads to a different final result. An example  is shown in Fig. \ref{fig:6}.

This is a fundamental difference between rotations and displacements, that commute, see Fig. \ref{fig:2}, because the composition of displacements is ruled by the sum of vectors:
\be
\bw=\bu+\bv=\bv+\bu
\ee
This difference, which is a major point in physics, comes from the difference in the operators: vectors for the displacements and tensors for the rotations.

Any rotation can be specified by the knowledge of three parameters. This can be easily seen  from Eq. (\ref{eq:rotation1}): the parameters are the three components of $\bw$, that are not independent because
\bes
w=|\bw|=\sqrt{w_1^2+w_2^2+w_3^2}=1
\ees
and by the amplitude angle $\phi$. The choice of the parameters by which to express a rotation is not unique. Besides the use of the Cartesian components of $\bw$ and $\phi$, cf. Eq. (\ref{eq:rot3}), other choices are possible, let us see three of them:
\begin{enumerate}[i.]
\item {\it Physical angles}: The rotation axis $\bw$ is given through its spherical coordinates $\psi$, the {\it longitude}, $0\leq\psi<2\pi$, and $\theta$, the {\it colatitude}, $0\leq\theta\leq\pi$, see Fig. \ref{fig:7}, the third parameter being the {\it rotation  amplitude} $\phi$. 
\begin{figure}[h]
	\begin{center}
         \includegraphics[scale=.45]{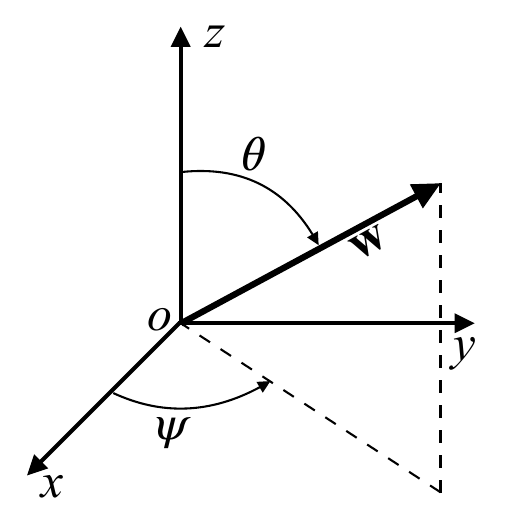}
	\caption{Physical angles.}
	\label{fig:7}
	\end{center}
\end{figure}
Then
\bes
\bw=(\sin\theta\cos\psi,\sin\theta\sin\psi,\cos\theta)\ \rightarrow\ \theta=\arccos w_3,\ \psi=\arctan\frac{w_2}{w_1},
\ees
and, Eq. (\ref{eq:rot3}),
\bes
\hspace{-10mm}\Rb=
\left[\begin{array}{ccc}
c\psi^2s\theta^2+c\phi(c\theta^2+s\psi^2s\theta^2) &s\psi c\psi s\theta^2(1-c\varphi)-c\theta s\phi &c\psi s\theta c\theta(1-c\varphi) +s\psi s\theta s\phi\\
s\psi c\psi s\theta^2(1-c\varphi)+c\theta s\phi & s\psi^2s\theta^2+c\phi(c\theta^2+c\psi^2s\theta^2) & s\psi s\theta c\theta(1-c\varphi) -c\psi s\theta s\phi\\
   c\psi s\theta c\theta(1-c\varphi) -s\psi s\theta s\phi &s\psi s\theta c\theta(1-c\varphi) +c\psi s\theta s\phi & c\theta^2+c\phi(c\psi^2s\theta^2+s\psi^2s\theta^2)
   \end{array}\right],
\ees
where $c\psi=\cos\psi,s\psi=\sin\psi,c\theta=\cos\theta,s\theta=\sin\theta,c\phi=\cos\phi$ and $s\phi=\sin\phi$.
We remark that all the components of $\Ro$ so expressed depend upon the first powers of the circular functions of $\phi$. Hence, for what is said above, with this representation of the rotations, the components of a rotated $r$-rank tensor depend upon the $r$th power of the circular functions of $\phi$, i.e. of the physical rotation, but not of $\psi$ nor of $\theta$. 

\item {\it Euler's angles}: In this case, the three parameters are the amplitude of three particular rotations into which the rotation is decomposed. Such parameters are the angles $\psi$, the {\it precession}, $\theta$, the {\it nutation}, and $\phi$, the {\it proper rotation}, see Fig. \ref{fig:8}
\begin{figure}[b]
	\begin{center}
         \includegraphics[scale=.35]{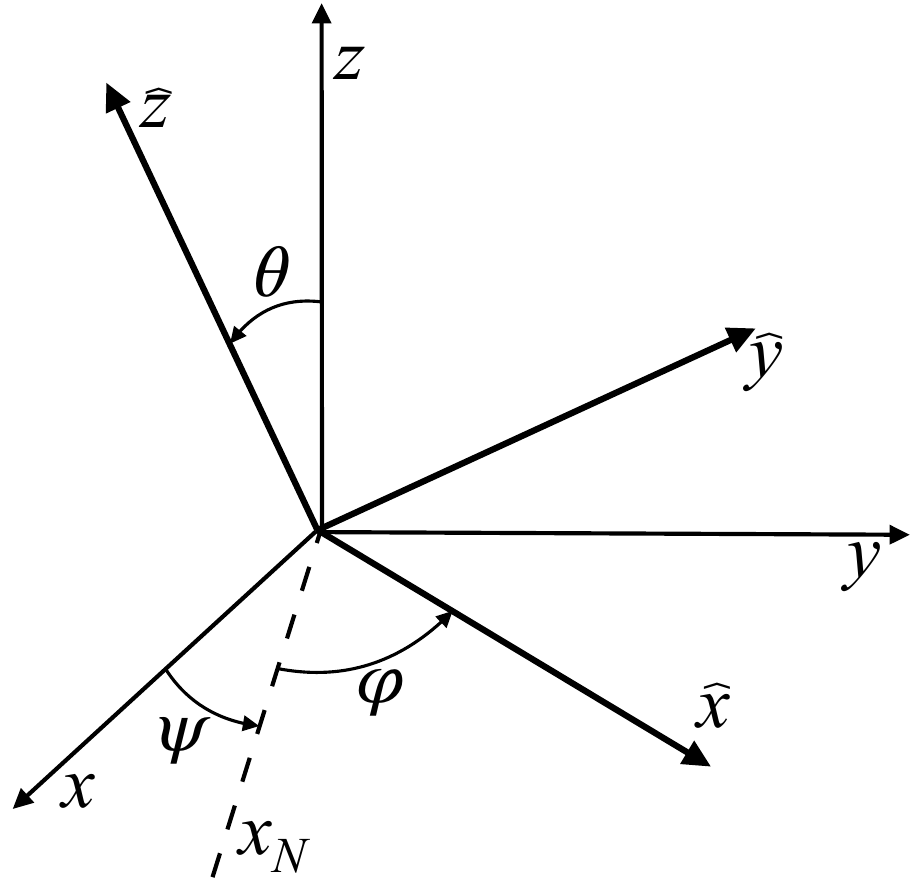}
	\caption{Euler's angles.}
	\label{fig:8}
	\end{center}
\end{figure}
These three  rotations are represented in Fig. \ref{fig:9}. The first one, of amplitude $\psi$, is made about $z$ to carry the axis $x$ onto the {\it knots line} $x_N$, the line perpendicular to both the axes $z$ and $\hat{z}$,  and $y$ onto $\overline{y}$; by Eq. (\ref{eq:rotcomponents}), in the frame $\{x,y,z\}$ it is
\bes
\Ro_\psi=
\left[\begin{array}{ccc}\cos\psi& -\sin\psi &0 \\\sin\psi & \cos\psi & 0 \\0 & 0 & 1\end{array}\right].
\ees
The second one, of amplitude $\theta$, is made about $x_N$ to carry $z$ onto $\hat{z}$; in the frame $\{x_N,\overline{y},z\}$, it is
\bes
\Ro_\theta=
\left[\begin{array}{ccc}1&0&0\\0&\cos\theta&-\sin\theta\\0& \sin\theta & \cos\theta \end{array}\right],
\ees
while in the frame $\{x,y,z\}$, 
\bes
\Ro_\theta^o=(\Ro_\psi^{-1})^\top\Ro_\theta\Ro_\psi^{-1}=\Ro_\psi\Ro_\theta\Ro_\psi^{\top}.
\ees
\begin{figure}[h]
	\begin{center}
         \includegraphics[scale=.5]{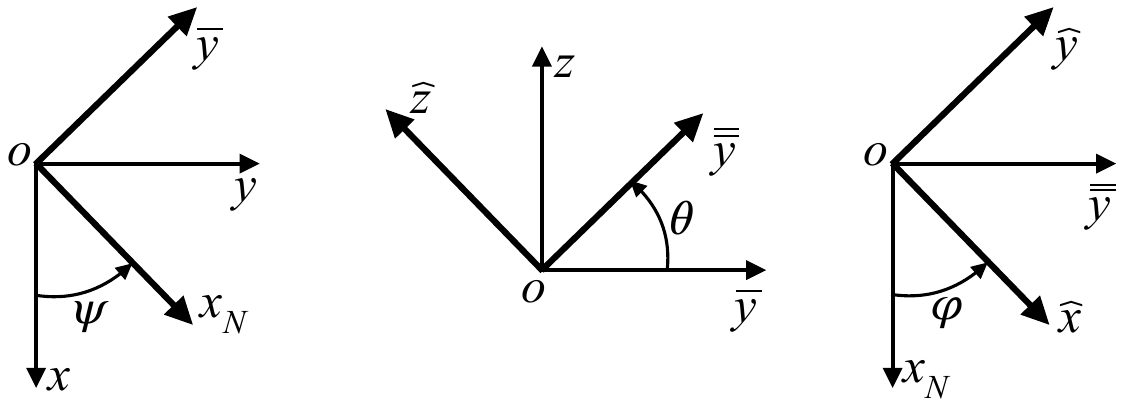}
	\caption{Euler's rotations, as seen from the respective axes of rotation.}
	\label{fig:9}
	\end{center}
\end{figure}

The last rotation, of amplitude $\phi$, is made about $\hat{z}$ to carry $x_N$ onto $\hat{x}$ and $\overline{\overline{y}}$ onto $\hat{y}$; in the frame $\{x_N,\overline{\overline{y}},\hat{z}\}$, it is
\bes
\Ro_\phi=\left[\begin{array}{ccc}\cos\varphi & -\sin\varphi &0 \\\sin\varphi & \cos\varphi & 0 \\0 & 0 & 1\end{array}\right],
\ees
while in $\{x,y,z\}$, 
\bes
\Ro_\phi^o=(\Ro_\psi^{-1})^\top(\Ro_\theta^{-1})^\top\Ro_\phi\Ro_\theta^{-1}\Ro_\psi^{-1}=\Ro_\psi\Ro_\theta\Ro_\phi\Ro_\theta^\top\Ro_\psi^{\top}.
\ees
Any vector $\bu$ is transformed, by the global rotation, into the vector
\bes
\hat{\bu}=\Ro\bu.
\ees
But we can also write 
\bes
\hat{\bu}=\Ro_\phi^o\overline{\overline{\bu}},
\ees
where $\overline{\overline{\bu}}$ is the vector transformed by the  rotation $\Ro_\theta^o$,
\bes
\overline{\overline{\bu}}=\Ro_\theta^o\overline{\bu},
\ees
and $\overline{\bu}$ is the vector transformed by the rotation $\Ro_\psi$:
\bes
\overline{\bu}=\Ro_\psi\bu.
\ees
Finally,
\bes
\hat{\bu}=\Ro\bu=\Ro_\phi^o\Ro_\theta^o\Ro_\psi\bu\ \rightarrow\ \Ro=\Ro_\phi^o\Ro_\theta^o\Ro_\psi,
\ees
i.e. the global rotation tensor is obtained composing, in the opposite order of execution of the rotations, the three tensors all expressed in the original basis. However,
\bes
\Ro=\Ro_\phi^o\Ro_\theta^o\Ro_\psi=\Ro_\psi\Ro_\theta\Ro_\phi\Ro_\theta^\top\Ro_\psi^\top\Ro_\psi\Ro_\theta\Ro_\psi^\top\Ro_\psi=\Ro_\psi\Ro_\theta \Ro_\phi,
\ees
i.e., the global rotation tensor is also equal to the composition of the three rotations, in the order of execution, if the three rotations are expressed in their own particular bases. This result is general, not bounded to the Euler's rotations nor to three rotations.

Performing the tensor multiplications, we  get
\bes
\Ro=\left[\begin{array}{ccc}
\cos\psi\cos\varphi -\sin\psi\sin\phi\cos\theta& -\cos\psi\sin\varphi-\sin\psi\cos\phi\cos\theta &\sin\psi\sin\theta \\
\sin\psi\cos\varphi +\cos\psi\sin\phi\cos\theta& -\sin\psi\sin\varphi+\cos\psi\cos\phi\cos\theta & -\cos\psi\sin\theta \\
\sin\phi\sin\theta & \cos\phi\sin\theta  & \cos\theta
\end{array}\right].
\ees
The components of a vector $\bu$ in the basis $\hat{\Ba}$ are then given by
\bes
\hat{\bu}=\Ro^\top\bu=\Ro_\phi^\top\Ro_\theta^\top\Ro_\psi^\top\bu,
\ees
and those of a second-rank tensor
\bes
\hat{\L}=\Ro^\top\L\Ro=\Ro_\phi^\top\Ro_\theta^\top\Ro_\psi^\top\L\Ro_\psi\Ro_\theta\Ro_\phi.
\ees

\item {\it Coordinate angles}: In this case, the rotation $\Ro$ is decomposed into three successive rotations $\alpha,\beta,\gamma$, respectively about the axes $x$, $y$ and $z$ of each rotation, i.e.
\bes
\Ro=\Ro_\alpha\Ro_\beta\Ro_\gamma
\ees
with
\bes
\Ro_\alpha=\left[\begin{array}{ccc}1&0&0\\0&\cos\alpha & -\sin\alpha  \\0&\sin\alpha& \cos\alpha\end{array}\right],
\Ro_\beta=\left[\begin{array}{ccc}\cos\beta&0 & -\sin\beta\\0&1&0  \\\sin\beta&0& \cos\beta\end{array}\right],
\Ro_\gamma=\left[\begin{array}{ccc}\cos\gamma & -\sin\gamma &0 \\\sin\gamma & \cos\gamma & 0 \\0 & 0 & 1\end{array}\right],
\ees
so finally
\bes
\Ro=\left[\begin{array}{ccc}
\cos\beta\cos\gamma&-\cos\beta\sin\gamma&-\sin\beta\\
\cos\alpha\sin\gamma-\sin\alpha\sin\beta\cos\gamma&\cos\alpha\cos\gamma+\sin\alpha\sin\beta\sin\gamma&-\sin\alpha\cos\beta\\
\sin\alpha\sin\gamma+\cos\alpha\sin\beta\cos\gamma&\sin\alpha\cos\gamma-\cos\alpha\sin\beta\sin\gamma&\cos\alpha\cos\beta
\end{array}\right].
\ees

\end{enumerate}

\begin{figure}[h]
	\begin{center}
         \includegraphics[scale=.6]{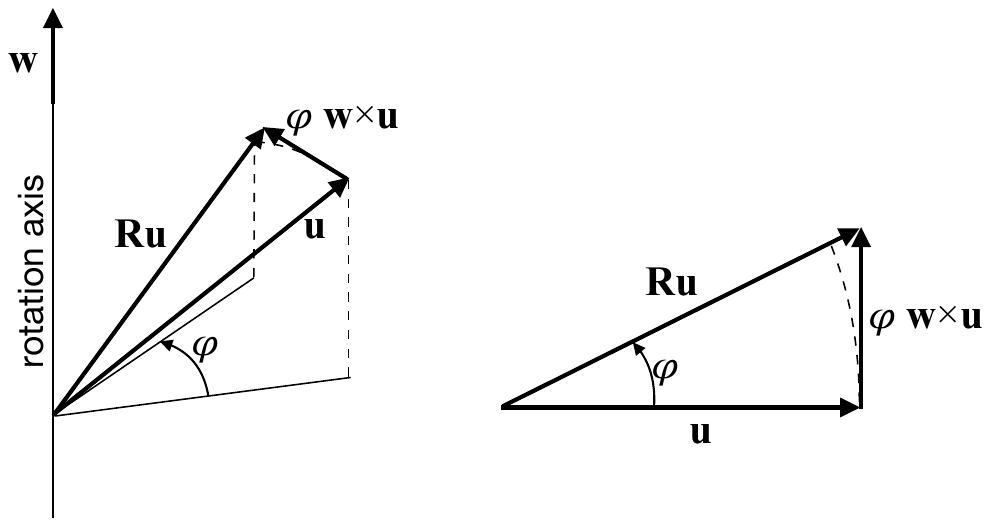}
	\caption{Small rotations.}
	\label{fig:10}
	\end{center}
\end{figure}

Let us now consider the case of {\it small rotations}, i.e. $|\phi|\rightarrow0$.  In such a case, 
\bes
\sin\phi\simeq\phi,\ \ 1-\cos\phi\simeq0
\ees
so that the Euler's theorem, Eq. (\ref{eq:rotation1}), becomes
\bes
\Rb\simeq\I+\phi\W,
\ees
i.e. in the small rotations approximation, any vector $\bu$ is transformed into 
\be
\label{eq:smallrot}
\Rb\bu\simeq(\I+\phi\W)\bu=\bu+\phi\bw\times\bu,
\ee
i.e. by a skew tensor and not by a rotation tensor. The term $(1-\cos\phi)\W^2\bu$ has disappeared, as it is a higher-order infinitesimal quantity, and the term $\phi\bw\times\bu$ is orthogonal to $\bu$. Because $\phi\rightarrow0$, the arc is approximated by its tangent, the vector $\phi\bw\times\bu$, see Fig. \ref{fig:10}.
Applying to Eq. (\ref{eq:smallrot}) the procedure already seen for the composition of finite amplitude rotations, we get
\bes
\besp
&\bu_1\ =\Rb_1\bu=(\I+\phi_1\W_1)\bu=\bu+\phi_1\bw_1\times\bu,\\
&\bu_{21}=\Rb_2\bu_1=(\I+\phi_2\W_2)\bu_1=\bu_1+\phi_2\bw_2\times\bu_1\\
&\hspace{5.5mm}=\bu+\phi_1\bw_1\times\bu+\phi_2\bw_2\times\bu\\
&\hspace{10mm}+\phi_1\phi_2\bw_2\times(\bw_1\times\bu).
\end{split}
\ees
If the order of the rotations is changed, the last term becomes $\phi_1\phi_2\bw_1\times(\bw_2\times\bu)$,  which is, in general, different from $\phi_1\phi_2\bw_2\times(\bw_1\times\bu)$:
  {Strictly speaking, also small rotations do not commute}\footnote{This can happen for some vectors, all the times that $\bw_1\cdot\bu=\bw_2\cdot\bu$, like for the case of a vector $\bu$ orthogonal to both $\bw_1$ and $\bw_2$; however, this is no more than a curiosity, it has no importance in practice.}.
 However, for small rotations, $\phi_1\phi_2$  is negligible with respect to $\phi_1$ and $\phi_2$: In this approximation, small rotations commute.
We remark that the approximation (\ref{eq:smallrot})  gives, for the displacements, a law that is quite similar to that of the velocities of the points of a rigid body:
\bes
\bv=\bv_0+\bom\times(p-o)
\ees
This is quite natural, because
\bes
\omega=\frac{d\phi}{dt},
\ees
i.e.  a small amplitude rotation can be seen as the rotation made with finite angular velocity $\bom$ in a small time interval $dt$.

\section{Reflexions}
\label{sec:symmetries}
Let us consider now  tensors  $\Sy\in Orth(\Ve)$ that are not a rotation, i.e. such that $\det\Sy=-1$. Let us call $\Sy$ an {\it improper rotation}. A particular improper rotation, whose all eigenvalues are equal to -1, is the {\it inversion} or {\it reflexion tensor}:
\bes
\Sy_I=-\I.
\ees
The effect of $\Sy_I$ is to transform any basis $\Ba$ into the basis $-\Ba$, i.e. with all the basis vectors changed in orientation (or, equivalently, to change the sign of all the components of a vector). In other words, $\Sy_I$ changes the orientation of the space. This is also the effect of any other improper rotation $\Sy$, that can be decomposed into a proper rotation $\Ro$ followed by the reflexion $\Sy_I$\footnote{The application of  Binet's theorem shows immediately that $\det\Sy=-1$, while $\Sy_I\Ro(\Sy_I\Ro)^\top=\Sy_I\Ro\Ro^\top\Sy_I^\top=-\I(-\I)^\top=\I$: The decomposition in Eq. (\ref{eq:decimprot}) actually gives an improper rotation.}:
\be
\label{eq:decimprot}
\Sy=\Sy_I\Ro.
\ee
Let $\n\in\S$; then
\be
\label{eq:sym2r}
\Sy_R=\I-2\n\otimes\n
\ee
is the tensor that operates the transformation of  symmetry with respect to a plane orthogonal to $\n$. In fact
\bes
\Sy_R\n=-\n, \ \ \Sy_R\mathbf{m}=\mathbf{m}\ \ \forall\mathbf{m}\in\Ve|\ \mathbf{m}\cdot\n=0.
\ees
$\Sy_R$ is an improper rotation; in fact, by Eq. (\ref{eq:propdyad0}), 
\bes
\besp
(\I-2\n\otimes\n)(\I-2\n\otimes\n)^\top&=(\I-2\n\otimes\n)(\I-2\n\otimes\n)\\
&=\I-2\n\otimes\n-2\n\otimes\n+4(\n\otimes\n)(\n\otimes\n)=\I,
\end{split}
\ees
while by the same definition of trace and through Eqs.   (\ref{eq:trilin4}) and (\ref{eq:detsum}),
\bes
\det(\I-2\n\otimes\n)=1-2\tr(\n\otimes\n)+4\frac{\tr^2(\n\otimes\n)-\tr(\n\otimes\n)(\n\otimes\n)}{2}-8\det(\n\otimes\n)=-1.
\ees
Let $\Sy=\Sy_I\Ro$ be an improper rotation. Then
\bes
\besp
(\Sy\bu)\times(\Sy\bv)&=(\Sy_I\Ro\bu)\times(\Sy_I\Ro\bv)=\det(\Sy_I\Ro)\left[(\Sy_I\Ro)^{-1}\right]^\top(\bu\times\bv)\\
&=\det\Sy_I\det\Ro(\Ro^{-1}\Sy_I^{-1})^\top(\bu\times\bv)=-(-\Ro^{-1}\I)^\top(\bu\times\bv)=\Ro(\bu\times\bv).
\end{split}
\ees
The transformation by $\Sy$ of any vector $\bu$ gives
\bes
\Sy\bu=\Sy_I\Ro\bu=-\Ro\bu,
\ees
i.e. it changes the orientation of the rotated vector; this is not the case when the same improper rotations transforms the vectors of a cross product: The rotated vector, result of the cross product, does not change of orientation, i.e. the cross product is insensitive to a reflexion. That is why, strictly speaking, the result of a cross product is not a vector, but a {\it pseudo-vector}: It behaves like vectors apart for the reflexions. For the same reason, a scalar result of a mixed product (scalar plus cross product of three vectors) is called a {\it pseudo-scalar} because in this case, the scalar result of the mixed product changes of sign under a reflexion, which can be checked easily.

\section{Polar decomposition}

\begin{teo}{\bf(Square root theorem).} Consider $\L\in Sym(\Ve)$ and positive definite. Then, there exists a unique tensor $\U\in Sym(\Ve)$ and positive definite such that
\bes
\L=\U^2.
\ees
\begin{proof}
Existence: Consider $\L,\U,\V\in Sym(\Ve)$ and positive definite, and 
\bes
\L=\omega_i\e_i\otimes\e_i
\ees
a spectral decomposition of $\L$,  $\omega_i>0\ \forall i$. Define $\U$ as
\bes
\U=\sqrt{\omega_i}\e_i\otimes\e_i;
\ees
then, by Eq. (\ref{eq:propdyad0})$_1$, we get
\bes
\U^2=\L.
\ees
Uniqueness: Suppose that also
\bes
\V^2=\L
\ees
and let $\e$ be an eigenvector of $\L$ corresponding to the (positive) eigenvalue $\omega$. Then, if $\lambda=\sqrt{\omega}$,
\bes
\bo=(\U^2-\omega\I)\e=(\U+\lambda\I)(\U-\lambda\I)\e,
\ees
and once we set 
\bes
\bv=(\U-\lambda\I)\e,
\ees
we get
\bes
\U\bv=-\lambda\bv\ \Rightarrow\ \bv=\bo\ \Rightarrow\ \U\e=\lambda\e
\ees
because $\U$ is positive definite and $-\lambda$ cannot be an eigenvalue of $\U$ because $\lambda>0$. In a similar  way,
\bes
\V\e=\lambda\e\ \Rightarrow\ \U\e=\V\e\ 
\ees
for every eigenvector $\e$ of $\L$. Because, based on the spectral theorem, it exists a basis of eigenvectors of $\L$, $\U=\V$.
\end{proof}
\end{teo}
We symbolically write that
\bes
\U=\sqrt{\L}.
\ees

For any $\F\in Lin(\Ve)$, both $\F\F^\top$ and $\F^\top\F$ clearly $\in Sym(\Ve)$. If in addition  $\det\F>0$, then 
\bes
\bu\cdot\F^\top\F\bu=(\F\bu)\cdot(\F\bu)\geq0
\ees
with the zero value  obtained $\iff\F\bu=\bo$ or, what is equivalent, because $\det\F>0\Rightarrow\F$ is invertible, $\iff\bu=\bo$. As a consequence, $\F^\top\F$ is positive definite. In a similar way, it can be proved that $\F\F^\top$ is also positive definite.

A particular tensor decomposition\footnote{This decomposition is fundamental to the theory of deformation of continuum bodies.}  is given by the

\begin{teo}{\bf (Polar decomposition theorem).}  $\forall\F\in Lin(\Ve)|\det\F>0$ exist, and are uniquely determined, two positive definite tensors $\U,\V\in Sym(\Ve)$   and a rotation $\Rb$ such that
\bes
{\F=\Rb\U=\V\Rb}.
\ees
\begin{proof} Uniqueness: 
Let $\F=\Ro\U$ be a {\it right polar decomposition} of $\F$; because $\Ro\in Orth(\Ve)^+$ and $\U\in Sym(\Ve)$, 
\bes
\F^\top\F=\U\Ro^\top\Ro\U=\U^2\Rightarrow\U=\sqrt{\F^\top\F}.
\ees
By the square-root theorem, tensor $\U$ is unique, and because
\bes
\Ro=\F\U^{-1},
\ees
$\Ro$ is unique too.

Now, let $\F=\V\Ro$ be a {\it left polar decomposition} of $\F$; by the same procedure, we get
\bes
\F\F^\top=\V^2\rightarrow\V=\sqrt{\F\F^\top},
\ees
so $\V$ is unique, and also,
\bes
\Ro=\V^{-1}\F.
\ees
Existence: let 
\bes
\U=\sqrt{\F^\top\F}
\ees
so $\U\in Sym(\Ve)$ and it is positive definite, and let
\bes
\Ro=\F\U^{-1}.
\ees
To prove that $\F=\Ro\U$ is a right polar decomposition, we just have to show that $\Ro\in Orth(\Ve)^+$. Since $\det\F>0$ and $\det\U>0$ (the latter because all the eigenvalues of $\U$ are strictly positive), by the theorem of Binet, also $\det\Ro>0$. Then,
\bes
\besp
\Ro^\top\Ro&=(\F\U^{-1})^\top(\F\U^{-1})=\U^{-1}\F^\top\F\U^{-1}\\
&=\U^{-1}\U^2\U^{-1}=\I\Rightarrow\Ro\in Orth(\Ve)^+.
\end{split}
\ees
Now, let
\bes
\V=\Ro\U\Ro^\top,
\ees
then $\V\in Sym(\Ve)$ and is positive definite, see Exercise 22, and 
\bes
\V\Ro=\Ro\U\Ro^\top\Ro=\Ro\U=\F,
\ees
which completes the proof.
\end{proof}
\end{teo}


\section{Exercises}
\begin{enumerate}
\item Prove that 
\bes
\L\bo=\bo\ \ \forall\L\in Lin(\Ve).
\ees
\item \label{ex:2chap2} Prove that, if a straight line $r$ has the direction of $\bu\in\S$, then the tensor giving the projection of a vector $\bv\in\Ve$ on $r$ is $\bu\otimes\bu$ (the {\it orthogonal projector}), while the one giving the projection on a direction orthogonal to $r$ is $\I-\bu\otimes\bu$ (the {\it complementary projector}), see Fig. \ref{fig:54}.
\begin{figure}[b]
\begin{center}
\includegraphics[scale=.5]{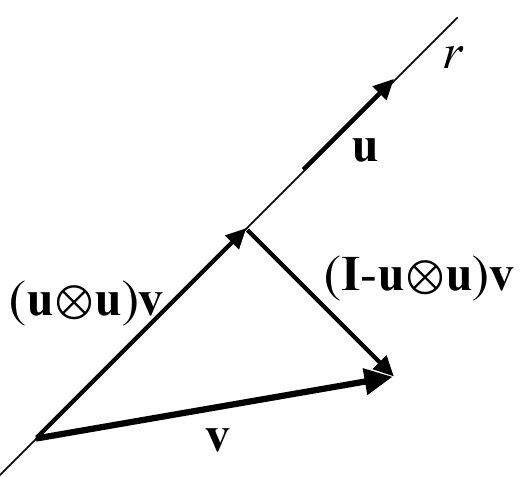}
\caption{Projected vectors.}
\label{fig:54}
\end{center}
\end{figure}
\item  \label{ex:eserc3chap2} For any $\alpha\in\R,\a,\b\in\Ve$ and $\A,\B\in Lin(\Ve)$, prove that
\bes 
(\alpha\A)^\top=\alpha\A^\top,\ \ (\A+\B)^\top=\A^\top+\B^\top,\ \ 
(\a\otimes\b)\A=\a\otimes(\A^\top\b).
\ees
\item Prove that 
\bes
\L+\O=\L\ \ \ \forall\L\in Lin(\Ve).
\ees
\item Prove that
\bes
\tr\I=3,\ \ \tr\O=0.
\ees
\item Prove that, $\forall\A,\B\in Lin(\Ve)$,
\bes
 \tr(\A\B)=\tr(\B\A).
\ees
\item \label{ex:2_5} Prove that, $\forall\L,\M,\N\in Lin(\Ve)$,
\bes
\L^\top\cdot\M^\top=\L\cdot\M,\ \ \L\M\cdot\N=\L\cdot\N\M^\top=\M\cdot\L^\top\N. 
\ees
\item Prove the assertions in Eq. (\ref{eq:propdyad0}).
\item \label{ex:2_12} Prove that any  form defined by a tensor $\L$ can be written as a scalar product of tensors:
\bes
\bv\cdot\L\bw=\L\cdot\bv\otimes\bw\ \ \forall\bv,\bw\in\Ve,\L\in Lin(\Ve).
\ees
\item \label{ex:2_6} Prove that $Sym(\Ve)$ and $Skw(\Ve)$ are orthogonal, i.e. prove that
\bes
\A\cdot\B=0\ \ \forall\A\in Sym(\Ve),\ \B\in Skw(\Ve).
\ees
\item For any $\L\in Lin(\Ve)$, prove that, if $\A\in Sym(\Ve)$, then
\bes
\A\cdot\L=\A\cdot\L^s,
\ees
while if $\B\in Skw(\Ve)$, then
\bes
\B\cdot\L=\B\cdot\L^a.
\ees
\item \label{ex:17} Let $\A,\B,\C,\D\in Lin(\Ve)$; prove that
\bes
\A\cdot(\B\C\D)=(\B^\top\A)\cdot(\C\D)=(\A\D^\top)\cdot(\B\C).
\ees
\item Prove that $\L\cdot\W=0\ \forall\W\in Skw(\Ve) \iff\L\in Sym(\Ve)$.
\item Express by components the second principal invariant $I_2$ of a tensor $\L$.
\item Prove that, if $\a=(a_1,a_2,a_3),\b=(b_1,b_2,b_3),\c=(c_1,c_2,c_3)$, then
\bes
\a\times\b\cdot\c=\det\left[\begin{array}{ccc}  a_1&   a_2&  a_3 \\ b_1 & b_2  & b_3  \\ c_1 & c_2  & c_3 \end{array}\right].
\ees
\item Prove the uniqueness of the inverse tensor.
\item \label{ex:detdyad} Show, using the Cartesian components, that all the dyads are singular.
\item Prove that if $\L$ is invertible and $\alpha\in\R-\{0\}$ then 
\bes
(\alpha\L)^{-1}=\alpha^{-1}\L^{-1}.
\ees
\item Prove that if $\W$ is the axial tensor of $\bw$, then
\bes
\W\W=-\frac{1}{2}|\W|^2(\I-\bw\otimes\bw).
\ees
\item Prove that for any two axial couples $\bw_1,\W_1$ and $\bw_2,\W_2$, we have:
\bes
\bw_1\cdot\bw_2=\frac{1}{2}\W_1\cdot\W_2.
\ees
\item Prove that $\forall\bu,\bv\in\Ve,\ \bu\times\bv=\bo\iff\bu\otimes\bv\in Sym(\Ve)$.
\item \label{ex:16} Let $\L\in Sym(\Ve)$ and positive definite and $\Ro\in Orth(\Ve)^+$, then prove that $\Ro\L\Ro^\top\in Sym(\Ve)$ and that it is positive definite.
\item Prove that the spectrum of  $\L^{sph}$ is composed of only
\bes
\lambda^{sph}=\frac{1}{3}\tr\L,
\ees
and that any $\bu\in\S$ is an eigenvector. 
\item Prove that the eigenvalues $\lambda^{dev}$ of $\L^{dev}$ are given by
\bes
\lambda^{dev}=\lambda-\lambda^{sph},
\ees
where $\lambda$ is an eigenvalue of $\L$.

\end{enumerate}

\chapter{Fourth rank tensors}
\label{ch:3}
\section{Fourth-rank tensors}
A {\it fourth-rank tensor} $\Lq$ is any linear application from $Lin(\Ve)$ to $Lin(\Ve)$:
\bes
\Lq:Lin(\Ve)\rightarrow Lin(\Ve)\ |\ \Lq(\alpha_i\A_i)=\alpha_i\Lq\A_i\ \forall\alpha_i\in\R,\ \A_i\in Lin(\Ve),\ i=1,...,n.
\ees
Defining the {\it sum of two fourth-rank tensors} as
\bes
\label{eq:tenssum4}
(\Lq_1+\Lq_2)\A=\Lq_1\A+\Lq_2\A\ \ \forall\A\in Lin(\Ve),
\ees
the {\it product of a scalar by a fourth-rank tensor} as
\bes
(\alpha\Lq)\A=\alpha(\Lq\A)\ \ \forall\alpha\in\R,\A\in Lin(\Ve)
\ees
and the {\it null fourth-rank tensor} $\Oq$ as the unique tensor such that
\bes
\Oq\A=\O\ \forall\A\in Lin(\Ve),
\ees
then the set of all the tensors $\Lq$ that operate on $Lin(\Ve)$ forms a vector space, denoted by $\Lq$in$(\Ve)$. We define the {\it  fourth-rank identity tensor} $\Iq$ as the unique tensor such that
\bes
\Iq\A=\A\ \ \forall \A\in Lin(\Ve).
\ees
It is apparent that the algebra of fourth-rank tensors is similar to that of second-rank tensors and, in fact, several operations with fourth-rank tensors can be introduced in almost the same way, in some sense shifting from $\Ve$ to $Lin(\Ve)$ the operations. However, the algebra of fourth-rank tensors is richer than that of the second-rank ones and some care must be taken.

In the following sections, we consider some of the operations that can be done with fourth-rank tensors.

\section{Dyads, tensor components}
For any couple of tensors $\A$ and $\B\in Lin(\Ve)$, the ({\it tensor}) {\it dyad} $\A\otimes\B$ is the fourth-rank tensor defined by
\bes
(\A\otimes\B)\L:=\B\cdot\L\ \A\ \ \forall\L\in Lin(\Ve).
\ees
The application defined above is actually a fourth-rank tensor because of the bilinearity of the scalar product of second-rank tensors. Applying this rule to the nine dyads of the basis $\Ba^2=\{\e_i\otimes\e_j,\ i,j=1,2,3\}$ of $Lin(\Ve)$ leads to the  introduction of the 81 fourth-rank tensors 
\bes
\e_i\otimes\e_j\otimes\e_k\otimes\e_l:=(\e_i\otimes\e_j)\otimes(\e_k\otimes\e_l)
\ees
that form a basis, $\Ba^4=\{\e_i\otimes\e_j\otimes\e_k\otimes\e_l,\ i,j=1,2,3\}$, for $\Lq$in$(\Ve)$. We remark hence that dim$(\Lq$in$(\Ve))=81$. A useful result is that
\be
\label{eq:t4dyadsondyad}
(\e_i\otimes\e_j\otimes\e_k\otimes\e_l)(\e_p\otimes\e_q)=(\e_k\otimes\e_l)\cdot(\e_p\otimes\e_q)(\e_i\otimes\e_j)=\delta_{kp}\delta_{lq}(\e_i\otimes\e_j).
\ee
 Any fourth-rank tensor  can be expressed as a linear combination ({the \it canonical decomposition}):
\bes
\Lq=L_{ijkl}\ \e_i\otimes\e_j\otimes\e_k\otimes\e_l,\ \ i,j=1,2,3,
\ees
where the $L_{ijkl}$s are the 81 {\it Cartesian components} of $\Lq$ with respect to $\Ba^4$. The $L_{ijkl}$s are defined by the operation
\bes
\besp
(\e_i\otimes\e_j)\cdot\Lq(\e_k\otimes\e_l)&=(\e_i\cdot\e_j)\cdot(L_{pqrs}\e_p\otimes\e_q\otimes\e_r\otimes\e_s)(\e_k\otimes\e_l)\\
&=(\e_i\otimes\e_j)\cdot(L_{pqrs}\delta_{rk}\delta_{sl}\e_p\otimes\e_q)\\
&=L_{pqrs}\delta_{rk}\delta_{sl}\delta_{ip}\delta_{jq}=L_{ijkl}.
\end{split}
\ees
The components of a tensor dyad can  be computed without any difficulty:
\bes
\besp
\A\otimes\B=(A_{ij}\e_i\otimes\e_j)\otimes (B_{kl}\e_k\otimes\e_l)=A_{ij}B_{kl}\e_i\otimes\e_j\otimes \e_k\otimes\e_l&\ \Rightarrow\\ 
(\A\otimes\B)_{ijkl}=A_{ij}B_{kl}&,
\end{split}
\ees
so that, in particular,
\bes
\label{eq:compdyad4}
((\a\otimes\b)\otimes(\c\otimes\d))_{ijkl}=
a_ib_jc_kd_l.
\ees
Concerning the identity of $\Lq$in$(\Ve)$,
\bes
\besp
I_{ijkl}=(\e_i\otimes\e_l)\cdot\Iq(\e_k\otimes\e_l)=(\e_i\otimes\e_j)\cdot(\e_k\otimes\e_l)=\e_i\cdot\e_k\e_j\cdot\e_l=\delta_{ik}\delta_{jl}&\ \Rightarrow\\
\Iq=\delta_{ik}\delta_{jl}(\e_i\otimes\e_l\otimes\e_k\otimes\e_l)&.
\end{split}
\ees
The components of  $\A\in Lin(\Ve)$, resulting from the application of  $\Lq\in\Lq$in$(\Ve)$ on  $\B\in Lin(\Ve)$, can now be easily calculated:
\be
\label{eq:Lu4}
\besp
\A=\Lq\B&=L_{ijkl}(\e_i\otimes\e_j\otimes\e_k\otimes\e_l)(B_{pq}\e_p\otimes\e_q)\\
&=L_{ijkl}B_{pq}\delta_{kp}\delta_{lq}(\e_i\otimes\e_j)\\
&=L_{ijkl}B_{kl}(\e_i\otimes\e_j)\ \Rightarrow\ A_{ij}= L_{ijkl}B_{kl}.
\end{split}
\ee
Moreover,
\bes
\besp
\label{eq:t4ondyad}
\Lq(\A\otimes\B)\C=\Lq((\A\otimes\B)\C)=\Lq(\B\cdot\C\A)=\B\cdot\C\ \Lq\A=((\Lq\A)\otimes\B)\C&\ \Rightarrow\\\Lq(\A\otimes\B)=(\Lq\A)\otimes\B&.
\end{split}
\ees
Using this result and Eq. (\ref{eq:t4dyadsondyad}), we can determine the components of a product of fourth-rank tensors:
\be
\label{eq:multt4}
\besp
\Aq\Bq&=A_{ijkl}(\e_i\otimes\e_j\otimes\e_k\otimes\e_l)B_{pqrs}(\e_p\otimes\e_q\otimes\e_r\otimes\e_s)\\
&=A_{ijkl}B_{pqrs}(\e_i\otimes\e_j\otimes\e_k\otimes\e_l)(\e_p\otimes\e_q)\otimes(\e_r\otimes\e_s)\\
&=A_{ijkl}B_{pqrs}[(\e_i\otimes\e_j\otimes\e_k\otimes\e_l)(\e_p\otimes\e_q)]\otimes(\e_r\otimes\e_s)\\
&=A_{ijkl}B_{pqrs}[\delta_{kp}\delta_{lq}(\e_i\otimes\e_j)]\otimes(\e_r\otimes\e_s)\\
&=A_{ijkl}B_{klrs}(\e_i\otimes\e_j\otimes\e_r\otimes\e_s)\ \Rightarrow\ (\Aq\Bq)_{ijrs}=A_{ijkl}B_{klrs}.
\end{split}
\ee

Depending upon four indices, a fourth-rank tensor $\Lq$ cannot be represented by a matrix; however, we will see in Section \ref{sec:kelvin} that a matrix representation of a fourth-rank tensor is still possible and that it is currently used in some cases, e.g. in elasticity.

\section{Conjugation product, transpose, symmetries}
For any two tensors $\A,\B\in Lin(\Ve)$ we call {\it conjugation product} the  tensor $\A\boxtimes\B\in\Lq$in$(\Ve)$ defined by the operation
\bes
(\A\boxtimes\B)\L:=\A\L\B^\top\ \forall\L\in Lin(\Ve).
\ees
As a consequence, for the dyadic tensors of $\Ba^2$, 
\be
\label{eq:t4_1}
(\e_i\otimes\e_j)\boxtimes(\e_k\otimes\e_l)=\e_i\otimes\e_k\otimes\e_j\otimes\e_l,
\ee
so that
\bes
(\A\boxtimes\B)_{ijkl}=A_{ik}B_{jl}.
\ees
Moreover, by the uniqueness of the identity $\Iq,\ \forall\A\in Lin(\Ve)$,
\bes
(\I\boxtimes\I)\A=\I\A\I^\top=\A\ \Rightarrow\ \Iq=\I\boxtimes\I.
\ees
The {\it transpose} of a fourth-rank tensor $\Lq$ is the unique tensor $\Lq^\top$ such that
\bes
\A\cdot(\Lq\B)=\B\cdot(\Lq^\top\A)\ \forall\A,\B\in Lin(\Ve).
\ees
By this definition, setting $\A=\e_i\otimes\e_j,\B=\e_k\otimes\e_l$ gives
\bes
(L^\top)_{ijkl}=L_{klij}.
\ees
A consequence  is that
\bes
\A\cdot(\Lq\B)=\B\cdot(\Lq^\top\A)=\A\cdot(\Lq^\top)^\top\B\ \Rightarrow\ (\Lq^\top)^\top=\Lq.
\ees
Moreover,  
\bes
\besp
\M\cdot(\A\otimes\B)^\top\L&=\L\cdot(\A\otimes\B)\M\\&=\L\cdot\A\M\cdot\B=\M\cdot(\B\A\cdot\L)\\
&=\M\cdot(\B\otimes\A)\L\ \Rightarrow \ (\A\otimes\B)^\top=\B\otimes\A,
\end{split}
\ees
while, cf. Exercise 7, Chapter 2,
\bes
\besp
\M\cdot(\A\boxtimes\B)^\top\L&=\L\cdot(\A\boxtimes\B)\M\\
&=\L\cdot\A\M\B^\top=\A^\top\L\cdot\M\B^\top=\M^\top\A^\top\L\cdot\B^\top\\
&=(\M^\top\A^\top\L)^\top\cdot(\B^\top)^\top=\L^\top\A\M\cdot\B=\A\M\cdot\L\B\\
&=\M\cdot\A^\top\L\B=\M\cdot(\A^\top\boxtimes\B^\top)\L\ \Rightarrow\\
&\ \ \ \  (\A\boxtimes\B)^\top=\A^\top\boxtimes\B^\top.
\end{split}
\ees
The property
\bes
(\Aq\Bq)^\top=\Bq^\top\Aq^\top
\ees
can be proved in the same manner used for the analogous property of the second-rank tensors.

A tensor $\Lq\in \Lq$in$(\Ve)$ is {\it symmetric} $\iff\Lq=\Lq^\top$. It is then evident that 
\bes
\Lq=\Lq^\top\Rightarrow\ L_{ijkl}=L_{klij},
\ees
relations that are  known as  {\it major symmetries}. There are 36 major symmetries  on the whole, so that a symmetric fourth-rank tensor has 45 independent components. Moreover,
\bes
\besp
&\A\boxtimes\B=(\A\boxtimes\B)^\top=\A^\top\boxtimes\B^\top\iff\A=\A^\top,\B=\B^\top,\\
&\A\otimes\B=(\A\otimes\B)^\top=\B\otimes\A\iff\B=\lambda\A,\ \lambda\in\R.
\end{split}
\ees
Let us now consider the case of a $\Lq\in\Lq$in$(\Ve)$ such that
\bes
\Lq\A=(\Lq\A)^\top\ \  \forall\A\in Lin(\Ve).
\ees
Then, by Eq. (\ref{eq:Lu4}),
\bes
L_{ijkl}=L_{jikl},
\ees
relations that are  called {\it left minor symmetries}: a tensor $\Lq$ having the  left minor symmetries has values in $Sym(\Ve)$.  On the whole, the number of left minor symmetries is 27. Finally, consider the case of a $\Lq\in\Lq$in$(\Ve)$ such that
\bes
\Lq\A=\Lq(\A^\top)\ \  \ \forall\A\in Lin(\Ve);
\ees
then, again by Eq. (\ref{eq:Lu4}), we get 
\bes
L_{ijkl}=L_{ijlk},
\ees
which are relations called {\it minor right-symmetries}, whose total number is also 27. It is immediate to recognize that if $\Lq$ has the minor right-symmetries, then
\bes
\Lq\W=\O\ \ \ \forall\W\in Skw(\Ve).
\ees
We say that {\it a tensor has  the minor symmetries} if it has both the right and left minor symmetries; the total number of minor symmetries is 45, because, as can be easily checked, some of the left and right minor symmetries are the same, so finally  a tensor with the minor symmetries has 36 independent components.

If $\Lq\in\Lq$in$(\Ve)$ has the major and minor symmetries, then the number of independent symmetry relations is actually 60 (some minor and major symmetries coincide), so in such a case $\Lq$ depends upon 21 independent components only. This is the case, for instance, of the classical elasticity tensor.

Finally, the 6 {\it Cauchy-Poisson symmetries}\footnote{The Cauchy-Poisson symmetries have played an important role in a celebrated diatribe of the XIXth century in elasticity, that between the so-called {\it rari-} and {\it muti-constant} theories.} are those of the type
\bes
L_{ijkl}=L_{ikjl}.
\ees
A tensor having the major, minor and Cauchy-Poisson symmetries is  said to be {\it completely symmetric}, i.e. swapping any couple of indices gives an identical component. In that case, the number of independent components is  only 15.

\section{Trace and scalar product of fourth-rank tensors}
We can introduce the scalar product between fourth-rank tensors in the same way we did for second-rank tensors. We first introduce the concept of {\it trace for fourth-rank tensors}, denoted by $\tr_4$, once again using the dyad (here, the tensor dyad):
\bes
\tr_4\A\otimes\B:=\A\cdot\B.
\ees
The easy proof that $\tr_4:\Lq$in$(\Ve)\rightarrow\R$ is a linear form is based upon the properties of the scalar product of second-rank tensors and it is left to the reader. An immediate result is that
\bes
\tr_4\A\otimes\B=A_{ij}B_{ij},
\ees
Then, using the canonical decomposition, we have that
\bes
\tr_4\Lq=\tr_4(L_{ijkl}(\e_i\otimes\e_j)\otimes(\e_k\otimes\e_l))=L_{ijkl}(\e_i\otimes\e_j)\cdot(\e_k\otimes\e_l)=L_{ijkl}\delta_{ik}\delta_{jl}=L_{ijij}
\ees
and that
\bes
\tr_4\Lq^\top=
\tr_4(L_{klij}(\e_i\otimes\e_j)\otimes(\e_k\otimes\e_l))=L_{klij}(\e_i\otimes\e_j)\cdot(\e_k\otimes\e_l)=L_{klij}\delta_{ik}\delta_{jl}=L_{ijij}=\tr_4\Lq.
\ees
Then, we define the {\it scalar product of fourth-rank tensors} as
\bes
\Aq\cdot\Bq:=\tr_4(\Aq^\top\Bq).
\ees
By the properties of $\tr_4$, the scalar product is a positive definite symmetric bilinear form:
\bes
\besp
&\alpha\Aq\cdot\beta\Bq=\tr_4(\alpha\Aq^\top\beta\Bq)=\alpha\beta\tr_4(\Aq^\top\Bq)=\alpha\beta\Aq\cdot\Bq,\\
&\Aq\cdot\Bq=\tr_4(\Aq^\top\Bq)=\tr_4(\Aq^\top\Bq)^\top=\tr_4(\Bq^\top\Aq)=\Bq\cdot\Aq,\\
&\Aq\cdot\Aq=\tr_4(\Aq^\top\Aq)=(\Aq^\top\Aq)_{ijij}=A_{klij}A_{klij}>0\ \forall\Aq\in{\Lq}in(\Ve), \Aq\cdot\Aq=0\iff\Aq=\Oq.
\end{split}
\ees
By components
\bes
\besp
\Aq\cdot\Bq&=\tr_4((A_{klij}\e_i\otimes\e_j\otimes\e_k\otimes\e_l)(B_{pqrs}\e_p\otimes\e_q\otimes\e_r\otimes\e_s))\\
&=\tr_4(A_{klij}B_{pqrs}\delta_{kp}\delta_{lq}(\e_i\otimes\e_j)\otimes(\e_r\otimes\e_s))\\
&=A_{klij}B_{pqrs}\delta_{kp}\delta_{lq}(\e_i\otimes\e_j)\cdot(\e_r\otimes\e_s)=A_{klij}B_{pqrs}\delta_{kp}\delta_{lq}\delta_{ir}\delta_{js}=A_{klij}B_{klij}.
\end{split}
\ees
The rule for computing the scalar product is hence always the same, as was already seen for vectors and second-rank tensors: All the indexes are to be saturated.

In complete analogy with vectors and second-rank tensors, we say that $\Aq$ is {\it orthogonal to} $\Bq\iff$
\bes
\Aq\cdot\Bq=0
\ees
and we define the {\it norm of } $\Lq$ as
\bes
|\Lq|:=\sqrt{\Lq\cdot\Lq}=\sqrt{\tr_4\Lq^\top\Lq}=\sqrt{L_{ijkl}L_{ijkl}}.
\ees

\section{Projectors, identities}
For the spherical part of any $\A\in Sym(\Ve)$ we can write
\bes
\A^{sph}:=\frac{1}{3}\tr\A\ \I=\frac{1}{3}\I\cdot\A\ \I=\frac{1}{3}(\I\otimes\I)\A=\Sq^{sph}\A,
\ees
where
\bes
\Sq^{sph}:=\frac{1}{3}\I\otimes\I
\ees
is the {\it spherical projector}, i.e. the fourth-rank tensor that extracts from any $\A\in Lin(\Ve)$ its spherical part. 
Moreover,
\bes
\A^{dev}:=\A-\A^{sph}=\Iq\A-\Sq^{sph}\A=\Dq^{dev}\A,
\ees
where
\bes
\Dq^{dev}:=\Iq-\Sq^{sph}
\ees
is the {\it deviatoric projector}, i.e. the fourth-rank tensor that extracts from any $\A\in Lin(\Ve)$ its deviatoric part.
It is worth noting that
\bes
\Iq=\Sq^{sph}+\Dq^{dev}.
\ees
Moreover, about the components of $\Sq^{sph}$,
\bes
\besp
S^{sph}_{ijkl}&=(\e_i\otimes\e_j)\cdot\frac{1}{3}(\I\otimes\I)(\e_k\otimes\e_l)=\frac{1}{3}\I\cdot(\e_i\otimes\e_j)\ \I\cdot(\e_k\otimes\e_l)\\
&=\frac{1}{3}\tr(\e_i\otimes\e_j)\tr(\e_k\otimes\e_l)=\frac{1}{3}\delta_{ij}\delta_{kl}\ \rightarrow\ 
\Sq^{sph}=\frac{1}{3}\delta_{ij}\delta_{kl}(\e_i\otimes\e_j\otimes\e_k\otimes\e_l).
\end{split}
\ees
We remark that
\bes
\Sq^{sph}=(\Sq^{sph})^\top.
\ees
We introduce now the tensor $\Iq^s$, {\it restriction of $\Iq$ to $\A\in Sym(\Ve)$}. It can be introduced as follows: $\forall \A\in Sym(\Ve)$
\bes
\A=\frac{1}{2}(\A+\A^\top),
\ees
and
\bes
\A=\Iq\A=\frac{1}{2}(\Iq\A+\Iq\A^\top)=\frac{1}{2}(I_{ijkl}A_{kl}+I_{ijkl}A_{lk})(\e_i\otimes\e_j\otimes\e_k\otimes\e_l);
\ees
because $\A=\A^\top$,  there is insensitivity to the swap of indexes $k$ and $l$, so
\bes
\A=\frac{1}{2}(I_{ijkl}A_{kl}+I_{ijlk}A_{lk})(\e_i\otimes\e_j\otimes\e_k\otimes\e_l)=\frac{1}{2}(\delta_{ik}\delta_{jl}+\delta_{il}\delta_{jk})A_{kl}(\e_i\otimes\e_j\otimes\e_k\otimes\e_l).
\ees
Then, if we admit the interchangeability of indexes $k$ and $l$, i.e. if we postulate the existence of the  minor right-symmetries for $\Iq$, then $\Iq=\Iq^s$, with
\bes
\Iq^s=\frac{1}{2}(\delta_{ik}\delta_{jl}+\delta_{il}\delta_{jk})(\e_i\otimes\e_j\otimes\e_k\otimes\e_l).
\ees
It is apparent that
\bes
I^s_{ijkl}=I^s_{klij},
\ees
i.e. $\Iq^s=(\Iq^s)^\top$, but also that
\bes
I^s_{ijkl}=\frac{1}{2}(\delta_{il}\delta_{jk}+\delta_{ik}\delta_{jl})=I^s_{jikl},
\ees
i.e., $\Iq^s$ has also the minor left-symmetries; in other words, $\Iq^s$ has the major and minor symmetries, like an elasticity tensor, while this is not the case for $\Iq$. In fact
\bes
I_{ijkl}=I_{jilk}=\delta_{ik}\delta_{jl}\neq\delta_{il}\delta_{jk}=I_{jikl}=I_{ijlk}.
\ees
Because $\Sq^{sph}$ and $\Dq^{dev}$ operate on $Sym(\Ve)$, it is immediate to recognize that it is also
\bes
\Dq^{dev}=\Iq^s-\Sq^{sph}\ \Rightarrow\ \Iq^s=\Sq^{sph}+\Dq^{dev}.
\ees
It is worth noting that
\bes
(\Dq^{dev})^\top=(\Iq^s-\Sq^{sph})^\top=(\Iq^s)^\top-(\Sq^{sph})^\top=\Iq^s-\Sq^{sph}=\Dq^{dev}.
\ees
We can now determine the components of $\Dq^{dev}$:
\bes
\besp
&D^{dev}_{ijkl}=I^s_{ijkl}-S^{sph}_{ijkl}=\frac{1}{2}(\delta_{ik}\delta_{jl}+\delta_{il}\delta_{jk})-\frac{1}{3}\delta_{ij}\delta_{kl}\ \rightarrow\\ &\Dq^{dev}=\left[\frac{1}{2}(\delta_{ik}\delta_{jl}+\delta_{il}\delta_{jk})-\frac{1}{3}\delta_{ij}\delta_{kl}\right](\e_i\otimes\e_j\otimes\e_k\otimes\e_l).
\end{split}
\ees
We remark that the result (\ref{eq:orthsphdev}) implies that $\Sq^{sph}$ and $\Dq^{dev}$ are {\it orthogonal projectors}, i.e. they project the same $\A\in Sym(\Ve)$ into two orthogonal subspaces of $\Ve$, $Sph(\Ve)$ and $Dev(\Ve)$.

The tensor $\Tq^{trp}\in\Lq$in$(\Ve)$ defined by the operation
\bes
\Tq^{trp}\A=\A^\top\ \ \ \forall\A\in Lin(\Ve),
\ees
is the {\it transposition projector}, whose components are
\bes
T^{trp}_{ijkl}=(\e_i\otimes\e_j)\cdot\Tq^{trp}(\e_k\otimes\e_l)=(\e_i\otimes\e_j)\cdot(\e_l\otimes\e_k)=\delta_{il}\delta_{jk}.
\ees

The following operation defines the {\it symmetry projector} $\Sq^{sym}\in\Lq$in$(\Ve)$:
\bes
\Sq^{sym}\A=\frac{1}{2}(\A+\A^\top)\ \forall\A\in Lin(\Ve),
\ees
while the {\it antisymmetry projector} $\Wq^{skw}\in\Lq$in$(\Ve)$ is defined by
\bes
\Wq^{skw}\A=\frac{1}{2}(\A-\A^\top)\ \forall\A\in Lin(\Ve).
\ees
Also $\Sq^{sym}$ and $\Wq^{skw}$ are orthogonal projectors, because they project the same $\A\in Lin(\Ve)$ into two orthogonal subspaces of $Lin(\Ve)$: $Sym(\Ve)$ and $Skw(\Ve)$, see Exercise 10, Chapter 2.

We prove now two properties of the projectors: $\forall\A\in Lin(\Ve)$,
\be
\label{eq:propproj1}
(\Sq^{sym}+\Wq^{skw})\A=\frac{1}{2}(\A+\A^\top)+\frac{1}{2}(\A-\A^\top)=\A=\Iq\A\ \Rightarrow\ \Sq^{sym}+\Wq^{skw}=\Iq.
\ee
Then,
\be
\label{eq:propproj2}
(\Sq^{sym}-\Wq^{skw})\A=\frac{1}{2}(\A+\A^\top)-\frac{1}{2}(\A-\A^\top)=\A^\top=\Tq^{trp}\A\ \Rightarrow\ \Sq^{sym}-\Wq^{skw}=\Tq^{trp}.
\ee

\section{Orthogonal conjugator}
For any $\U\in Orth(\Ve)$, we define its {\it orthogonal conjugator} $\Uq\in\Lq$in$(\Ve)$ as
\bes
\Uq:=\U\boxtimes\U.
\ees
\begin{teo} {\bf (Orthogonality of $\Uq$).} The orthogonal conjugator is an orthogonal tensor of $\Lq$in$(\Ve)$, i.e. it preserves the scalar product between tensors:
\bes
\Uq\A\cdot\Uq\B=\A\cdot\B\ \ \forall\A,\B\in Lin(\Ve).
\ees
\begin{proof}
By the assertion in Exercise 12 of Chapter 2 
and because $\U\in Orth(\Ve)$, we have
\bes
\besp
\Uq\A\cdot\Uq\B&=(\U\boxtimes\U)\A\cdot(\U\boxtimes\U)\B=\U\A\U^\top\cdot\U\B\U^\top\\
&=\U^\top\U\A\U^\top\cdot\B\U^\top=
\A\U^\top\cdot\B\U^\top=\A\U^\top\U\cdot\B=\A\cdot\B.
\end{split}
\ees
\end{proof}
\end{teo}
Just as for tensors of $Orth(\Ve)$, we also have
\bes
\Uq\Uq^\top=\Uq^\top\Uq=\Iq.
\ees
In fact, see the assertion in Exercise 4:
\be
\label{eq:unitary4}
\Uq\Uq^\top=(\U\boxtimes\U)(\U^\top\boxtimes\U^\top)=\U\U^\top\boxtimes\U\U^\top=\I\boxtimes\I=\Iq.
\ee
The orthogonal conjugators also have  some properties in relation with projectors:
\begin{teo} $\Sq^{sph}$ is unaffected by any orthogonal conjugator, while $\Dq^{dev}$ commutes with any orthogonal conjugator.
\begin{proof}
For any $\L\in Sym(\Ve)$ and $\U\in Orth(\Ve)$,
\bes
\besp
\Uq\Sq^{sph}\L&=(\U\boxtimes\U)\left(\frac{1}{3}\I\otimes\I\right)\L=\frac{1}{3}(\tr\L)(\U\boxtimes\U)\I=\frac{1}{3}(\tr\L)\U\I\U^\top\\
&=\frac{1}{3}(\tr\L)\I=\frac{1}{3}\I\cdot\L\ \I=\frac{1}{3}(\I\otimes\I)\L=\Sq^{sph}\L.
\end{split}
\ees
Moreover,
\bes
\besp
\Sq^{sph}\Uq\L&=\left(\frac{1}{3}\I\otimes\I\right)(\U\boxtimes\U)\L=\frac{1}{3}(\I\otimes\I)(\U\L\U^\top)=\frac{1}{3}(\I\cdot\U\L\U^\top)\I\\
&=\frac{1}{3}\tr(\U\L\U^\top)\I=\frac{1}{3}\tr(\U^\top\U\L)\I=\frac{1}{3}(\tr\L)\I=\frac{1}{3}\I\cdot\L\I=\frac{1}{3}(\I\otimes\I)\L=\Sq^{sph}\L.
\end{split}
\ees
Thus, we have proved that
\bes
\Sq^{sph}\Uq=\Uq\Sq^{sph}=\Sq^{sph},
\ees
i.e. that the spherical projector $\Sq^{sph}$ is unaffected by any orthogonal conjugator. Furthermore
\bes
\Dq^{dev}\Uq\L=(\Iq^s-\Sq^{sph})\Uq\L=\Iq^s\Uq\L-\Sq^{sph}\Uq\L=\Uq\L-\Sq^{sph}\L=(\Uq-\Sq^{sph})\L
\ees
and
\bes
\Uq\Dq^{dev}\L=\Uq(\Iq^s-\Sq^{sph})\L=\Uq\Iq^s\L-\Uq\Sq^{sph}\L=\Uq\L-\Sq^{sph}\L=(\Uq-\Sq^{sph})\L,
\ees
so that
\bes
\Dq^{dev}\Uq=\Uq\Dq^{dev}.
\ees
\end{proof}
\end{teo}

\section{Rotations and symmetries}
We ponder now how to rotate a fourth-rank tensor, i.e., what are the components of 
\bes
\Lq=L_{ijkl}\e_i\otimes\e_j\otimes\e_k\otimes\e_l
\ees
in a basis $\Ba'=\{\e_1',\e_2',\e_3'\}$ obtained rotating the basis $\Ba=\{\e_1,\e_2,\e_3\}$ by the rotation $\R=R_{ij}\e_i\otimes\e_j, \R\in Orth(\Ve)^+$. The procedure is exactly the same already followed for vectors and second-rank tensors:
\bes
\besp
\Lq&=L_{ijkl}\e_i\otimes\e_j\otimes\e_k\otimes\e_l=L_{ijkl}R^\top_{pi}\e_p'\otimes R^\top_{qj}\e_q'\otimes R^\top_{rk}\e_r'\otimes R^\top_{sl}\e_s'\\
&=R^\top_{pi}R^\top_{qj}R^\top_{rk}R^\top_{sl}L_{ijkl}\e_p'\otimes\e_q'\otimes\e_r'\otimes\e_s',
\end{split}
\ees
i.e.
\bes
L_{pqrs}'=R^\top_{pi}R^\top_{qj}R^\top_{rk}R^\top_{sl}L_{ijkl}.
\ees
We see clearly that the components of $\Lq$ in the basis $\Ba'$ are a linear combination of those in $\Ba$, the coefficients of the linear combination being fourth powers of the director cosines, the $R_{ij}$s. 
The introduction of the orthogonal conjugator\footnote{Here the symbol $\R$ indicates the orthogonal conjugator of $\Ro$, not the set of real numbers.} of the rotation $\Ro$,
\bes
\R=\Ro\boxtimes\Ro,
\ees
allows us to give a compact expression for the rotation of second- and fourth-rank tensors (for completeness we recall also that of a vector $\bw$); 
\bes
\begin{array}{c}
\bw'=\Ro^\top\bw,\medskip\\
\L'=\Ro^\top\L\Ro=(\Ro^\top\boxtimes\Ro^\top)\L=\R^\top\L,\medskip\\
\Lq'=(\Ro^\top\boxtimes\Ro^\top)\Lq(\Ro\boxtimes\Ro)=\R^\top\Lq\R.
\end{array}
\ees
Checking  the above relations with the orthogonal conjugator $\R$ is left to the reader. 
It is worth noting that, actually, these transformations are valid not only for $\Ro\in Orth(\Ve)^+$, but more generally for any $\U\in Orth(\Ve)$, i.e. also for symmetries.

If  $\U$  denotes the tensor of change of basis under any orthogonal transformation, i.e. if we put $\U=\Ro^\top$ for the rotations, then the above relations become
\be
\label{eq:orthtransfor}
\begin{array}{c}
\bw'=\U\bw,\medskip\\
\L'=\U\L\U^\top=(\U\boxtimes\U)\L=\Uq\L,\medskip\\
\Lq'=(\U\boxtimes\U)\Lq(\U\boxtimes\U)^\top=\Uq\Lq\Uq^\top.
\end{array}
\ee

Finally, we say that $\L\in Lin(\Ve)$ or $\Lq\in\Lq$in$(\Ve)$ is {\it invariant under an orthogonal transformation} $\U$ if
\bes
\U\L\U^\top=\L,\ \ \ \Uq\Lq\Uq^\top=\Lq;
\ees
right multiplying both terms by $\U$ or by $\Uq$ and through  Eq. (\ref{eq:unitary4}), we get that $\L$ or $\Lq$ are invariant under $\U\iff$
\bes
\U\L=\L\U,\ \ \ \Uq\Lq=\Lq\Uq,
\ees
i.e.  $\iff\L$ and $\U$, or $\Lq$ and $\Uq$ commute. This relation allows, for example, the analysis of material symmetries in anisotropic elasticity.

If a tensor is invariant under {\it any} orthogonal transformation, i.e. if the previous equations hold true $\forall\U\in Orth(\Ve)$,
then the tensor is said to be {\it isotropic}. 
A general result\footnote{Actually, this is  quite a famous result in classical elasticity, the {\it Lamé's equation}, defining an isotropic elastic material.} is that a fourth-rank tensor $\Lq$ is isotropic $\iff$ there exist two scalar functions $\lambda,\mu$ such that
\bes
\Lq\A=2\mu \A+\lambda\tr\A\ \I\ \ \ \forall\A\in Sym(\Ve).
\ees 
The reader is referred to the book of Gurtin (see the references) for the proof of this result and  for a deeper insight into isotropic functions.

\section{The Kelvin formalism}
\label{sec:kelvin}
As already mentioned, though fourth-rank tensors cannot be organized in and represented by a matrix, nevertheless a matrix formalism for these operators exists. 
Such formalism is due to Kelvin\footnote{W. Thomson (Lord Kelvin):  Elements of a mathematical theory of elasticity. {\it Philos. Trans. R. Soc.}, 146, 481-498, 1856. Later, Voigt (W. Voigt: {\it Lehrbuch der Kristallphysik}. B. G. Taubner, Leipzig, 1910) gave another, similar matrix formalism for tensors, more widely known than the Kelvin one, but less effective.} and it is strictly related to the theory of elasticity, i.e. it concerns the Cauchy's stress tensor $\bsig$, the strain tensor $\beps$ and the elasticity tensor $\Ey$. The relation between $\bsig$ and $\beps$ is given by the celebrated (generalized) {\it Hooke's law}:
\bes
\bsig=\Ey\beps.
\ees
Both $\bsig,\beps\in Sym(\Ve)$ while $\Ey=\Ey^\top$ and it has also the minor symmetries, so $\Ey$ has just 21 independent components\footnote{Actually, the Kelvin formalism can also be extended without major difficulties  to tensors that do not possess all the symmetries.}. 
In the Kelvin formalism, the six independent components of $\bsig$ and $\beps$ are organized into column vectors and renumbered as follows
\bes
\{\sigma\}=\left\{
\begin{array}{c}
\sigma_1=\sigma_{11}\\
\sigma_2=\sigma_{22}\\
\sigma_3=\sigma_{33}\\
\sigma_4=\sqrt{2}\sigma_{23}\\
\sigma_5=\sqrt{2}\sigma_{31}\\
\sigma_6=\sqrt{2}\sigma_{12}
\end{array}
\right\},\ \ \
\{\varepsilon\}=\left\{
\begin{array}{c}
\varepsilon_1=\varepsilon_{11}\\
\varepsilon_2=\varepsilon_{22}\\
\varepsilon_3=\varepsilon_{33}\\
\varepsilon_4=\sqrt{2}\varepsilon_{23}\\
\varepsilon_5=\sqrt{2}\varepsilon_{31}\\
\varepsilon_6=\sqrt{2}\varepsilon_{12}
\end{array}
\right\}.
\ees
The elasticity tensor $\Ey$ is reduced to a $6\times6$ matrix $[E]$ as a  consequence of the minor symmetries induced by the symmetry of $\bsig$ and $\beps$; this matrix is symmetric because $\Ey=\Ey^\top$:
\bes
\hspace{-2cm}[E]=\left[
\begin{array}{cccccc}
E_{11}=E_{1111}&E_{12}=E_{1122}&E_{13}=E_{1133}&E_{14}=\sqrt{2}E_{1123}&E_{15}=\sqrt{2}E_{1131}&E_{16}=\sqrt{2}E_{1112}\\
E_{12}=E_{1122}&E_{22}=E_{2222}&E_{23}=E_{2233}&E_{24}=\sqrt{2}E_{2223}&E_{25}=\sqrt{2}E_{2231}&E_{26}=\sqrt{2}E_{2212}\\
E_{13}=E_{1133}&E_{23}=E_{2233}&E_{33}=E_{3333}&E_{34}=\sqrt{2}E_{3323}&E_{35}=\sqrt{2}E_{3331}&E_{36}=\sqrt{2}E_{3312}\\
E_{14}=\sqrt{2}E_{1123}&E_{24}=\sqrt{2}E_{2223}&E_{34}=\sqrt{2}E_{3323}&E_{44}=2E_{2323}&E_{45}=2E_{2331}&E_{46}=2E_{2312}\\
E_{15}=\sqrt{2}E_{1131}&E_{25}=\sqrt{2}E_{2231}&E_{35}=\sqrt{2}E_{3331}&E_{45}=2E_{2331}&E_{55}=2E_{3131}&E_{56}=2E_{3112}\\
E_{16}=\sqrt{2}E_{1112}&E_{26}=\sqrt{2}E_{2212}&E_{36}=\sqrt{2}E_{3312}&E_{46}=2E_{2312}&E_{56}=2E_{3112}&E_{66}=2E_{1212}
\end{array}
\right].
\ees
In this way, the matrix product
\be
\label{eq:hookematrix}
\{\sigma\}=[E]\{\eps\}
\ee
is equivalent to the tensor form of the Hooke's law and all the operations can be done by the aid of classical  matrix algebra\footnote{Mehrabadi and Cowin have shown that the Kelvin formalism transforms second- and fourth-rank tensors on $\R^3$ into vectors and second-rank tensors on $\R^6$ (M. M. Mehrabadi, S. C. Cowin: {Eigentensors of linear anisotropic elastic materials}. {\it Q. J. Mech. Appl. Math.}, 43, 15-41, 1990).}, e.g. the computation of the inverse of $\Ey$, the {\it compliance tensor}.

An important operation is the expression of tensor $\Uq$ in Eq. (\ref{eq:orthtransfor}) in the Kelvin formalism; some tedious but straightforward passages give the result:
\bes
\label{eq:matrot}
\hspace{-6mm}[U]=
\begin{small}
\left[\begin{array}{cccccc}
U^2_{11} &U^2_{12} & U^2_{13}  & \sqrt{2}U_{12}U_{13} & \sqrt{2}U_{13}U_{11} & \sqrt{2}U_{11}U_{12} \\
U^2_{21}  &U^2_{22}  & U^2_{23}  & \sqrt{2}U_{22}U_{23} & \sqrt{2}U_{23}U_{21} & \sqrt{2}U_{21}U_{22}\\
U^2_{31}  & U^2_{32} & U^2_{33} & \sqrt{2}U_{32}U_{33} & \sqrt{2}U_{33}U_{31} & \sqrt{2}U_{31}U_{32} \\
\sqrt{2}U_{21}U_{31} & \sqrt{2}U_{22}U_{32} & \sqrt{2}U_{23}U_{33} & 
U_{23}U_{32}+U_{22}U_{33} & U_{33}U_{21}+U_{31}U_{23} & U_{31}U_{22}+U_{32}U_{21} \\
\sqrt{2}U_{31}U_{11} & \sqrt{2}U_{32}U_{12} & \sqrt{2}U_{33}U_{13} & 
U_{32}U_{13}+U_{33}U_{12} & U_{31}U_{13}+U_{33}U_{11} & U_{31}U_{12}+U_{32}U_{11} \\
 \sqrt{2}U_{11}U_{21} & \sqrt{2}U_{12}U_{22} & \sqrt{2}U_{13}U_{23} & 
U_{12}U_{23}+U_{13}U_{22} & U_{11}U_{23}+U_{13}U_{21} & U_{11}U_{22}+U_{12}U_{21}\end{array}\right]
 \end{small}
\ees
With some work, it can be checked that
\bes
[U][U]^\top=[U]^\top[U]=[I],
\ees
i.e. that $[U]$ is an orthogonal matrix in $\R^6$. Of course,
\bes
[R]=[U]^\top
\ees
is the matrix that in the Kelvin formalism represents the tensor $\Ro=\U^\top$. The change of basis for $\bsig$ and $\beps$ are hence done through the relations
\bes
\{\sigma'\}=[U]\{\sigma\},\ \ \ \{\eps'\}=[U]\{\eps\},
\ees
which applied to Eq. (\ref{eq:hookematrix}) give
\bes
\{\sigma\}=[E]\{\eps\}\ \rightarrow\ [U]^\top\{\sigma'\}=[E][U]^\top\{\eps'\}\ \rightarrow\ \{\sigma'\}=[U][E][U]^\top\{\eps'\}
\ees
i.e. in the basis $\Ba'$
\bes
\{\sigma'\}=[E']\{\eps'\},
\ees
where 
\bes
[E']=[U][E][U]^\top=[R]^\top[E][R]
\ees
is the matrix representing $\Ey$ in $\Ba'$ in the Kelvin formalism.
Though it is possible to give the expression of the components of $[E']$, they are so long that they are omitted here.

\section{The polar formalism for plane tensors}
The Cartesian representation of tensors makes use of quantities that are basis-dependent, and the change of basis implies algebraic transformations rather complicate. 
The question of representing tensors using other quantities than Cartesian components is hence of importance. In particular, it should be interesting to represent a tensor making use of only invariants of the tensor itself and of angles, the simplest geometrical way to determine a direction.

In the case of plane tensors this has been done by Verchery\footnote{G. Verchery: {\it Les invariants des tenseurs d'ordre 4 du type de l'élasticité}, Proc. Colloque EUROMECH 115, 1979.} who introduced the so-called {\it polar formalism}. This is basically a mathematical technique to find the invariants of a tensor of any rank. Here, we give just a short insight into the polar formalism of fourth-rank tensors of the elastic type, i.e. having the minor and major symmetries, omitting the proof of the results\footnote{A detailed presentation of the method can be found in P. Vannucci: {\it Anisotropic elasticity}, Springer, 2018.}.

The Cartesian components of a plane fourth-rank elasticity-type tensor $\Tq$ in a frame rotated through an angle $\theta$ can be expressed as
\bes
\besp
&T_{1111}=T_0+2T_1+R_0\cos4(\Phi_0-\theta)+4R_1\cos2(\Phi_1-\theta),\\
&T_{1112}=R_0\sin4(\Phi_0-\theta)+2R_1\sin2(\Phi_1-\theta),\\
&T_{1122}=-T_0+2T_1-R_0\cos4(\Phi_0-\theta),\\
&T_{1212}=T_0-R_0\cos4(\Phi_0-\theta),\\
&T_{1222}=-R_0\sin4(\Phi_0-\theta)+2R_1\sin2(\Phi_1-\theta),\\
&T_{2222}=T_0+2T_1+R_0\cos4(\Phi_0-\theta)-4R_1\cos2(\Phi_1-\theta).
\end{split}
\ees
In the above equations, $T_0,T_1,R_0,R_1$ are tensor invariants, with all of them non negative, while $\Phi_0$ and $\Phi_1$ are angles whose difference, $\Phi_0-\Phi_1$, is also a tensor invariant, so fixing one of the two polar angles corresponds to fixing a frame. 
In particular,  the tensor invariants have a direct physical meaning (e.g., for the elasticity tensor,  they are linked to  material symmetries and to  strain energy decomposition). We remark also that the change of frame is extremely simple in the polar formalism: It is sufficient to subtract the angle $\theta$ formed by the new frame from the two polar angles. 

The Cartesian expression of the polar invariants can be found inverting the previous expressions:
\bes
\besp
&T_0=\frac{1}{8}(T_{1111}-2T_{1122}+4T_{1212}+T_{2222}),\\
&T_1=\frac{1}{8}(T_{1111}+2T_{1122}+T_{2222}),\\
&R_0=\frac{1}{8}\sqrt{(T_{1111}-2T_{1122}-4T_{1212}+T_{2222})^2+16(T_{1112}-T_{1222})^2},\\
&R_1=\frac{1}{8}\sqrt{(T_{1111}-T_{2222})^2+4(T_{1112}+T_{1222})^2},\\
&\tan4\Phi_0=\frac{4(T_{1112}-T_{1222})}{T_{1111}-2T_{1122}-4T_{1212}+T_{2222}},\\
&\tan2\Phi_1=\frac{2(T_{1112}+T_{1222})}{T_{1111}-T_{2222}}.
\end{split}
\ees

\section{Exercises}
\begin{enumerate}
\item Prove Eq. (\ref{eq:t4_1}).
\item Prove that
\bes
(\Aq\Bq)^\top=\Bq^\top\Aq^\top.
\ees
\item Prove that
\bes
\A\otimes\B\Lq=\A\otimes\Lq^\top\B.
\ees
\item \label{ex:3_4} Prove  that
\bes
(\A\boxtimes\B)(\C\boxtimes\D)=\A\C\boxtimes\B\D.
\ees
\item Prove Eq. (\ref{eq:multt4}) using the result of the previous exercise.
\item Prove that
\bes
(\A\otimes\B)(\C\boxtimes\D)=\A\otimes((\C^\top\boxtimes\D^\top)\B).
\ees
\item Prove that
\bes
(\A\boxtimes\B)(\C\otimes\D)=((\A\boxtimes\B)\C)\otimes\D.
\ees
\item Let $\p\in\S$ and $\P=\p\otimes\p$, then prove that
\bes
\P\boxtimes\P=\P\otimes\P.
\ees
\item Prove that, $\forall\Aq\in\Lq$in$(\Ve)$,
\bes
\Iq\Aq=\Aq\Iq=\Aq.
\ees
\item Show that
\bes
(\A\otimes\B)\cdot(\C\otimes\D)=\A\cdot\C\ \B\cdot\D.
\ees
\item Show that
\bes
\Sq^{sph}=\frac{\I}{|\I|}\otimes\frac{\I}{|\I|}.
\ees
\item Show that 
\bes
dim(Sph(\Ve))=1,\ \ dim(Dev(\Ve))=5.
\ees
\item Show the following properties of $\Sq^{sph}$ and $\Dq^{dev}$:
\bes
\begin{array}{c}
\Sq^{sph}\Sq^{sph}=\Sq^{sph},\medskip\\
\Dq^{dev}\Dq^{dev}=\Dq^{dev},\medskip\\
\Sq^{sph}\Dq^{dev}=\Dq^{dev}\Sq^{sph}=\Oq.
\end{array}
\ees
\item Prove the results in Eqs. (\ref{eq:propproj1}) and (\ref{eq:propproj2}) using the components.
\item Show that
\bes
\begin{array}{c}
\Sq^{sph}\cdot\Sq^{sph}=1,\medskip\\
\Dq^{dev}\cdot\Dq^{dev}=5,\medskip\\
\Sq^{sph}\cdot\Dq^{dev}=0.
\end{array}
\ees
\item Make explicit the orthogonal conjugator $\Sq_R$ of the tensor $\Sy_R$ in Eq. (\ref{eq:sym2r}).
\item Using the polar formalism, it can be proved that the material symmetries conditions in plane elasticity are all condensed into the equation
\bes
R_0R_1\sin4(\Phi_0-\Phi_1)=0;
\ees
determine the different types of possibles elastic symmetries. 

\end{enumerate}

\chapter{Tensor analysis: curves}
\label{ch:4}

\section{Curves of points, vectors and tensors}
The scalar products in $\Ve,Lin(\Ve)$ and $\Lq$in$(\Ve)$ allow us to define a {\it norm}, the {\it Euclidean norm}, so they automatically endow these spaces with  a {\it metric}, i.e. we are able to measure and calculate a distance between two elements of such a space and in $\Eu$. This allows us to generalise the concepts of continuity and differentiability already known in $\R$, whose definition intrinsically makes use of a distance between real quantities.

Let $\pi_n=\{p_n\in\Eu,n\in\Nq\}$ be a sequence of points in $\Eu$. We say that $\pi_n$ {\it converges to} $p\in\Eu$ if
\bes
\lim_{n\rightarrow\infty} d(p_n-p)=0.
\ees
A similar definition can be given for sequences of vectors or tensors of any rank. Through this definition of convergence we can now make the concepts of continuity and of curve precise.

Let $[a,b]$ be an interval of $\R$; the function 
\bes
p=p(t):[a,b]\rightarrow\Eu
\ees
 is {\it continuous} at $t\in[a,b]$ if for each sequence $\{t_n\in[a,b],n\in\Nq\}$ that converges to $t$, the sequence $\pi_n$ defined by $p_n=p(t_n)\ \forall n\in\Nq$ converges to $p(t)\in\Eu$. The function $p=p(t)$ is a {\it curve in }$\Eu\iff$ it is continuous $\forall t\in[a,b]$. In the same way  we can define  curves of vectors and  tensors:
\bes
\begin{array}{c}
\bv=\bv(t):[a,b]\rightarrow\Ve,\medskip\\
\L=\L(t):[a,b]\rightarrow Lin(\Ve),\medskip\\
\Lq=\Lq(t):[a,b]\rightarrow {\Lq}\mathrm{in}(\Ve).
\end{array}
\ees
Mathematically, a curve is a function that lets correspond to  a real value $t$ (the {\it parameter}) in a given interval, an element of a space: $\Eu,\Ve, Lin(\Ve)$ or $\Lq(\Ve)$.

\section{Differention of curves}
\label{sec:diffcurves}
Let $\bv=\bv(t):[a,b]\rightarrow\Ve$ be a curve of vectors and $g=g(t):[a,b]\rightarrow\R$ a scalar function. We say that $\bv$ is {\it of the order o with respect to g in} $t_0\iff$
\bes
\lim_{t\rightarrow t_0}\frac{|\bv(t)|}{|g(t)|}=0,
\ees
and we write
\bes
\bv(t)=o(g(t))\ \mathrm{for}\ t\rightarrow t_0.
\ees
A similar definition can be given for a curve of tensors of any rank. We then say that the curve $\bv$ is {\it differentiable} in $t_0\in]a,b[\iff\exists \bv'\in\Ve$ such that
\bes
\bv(t)-\bv(t_0)=(t-t_0)\bv'(t_0)+o(t-t_0).
\ees
We call $\bv'(t_0)$ the {\it derivative} \footnote{The derivative is also written as $\dfrac{d\bv}{dt},\ \bv_{,t}$ or also as $\dot{\bv}$, with the last symbol  usually reserved, in physics, to the case where $t$ is the time. For the sake of brevity, we omit to indicate the derivative of $\bv$ at $t_0$ as $\bv'(t_0)$, writing simply $\bv'$.} $of\ \bv$ at $t_0$. Applying  the definition of derivative to $\bv'$ we define the {\it second derivative} $\bv''$ of $\bv$ and recursively all the derivatives of higher orders. We say that $\bv$ is of {\it class} C$^n$ if it is continuous with its derivatives up to the order $n$; if $n\geq1$, $\bv$ is said to be {\it smooth}. A curve $\bv(t)$ of class C$^n$ is said to be {\it regular} if $\bv'\neq\bo\ \forall t$. 
Similar definitions can be given for curves in $\Eu, Lin(\Ve)$ and $\Lq$in$(\Ve)$, thus defining derivatives of points and tensors. We remark that the derivative of a curve in $\Eu$, defined as a difference of points, is a curve in $\Ve$ (we say, in short,  that {\it the derivative of a point is a vector}). For what concerns tensors, the derivative of a tensor of rank $r$ is a tensor of the same rank.

Let $\bu,\bv$ be curves in $\Ve$, $\L,\M$ curves in $Lin(\Ve)$, $\Lq,\Mq$ curves in $\Lq$in$(\Ve)$ and $\alpha$ a scalar function, all of them defined  and at least of class C$^1$ on $[a,b]$. The same definition of derivative of a curve gives the following results, whose proof is let to the reader:\\
\bes
\begin{array}{c}
(\bu+\bv)'={\bu}'+{\bv}',\medskip\\
(\alpha\bv)'={\alpha}'\bv+\alpha{\bv}',\medskip\\
(\bu\cdot\bv)'={\bu}'\cdot\bv+\bu\cdot{\bv}',\medskip\\
(\bu\times\bv)'={\bu}'\times\bv+\bu\times{\bv}',\medskip\\
(\bu\otimes\bv)'={\bu}'\otimes\bv+\bu\otimes{\bv}',\medskip\\
(\L+\M)'={\L}'+{\M}',\medskip\\
(\alpha\L)'={\alpha}'\L+\alpha{\L}',\medskip\\
(\L\bv)'={\L}'\bv+\L{\bv}',\medskip\\
(\L\M)'={\L}'\M+\L{\M}',\medskip\\
(\L\cdot\M)'={\L}'\cdot\M+\L\cdot{\M}',\medskip\\
\end{array}
\ees

\bes
\begin{array}{c}
(\L\otimes\M)'={\L}'\otimes\M+\L\otimes{\M}',\medskip\\
(\L\boxtimes\M)'={\L}'\boxtimes\M+\L\boxtimes{\M}',\medskip\\
(\Lq+\Mq)'={\Lq}'+{\Mq}',\medskip\\
(\alpha\Lq)'={\alpha}'\Lq+\alpha{\Lq}',\medskip\\
(\Lq\L)'={\Lq}'\L+\Lq{\L}',\medskip\\
(\Lq\Mq)'={\Lq}'\Mq+\Lq{\Mq}',\medskip\\
(\Lq\cdot\Mq)'={\Lq}'\cdot\Mq+\Lq\cdot{\Mq}'.
\end{array}
\ees
We remark that the derivative of any kind of product is made according to the usual rule of the derivative of a product of functions.

\begin{figure}[h]
\begin{center}
\includegraphics[scale=1]{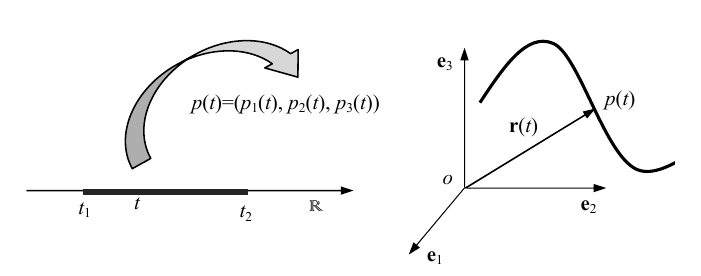}
\caption{Mapping of a curve of points.}
\label{fig:f11}
\end{center}
\end{figure}

Let $\mathcal{R}=\{o; \Ba\}$ be a reference frame of the euclidean space $\mathcal{E}$, composed of an {\it origin o} and a basis $\Ba=\{\e_1,\e_2,\e_3\}$ of $\Ve, \e_i\cdot\e_j=\delta_{ij}\ \forall i,j=1,2,3$  and let us consider a point $p(t)=(p_1(t),p_2(t),p_3(t))$. If the three {\it coordinates} $p_i(t)$ are three continuous functions  over the interval $[t_1,t_2]\in\R$, then, by the definition given above, the mapping $p(t):[t_1,t_2]\rightarrow\mathcal{E}$ is a  curve in $\mathcal{E}$ and the
 equation
\bes
p(t)=(p_1(t),p_2(t),p_3(t))\ \rightarrow\ \left\{
\begin{array}{l}
p_1=p_1(t)\\
p_2=p_2(t)\\
p_3=p_3(t)
\end{array}
\right.
\ees
is the {\it parametric point equation of the curve}: To each value of $t\in[t_1,t_2]$ it corresponds a point of the curve in $\mathcal{E}$, see Fig. \ref{fig:f11}.

The vector function $\gr{r}(t)=p(t)-o$ is the {\it position vector} of point $p$ in $\mathcal{R}$; the equation
\bes
\gr{r}(t)=r_i(t)\e_i=r_1(t)\gr{e}_1+r_2(t)\gr{e}_2+r_3(t)\gr{e}_3\ \rightarrow\ \left\{
\begin{array}{l}
r_1=r_1(t)\\
r_2=r_2(t)\\
r_3=r_3(t)
\end{array}
\right.
\ees
is the {\it parametric vector equation} of the curve: To each value of $t\in[t_1,t_2]$ there corresponds a vector of $\mathcal{V}$ that determines a point of the curve in $\mathcal{E}$ through the operation $p(t)=o+\gr{r}(t)$.

 Similarly,  if the components $L_{ij}(t)$ are continuous functions of a parameter $t$, the mapping $\gr{L}(t):[t_1,t_2]\rightarrow Lin(\mathcal{V})$ defined by
\bes
\gr{L}(t)=L_{ij}(t)\gr{e}_i\otimes\gr{e}_j,\ \ i,j=1,2,3,
\ees
is a  curve of tensors. In the same way we can give a curve of fourth-rank tensors $\Lq(t):[t_1,t_2]\rightarrow \Lq$in$(\Ve)$ by
\bes
\Lq(t)=L_{ijkl}(t)\e_i\otimes\e_j\otimes\e_k\otimes\e_l,\ \ i,j,k,l=1,2,3.
\ees

It is worth noting that the choice of the parameter is {\it not unique}: The equation $p=p[\tau(t)]$ still represents the same curve $p=p(t)$, through the {\it change of parameter} $\tau=\tau(t)$.

\begin{figure}[h]
\begin{center}
\includegraphics[scale=.9]{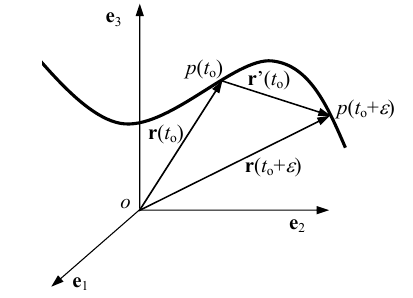}
\caption{Derivative of a curve.}
\label{fig:f12}
\end{center}
\end{figure}

The definition given above for the derivative  of a curve of points $p=p(t)$ in $t=t_0$ is equivalent to the following one\footnote{This is true also for the derivatives of vector or tensor curves.} (probably more familiar to the reader)
\bes
\frac{dp(t)}{dt}=\lim_{\eps\rightarrow0}\frac{p(t_0+\eps)-p(t_0)}{\eps},
\ees
represented in Fig. \ref{fig:f12}, where it is apparent that  $\r'(t)=\dfrac{dp(t)}{dt}$ is a vector.

An important case is that of a vector $\gr{v}(t)$ whose norm $v(t)$ is constant $\forall t$:
\be
\label{eq:constvect}
(v^2)'=(\gr{v}\cdot\gr{v})'=\gr{v}'\cdot\gr{v}+\gr{v}\cdot\gr{v}'=2\gr{v}'\cdot\gr{v}=0:
\ee
the derivative of such a vector is  orthogonal to it $\forall t$. The contrary is also true, as is immediately apparent.

Finally, using the above rules and assuming that the reference frame $\mathcal{R}$ is independent of $t$, we get easily that
\be
\label{eq:derivcomp}
\begin{split}
&p'(t)=p'_i(t)\ \gr{e}_i,\\
&\gr{v}'(t)=v'_i(t)\ \gr{e}_i,\\
&\gr{L}'(t)=L'_{ij}(t)\ \gr{e}_i\otimes\gr{e}_j,\\
&\Lq'(t)=L'_{ijkl}(t)\ \e_i\otimes\e_j\otimes\e_k\otimes\e_l,
\end{split}
\ee
i.e., the derivative of a curve of points, vectors or tensors is simply calculated differentiating the coordinates of the components. Using this result, it is immediate to prove that
\bes
\begin{array}{c}
(\L^\top)'={\L}'^\top,\medskip\\
(\Lq^\top)'={\Lq}'^\top,
\end{array}
\ees
while for any invertible tensor $\L$ it is (we state the following results without proof\footnote{The interested reader can find these proofs in the text by Gurtin, see the suggested texts.}.)
\bes
\begin{array}{c}
(\L^{-1})'=-\L^{-1}\L'\L^{-1},\medskip\\
(\det\L)'=\det\L\ \tr(\L'\L^{-1})=\det\L\ {\L^\top}'\cdot\L^{-1}=\det\L\ \L'\cdot \L^{-\top}.
\end{array}
\ees

Let $\Q(t):\R\rightarrow Orth(\Ve)^+$ a differentiable function. We call {\it spin tensor} the tensor $\Sy(t)$ defined as
\bes
\Sy(t):=\Q'(t)\Q^\top(t).
\ees
Then, we have the following\footnote{The spin tensor and the following result are of importance in kinematics: If $t$ is time and $\Q(t)\in Orth{\Ve}^+$, then the axial vector of $\Sy(t)$ is $\bom(t)$, the {\it angular velocity}.}
\begin{teo} \bf{(Characterization of the spin tensor).}
$\Sy(t)\in Skw(\Ve)\ \forall t\in\R$.
\begin{proof}
As $\Q(t)\in Orth(\Ve)^+\  \forall t$, then
\bes
 \Q\Q^\top=\I\Rightarrow(\Q\Q^\top)'=\Q'\Q^\top+\Q{\Q^\top}'=\I'=\O\\ \Rightarrow\Q{\Q^\top}'=-\Q'\Q^\top
 \ees
 so
 \bes
 \Sy^\top=(\Q'\Q^\top)^\top=\Q{\Q^\top}'=-\Q'\Q^\top=-\Sy.
 \ees
\end{proof}
\end{teo}

\section{Integral of a curve of vectors and  length of a curve}
\begin{figure}[h]
\begin{center}
\includegraphics[scale=.9]{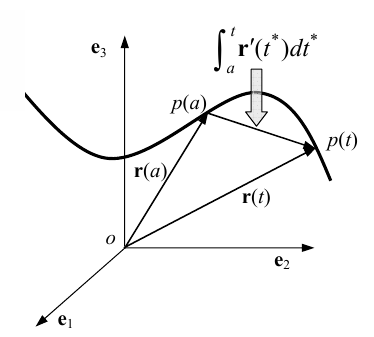}
\caption{Integral of a vector curve.}
\label{fig:f13}
\end{center}
\end{figure}

We define {\it integral of a curve of vectors} $\r(t)$ between $a$ and $b \in[t_1,t_2]$ the curve that is obtained by integrating each component of the curve:
\bes
\begin{split}
&\int_a^b\gr{r}(t)\ dt=\int_a^br_i(t)\ dt\ \gr{e}_i.
\end{split}
\ees
If the curve is regular, we can generalize the second fundamental theorem of the integral calculus
\bes
\gr{r}(t)=\gr{r}(a)+\int_a^t\gr{r}'(t^*)\ dt^*.
\ees
Because
\bes
\gr{r}(t)=p(t)-o,\ \ \ \gr{r}'(t)=(p(t)-o)'=p'(t),
\ees
we  also get
\bes
p(t)=p(a)+\int_a^tp'(t^*)\ dt^*.
\ees
The integral of a vector function is the generalization of the vector sum, see Fig. \ref{fig:f13}.

Let $\r(t):[a,b]\rightarrow\Eu$ be a regular curve, $\sigma$ a partition of $[a,b]$ of the type $a=t_0<t_1<...<t_n=b$, and 
\bes
\sigma_{max}={\max_{i=1,...,n}}|t_i-t_{i-1}|.
\ees
The length $\ell_\sigma$ of the polygonal line whose vertices are the points $\r(t_i)$  is hence:
\bes
\ell_\sigma=\sum_{i=1}^n{|\r(t_i)-\r(t_{i-1})|}.
\ees
We define {\it length of the curve} $\r(t)$ the (positive) number
\bes
\ell:=\sup_{\sigma}\ell_\sigma.
\ees
\begin{teo} Let $\gr{r}(t):[a,b]\Rightarrow\Eu$ be a regular curve, then
\bes
\ell=\int_a^b|\gr{r}'(t)|dt.
\ees
\begin{proof}
By the fundamental theorem of calculus,
\bes
\r(t_i)-\r(t_{i-1})=\int_{t_{i-1}}^{t_i}\r'(t)dt,
\ees
so that, using  Minkowski's inequality, 
\bes
|\r(t_i)-\r(t_{i-1})|=\left|\int_{t_{i-1}}^{t_i}\r'(t)dt\right|\leq\int_{t_{i-1}}^{t_i}|\r'(t)|dt,
\ees
whence
\be
\label{eq:longueur1}
\ell\leq\int_a^b|\r'(t)|dt.
\ee
Because $\r'(t)$ is continuous on $[a,b]$, $\forall \eps>0\ \exists\delta>0$ such that $|t-\overline{t}|<\delta\Rightarrow|\r'({t})-\r'(\overline{t})|<\eps$. Let $t\in[t_{i-1},t_i]$ and $\sigma_{max}<\delta$, which is always possible by the choice of the partition $\sigma$; again by the Minkowski's inequality,
\bes
|\r'(t)|\leq|\r'(t)-\r'(t_i)|+|\r'(t_i)|<\eps+|\r'(t_i)|,
\ees
whence
\bes
\besp
\int_{t_{i-1}}^{t_i}|\r'(t)|dt&<\int_{t_{i-1}}^{t_i}|\r'(t_i)|dt+\eps(t_i-t_{i-1})=\left|\int_{t_{i-1}}^{t_i}\r'(t_i)dt\right|+\eps(t_i-t_{i-1})\\
&\leq\left|\int_{t_{i-1}}^{t_i}\r'(t)dt\right|+\left|\int_{t_{i-1}}^{t_i}(\r'(t_i)-\r'(t))dt\right|+\eps(t_i-t_{i-1})\\
&\leq|\r(t_i)-\r(t_{i-1})|+2\eps(t_i-t_{i-1}).
\end{split}
\ees
Summing up over all the intervals $[t_{i-1},t_i]$ we get
\bes
\int_a^b|\r'(t)|dt\leq\ell_\sigma+2\eps(b-a)\leq\ell+2\eps(b-a),
\ees
and because $\eps$ is arbitrary, 
\bes
\int_a^b|\r'(t)|dt\leq\ell,
\ees
which by Eq. (\ref{eq:longueur1}) implies the thesis.
\end{proof}
\end{teo}
Let $t=f(\tau):[c,d]\rightarrow[a,b]$ be a bijective function that operates the change of parameter from $t$ to $\tau$. If $\r_t(t):[a,b]\rightarrow\Ve$ is a parametric equation of a curve, $\r_\tau:[c,d]\rightarrow\Ve$ is a {\it re-parameterization} of the same curve. We then have the following
\begin{teo}
The length of a curve does not depend upon its parameterization.
\begin{proof}
Let $\r_t(t):[a,b]\rightarrow\Eu$ be a regular curve and $t=f(\tau):[c,d]\rightarrow[a,b]$ be a change of parameter; then $dt=f'(\tau)d\tau$ and
\bes
\ell=\int_a^b|\r_t'(t)|dt=\int_c^d|\r_t'(f(\tau))f'(\tau)|d\tau=\int_c^d|\r_\tau'(\tau)|d\tau.
\ees
\end{proof}
\end{teo}
A simple way to determine a point $p(t)$ on a curve is to fix a point $p_0$ on the curve and to measure the length $s(t)$ of the arc of curve between $p_0=p(t=0)$ and $p(t)$. This length $s(t)$ is called  {\it curvilinear abscissa}\footnote{The curvilinear abscissa is also called {\it arc-length} or {\it natural parameter}.}:
\be
\label{eq:esse}
s(t)=\int_{0}^t|\gr{r}'(t^*)|dt^*=\int_{0}^t|(p(t^*)-o)'|dt^*.
\ee
From Eq. (\ref{eq:esse}) we get
\bes
\frac{ds}{dt}=|\gr{r}'(t)|>0,
\ees
so that $s(t)$ is an increasing function of $t$ and the length of an infinitesimal arc is
\bes
ds=\sqrt{dr_1^2+dr_2^2+dr_3^2}.
\ees
For a plane curve $y=f(x)$, we can always put $t=x$, which gives the parametric equation
\bes
p(t)=(t,f(t)),
\ees
or in vector form
\bes
\gr{r}(t)=t\ \gr{e}_1+f(t)\ \gr{e}_2,
\ees
from which we obtain
\be
\label{eq:ds}
\frac{ds}{dt}=|\gr{r}'(t)|=|p'(t)|=\sqrt{1+f'^2(t)},
\ee
which gives the length of a plane curve between $t=x_0$ and $t=x$ as a function of the abscissa $x$:
\bes
s(x)=\int_{x_0}^x\sqrt{1+f'^2(t)}dt.
\ees

\section{The Frenet-Serret basis}
We define the {\it tangent vector} $\btau(t)$ to a regular curve $p=p(t)$ as the vector
\bes
\label{eq:tau}
\btau(t):=\dfrac{p'(t)}{|p'(t)|}.
\ees
By the definition of the derivative, this unit vector is always oriented as the increasing values of $t$; hence, the straight line tangent to the curve in $p_0=p(t_0)$ has the equation
\bes
q(\bar{t})=p(t_0)+\bar{t}\ \btau(t_0).
\ees
If the curvilinear abscissa $s$ is chosen as parameter for the curve, through the change of parameter $s=s(t)$ we get
\be
\label{eq:tangentparams}
\btau(t)=\dfrac{p'(t)}{|p'(t)|}=\dfrac{p'[s(t)]}{|p'[s(t)]|}=\frac{1}{s'(t)}\frac{dp(s)}{ds}\frac{ds(t)}{dt}=\frac{dp(s)}{ds}\ \rightarrow \btau(s)=p'(s).
\ee
So, if the parameter of the curve is $s$, the derivative of the curve is  $\btau$, i.e. it is automatically a unit vector. The above equation, in addition, shows that the change of parameter does not change the direction of the tangent, because it is only a scalar, the derivative of the parameter's change, that multiplies the vector. Nevertheless, generally speaking, a change of parameter can change the orientation of the curve.

Because the norm of $\btau$ is constant, its derivative is a vector orthogonal to $\btau$, see Eq. (\ref{eq:constvect}). That is why we call {\it principal normal vector} to a curve the unit vector
\be
\label{eq:normal}
\bnu(t):=\frac{\btau'(t)}{|\btau'(t)|}.
\ee
$\bnu$ is defined only on the points of the curve where $\btau'\neq\gr{o}$, which implies that $\bnu$ is not defined on the points of a straight line. This simply means that there is not, among the infinite unit normal vectors to a straight line, a normal with special properties, a {\it principal} one,  linked to $\btau$ in a unique way.

Unlike $\btau$, whose orientation changes with the choice of the parameter, $\bnu$ is an {\it intrinsic} local characteristic of the curve: {\it It is not affected by the choice of the parameter}. In fact, by its same definition, $\bnu$ does not depend upon the reference frame; then, because the direction of $\btau$ is also independent upon the parameter's choice,  the only factor that could affect $\bnu$ is the orientation of the curve, which depends upon the parameter. But a change in the orientation affects, in (\ref{eq:normal}), both $\btau$ and the sign of the increment $dt$, so that $\btau'(t)=d\btau/dt$ does not change, nor does $\bnu$, which is hence an intrinsic property of the curve.

\begin{figure}[h]
\begin{center}
\includegraphics[scale=.8]{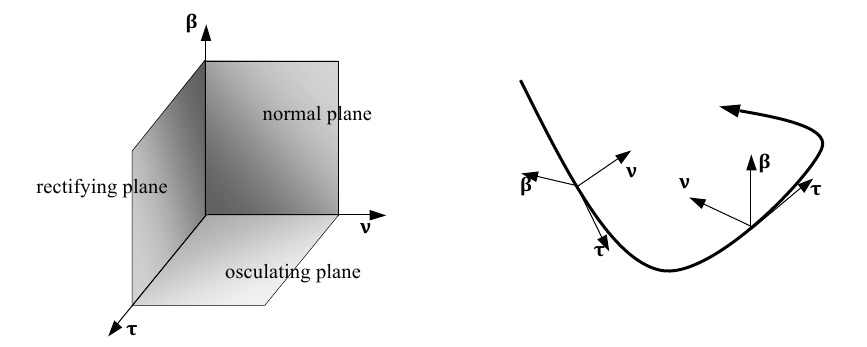}
\caption{The Frenet-Serret basis.}
\label{fig:f14}
\end{center}
\end{figure}

The vector 
\bes
\bb(t):=\btau(t)\times\bnu(t)
\ees
is called the {\it binormal vector}; by construction, it is orthogonal to $\btau$ and $\bnu$ and it is a unit vector. In addition, it is evident that
\bes
\btau\times\bnu\cdot\bb=1,
\ees
so the set $\{\btau,\bnu,\bb\}$ forms a positively oriented othonormal basis that can be defined at any regular point of a curve with $\btau'\neq\gr{o}$. Such a basis is called the {\it Frenet-Serret local basis}, local in the sense that it changes with the position along the curve. The plane $(\btau,\bnu)$ is the {\it osculating plane}, the plane $(\bnu,\bb)$ the {\it normal plane} and the plane $(\bb,\btau)$ the {\it rectifying plane}, see Fig. \ref{fig:f14}.
The osculating  plane is particularly important: If we consider a plane passing through three nonaligned points of the curve, when these points become closer and closer, still remaining on the curve, the plane tends to the osculating plane: The osculating plane at a point of a curve is hence the plane that better approaches the curve near the point.
A plane curve is entirely contained in the osculating plane, which is fixed. 

The principal normal $\bnu$ is always oriented toward the part  of the space, with respect to the rectifying plane, where the curve is; in particular, for a plane curve, $\bnu$ is always directed toward the concavity of the curve.
To show that, it is sufficient to prove that the vector $p(t+\eps)-p(t)$ forms with $\bnu$ an angle $\psi\leq\pi/2$, i.e. that $(p(t+\eps)-p(t))\cdot\bnu\geq0$. In fact,
\bes
\begin{split}
&p(t+\eps)-p(t)=\eps\ p'(t)+\frac{1}{2}\eps^2p''(t)+o(\eps^2)\ \Rightarrow\\
&(p(t+\eps)-p(t))\cdot\bnu=\frac{1}{2}\eps^2p''(t)\cdot\bnu+o(\eps^2),
\end{split}
\ees
but
\bes
p''(t)\cdot\bnu=(\btau'|p'|+\btau|p'|')\cdot\bnu=(|\btau'||p'|\bnu+\btau|p'|')\cdot\bnu=|\btau'||p'|,
\ees
so that, to within infinitesimal quantities of order $o(\eps^2)$, we obtain
\bes
(p(t+\eps)-p(t))\cdot\bnu=\frac{1}{2}\eps^2|\btau'||p'|\geq0.
\ees

\section{Curvature of a curve}
\begin{figure}[h]
\begin{center}
\includegraphics[scale=1]{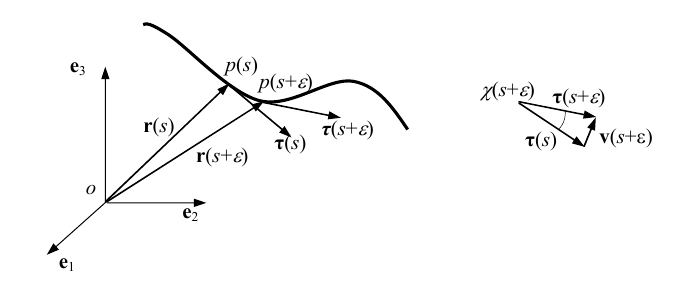}
\caption{Curvature of a curve.}
\label{fig:f15}
\end{center}
\end{figure}

It is important, in several situations, to evaluate how much a curve moves away from a straight line, in the neighborhood of a point. To do that, we calculate the angle formed by the tangents at two close points, determined by the curvilinear abscissae $s$ and $s+\eps$, and we measure the angle $\chi(s,\eps)$ that they form, see Fig. \ref{fig:f15}.

We then define {\it curvature} of the curve in $p=p(s)$ as the limit
\bes
c(s)=\lim_{\eps\rightarrow0}\left|\frac{\chi(s,\eps)}{\eps}\right|.
\ees
The curvature is hence a non-negative scalar that measures the rapidity of variation in the direction of the curve per unit length of the curve (that is why $c(s)$ is defined as a function of the curvilinear abscissa); by its same definition, the curvature is an {\it intrinsic property} of the curve, i.e. independent of the parameter's choice. For a straight line, the curvature is everywhere identically null.

The curvature is linked to the second derivative of the curve; referring to Fig. \ref{fig:f15}, it is
\bes
\begin{split}
c(s)&=\lim_{\eps\rightarrow0}\left|\frac{\chi(s,\eps)}{\eps}\right|=\lim_{\eps\rightarrow0}\left|\frac{\sin\chi(s,\eps)}{\eps}\right|=\lim_{\eps\rightarrow0}\left|\frac{2}{\eps}\sin\frac{\chi(s,\eps)}{2}\right|\\
&=\lim_{\eps\rightarrow0}\left|\frac{\gr{v}(s,\eps)}{\eps}\right|=\lim_{\eps\rightarrow0}\left|\frac{\btau(s+\eps)-\btau(s)}{\eps}\right|=|\btau'(s)|=|p''(s)|.
\end{split}
\ees

Another formula for the calculation of $c(s)$ can be obtained if we consider that
\bes
\frac{d\btau[s(t)]}{dt}=\frac{d\btau}{ds}\frac{ds}{dt}=\frac{d\btau}{ds}|p'(t)|\ \rightarrow\ \frac{d\btau}{ds}=\frac{1}{|p'(t)|}\frac{d\btau}{dt},
\ees
so that
\be
\label{eq:curv2}
c(s)=|\btau'(s)|=\frac{1}{|p'(t)|}\left|\frac{d\btau}{dt}\right|=\frac{|\btau'(t)|}{|p'(t)|}.
\ee

A better formula can be obtained using the {\it complementary projector} onto $\btau$, i.e. the tensor $\I-\btau\otimes\btau$, introduced in Exercise 2.2:
\bes
\begin{split}
\frac{d\btau}{ds}&=\frac{1}{|p'(t)|}\frac{d\btau}{dt}=\frac{1}{|p'(t)|}\frac{d}{dt}\frac{p'(t)}{|p'(t)|}=
\frac{1}{|p'|}\frac{p''|p'|-p'\dfrac{p''\cdot p'}{|p'|}}{|p'|^2}\\
&=\frac{p''-\btau\ p''\cdot\btau}{|p'|^2}=(\gr{I}-\btau\otimes\btau)\frac{p''}{|p'|^2}.
\end{split}
\ees
Consequently,
\bes
c(s)=\left|\frac{d\btau(s)}{ds}\right|=\frac{1}{|p'|^2}|(\gr{I}-\btau\otimes\btau)p''|.
\ees
Now, we use Eq. (\ref{eq:normW}) with $\bw=\btau$; denoting by $\W_\tau$ the axial tensor of $\btau$, then
\bes
\W_\tau\W_\tau=-\frac{1}{2}|\W_\tau|^2(\gr{I}-\btau\otimes\btau),
\ees
whence
\bes
\gr{I}-\btau\otimes\btau=-2\frac{\W_\tau\W_\tau}{|\W_\tau|^2}=-\W_\tau\W_\tau,
\ees
because if $\btau=(\tau_1,\tau_2,\tau_3)$, then
\bes
\begin{split}
|\W_\tau|^2&=\W_\tau\cdot\W_\tau=\left[\begin{array}{ccc}0 & -\tau_3 & \tau_2 \\\tau_3& 0 & -\tau_1 \\-\tau_2 & \tau_1 & 0\end{array}\right]\cdot\left[\begin{array}{ccc}0 & -\tau_3 & \tau_2 \\\tau_3& 0 & -\tau_1 \\-\tau_2 & \tau_1 & 0\end{array}\right]\\
&=2(\tau_1^2+\tau_2^2+\tau_3^2)=2.
\end{split}
\ees
So, because $\W_\tau\in Skw(\Ve)$,
\bes
\W_\tau\ \gr{u}=\btau\times\gr{u}\ \ \forall \gr{u}\in\mathcal{V}.
\ees
Finally, using Eq. (\ref{eq:doublecrossprod}), the orthogonality property of cross product, Eq. (\ref{eq:orthprop}) and  Eq. (\ref{eq:formulevectprod}), we get
\bes
\begin{split}
|(\gr{I}-\btau\otimes\btau)p''|&=|-\W_\tau\W_\tau p''|=|-\W_\tau(\btau\times p'')|=|-\btau\times(\btau\times p'')|\\
&=|\btau\times(\btau\times p'')|=|\btau\times p''|=\frac{|p'\times p''|}{|p'|},
\end{split}
\ees
so that, finally,
\be
\label{eq:curv1}
c=\frac{|p'\times p''|}{|p'|^3}.
\ee

Applying this last formula to a plane curve $p(t)=(x(t),y(t))$, we get
\bes
c=\frac{|x'y''-x''y'|}{(x'^2+y'^2)^\frac{3}{2}},
\ees
and if the curve is given in the form $y=y(x)$, so that the parameter $t=x$, then we obtain
\bes
c=\frac{|y''|}{(1+y'^2)^\frac{3}{2}}.
\ees

This last formula shows that if $|y'|\ll1$, then
\bes
c\simeq|y''|.
\ees
This result is fundamental in the linearized (infinitesimal) theory of beams, plates and shells. 

\section{The Frenet-Serret formulae}
From Eq. (\ref{eq:normal}) for $t=s$ and Eq. (\ref{eq:curv2}), we get
\be
\label{eq:fs1}
\frac{d\btau}{ds}=c\ \bnu
\ee
which is the {\it first Frenet-Serret Formula}, giving the variation in $\btau$ per unit length of the curve. Such a variation is a vector whose norm is the curvature and the  direction is that of $\bnu$. We remark that, because $t=s$, by Eq. (\ref{eq:tangentparams}) it is also
\be
\label{eq:frenetserret1avecs}
p''(s)=c(s)\bnu(s).
\ee

Let us now consider the variation in $\bb$ per unit length of the curve; because $\bb$ is a unit vector, we have
\bes
\frac{d\bb}{ds}\cdot\bb=0,
\ees
and 
\bes
\bb\cdot\btau=0\ \Rightarrow\ \frac{d(\bb\cdot\btau)}{ds}=\frac{d\bb}{ds}\cdot\btau+\bb\cdot\frac{d\btau}{ds}=0.
\ees
Through Eq. (\ref{eq:fs1}) and because $\bb\cdot\bnu=0$ we get
\bes
\frac{d\bb}{ds}\cdot\btau=-c\ \bb\cdot\bnu=0,
\ees
so that $\dfrac{d\bb}{ds}$ is necessarily parallel to $\bnu$. We then set
\bes
\label{eq:fs2}
\frac{d\bb}{ds}=\vartheta\bnu,
\ees
which is the {\it second Frenet-Serret formula}. The scalar $\vartheta(s)$ is called the {\it torsion of the curve} in $p=p(s)$. So, we see that the variation in $\bb$ per unit length is a vector parallel to $\bnu$ and proportional to the torsion of the curve.

We can now find the variation in $\bnu$ per unit length of the curve:
\bes
\frac{d\bnu}{ds}=\frac{d(\bb\times\btau)}{ds}=\frac{d\bb}{ds}\times\btau+\bb\times\frac{d\btau}{ds}=\vartheta\ \bnu\times\btau+c\ \bb\times\bnu,
\ees
so finally
\bes
\label{eq:fs3}
\frac{d\bnu}{ds}=-c\ \btau-\vartheta\ \bb,
\ees
which is the {\it third Frenet-Serret formula}: The variation in $\bnu$ per unit length of the curve is a vector of the rectifying plane.

The three formulae of Frenet-Serret (discovered independently by J. F. Frenet in 1847 and by J. A. Serret in 1851) can be condensed in the symbolic matrix product
\bes
\left\{\begin{array}{c}\btau' \\\bnu' \\\bb'\end{array}\right\}=
\left[\begin{array}{ccc}0 & c & 0 \\-c & 0 & -\vartheta \\0 & \vartheta & 0\end{array}\right]
\left\{\begin{array}{c}\btau \\\bnu \\\bb\end{array}\right\}.
\ees
The matrix in the equation above is called the {\it matrix of Cartan}, and it is skew. 

\section{The torsion of a curve}
We have introduced the torsion of a curve in the previous section, with the second formula of Frenet-Serret.
The torsion measures the deviation of a curve from flatness: If a curve is planar, it belongs to the osculating plane and $\bb$, which is perpendicular to the osculating pane, is hence a constant vector. So, its derivative is  null, and by the Frenet-Serret second formula, $\vartheta=0$.

Conversely, if $\vartheta=0$ everywhere, $\bb$ is a constant vector and hence the osculating plane does not change and the curve is planar. So we have that {\it a curve is planar if and only if the torsion is null $\forall p(s)$}.

Using the Frenet-Serret formulae in the expression of $p'''(s)$, we get a formula for the torsion:
\bes
\begin{split}
&p'(t)=|p'|\btau=\frac{dp}{ds}\frac{ds}{dt}=s'\btau\ \Rightarrow\ |p'|=s'\ \rightarrow\\
&p''(t)=s''\btau+s'\btau'=s''\btau+s'^2\frac{d\btau}{ds}=s''\btau+c\ s'^2\bnu\ \rightarrow\\
&p'''(t)=s'''\btau+s''\btau'+(c\ s'^2)'\bnu+c\ s'^2\bnu'\\
&\hspace{9.5mm}=s'''\btau+s''s'\frac{d\btau}{ds}+(c\ s'^2)'\bnu+c\ s'^3\frac{d\bnu}{ds}\\
&\hspace{9.5mm}=s'''\btau+s''s'c\bnu+(c\ s'^2)'\bnu-c\ s'^3(c\btau+\vartheta\bb)\\
&\hspace{9.5mm}=(s'''-c^2s'^3)\btau+(s''s'c+c's'^2+2c\ s's'')\bnu-c\ s'^3\vartheta\bb,
\end{split}
\ees
whence, through Eq. (\ref{eq:curv1}), 
\bes
\begin{split}
p'\times p''\cdot p'''&=s'\btau\times(s''\btau+c\ s'^2\bnu)\cdot[(s'''-c^2s'^3)\btau\\
&+(s''s'c+c's'^2+2c\ s's'')\bnu-c\ s'^3\vartheta\bb]\\
&=-c^2s'^6\vartheta=-c^2|p'|^6\vartheta=-\frac{|p'\times p''|^2}{|p'|^6}|p'|^6\vartheta,
\end{split}
\ees
so that, finally, 
\bes
\vartheta=-\frac{p'\times p''\cdot p'''}{|p'\times p''|^2}.
\ees
We remark that, while the curvature is linked to the second derivative of the curve, the torsion is also a function  of the third derivative.

Unlike curvature, which is intrinsically positive, the torsion can be negative. 
In fact, again using the Frenet-Serret formulae,
\bes
\begin{split}
p(s+\eps)-p(s)&=\eps\ p'+\frac{1}{2}\eps^2p''+\frac{1}{6}\eps^3p'''+o(\eps^3)\\
&=\eps\btau+\frac{1}{2}\eps^2c\bnu+\frac{1}{6}\eps^3(c\bnu)'+o(\eps^3)\\
&=\eps\btau+\frac{1}{2}\eps^2c\bnu+\frac{1}{6}\eps^3(c'\bnu-c^2\btau-c\ \vartheta\bb)+o(\eps^3)\\
&\Rightarrow(p(s+\eps)-p(s))\cdot\bb=-\frac{1}{6}\eps^3c\ \vartheta+o(\eps^3).
\end{split}
\ees

\begin{figure}[h]
\begin{center}
\includegraphics[scale=.9]{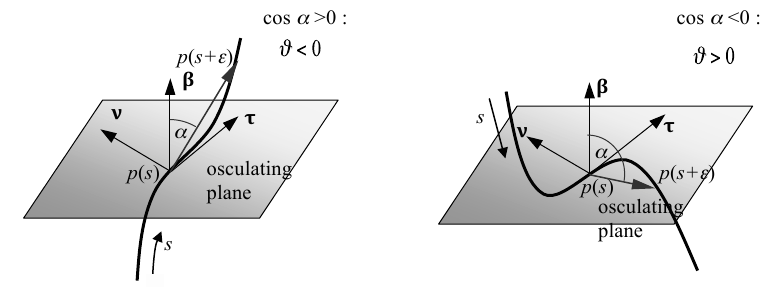}
\caption{Torsion of a curve.}
\label{fig:f16}
\end{center}
\end{figure}

The above dot product determines whether the point $p(s+\eps)$ is located, with respect to the osculating plane, on the side of $\bb$ or on the opposite one, see Fig. \ref{fig:f16}: If following the curve for increasing values of $s,\ \eps>0$, the point passes into the semi-space of $\bb$ from the opposite one, because $1/6\ c\ \eps^3>0$, it will be $\vartheta<0$, while in the opposite case it will be $\vartheta>0$.

This result is intrinsic, i.e. it does not depend upon the choice of the parameter, hence of the positive orientation of the curve; in fact, $\bnu$ is intrinsic, but changing the orientation of the curve, $\btau$, and hence $\bb$, change in orientation.

 \section{Osculating sphere and circle}
 The {\it osculating sphere}\footnote{The word osculating comes from the latin word {\it osculo}, which means to kiss; an osculating sphere or circle or plane is a geometric object that is very close to the curve, as close as two  lovers are in a kiss.} to a curve at a point $p$ is a sphere to which the curve tends to adhere in the neighborhood of $p$. Mathematically, if $q_s$ is the center of the sphere relative to the point $p(s)$, then
 \bes
 |p(s+\eps)-q_s|^2= |p(s)-q_s|^2+o(\eps^3).
 \ees
 Using this definition, discarding the terms of order $o(\eps^3)$ and using the Frenet-Serret formulae, we get:
 \bes
 \begin{split}
  |p(s+\eps)-q_s|^2&=|p(s)-q_s+\eps p'+\frac{1}{2}\eps^2p''+\frac{1}{6}\eps^3p'''+o(\eps^3)|^2\\
  &=|p(s)-q_s+\eps \btau+\frac{1}{2}\eps^2c\ \bnu+\frac{1}{6}\eps^3(c\bnu)'+o(\eps^3)|^2\\
  &=|p(s)-q_s|^2+2\eps(p(s)-q_s)\cdot\btau+\eps^2+\eps^2c(p(s)-q_s)\cdot\bnu\\
  &+\frac{1}{3}\eps^3(p(s)-q_s)\cdot(c'\bnu-c^2\btau-c\ \vartheta\bb)+o(\eps^3),
 \end{split}
 \ees
 which gives
 \bes
 \begin{split}
 &(p(s)-q_s)\cdot\btau=0,\\
  &(p(s)-q_s)\cdot\bnu=-\frac{1}{c}=-\rho,\\
   &(p(s)-q_s)\cdot\bb=-\frac{c'}{c^2\vartheta}=\frac{\rho'}{\vartheta},
 \end{split}
 \ees
 and finally
 \be
 \label{eq:spherecenter}
 q_s=p+\rho\ \bnu-\frac{\rho'}{\vartheta}\bb,
 \ee
so the center of the sphere belongs to the normal plane; the sphere is not defined for a plane curve. The quantity $\rho$ is the {\it radius of curvature} of the curve, defined as
\bes
\rho=\frac{1}{c}.
\ees

The radius of the osculating sphere is
\bes
\rho_s=|p-q_s|=\sqrt{\rho^2+\left(\frac{\rho'}{\vartheta}\right)^2}.
\ees

\begin{figure}[h]
\begin{center}
\includegraphics[scale=.5]{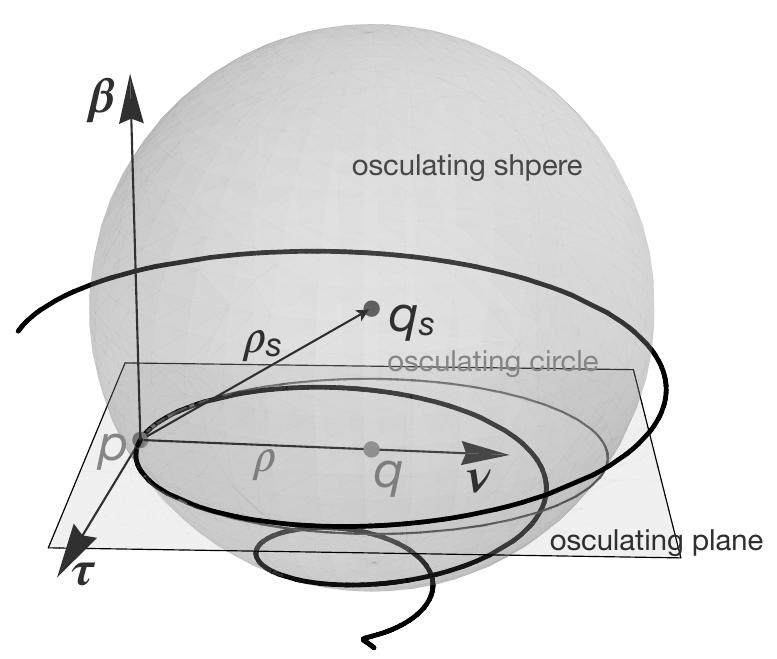}
\caption{Osculating plane, circle and sphere for a point $p$ of a conical helix.}
\label{fig:f17}
\end{center}
\end{figure}

The intersection between the osculating sphere and the osculating plane at a same point $p$ is the {\it osculating circle}. This circle has the property of sharing the same tangent in $p$ with the curve and its radius is the radius of curvature, $\rho$. From Eq. (\ref{eq:spherecenter}) we get the position of the osculating circle center $q$:
\be
\label{eq:cco}
q=p+\rho\ \bnu.
\ee
An example can be seen in Fig. \ref{fig:f17}, where the osculating plane, circle and sphere are shown for a point $p$ of a conical helix. 

The osculating circle is a diametral circle of the osculating sphere only when $q=q_s$, so if and only if
\bes
\frac{\rho'}{\vartheta}=-\frac{c'}{c^2\vartheta}=0,
\ees
i.e. when the curvature is constant.

\section{Evolute, involute and  envelopes of plane curves}
For any plane curve $\bg(s)$, the center of the osculating circle  $q$ describes a curve $\bd(\sigma)$ that is called the {\it evolute} of $\bg(s)$ ($s$ and $\sigma$ are curvilinear abscissae).  A point $q$ of the evolute is then given by Eq. (\ref{eq:cco}). We call {\it involute} of a curve $\bg(s)$ a curve $\bmu(\sigma)$ whose evolute is $\bg(s)$. 
We call {\it envelope} of a family of plane curves $\bphi(s,\kappa),\kappa\in\mathbb{R}$ being a parameter, a curve that is tangent, in each of  its points, to the curve of $\bphi(s,\kappa)$ passing through that point.

Let us consider the evolute $\bd(\sigma)$ of a curve $\bg(s)$; the tangent to $\bd(\sigma)$ is the vector, cf. Eq. (\ref{eq:cco}),
\bes
\btau_\delta=\frac{dq}{d\sigma}=\frac{dq}{ds}\frac{ds}{d\sigma}.
\ees
But, cf. again  Eq. (\ref{eq:cco}) and the Frenet-Serret formulae,
\bes
\frac{dq}{ds}=\frac{dp}{ds}+\frac{d\rho}{ds}\bnu+\rho\frac{d\bnu}{ds}=\btau+\frac{d\rho}{ds}\bnu-\rho\ c\ \btau=\frac{d\rho}{ds}\bnu,
\ees
so
\bes
\btau_\delta=\frac{dq}{d\sigma}=\frac{d\rho}{ds}\frac{ds}{d\sigma}\bnu.
\ees
Because
\bes
\left|{\frac{dq}{d\sigma}}\right|=|\btau_\delta|=|\bnu|=1,
\ees
then
\bes
\frac{d\rho}{ds}\frac{ds}{d\sigma}=1\ \Rightarrow\ \frac{d\rho}{ds}=\frac{d\sigma}{ds}
\ees
and
\bes
\btau_\delta=\bnu.
\ees
The evolute, $\bd(\sigma)$, of $\bg(s)$ is hence the envelope of its principal normals $\bnu(s)$. 

\begin{figure}[h]
\begin{center}
\includegraphics[width=.54\textwidth]{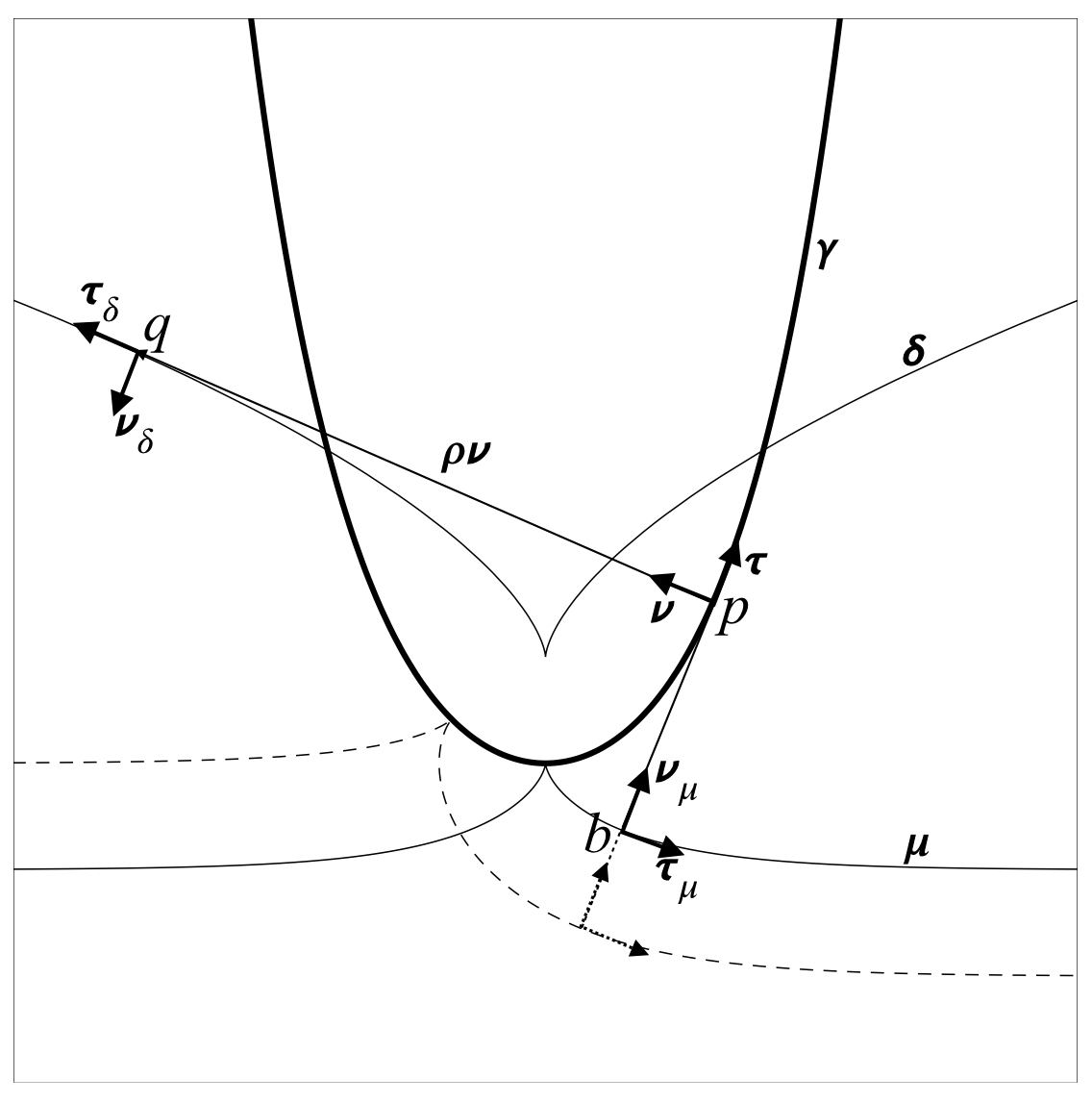}
\caption{Evolute, $\bd$, and involutes for $a=0$, denoted by $\bmu$, and $a=1$, dashed, of a catenary $\bg$.}
\label{fig:51}
\end{center}
\end{figure}

This result helps us in finding the equation of the involute $\bmu(\sigma)$ of a curve $\bg(s)$; let   $p=p(s)$ be a point of $\bg(s)$; then, if $b\in\bmu(\sigma)$, it must be that
\bes
(b-p)\cdot\bnu=0
\ees
where $\bnu$ is the principal normal to $\bg(s)$ in $p$, because $\bg(s)$ is the evolute of $\bmu(\sigma)$, which implies, for the last result, that $\btau=\bnu_\mu$, with $\btau$ the tangent to $\bg(s)$ in $p$ and $\bnu_\mu$ the principal normal to $\bmu(\sigma)$ in $b$, see Fig. \ref{fig:51}.

Therefore, 
\bes
b(s)-p(s)=f(s)\btau(s)\ \rightarrow\ b(s)=p(s)+f(s)\btau(s),
\ees
with $f=f(s)$ a scalar function of $s$; to remark that $b=b(s)$, i.e. the arc-length $s$ of $\bg(s)$ is the parameter also for $\bmu(s)$, but in general $\sigma\neq s$. Upon differentiation, we get
\bes
b'(s)=(1+f'(s))\btau(s)+f(s)c(s)\bnu(s).
\ees
Then, because $b'(s)=|b'(s)|\btau_\mu$ is orthogonal to $\bnu_\mu=\btau$,  it is parallel to $\bnu$, so it must be that
\bes
1+f'(s)=0\ \Rightarrow\ f(s)=a-s,\ \ \ a\in\R.
\ees
Finally, the equation of the involute $\bmu(s)$ to $\bg(s)$ is
\bes
b(s)=p(s)+(a-s)\btau(s),
\ees
and we remark that the involute is not unique.

\section{The theorem of Bonnet}
The  curvature, $c(s)$, and the torsion, $\vartheta(s)$, are the only differential parameters that completely describe a curve. In other words, given two functions $c(s)$ and $\vartheta(s)$, then a curve exists with such a curvature and torsion (we remark that there are no conditions bounding these parameters). This is proved in the following
\begin{teo}{\bf(Bonnet's theorem).} Given two scalar functions $c(s)\in$C$^1$ and $\vartheta(s)\in$C$^0$, there always exists a unique  curve $\bg\in$C$^3$ whose curvilinear abscissa is $s$, curvature $c(s)$ and torsion $\vartheta(s)$.
\begin{proof}
Let 
\bes
\mathbf{e}=\left(\begin{array}{c} \btau\\\bnu\\\bb\end{array}\right)
\ees
be the column vector whose elements are the vectors of the Frenet-Serret basis. Then
\be
\label{eq:bonnet1}
\frac{d\mathbf{e}(s)}{ds}=\mathbf{C}(s)\mathbf{e}(s)
\ee
with 
\bes
\mathbf{C}(s)=\left[
\begin{array}{ccc}
0&c(s)&0\\-c(s)&0&-\vartheta(s)\\
0&\vartheta(s)&0
\end{array}\right]
\ees
 the matrix of Cartan. Adding the initial condition 
\bes
\mathbf{e}(0)=\left(\begin{array}{c} \eu\\\ed\\\et\end{array}\right)
\ees
we have a Cauchy problem for the basis $\mathbf{e}(0)$. As known, such a problem admits a unique solution, i.e. we can associate to $c(s)$ and $\vartheta(s)$ a family of bases $\mathbf{e}(s)$ (that are orthonormal because if one of them were not so, the Cartan's matrix should not be skew). 
Call $\btau(s)$ the first vector of the basis $\mathbf{e}(s)$ and define the function 
\bes
p(s):=p_0+\int_0^s\btau(s^*)ds^*;
\ees
$p(s)$ is the curve looked for (it depends upon an arbitrary point $p_0$, i.e. upon an inessential rigid displacement). In fact, because $|\btau|=1$,  $s$ is the curvilinear abscissa of the curve. Then, it is sufficient to write the Frenet-Serret equations identifying them with the system (\ref{eq:bonnet1}).
\end{proof}
\end{teo}

\section{Canonic equations of a curve}
We call {\it the canonic equations} of a curve at a point $p_0$ the equations of the curve referred to the Frenet-Serret basis in $p_0$. For this purpose, we expand the curve in a Taylor series of initial point $p_0$:
\bes
p(s)=p_0+s\ p'(0)+\frac{1}{2}s^2p''(0)+\frac{1}{6}s^3p'''(0)+o(s^3).
\ees
In the Frenet-Serret basis,
\bes
p'(0)=\btau(0),\ \ p''(0)=c(0)\bnu(0),\ \ p'''(0)=\left.\dfrac{dc\bnu}{ds}\right|_{s=0}=c'(0)\bnu(0)-c^2(0)\btau(0)-c(0)\vartheta(0)\bb(0),
\ees
so
\bes
p(s)=p_0+s\ \btau(0)+\frac{1}{2}s^2c(0)\bnu(0)+\dfrac{1}{6}s^3(-c^2(0)\btau(0)+c'(0)\bnu(0)-c(0)\vartheta(0)\bb(0))+o(s^3).
\ees
The coordinates of a point $p(s)$ close to $p_0$, in the basis $\{\btau(0),\bnu(0),\bb(0)\}$, are hence
\bes
\besp
&p_1(s)=s-\frac{1}{6}c^2(0)s^3+o(s^3),\\
&p_2(s)=\frac{1}{2}c(0)s^2+\dfrac{1}{6}c'(0)s^3+o(s^3),\\
&p_3(s)=-\frac{1}{6}c(0)\vartheta(0)s^3+o(s^3).
\end{split}
\ees
The projections of the curve onto the planes of the Frenet-Serret basis  hence have, close to $p_0$ (i.e. retaining the first non null term in the expressions above), the following equations:
\begin{itemize}
\item On the osculating plane 
\bes
\left\{\begin{array}{l}p_1(s)=s,\\p_2(s)=\dfrac{1}{2}c(0)s^2,\end{array}\right.
\ees
or, eliminating $s$,
\bes
p_2=\dfrac{1}{2}c(0)p_1^2,
\ees
which is the equation of a parabola.
\item On the rectifying plane
\bes
\left\{\begin{array}{l}p_1(s)=s,\\p_3(s)=-\dfrac{1}{6}c(0)\vartheta(0)s^3,\end{array}\right.
\ees
or, eliminating $s$,
\bes
p_3=-\frac{1}{6}c(0)\vartheta(0)p_1^3,
\ees
which is the equation of a cubic parabola.
\item On the normal plane 
\bes
\left\{\begin{array}{l}p_2(s)=\dfrac{1}{2}c(0)s^2,\\p_3(s)=-\dfrac{1}{6}c(0)\vartheta(0)s^3,\end{array}\right.
\ees
or, eliminating $s$,
\bes
p_3^2=\frac{2}{9}\frac{\vartheta^2(0)}{c(0)}p_2^3,
\ees
which is the equation of a semicubic parabola, with a cusp at the origin, hence a singular point, though the curve $p(s)$ is regular.
\end{itemize}

\section{Exercises}
\begin{enumerate}
\item Using the same definition of derivative of a curve, prove the relations in Sect. (\ref{sec:diffcurves}).
\item Prove the relations in Eq. (\ref{eq:derivcomp}).
\item The curve whose polar equation is 
\begin{equation*}
r=a\ \theta,\ \ a\in\mathbb{R},
\end{equation*}
is an {\it Archimedes' spiral}, Fig. \ref{fig:60} a). Find its curvature $c(\theta)$ and its length $\ell(\theta)$, and prove that any straight line passing through the origin is divided by the spiral in segments of constant length $2\pi\ a$ (that is why the Archimede's spiral is used to record disks).
\begin{figure}[h]
\begin{center}
\includegraphics[width=.6\textwidth]{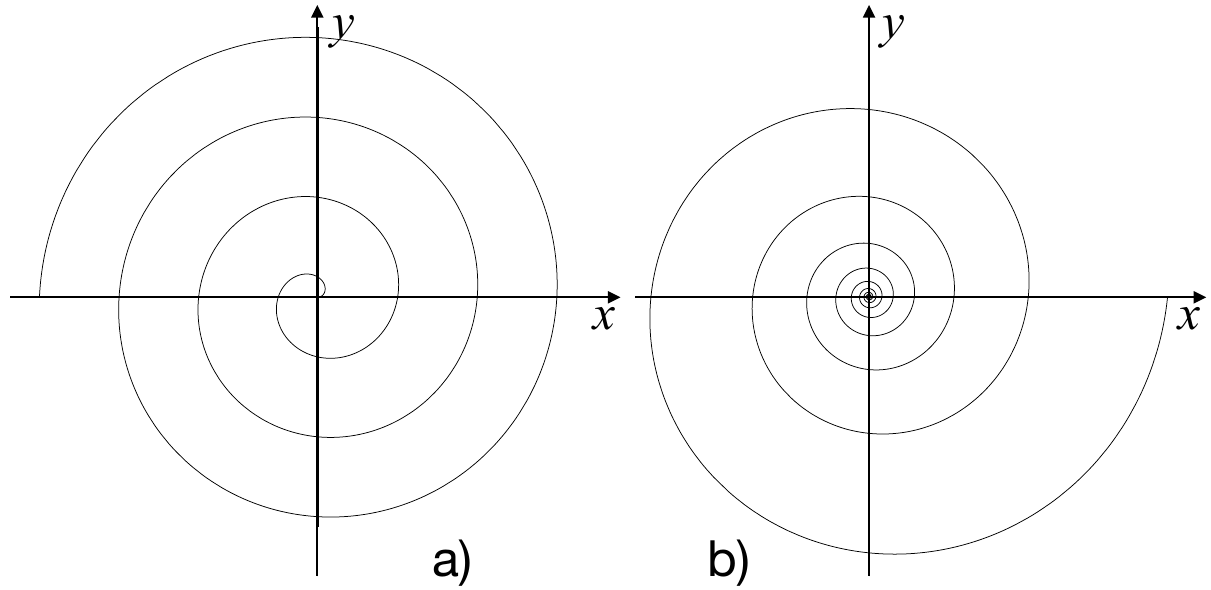}
\caption{The Archimedes', a), and logarithmic, b), spirals.}
\label{fig:60}
\end{center}
\end{figure}

\item The curve whose polar equation is 
\begin{equation*}
r=a\ \mathrm{e}^{b \theta},\ \ a,b\in\mathbb{R},
\end{equation*}
is the {\it logarithmic spiral}. Prove that the origin is an asymptotic point of the curve, find its curvature $c(\theta)$ and its length  $\ell(\theta)$, and show that the length of the segments in which a straight line by the origin is divided by two consecutive intersections with the spiral varies as a geometrical progression. Then, prove  its {\it equiangular property}: The angle $\alpha$ between $p(\theta)-o$ and $\btau(\theta)$ is constant. Finally, show that the evolute of the logarithmic spiral is a logarithmic spiral itself (an hence that its involute is still a logarithmic spiral, that's why Jc. Bernoulli coined for this curve the Latin sentence {\it eadem mutata resurgo}.)

\item The curve whose parametric equation is 
\begin{equation*}
p(\theta)=a(\cos\theta+\theta\sin\theta)\gr{e}_1+a(\sin\theta-\theta\cos\theta)\gr{e}_2
\end{equation*}
with  $\theta$ the angle formed by $p(\theta)-o$ with the axis $x_1$, is the {\it involute of the circle}, Fig. \ref{fig:61}. Find its curvature $c(\theta)$ and its length $\ell(\theta)$, and prove that its evolute is  exactly the circle of center $o$ and radius $a$ (that is why the involute of the circle is used to profile gears).
\begin{figure}[h]
\begin{center}
\includegraphics[width=.35\textwidth]{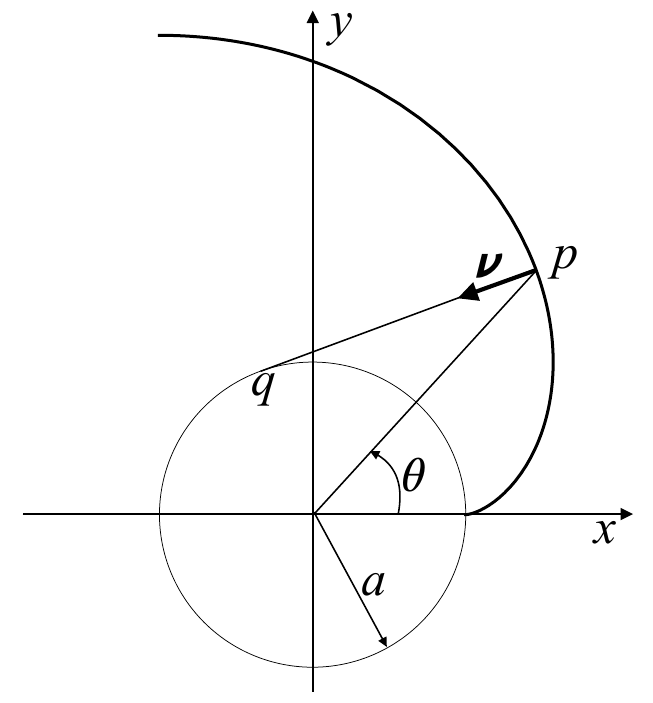}
\caption{The involute of the circle and its evolute, the circle.}
\label{fig:61}
\end{center}
\end{figure}

\item The curve whose parametric equation is 
\begin{equation*}
p(\theta)=a\cos\omega\theta\gr{e}_1+a\sin\omega\theta\gr{e}_2+b\omega\theta\gr{e}_3
\end{equation*}
is a {\it circular helix}, i.e. a {\it helix} that winds on a circular cylinder of radius $a$, Fig. \ref{fig:62}. Show that the angle $\phi$ formed by the helix and any generatrix of the cylinder is constant (a property that defines a helix in the general case). Then, find its length  $\ell(\theta)$, its curvature $c(\theta)$, torsion $\vartheta(\theta)$, and the pitch $d$, i.e., the distance between two successive intersections of the helix with a generatrix of the cylinder. 
Prove then the {\it Bertrand's theorem}: A curve is a cylindrical helix if and only if the ratio $c/\vartheta=const.$ 
Finally, prove that for the above circular helix there are two constants $A$ and $B$ such that
\begin{equation*}
p'\times p''=A\gr{u}(\theta)+B\gr{e}_3,
\end{equation*}
with
\begin{equation*}
\gr{u}=\sin\omega\theta\gr{e}_1-\cos\omega\theta\gr{e}_2;
\end{equation*}
then, find  $A$ and $B$.
\begin{figure}[h]
\begin{center}
\includegraphics[width=.35\textwidth]{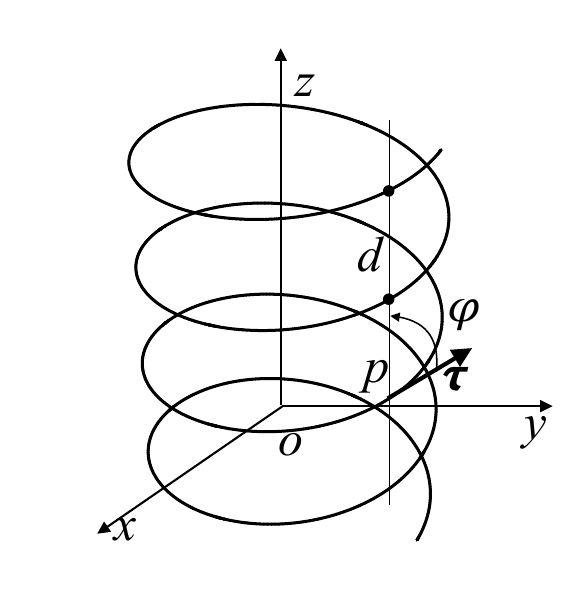}
\caption{The circular helix.}
\label{fig:62}
\end{center}
\end{figure}

\item Find the equation of the {\it cycloid}, i.e. of the curve that is the trace of a point of a circle of radius $r$ rolling without slipping on a horizontal axis, see Fig. \ref{fig:58}. Calculate the length of the cycloid for a complete round of the circle, determine its curvature, and show that the evolute of the cycloid is the cycloid itself (Huygens, 1659).
\begin{figure}[h]
\begin{center}
\includegraphics[width=.6\textwidth]{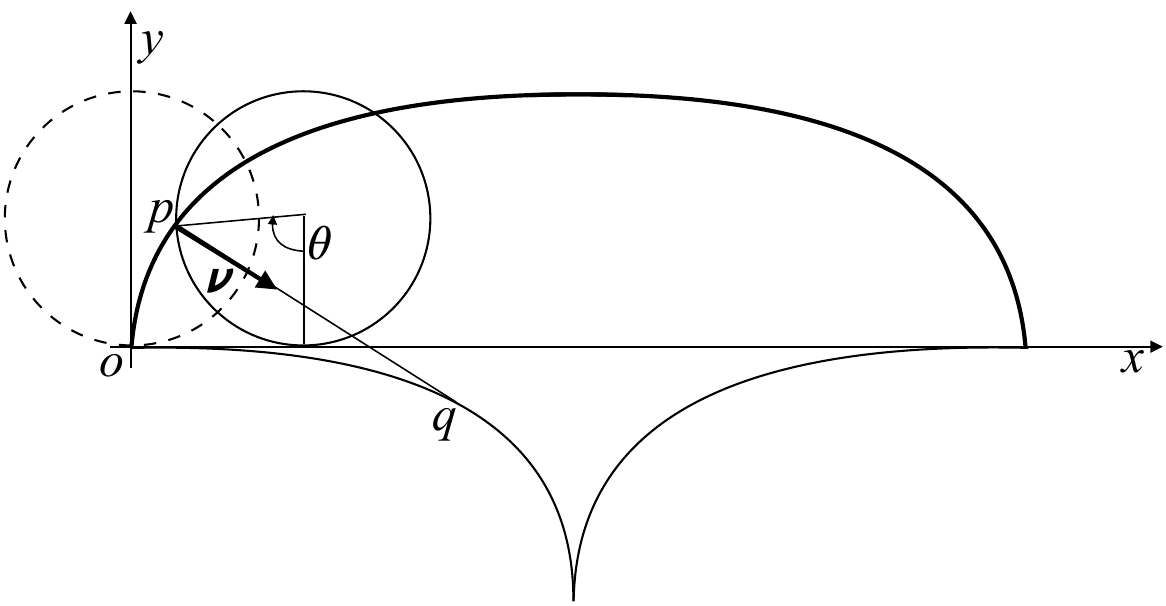}
\caption{The cycloid and its evolute.}
\label{fig:58}
\end{center}
\end{figure}

\item The planar curve whose parametric equation is 
\bes
p(t)=t\eu+\cosh t\ed
\ees
is the {\it catenary} (Jc. Bernoulli, 1690; Jn. Bernoulli, Leibniz, Huygens, 1691). It is the equilibrium curve of a heavy, perfectly flexible, and inextensible cable. Calculate the curvature of the catenary and the equation of its evolute and of its involutes (see Fig. \ref{fig:51}).

\item The planar curve whose parametric equation is 
\bes
p(t)=\left(\cos t+\ln\tan\frac{t}{2}\right)\eu+\sin t\ed
\ees
is the {\it tractrix} (Perrault, 1670; Newton, 1676; Huygens, 1693). This is the curve along which an object moves, under the influence of friction, when pulled on a horizontal plane by a line segment attached to a tractor that moves at a right angle to the initial line between the object and the puller at an infinitesimal speed, see Fig. \ref{fig:59}. Show that the length of the tangent to the tractrix between the points on the tractrix itself and the axis $x$ is constant $\forall t$, calculate the length of the curve between $t_1$ and $t_2$, calculate the curvature of the tractrix, and finally show that its evolute is the catenary.
\begin{figure}[h]
\begin{center}
\includegraphics[width=.6\textwidth]{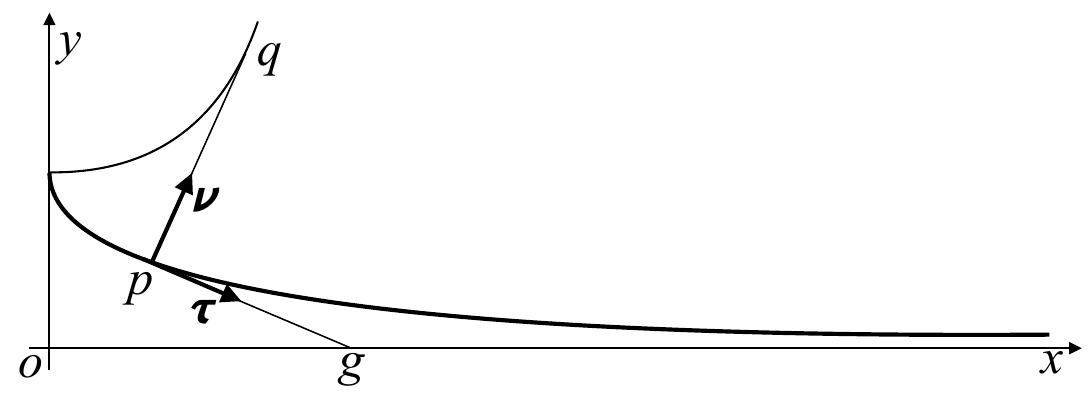}
\caption{The tractrix and its evolute.}
\label{fig:59}
\end{center}
\end{figure}

\item For the curve whose cylindrical equation is
\begin{equation*}
\left\{
\begin{split}
&r=1,\\
&z=\sin\theta
\end{split}
\right.
\end{equation*}
find the highest curvature and determine whether or not it is planar.

\item Let $p=p(t)$ be the path of a moving particle of masse $m$, with $t$ being the time. Define the velocity and the acceleration of $p$ as, respectively, the first and the second derivative of $p$ with respect to $t$. Decompose these two vectors in the Frenet-Serret basis and interpret physically the result. Recalling the second Newton's principle of mechanics, what about the forces on $p$?
\end{enumerate}

\chapter{Tensor analysis: fields}
\label{ch:5}

\section{Scalar, vector and tensor fields}
Let $\Omega\subset\Eu$ and $\f:\Omega\rightarrow\Ve$. We say that $\f$ is {\it continuous at} $p\in\Omega\iff\forall$ sequence $\pi_n=\{p_n\in\Omega,n\in\Nq\}$ that converges to $p\in\Eu$, the sequence $\{\bv_n=\f(p_n),n\in\Nq\}$ converges to $\f(p)$ in $\Ve$. The function $\f(p):\Omega\rightarrow\Ve$ is a {\it vector field on} $\Omega$ if it is continuous at each $p\in\Omega$. In the same way we can define a {\it scalar field} $\varphi(p):\Omega\rightarrow\R$ and a {\it tensor field}, $\L(p):\Omega\rightarrow Lin(\Ve)$.

A {\it deformation} is any continuous and bijective function $f(p):\Omega\rightarrow\Eu$, i.e. any transformation of a region $\Omega\subset\Eu$ into another region of $\Eu$; bijectivity imposes that to any point $p\in\Omega$ corresponds one and only one point in the transformed region, and vice-versa, which is the  mathematical condition expressing the physical constraint of mass conservation. 

Finally, the basic difference between fields/deformations and curves, is that a field or a deformation is defined over a subset of $\Eu$, not of $\R$. In practice, this implies that the components of the field/deformation are functions of three variables, the coordinates $x_i$ of a point $p\in\Omega$.

\section{Differentiation of fields, differential operators}
\label{sec:difffields}
Let $\psi(p)$ be a scalar or vector or tensor field or also a deformation; we define {\it directional derivative} of $\psi(p)$ in the direction of $\e\in\S$ the limit
\bes
\frac{d\psi(p)}{d\e}:=\lim_{\alpha\rightarrow0}\frac{\psi(p+\alpha\e)-\psi(p)}{\alpha}, \ \ \ \alpha\in\R.
\ees
The directional derivative measures the rate of variation of $\psi(p)$ in the direction of $\e$. In the particular case of $\e=\ei,\ i=1,2,3$, i.e. of the directions of the basis $\{\eu,\ed,\et\}$ of $\Ve$, then 
\bes
\frac{d\psi(p)}{d\ei}=\lim_{\alpha\rightarrow0}\frac{\psi(p+\alpha\ei)-\psi(p)}{\alpha}
\ees
is the {\it partial derivative of $\psi$ with respect to $x_i$}; e.g., if $i=1$, then
\bes
\frac{d\psi(p)}{d\eu}=\lim_{\alpha\rightarrow0}\frac{\psi(x_1+\alpha,x_2,x_3)-\psi(x_1,x_2,x_3)}{\alpha}.
\ees
The partial derivative with respect to $x_i$ is usually indicated as $\dfrac{\partial\psi}{\partial x_i}$ or also as $\psi_{,i}$.

Let $\bv(p):\Omega\rightarrow\Ve$; we say that $\bv$ is {\it differentiable in} $p_0\in\Omega\iff\exists\ \mathrm{grad}\bv\in Lin(\Ve)$ such that
\bes
\bv(p_0+\bu)=\bv(p)+\mathrm{grad}\bv(p)\ \bu+o(u)
\ees
when $\bu\rightarrow\bo$. If $\bv$ is differentiable $\forall p\in\Omega,\mathrm{grad}\bv$ defines a tensor field on $\Omega$ called the {\it gradient of} $\bv$. It is also possible to define higher order differential operators, using higher order tensors, but this will not be done here. If $\bv$ is continuous with $\mathrm{grad}\bv\ \forall p\in\Omega$, then $\bv$ is of class C$^1$ ({\it smooth}).

Let $\bv$ be a vector field of class C$^1$ on $\Omega$. Then, the {\it divergence} of $\bv$ is the scalar field defined by 
\bes
\div\bv:=\tr(\mathrm{grad}\bv),
\ees
while $\curl\bv$ is the unique vector field  that satisfies the relation
\bes
(\mathrm{grad}\bv-\mathrm{grad}\bv^\top)\bu=(\curl\bv)\times\bu\ \ \forall\bu\in\Ve.
\ees
The {\it divergence of a tensor field} $\L$ is the unique vector field $\div\L$ that satisfies
\bes
(\div\L)\cdot\bu=\div(\L^\top\bu)\ \ \forall\bu=const.\in\Ve.
\ees
Let $\varphi(p):\Omega\rightarrow\R$ be a scalar field over $\Omega$. Similar to the case of  vector fields, we say that $\varphi$ is {\it differentiable at} $p_0\in\Omega\iff\exists\ \mathrm{grad}\varphi\in\Ve$ such that
\bes
\varphi(p+\bu)=\varphi(p)+\mathrm{grad}\varphi(p)\cdot\bu+o(u)
\ees
when $\bu\rightarrow\bo$. If $\varphi$ is differentiable $\forall p\in\Omega,\mathrm{grad}\varphi$ defines a vector field on $\Omega$ called the {\it gradient of} $\varphi$. If $\mathrm{grad}\varphi$ is differentiable, its gradient is the tensor $\mathrm{grad}^{II}\varphi$ called {\it second gradient} or {\it Hessian}. It is immediate to show that under continuity assumption,
\bes
\mathrm{grad}^{II}\varphi=(\mathrm{grad}^{II}\varphi)^\top.
\ees
A {\it level set} of a scalar field $\varphi(p)$ is the set $\S_L$ such that
\bes
\varphi(p)=const.\ \ \forall p\in\S_L.
\ees
Considering hence two points $p$ and $ p+\bu$ of the same $\S_L$,  then by the  definition of differentiability of $\varphi(p)$ itself, we see that $\mathrm{grad}\varphi$ is a vector that is orthogonal to 
$\S_L$ at $p$. The curves of $\Eu$ that are tangent to $\mathrm{grad}\varphi\ \forall p\in\Omega$ are the {\it streamlines of } $\varphi$; they have the property to be orthogonal to any $\S_L$ of $\phi\ \forall p\in\Omega$.

$\mathrm{grad}\phi$ allows to calculate the  directional derivative of $\phi$ along any direction $\n\in\S$ as
\bes
\frac{d\phi}{d\n}=\mathrm{grad}\phi\cdot\n.
\ees
The highest variation of $\phi$ is hence in the direction of $\mathrm{grad}\phi$, and $|\mathrm{grad}\phi|$ is the value of this variation; we remark also that $\mathrm{grad}\phi$ is a vector directed as the increasing values of $\phi$.

Similarly, for a vector field $\bv$ the  directional derivative along  any direction $\n\in\S$ can be computed as
\bes
\frac{d\bv}{d\n}=\mathrm{grad}\bv\ \n.
\ees

Let $\psi$ be a scalar of vector field of class C$^2$ at least. Then, the {\it laplacian $\Delta\psi$ of} $\psi$ is defined by
\bes
\Delta\psi:=\div(\mathrm{grad}\psi).
\ees
By the linearity of the trace, and hence of the divergence, we see easily that the laplacian of a vector field is the vector field whose components are the laplacian of each corresponding component of the field. A field is said to be {\it harmonic} on $\Omega$ if its laplacian is null $\forall p\in\Omega$.

The definitions given above for differentiable field, gradient and class C$^1$ can be repeated {\it verbatim} for a deformation $f(p):\Omega\rightarrow\Eu$.

\section{Properties of the differential operators}
The differential operators, gradient, divergence, curl and laplacian, have some interesting properties, that are useful for calculations; they are introduced in this section.
\begin{teo} \label{teo:propgrad} {\bf{(Gradient of products).}} Let $\phi,\psi$  be scalar  and $\bu,\bv,\bw$ be vector fields, with all of them differentiable. Then:
\be
\label{eq:gradproducts}
\begin{array}{ll}
\mathrm{i)}& \mathrm{grad}(\phi\psi)=\phi\ \mathrm{grad}\psi+\psi\ \mathrm{grad}\phi,\medskip\\
\mathrm{ii)} &\mathrm{grad}(\phi\bv)=\phi\ \mathrm{grad}\bv+\bv\otimes\mathrm{grad}\phi,\medskip\\
\mathrm{iii)}& \mathrm{grad}(\bv\cdot\bw)=(\mathrm{grad}\bw)^\top\bv+(\mathrm{grad}\bv)^\top\bw.
\end{array}
\ee
\begin{proof} The proof is based upon the definition of gradient itself\footnote{For the sake of brevity, we omit to indicate the point $p$, e.g. we simply write $\phi$ for $\phi(p)$, $\grad\phi$ for $\grad\phi(p)$ etc.}:

i) \hspace{37mm}$(\phi\psi)(p+\bu)=\phi\psi+\grad(\phi\psi)\cdot\bu+o(u),$

but  also
\bes
\besp
(\phi\psi)(p+\bu)&=\phi(p+\bu)\psi(p+\bu)=(\phi+\grad\phi\cdot \bu+o(u))(\psi+\grad\psi\cdot \bu+o(u))\\
&=\phi\psi+\phi\ \grad\psi\cdot \bu+\psi\ \grad\phi\cdot \bu+o(u)\\
&=\phi\psi+(\phi\ \grad\psi+\psi\ \grad\phi)\cdot\bu+o(u),
\end{split}
\ees
so by comparison
\bes
\grad(\phi\psi)=\phi\grad\psi+\psi\grad\phi.
\ees
ii) in the same way 
\bes
(\phi\bv)(p+\bu)=\phi\bv+\grad(\phi\bv)\bu+o(u)=\phi\bv+\grad(\phi\bv)\bu+o(u),
\ees
but also
\bes
\besp
(\phi\bv)(p+\bu)&=\phi(p+\bu)\bv(p+\bu)=(\phi+\grad\phi\cdot \bu+o(u))(\bv+\grad\bv\ \bu+o(u))\\
&=\phi\bv+\phi\grad\bv\ \bu+\grad\phi\cdot \bu\ \bv+o(u)\\
&=\phi\bv+(\phi\ \grad\bv+\bv\otimes\grad\bv)\bu+o(u),
\end{split}
\ees
so comparing the two results we get
\bes
\grad(\phi\bv)=\phi\ \grad\bv+\bv\otimes\grad\bv.
\ees
iii) in the same way
\bes
(\bv\cdot\bw)(p+\bu)=\bv\cdot\bw+\grad(\bv\cdot\bw)\cdot\bu+o(u),
\ees
but also
\bes
\besp
(\bv\cdot\bw)(p+\bu)&=\bv(p+\bu)\cdot\bw(p+\bu)=(\bv+\grad\bv\ \bu+o(u))\cdot(\bw+\grad\bw\ \bu+o(u))\\
&=\bv\cdot\bw+\bv\cdot(\grad\bw\ \bu)+(\grad\bv\ \bu)\cdot\bw+o(u)\\
&=\bv\cdot\bw+((\grad\bw)^\top \bv+(\grad\bv)^\top\bw)\cdot\bu+o(u),
\end{split}
\ees
whence, by comparison of the two results,
\bes
\grad(\bv\cdot\bw)=(\grad\bw)^\top \bv+(\grad\bv)^\top\bw.
\ees
\end{proof}
\end{teo}
Another important result\footnote{This result is fundamental to fluid mechanics, as it allows us to get an interesting form of the Navier-Stokes equations.}, relating the gradient and the curl of a vector field, is the following theorem.
\begin{teo} \label{teo:rotoremecaflu}  If $\bv$ is a differentiable vector field, then
\bes
(\mathrm{grad}\bv)\bv=(\curl\bv)\times\bv+\dfrac{1}{2}\mathrm{grad}\bv^2.
\ees
\begin{proof}
\bes
\besp
(\curl\bv)\times\bv&=(\grad\bv-(\grad\bv)^\top)\bv=(\grad\bv)\bv-(\grad\bv)^\top\bv\\
&=(\grad\bv)\bv-\dfrac{1}{2}((\grad\bv)^\top\bv+(\grad\bv)^\top\bv),
\end{split}
\ees
and by property iii) of the previous theorem,
\bes
(\grad\bv)^\top\bv+(\grad\bv)^\top\bv=\grad(\bv\cdot\bv)=\grad\bv^2,
\ees
so that
\bes
(\curl\bv)\times\bv=(\grad\bv)\bv-\dfrac{1}{2}\grad\bv^2,
\ees
whence we obtain  the thesis.
\end{proof}
\end{teo}
The proof of the following properties of the gradient are left to the reader as an exercise:
\be
\label{eq:propgrad}
\begin{array}{c}
\mathrm{grad}(\bv\cdot\bw)=(\mathrm{grad}\bw)\bv+(\mathrm{grad}\bv)\bw+\bv\times\curl\bw+\bw\times\curl\bv,\medskip\\
\mathrm{grad}(\bu\cdot\bv\ \bw)=(\bu\cdot\bv)\mathrm{grad}\bw+(\bw\otimes\bu)\mathrm{grad}\bv+(\bw\otimes\bv)\mathrm{grad}\bu,\medskip\\
\mathrm{grad}\bv\cdot\mathrm{grad}\bv^\top=\div((\mathrm{grad}\bv)\bv-(\div\bv)\bv)+(\div\bv)^2.
\end{array}
\ee

\begin{teo}
\label{teo:proddiv}
{\bf{(Divergence of products).}} Let $\phi,\bu,\bv,\bw,\L$ be differentiable scalar, vector or tensor fields. Then:
\bes
\begin{array}{ll}
\mathrm{i)}&\div(\phi\bv)=\phi\div\bv+\bv\cdot\mathrm{grad}\phi,\medskip\\
\mathrm{ii)}&\div(\bv\otimes\bw)=\bv\div\bw+(\mathrm{grad}\bv)\bw,\medskip\\
\mathrm{iii)}&\div(\phi\L)=\phi\div\L+\L\mathrm{grad}\phi,\medskip\\
\mathrm{iv)}&\div(\L^\top\bv)=\L\cdot\mathrm{grad}\bv+\bv\cdot\div\L,\medskip\\
\mathrm{v)}&\div(\bv\times\bw)=\bw\cdot\curl\bv-\bv\cdot\curl\bw.
\end{array}
\ees
\begin{proof}
i) Using the definition of divergence and property ii) of Theorem \ref{teo:propgrad}, we get
\bes
\besp
\div(\phi\bv)&=\tr(\grad(\phi\bv))=\tr(\phi\ \grad\bv+\bv\otimes\grad\phi)\\
&=\phi\ \tr(\grad\bv)+\tr(\bv\otimes\grad\phi)=\phi\ \div\bv+\bv\cdot\grad\phi.
\end{split}
\ees
ii) By the definition of divergence of a tensor, $\forall \a=const.\in\Ve$, and using the previous property along with property iii) of Theorem \ref{teo:propgrad}:
\bes
\besp
\div(\bv\otimes\bw)\cdot\a&=\div((\bv\otimes\bw)^\top\a)=\div(\bw\otimes\bv\ \a)=\div(\bv\cdot\a\ \bw)\\&=\bv\cdot\a\ \div\bw+\bw\cdot\grad(\a\cdot\bv)\\
&=\div\bw\ \bv\cdot\a+\bw\cdot(\grad\bv)^\top\a+\bw\cdot(\grad\a)^\top\bv\\
&=(\bv\div\bw+\grad\bv\ \bw)\cdot\a.
\end{split}
\ees
iii) By the definition of divergence of a tensor, $\forall \a=const.\in\Ve$, and using the  property i) along with the iii) of Theorem \ref{teo:propgrad}:
\bes
\besp
\div(\phi\L)\cdot\a&=\div((\phi\L)^\top\a)=\div(\phi\L^\top\a)=\phi\div(\L^\top\a)+\L^\top\a\cdot\grad\phi\\
&=\phi\div\L\cdot\a+\a\cdot\L\grad\phi=(\phi\div\L+\L\grad\phi)\cdot\a.
\end{split}
\ees
iv) By the definition of divergence of a tensor and using the  previous property:
\bes
\besp
\div(\L^\top\bv)&=\div(\L^\top v_j\ej)=\div((v_j\L^\top)\ej)=\div(v_j\L^\top)^\top\cdot\ej\\
&=\div(v_j\L)\cdot\ej=v_j\div\L\cdot\ej+\L\ \grad v_j\cdot\ej\\&=\bv\cdot\div\L+(L_{pq}\ep\otimes\eq(\grad v_j)_m\e_m)\cdot\ej\\
&=\bv\cdot\div\L+L_{pq}v_{j,m}\delta_{qm}\delta_{jp}=\bv\cdot\div\L+\L\cdot\grad\bv.
\end{split}
\ees

v) This property can be proved making use of the expression of the cross product with the Ricci's alternator, given in Eq.  (\ref{eq:prodvectricci}):
\bes
\curl\bv=\epsilon_{ijk}v_{k,j}\ei 
\ees 
and
\bes
\bv\times\bw=\epsilon_{ijk}v_jw_k\ei,
\ees
whence
\bes
\div(\bv\times\bw)=\div(\epsilon_{ijk}v_jw_k\ei)=\epsilon_{ijk}(v_jw_k)_{,i}=\epsilon_{ijk}v_{j,i}w_k+\epsilon_{ijk}w_{k,i}v_j.
\ees
Moreover,
\bes
\bw\cdot\curl\bv=w_m\e_m\cdot\epsilon_{pqr}v_{r,q}\ep=\epsilon_{pqr}v_{r,q}w_m\delta_{pm}=\epsilon_{pqr}v_{r,q}w_p=\epsilon_{qrp}v_{r,q}w_p
\ees
and
\bes
\bv\cdot\curl\bw=v_m\e_m\cdot\epsilon_{pqr}w_{r,q}\ep=\epsilon_{pqr}w_{r,q}v_m\delta_{pm}=\epsilon_{pqr}w_{r,q}v_p=-\epsilon_{qpr}w_{r,q}v_p,
\ees
so finally, comparing the last three results (all the subscripts are dummy indexes, so their denomination is inessential),
\bes
\div(\bv\times\bw)=\bw\cdot\curl\bv-\bv\cdot\curl\bw.
\ees
\end{proof}
\end{teo}
The divergence has also the following properties
\be
\label{eq:divproperties}
\begin{array}{c}
\div(\mathrm{grad}\bv^\top)=\mathrm{grad}(\div\bv),\medskip\\
\div((\mathrm{grad}\bv)\bv)=\mathrm{grad}\bv\cdot\mathrm{grad}\bv^\top+\bv\cdot\mathrm{grad}(\div\bv),\medskip\\
\div(\phi\L\bv)=\phi\L^\top\cdot\mathrm{grad}\bv+\phi\bv\cdot\div\L^\top+\L\bv\cdot\mathrm{grad}\phi,
\end{array}
\ee
whose proof is a good exercise for the reader.

The relations of gradient and divergence with the curl are given by the following theorem.
\begin{teo}
Let $\phi$ and $\bv$  be scalar and vector fields of class C$^2$; then
\bes
\begin{array}{ll}
\mathrm{i)}&\div(\curl\bv)=0,\medskip\\
\mathrm{ii)}&\curl(\mathrm{grad}\phi)=\bo.
\end{array}
\ees
\begin{proof}
i) Using again the Ricci's alternator to represent the cross product,
\bes
\besp
\div(\curl\bv)&=\div(\epsilon_{ijk}v_{k,j}\ei)=\epsilon_{ijk}v_{k,j}\div\ei+\epsilon_{ijk}v_{k,ji}=\epsilon_{ijk}v_{k,ji}\\&=v_{3,21}+v_{1,32}+v_{2,13}-v_{2,31}-v_{3,12}-v_{1,23}=0.
\end{split}
\ees
ii) In a similar manner,
\bes
\curl(\grad\phi)=\epsilon_{ijk}\phi_{,kj}\ei=\phi_{,32}+\phi_{,13}+\phi_{,21}-\phi_{,23}-\phi_{,31}-\phi_{,12}=0.
\ees
\end{proof}
\end{teo}
The following theorem gives an interesting relation between the curl of a vector and the divergence of its axial tensor.
\begin{teo} 
\label{teo:curlaxvect}
{\bf{(Curl of an axial vector).}} Let $\bw$ be a differentiable vector field and $\W$ its axial tensor field. Then,
\bes
\curl\bw=-\div\W.
\ees
\begin{proof}
Using properties iv) and v) of Theorem \ref{teo:proddiv}  and because $\W=-\W^\top, \forall \a=const.\in\Ve$ we get
\bes
\begin{array}{c}
\div(\bw\times\a)=\a\cdot\curl\bw-\bw\cdot\curl\a=\a\cdot\curl\bw,\medskip\\
\div(\W\a)=\div(-\W^\top\a)=-\W\cdot\grad\a-\a\cdot\div\W=-\a\cdot\div\W.
\end{array}
\ees
Now, because $\forall \a,  \bw\times\a=\W \a\Rightarrow\div(\bw\times\a)=\div(\W\a)$, we get the thesis.
\end{proof}
\end{teo}
The way the curl of a curl\footnote{This relation is useful in fluid mechanics, for writing the vorticity equation.} is computed is given by the following theorem.
\begin{teo} {\bf{(Curl of a curl).}} Let $\bv$be  a  vector field of class $\geq$C$^{2}$. Then,
\bes
\curl(\curl\bv)=\mathrm{grad}(\div\bv)-\Delta\bv.
\ees
\begin{proof}
Using properties iv) and v) of Theorem \ref{teo:proddiv}, along with the first of Eq.  (\ref{eq:divproperties}),   $\forall \a=const.\in\Ve$ we get
\bes
\div((\curl\bv)\times\a)=\a\cdot\curl(\curl\bv)-\curl\bv\cdot\curl\a=\a\cdot\curl(\curl\bv),
\ees
and by the definition of curl and  laplacian
\bes
\besp
\div((\curl\bv)\times\a)&=\div((\grad\bv-(\grad\bv)^\top)\a)=\div(\grad\bv\ \a)-\div((\grad\bv)^\top\a)\\
&=(\grad\bv)^\top\cdot\grad\a+\a\cdot\div(\grad\bv)^\top-\grad\bv\cdot\grad\a-\a\cdot\div(\grad\bv)\\
&=\a\cdot(\div(\grad\bv)^\top-\div(\grad\bv))=\a\cdot(\grad(\div\bv)-\Delta\bv),
\end{split}
\ees
whence, by comparison,
\bes
\curl(\curl\bv)=\grad(\div\bv)-\Delta\bv.
\ees
\end{proof}
\end{teo}
The proof of the following properties of the curl are can be done using the above results and it is a good exercise:
\be
\label{eq:propcurl}
\begin{array}{c}
\curl(\phi\bv)=\phi\curl\bv+\mathrm{grad}\phi\times\bv,\medskip\\
\curl(\bv\times\bw)=(\mathrm{grad}\bv)\bw-(\mathrm{grad}\bw)\bv+\bv\div\bw-\bw\div\bv.
\end{array}
\ee

Finally,   we have a theorem also for the laplacian of a product.
\begin{teo}
{\bf{(Laplacian of products).}} Let $\phi,\psi,\bu,\bv$ be scalar and vector fields of class $\geq$C$^2$. Then:
\bes
\begin{array}{ll}
\mathrm{i)}&\Delta(\phi\psi)=2\mathrm{grad}\phi\cdot\mathrm{grad}\psi+\phi\Delta\psi+\psi\Delta\phi,\medskip\\
\mathrm{ii)}&\Delta(\bv\cdot\bw)=2\mathrm{grad}\bv\cdot\mathrm{grad}\bw+\bv\cdot\Delta\bw+\bw\cdot\Delta\bv.
\end{array}
\ees
\begin{proof} i) Using properties i) of Theorems \ref{teo:propgrad}  and \ref{teo:proddiv}, we get
\bes
\besp
\Delta(\phi\psi)&=\div(\grad(\phi\psi))=\div(\phi\ \grad\psi+\psi\ \grad\phi)\\
&=\div(\phi\ \grad\psi)+\div(\psi\ \grad\phi)\\
&=\phi\ \div(\grad\psi)+\grad\psi\cdot\grad\phi+\psi\ \div(\grad\phi)+\grad\phi\cdot\grad\psi\\
&=2\grad\phi\cdot\grad\psi+\phi\ \Delta\psi+\psi\ \Delta\phi.
\end{split}
\ees
ii) Using properties iii) of Theorem \ref{teo:propgrad} and iv) of Theorem \ref{teo:proddiv}, we obtain
\bes
\besp
\Delta(\bv\cdot\bw)&=\div(\grad(\bv\cdot\bw))=\div((\grad\bw)^\top\bv+(\grad\bv)^\top\bw)\\
&=\div((\grad\bw)^\top\bv)+\div((\grad\bv)^\top\bw)\\
&=\grad\bw\cdot\grad\bv+\bv\cdot\div(\grad\bw)+\grad\bv\cdot\grad\bw+\bw\cdot\div(\grad\bv)\\
&=2\grad\bv\cdot\grad\bw+\bv\cdot\Delta\bw+\bw\cdot\Delta\bv.
\end{split}
\ees
\end{proof}
\end{teo}

\section{Theorems on fields}
We recall here some classical theorems on fields and operators. 
\begin{teo}{\bf(Harmonic fields).} If $\bv(p)$ is a vector field of class $\geq\mathrm{C}^2$ such that
\bes
\div\bv=0,\ \ \ \curl\bv=\bo,
\ees
then $\bv$ is harmonic: $\Delta\bv=\bo$.
\begin{proof}
By the definition of curl,
\bes
\curl\bv=\bo\Rightarrow\ \grad\bv-(\grad\bv)^\top=\bo\Rightarrow\div(\grad\bv-(\grad\bv)^\top)=\bo,
\ees
and through Eq.  (\ref{eq:divproperties})$_1$, the definition of laplacian and because by hypothesis $\div\bv=0$, we have
\bes
\div(\grad\bv-(\grad\bv)^\top)=\Delta\bv-\grad(\div\bv)=\Delta\bv.
\ees
\end{proof}
\end{teo}
We state now without proof a  lemma\footnote{In the following, $\partial\Omega$ indicates the boundary of $\Omega$.} that, basically, allows us to transform a volume integral on a domain $\Omega$ to a surface integral on the boundary surface $\partial\Omega$.
\begin{teo}{\bf{(Divergence lemma).}} Let $\bv(p)$ be a vector field of class $\geq\mathrm{C}^1$ on a regular region $\Omega\subset\Eu$. Then,
\bes
\int_{\partial\Omega}\bv\otimes\n\ dA=\int_\Omega\mathrm{grad}\bv\ dV.
\ees
\end{teo}
This lemma is fundamental for proving the three forms of the Gauss theorem, which is of the paramount importance in many fields of mathematical physics.
\begin{teo}{\bf{(Divergence or Gauss theorem).}} Let $\phi,\bv,\L$ be, respectively, a scalar, vector and tensor field of class $\geq\mathrm{C}^1$ on a regular region $\Omega\subset\Eu$ . Then:
\bes
\begin{split}
&\mathrm{i)}\ \ \int_{\partial\Omega}\phi\n\ dA=\int_\Omega\mathrm{grad}\phi\ dV,\medskip\\
&\mathrm{ii)}\ \int_{\partial\Omega}\bv\cdot\n\ dA=\int_\Omega\div\bv\ dV,\medskip\\
&\mathrm{iii)}\int_{\partial\Omega}\L\n\ dA=\int_\Omega\div\L\ dV.
\end{split}
\ees
\begin{proof}
i) $\forall\a=const.\in\Ve$, by the lemma of divergence,
\bes
\int_{\Omega}\grad(\phi\a)dV=\int_{\partial\Omega}\phi\a\otimes\n\ dA=\a\otimes\int_{\partial\Omega}\phi\n\ dA,
\ees
but also, by  ii) of Theorem \ref{teo:propgrad},
\bes
\int_{\Omega}\grad(\phi\a)dV=\int_\Omega(\phi\ \grad\a+\a\otimes\grad\phi)dV=\a\otimes\int_\Omega\grad\phi\ dV,
\ees
whence, by comparison,
\bes
\int_{\partial\Omega}\phi\n\ dA=\int_\Omega\grad\phi\ dV.
\ees
ii) Again by the divergence lemma,
\bes
\tr\int_\Omega\grad\bv\ dV=\tr\int_{\partial\Omega}\bv\otimes\n\ dA=\int_{\partial\Omega}\tr(\bv\otimes\n)dA=\int_{\partial\Omega}\bv\cdot\n\ dA,
\ees
but also
\bes
\tr\int_\Omega\grad\bv\ dV=\int_\Omega \tr(\grad\bv)dV=\int_\Omega\div\bv\ dV,
\ees
so, by comparison,
\bes
\int_{\partial\Omega}\bv\cdot\n\ dA=\int_\Omega\div\bv\ dV.
\ees
iii) $\forall\a=const.\in\Ve$, by the lemma of divergence, property iv) of Theorem \ref{teo:proddiv} and  ii) just proved,
\bes
\int_\Omega\div(\L^\top\a)dV=\int_{\partial\Omega}(\L^\top\a)\cdot\n\ dA=\int_{\partial\Omega}\a\cdot\L\n\ dA=\a\cdot\int_{\partial\Omega}\L\n\ dA,
\ees
but also
\bes
\int_\Omega\div(\L^\top\a)dV=\int_\Omega(\div\L)\cdot\a+\L\cdot\grad\a\ dV=\a\cdot\int_\Omega\div\L\ dV,
\ees
so, once more by comparison,
\bes
\int_{\partial\Omega}\L\n\ dA=\int_\Omega\div\L\ dV.
\ees
\end{proof}
\end{teo}
The following identities follow directly from the Gauss theorem:
\be
\label{eq:propthgauss}
\besp
&\int_{\partial\Omega}\bv\cdot\L\n\ dA=\int_\Omega(\bv\cdot\div\L+\L\cdot\grad\bv)dV,\\
&\int_{\partial\Omega}(\L\n)\otimes\bv\ dA=\int_\Omega((\div\L)\otimes\bv+\L(\grad\bv^\top))dV,\\
&\int_{\partial\Omega}(\bw\cdot\n)\bv\ dA=\int_\Omega(\bv\div\bw+(\grad\bv)\bw)dV.
\end{split}
\ee
A direct consequence of the Gauss theorem is the following result.
\begin{teo}{\bf{(Flux theorem).}} Let $\bv(p)$ be a vector field of class $\geq\mathrm{C}^1$ on an open subset $\mathtt{R}$ of $\Eu$. Then,
\bes
\div\bv=0\iff\int_{\partial\Omega}\bv\cdot\n\ dA=0\ \ \ \forall\Omega\subset\mathtt{R}.
\ees
\begin{proof}
It immediately follows from the  ii) of the Gauss theorem.
\end{proof}
\end{teo}
Another consequence of the Gauss Theorem is the next theorem.
\begin{teo}{\bf{(Curl theorem).}} Let $\bv(p)$ be a vector field of class $\geq\mathrm{C}^1$ on a regular region $\Omega\subset\Eu$; then
\bes
\int_{\partial\Omega}\n\times\bv\ dA=\int_\Omega\curl\bv\ dV.
\ees
\begin{proof} If $\V$ is the axial tensor of $\bv$, by Theorem \ref{teo:curlaxvect} and  iii) of the Gauss theorem, 
\bes
\int_{\partial\Omega}\n\times\bv\ dA=-\int_{\partial\Omega}\bv\times\n\ dA=-\int_{\partial\Omega}\V\n\ dA=-\int_\Omega\div\V\ dV=\int_\Omega\curl\bv\ dV.
\ees
\end{proof}
\end{teo}
The following classical theorems on fields are recalled here without proof.
\begin{teo}{\bf{(Potential theorem).}} Let $\bv(p)$ be a vector field of class $\geq\mathrm{C}^1$ on a simply connected region $\Omega\subset\Eu$. Then,
\bes
\curl\bv=\bo\iff\bv=\mathrm{grad}\phi
\ees
with $\phi(p)$ the {\it potential}, a scalar field of class $\geq\mathrm{C}^2$.
\end{teo}

\begin{figure}[h]
\begin{center}
\includegraphics[scale=.6]{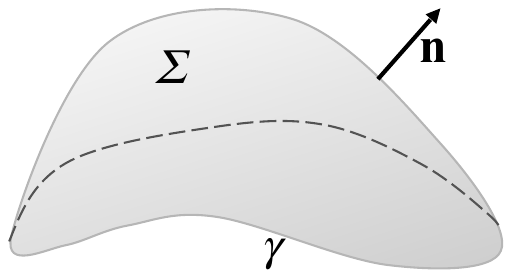}
\caption{Scheme for the Stokes theorem.}
\label{fig:f18}
\end{center}
\end{figure}

\begin{teo}{\bf{(Stokes theorem).}} Let $\bv(p)$ be a vector field of class $\geq\mathrm{C}^1$ on a regular region $\Omega\subset\Eu$,  $\Sigma$ an open surface whose support is the closed line $\gamma$ and  $\n\in\S$ the   normal to $\Sigma$, see Fig. \ref{fig:f18}. Then,
\bes
\oint_\gamma\bv\cdot d\ell=\int_\Sigma\curl\bv\cdot\n\ dA.
\ees
\end{teo}
The parametric equation of $\gamma$ must be chosen in such a way that
\bes
p'(t_1)\times p'(t_2)\cdot\n>0\ \ \ \forall t_2>t_1.
\ees
\begin{teo}{\bf{(Green's formula).}} Let $\phi(p),\psi(p)$ be two scalar fields of class $\geq\mathrm{C}^2$ on a regular region $\Omega\subset\Eu$, with $\n\in\S$ the normal to $\partial\Omega$. Then,
\bes
\int_{\partial\Omega}\left(\psi\frac{d\phi}{d\n}-\phi\frac{d\psi}{d\n}\right)dA=\int_\Omega(\psi\ \Delta\phi-\phi\ \Delta\psi)dV.
\ees
\end{teo}

\section{Differential operators in Cartesian coordinates}
\label{sec:diffopercart}
The Cartesian expression 
of the differential operators can be found without difficulty by applying the properties of such operators shown previously and considering that the vectors of a Cartesian basis are fixed. The final result is\footnote{In what follows, and also in the following sections, $f,\bv,\L$ are, respectively,  scalar, vector and tensor fields.},
\be
\label{eq:cartcoorddiffop}
\besp
&\mathrm{grad}f=f_{,i}\ \ei,\\
&\mathrm{grad}\bv=v_{i,j}\ei\otimes\ej,\\
&\div\bv=v_{i,i},\\
&\div\L=L_{ij,j}\ei,\\
&\Delta f=f_{,ii},\\
&\Delta\bv=\Delta v_i\ei=v_{i,jj}\ei,\\
&\curl\bv=(v_{3,2}-v_{2,3})\eu+(v_{1,3}-v_{3,1})\ed+(v_{2,1}-v_{1,2})\et.
\end{split}
\ee
The so-called {\it operator nabla} $\nabla$,
\be
\label{eq:operatorenabla}
\nabla:=\frac{\partial\cdot}{\partial x_i}\e_i=\frac{\partial\cdot}{\partial x_1}\e_1+\frac{\partial\cdot}{\partial x_2}\e_2+\frac{\partial\cdot}{\partial x_3}\e_3
\ee
is often used to indicate the differential operators:
\bes
\besp
&\mathrm{grad}f=\nabla f,\\
&\div\bv=\nabla\cdot\bv,\\
&\curl\bv=\nabla\times\bv,\\
&\Delta f=\nabla^2f.
\end{split}
\ees

\section {Differential operators in cylindrical coordinates}
\label{sec:cylcoord}
\begin{figure}[h]
\begin{center}
\includegraphics[scale=.7]{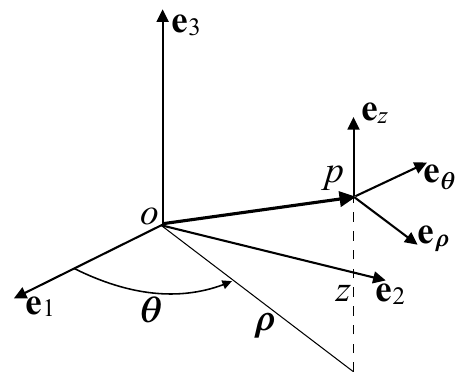}
\caption{Cylindrical coordinates.}
\label{fig:f19}
\end{center}
\end{figure}

The cylindrical coordinates $\rho,\theta,z$ of a point $p$, whose Cartesian coordinates in the (fixed) frame $\{o;\e_1,\e_2,\e_3\}$ are $p=(x_1,x_2,x_3)$, are shown in Fig. \ref{fig:f19}. They are related together by
\be
\label{eq:cylindcoord1}
\besp
&\rho=\sqrt{x_1^2+x_2^2},\\
&\theta=\arctan\dfrac{x_2}{x_1},\\
&z=x_3,
\end{split}
\ee
or conversely
\be
\label{eq:cylindcoord2}
\besp
&x_1=\rho\cos\theta,\\
&x_2=\rho\sin\theta,\\
&x_3=z.
\end{split}
\ee
To notice that $\rho\geq0$ and that the {\it anomaly} $\theta$ is bounded by $0\leq\theta<2\pi$.
A vector $p-o=x_i\ei$, in the cylindrical basis is expressed has
\bes
p-o=\rho\e_\rho+z\ez
\ees 
and the rotation tensor transforming the Cartesian basis $\{\eu,\ed,\et\}$ into the cylindrical one, $\{\e_\rho,\e_\theta,\ez\}$ is
\bes
\Q=\left[\begin{array}{ccc}\cos\theta & -\sin\theta & 0 \\\sin\theta & \cos\theta & 0 \\0 & 0 & 1\end{array}\right],
\ees
so the relations between the vectors of the Cartesian and the cylindrical bases are
\bes
\besp
&\e_\rho=\cos\theta\eu+\sin\theta\ed,\\
&\e_\theta=-\sin\theta\eu+\cos\theta\ed,\\
&\ez=\et,
\end{split}
\ees
and viceversa
\be
\label{eq:vectcartcylind}
\besp
&\eu=\cos\theta\e_\rho-\sin\theta\e_\theta,\\
&\ed=\sin\theta\e_\rho+\cos\theta\e_\theta,\\
&\et=\ez.
\end{split}
\ee
The question is: How can we express the differential operators in the (moving) frame $\{p;\e_\rho,\e_\theta,\e_z\}$? To this end, we can proceed as follows: From Eq.  (\ref{eq:cylindcoord1}) 
\be
\label{eq:derivpartscal}
f_{,i}=f_{,\rho}\frac{\partial\rho}{\partial x_i}+f_{,\theta}\frac{\partial\theta}{\partial x_i}+f_{,z}\frac{\partial z}{\partial x_i}\rightarrow\ \left\{
\besp
&f_{,1}=f_{,\rho}\frac{x_1}{\rho}-f_{,\theta}\frac{x_2}{\rho^2},\\
&f_{,2}=f_{,\rho}\frac{x_2}{\rho}+f_{,\theta}\frac{x_1}{\rho^2},\\
&f_{,3}=f_{,z}.
\end{split}
\right.
\ee
So, by Eqs.  (\ref{eq:cartcoorddiffop})$_1$ and (\ref{eq:vectcartcylind}),
\bes
\grad f=f_i\ei=\left(f_{,\rho}\frac{x_1}{\rho}-f_{,\theta}\frac{x_2}{\rho^2}\right)\hspace{-.8mm}(\cos\theta\e_\rho-\sin\theta\e_\theta)+\left(f_{,\rho}\frac{x_2}{\rho}+f_{,\theta}\frac{x_1}{\rho^2}\right)\hspace{-.8mm}(\sin\theta\e_\rho+\cos\theta\e_\theta)+f_{,z}\ez.
\ees
Finally, by Eq.  (\ref{eq:cylindcoord2}) and through some standard operations, we obtain
\bes
\grad f=f_{,\rho}\e_\rho+\frac{1}{\rho}f_{,\theta}\ \e_\theta+f_{,z}\ez.
\ees
The gradient of a vector field $\bv$ can be obtained in a similar way. If we denote by $\bv_{Cart}$ the vector $\bv$ expressed by its Cartesian components $(v_1,v_2,v_3)$ and by $\bv_{cyl}$ the same vector  expressed through the cylindrical ones, $\{v_\rho,v_\theta,v_z\}$, then, cf. Section \ref{sec:changeofbasis},
\be
\label{eq:cartcompvectcyl}
\bv_{Cart}=\Q\bv_{cyl}\rightarrow\left\{
\besp
&v_1=v_\rho\cos\theta-v_\theta\sin\theta,\\
&v_2=v_\rho\sin\theta+v_\theta\cos\theta,\\
&v_3=v_z.
\end{split}
\right.
\ee
Applying Eq.  (\ref{eq:derivpartscal}) to these components, we get
\be
\label{eq:derivcompvectfield}
\besp
&v_{i,1}=v_{i,\rho}\frac{x_1}{\rho}-v_{i,\theta}\frac{x_2}{\rho^2},\\
&v_{i,2}=v_{i,\rho}\frac{x_2}{\rho}+v_{i,\theta}\frac{x_1}{\rho^2},\\
&v_{i,3}=v_{i,z}.
\end{split}
\ee
Injecting these expressions into Eq.  (\ref{eq:cartcoorddiffop})$_2$ for the $v_{i,j}$s and the (\ref{eq:vectcartcylind})  for the $\ei$s, gives finally\footnote{Though straightforward, the details of the calculations for this formula, as for the following ones, are particularly long and tedious, and for this reason they are omitted here; however, they are a very good exercise for the reader.}
\bes
\besp
\grad\bv&=v_{\rho,\rho}(\e_\rho\otimes\e_\rho)+\dfrac{1}{\rho}(v_{\rho,\theta}-v_\theta)(\e_\rho\otimes\e_\theta)+v_{\rho,z}(\e_\rho\otimes\ez)\\
&+v_{\theta,\rho}(\e_\theta\otimes\e_\rho)+\dfrac{1}{\rho}(v_{\theta,\theta}+v_\rho)(\e_\theta\otimes\e_\theta)+v_{\theta,z}(\e_\theta\otimes\ez)\\
&+v_{z,\rho}(\ez\otimes\e_\rho)+\dfrac{1}{\rho}v_{z,\theta}(\ez\otimes\e_\theta)+v_{z,z}(\ez\otimes\ez),
\end{split}
\ees
or, in matrix form,
\bes
\grad\bv=\left[
\begin{array}{ccc}
v_{\rho,\rho}&\dfrac{1}{\rho}(v_{\rho,\theta}-v_\theta)&v_{\rho,z}\smallskip\\
v_{\theta,\rho}&\dfrac{1}{\rho}(v_{\theta,\theta}+v_\rho)&v_{\theta,z}\smallskip\\
v_{z,\rho}&\dfrac{1}{\rho}v_{z,\theta}&v_{z,z}
\end{array}
\right],
\ees
By the definition of divergence, we get immediately
\be
\label{eq:divvectcyl}
\div\bv=v_{\rho,\rho}+\dfrac{1}{\rho}(v_{\theta,\theta}+v_\rho)+v_{z,z}.
\ee
Now, from Eq.  (\ref{eq:cartcoorddiffop}), we see that $\div\L$ is the vector whose components are the divergence of the rows of the matrix representing $\L$. So, we need first to calculate the Cartesian components of $\L$ as functions of the cylindrical ones, cf. Section \ref{sec:changeofbasis}:
\bes
\L_{Cart}=\Q\L_{cyl}\Q^\top\rightarrow\left\{
\begin{array}{l}
L_{11}=-\sin\theta (L_{\rho \theta } \cos\theta-L_{\theta \theta } \sin\theta)+\cos\theta (L_{\rho \rho } \cos\theta-L_{\theta \rho } \sin\theta),\medskip\\
L_{12}= \cos\theta (L_{\rho \theta } \cos\theta-L_{\theta \theta } \sin\theta)+\sin\theta (L_{\rho \rho } \cos\theta-L_{\theta \rho } \sin\theta),\medskip\\
L_{13}=L_{\rho z} \cos\theta-L_{\theta z} \sin\theta,\medskip\\
L_{21}= -\sin\theta (L_{\theta \theta } \cos\theta+L_{\rho \theta } \sin\theta)+\cos\theta (L_{\theta \rho } \cos\theta+L_{\rho \rho } \sin\theta),\medskip\\
L_{22}= \cos\theta (L_{\theta \theta } \cos\theta+L_{\rho \theta } \sin\theta)+\sin\theta (L_{\theta \rho } \cos\theta+L_{\rho \rho } \sin\theta),\medskip\\
L_{23}=L_{\theta z} \cos\theta+L_{\rho z} \sin\theta,\medskip\\
L_{31}= L_{z\rho } \cos\theta-L_{z\theta } \sin\theta,\medskip\\
L_{32}=L_{z\theta } \cos\theta+L_{z\rho } \sin\theta,\medskip\\
L_{33}=L_{zz}.
\end{array}
\right.
\ees
 Then, applying Eqs.  (\ref{eq:vectcartcylind}) and (\ref{eq:divvectcyl}) in Eq.  (\ref{eq:cartcoorddiffop})$_3$ for the  vectors $\bv_i=(L_{i1},L_{i2},L_{i3}),\ i=1,2,3$, we get, through long but standard passages and after putting $\theta=0$ in order to obtain the components of $\div\L$ in the basis $\{\e_\rho,\e_\theta,\ez\}$,
\bes
\besp
\div\L&=\left(
\dfrac{1}{\rho}((\rho L_{\rho\rho})_{,\rho}+L_{\rho\theta,\theta}-L_{\theta\theta})+L_{\rho z,z}\right)\e_\rho\\
&+\left(L_{\theta\rho,\rho}+\dfrac{1}{\rho}(L_{\theta\theta,\theta}+L_{\rho\theta}+L_{\theta\rho})+L_{\theta z,z}\right)\e_\theta\\
&+\left(\dfrac{1}{\rho}((\rho L_{z\rho})_{,\rho}+L_{z\theta,\theta})+L_{zz,z}\right)\ez.
\end{split}
\ees
To obtain $\Delta f=f,_{ii}$, we need to apply twice Eq.  (\ref{eq:derivpartscal}), which gives
\bes
\besp
f_{,11}&=\left(f_{,\rho}\frac{x_1}{\rho}-f_{,\theta}\frac{x_2}{\rho^2}\right)_{,1}=f_{,\rho 1}\frac{x_1}{\rho}+f_{,\rho}\frac{\rho-x_1\rho_{,1}}{\rho^2}-f_{,\theta 1}\frac{x_2}{\rho^2}+f_{,\theta}\frac{2x_2\rho\rho_{,1}}{\rho^4}\\
&=\left(f_{,\rho\rho}\frac{x_1}{\rho}-f_{,\rho\theta}\frac{x_2}{\rho^2}\right)\frac{x_1}{\rho}+f_{,\rho}\frac{\rho^2-x_1^2}{\rho^3}-\left(f_{,\rho\theta}\frac{x_1}{\rho}-f_{,\theta\theta}\frac{x_2}{\rho^2}\right)\frac{x_2}{\rho^2}+f_{,\theta}\frac{2x_1x_2}{\rho^4}\\
&=f_{,\rho\rho}\cos^2\theta-2f_{,\rho\theta}\frac{\sin\theta\cos\theta}{\rho}+f_{,\rho}\frac{\sin^2\theta}{\rho}+f_{,\theta\theta}\frac{\sin^2\theta}{\rho^2}+2f_{,\theta}\frac{\sin\theta\cos\theta}{\rho^2},
\end{split}
\ees 
\bes
\besp
f_{,22}&=\left(f_{,\rho}\frac{x_2}{\rho}+f_{,\theta}\frac{x_1}{\rho^2}\right)_{,2}=f_{,\rho 2}\frac{x_2}{\rho}+f_{,\rho}\frac{\rho-x_2\rho_{,2}}{\rho^2}+f_{,\theta 2}\frac{x_1}{\rho^2}-f_{,\theta}\frac{2x_1\rho\rho_{,2}}{\rho^4}\\
&=\left(f_{,\rho\rho}\frac{x_2}{\rho}+f_{,\rho\theta}\frac{x_1}{\rho^2}\right)\frac{x_2}{\rho}+f_{,\rho}\frac{\rho^2-x_2^2}{\rho^3}+\left(f_{,\rho\theta}\frac{x_2}{\rho}+f_{,\theta\theta}\frac{x_1}{\rho^2}\right)\frac{x_1}{\rho^2}-f_{,\theta}\frac{2x_1x_2}{\rho^4}\\
&=f_{,\rho\rho}\sin^2\theta+2f_{,\rho\theta}\frac{\sin\theta\cos\theta}{\rho}+f_{,\rho}\frac{\cos^2\theta}{\rho}+f_{,\theta\theta}\frac{\cos^2\theta}{\rho^2}-2f_{,\theta}\frac{\sin\theta\cos\theta}{\rho^2},
\end{split}
\ees 
\bes
f_{,33}=f_{,zz}.
\ees
Then, adding these terms together, we finally have
\bes
\Delta f=\frac{1}{\rho}(\rho f_{,\rho})_{,\rho}+\frac{1}{\rho^2}f_{,\theta\theta}+f_{,zz}.
\ees
The laplacian $\Delta\bv$ of a vector field $\bv$, Eq.  (\ref{eq:cartcoorddiffop})$_6$, can be obtained by following the same steps for  each one of the components in Eq.  (\ref{eq:cartcompvectcyl}), which gives
\bes
\besp
\Delta\bv&=
\left(\dfrac{1}{\rho}(\rho v_{\rho,\rho})_{,\rho}+\dfrac{1}{\rho^2}v_{\rho,\theta\theta}+v_{\rho,zz}-\dfrac{1}{\rho^2}(v_\rho+2v_{\theta,\theta})\right)\e_\rho\\
&+\left(\dfrac{1}{\rho}(\rho v_{\theta,\rho})_{,\rho}+\dfrac{1}{\rho^2}v_{\theta,\theta\theta}+v_{\theta,zz}-\dfrac{1}{\rho^2}(v_\theta-2v_{\rho,\theta})\right)\e_\theta\\
&+\left(\dfrac{1}{\rho}(\rho v_{z,\rho})_{,\rho}+\dfrac{1}{\rho^2}v_{z,\theta\theta}+v_{z,zz}\right)\ez.
\end{split}
\ees
Finally, injecting Eqs.  (\ref{eq:vectcartcylind}), (\ref{eq:cartcompvectcyl}) and (\ref{eq:derivcompvectfield})  into Eq.  (\ref{eq:cartcoorddiffop})$_7$ gives
\be
\label{eq:rotorecoordcylind}
\curl\bv=\left(\dfrac{1}{\rho}v_{z,\theta}-v_{\theta,z}\right)\e_\rho+(v_{\rho,z}-v_{z,\rho})\e_\theta+\left(\dfrac{1}{\rho}((\rho v_\theta)_{,\rho}-v_{\rho,\theta})\right)\ez.
\ee

\section{Differential operators in spherical coordinates}
\label{sec:sphercoord}
\begin{figure}[b]
\begin{center}
\includegraphics[scale=.7]{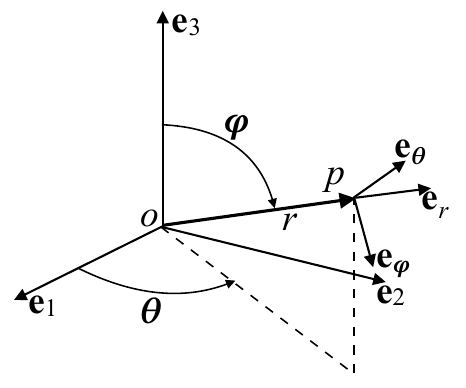}
\caption{Spherical coordinates.}
\label{fig:f20}
\end{center}
\end{figure}

The spherical coordinates $r,\phi,\theta$ of a point $p$, whose Cartesian coordinates in the (fixed) frame $\{o;\e_1,\e_2,\e_3\}$ are $p=(x_1,x_2,x_3)$, are shown in Fig. \ref{fig:f20}. They are related together by
\bes
\besp
&r=\sqrt{x_1^2+x_2^2+x_3^2},\\
&\phi=\arctan\dfrac{\sqrt{x_1^2+x_2^2}}{x_3},\\
&\theta=\arctan\dfrac{x_2}{x_1},
\end{split}
\ees
or conversely
\bes
\besp
&x_1=r\cos\theta\sin\phi,\\
&x_2=r\sin\theta\sin\phi,\\
&x_3=r\cos\phi.
\end{split}
\ees
We note that $r\geq0$ and that the {\it anomaly} $\theta$ is bounded by $0\leq\theta<2\pi$ while the {\it colatitude} $\phi$ by $0\leq\phi\leq\pi$.

The procedure to determine the expression of the differential operators in spherical coordinates, i.e. in  the (moving) frame $\{p;\e_r,\e_\phi,\e_\theta\}$, is identical to that used for the cylindrical coordinates, but the analytical developments are even more complicated and long, so they are omitted here and only the final formulae are given in the following:
\bes
\grad f=f_{,r}\e_r+\dfrac{1}{r}f_{,\phi}\e_\phi+\dfrac{1}{r\sin\phi}f_{,\theta}\e_\theta,
\ees
\bes
\besp
\grad\bv&=
v_{r,r}\e_r\otimes\e_r+\dfrac{1}{r}(v_{r,\phi}-v_\phi)\e_r\otimes\e_\phi+\dfrac{1}{r}\left(\dfrac{1}{\sin\phi}v_{r,\theta}-v_\theta\right)\e_r\otimes\e_\theta\\
&+v_{\phi,r}\e_\phi\otimes\e_r+\dfrac{1}{r}(v_{\phi,\phi}+v_r)\e_\phi\otimes\e_\phi+\dfrac{1}{r}\left(\dfrac{1}{\sin\phi}v_{\phi,\theta}-v_\theta\cot\phi\right)\e_\phi\otimes\e_\theta\\
&+v_{\theta,r}\e_\theta\otimes\e_r+\dfrac{1}{r}v_{\theta,\phi}\e_\theta\otimes\e_\phi+\dfrac{1}{r}\left(\dfrac{1}{\sin\phi}v_{\theta,\theta}+v_r+v_\phi\cot\phi\right)\e_\theta\otimes\e_\theta,
\end{split}
\ees
or, in matrix form,
\bes
\grad\bv=\left[
\begin{array}{ccc}
v_{r,r}&\dfrac{1}{r}(v_{r,\phi}-v_\phi)&\dfrac{1}{r}\left(\dfrac{1}{\sin\phi}v_{r,\theta}-v_\theta\right)\\
v_{\phi,r}&\dfrac{1}{r}(v_{\phi,\phi}+v_r)&\dfrac{1}{r}\left(\dfrac{1}{\sin\phi}v_{\phi,\theta}-v_\theta\cot\phi\right)\\
v_{\theta,r}&\dfrac{1}{r}v_{\theta,\phi}&\dfrac{1}{r}\left(\dfrac{1}{\sin\phi}v_{\theta,\theta}+v_r+v_\phi\cot\phi\right)
\end{array}
\right],
\ees
\bes
\div\bv=\frac{1}{r^2}(r^2v_r)_{,r}+\frac{1}{r\sin\phi}((v_\phi\sin\phi)_{,\phi}+v_{\theta,\theta}),
\ees
\bes
\besp
\div\L&=
\left(\dfrac{1}{r^2}(r^2L_{rr})_{,r}+\dfrac{1}{r}L_{r\phi,\phi}+\dfrac{1}{r\sin\phi}L_{r\theta,\theta}-\dfrac{L_{\phi\phi}+L_{\theta\theta}}{r}+\dfrac{\cot\phi}{r}L_{r\phi}\right)\e_r\\
&+\left(\dfrac{1}{r^2}(r^2L_{\phi r})_{,r}+\dfrac{1}{r}L_{\phi\phi,\phi}+\dfrac{1}{r\sin\phi}L_{\phi\theta,\theta}+\dfrac{1}{r}L_{r\phi}+\dfrac{\cot\phi}{r}(L_{\phi\phi}-L_{\theta\theta})\right)\e_\phi\\
&+\left(\dfrac{1}{r^2}(r^2L_{\theta r})_{,r}+\dfrac{1}{r}L_{\theta\phi,\phi}+\dfrac{1}{r\sin\phi}L_{\theta\theta,\theta}+\dfrac{1}{r}L_{r\theta}+\dfrac{\cot\phi}{r}(L_{\phi\theta}+L_{\theta\phi})\right)\e_\theta,
\end{split}
\ees
\bes
\Delta f=\frac{1}{r^2}(r^2f_{,r})_{,r}+\frac{1}{r^2\sin\phi}\left(\frac{f_{,\theta\theta}}{\sin\phi}+(f_{,\phi}\sin\phi)_{,\phi}\right),
\ees
\bes
\besp
\Delta\bv\hspace{-1mm}&=\hspace{-1mm}
\left(v_{r,rr}+\dfrac{2v_{r,r}}{r}+\dfrac{v_{r,\phi\phi}-2v_{\phi,\phi}}{r^2}+\dfrac{v_{r,\phi}-2v_\phi}{r^2\tan\phi}+\dfrac{1}{r^2\sin\phi}\left(\dfrac{v_{r,\theta\theta}}{\sin\phi}-2v_{\theta,\theta}\right)-\dfrac{2v_r}{r^2}\right)\e_r\\
&+\left(v_{\phi,rr}\hspace{-1mm}+\hspace{-1mm}\dfrac{2v_{\phi,r}}{r}\hspace{-1mm}+\hspace{-1mm}\dfrac{v_{\phi,\phi\phi}\hspace{-1mm}+\hspace{-1mm}2v_{r,\phi}}{r^2}\hspace{-1mm}+\hspace{-1mm}\dfrac{v_{\phi,\phi}\hspace{-1mm}-\hspace{-1mm}v_\phi\cot\phi}{r^2\tan\phi}\hspace{-1mm}+\hspace{-1mm}\dfrac{1}{r^2\sin^2\phi}\left(v_{\phi,\theta\theta}\hspace{-1mm}-\hspace{-1mm}2v_{\theta,\theta}\cos\phi\right)\hspace{-1mm}-\hspace{-1mm}\dfrac{v_\phi}{r^2}\right)\e_\phi\\
&+\left(v_{\theta,rr}\hspace{-1mm}+\hspace{-1mm}\dfrac{2v_{\theta,r}}{r}\hspace{-1mm}+\hspace{-1mm}\dfrac{v_{\theta,\phi\phi}}{r^2}\hspace{-1mm}+\hspace{-1mm}\left(v_{\theta,\phi}\hspace{-1mm}+\hspace{-1mm}\dfrac{2v_{\phi,\theta}}{\sin\phi}\right)\dfrac{1}{r^2\tan\phi}\hspace{-1mm}+\hspace{-1mm}\dfrac{1}{r^2\sin\phi}\left(\dfrac{v_{\theta,\theta\theta}}{\sin\phi}\hspace{-1mm}+\hspace{-1mm}2v_{r,\theta}\right)\hspace{-1mm}-\hspace{-1mm}\dfrac{v_\theta}{r^2\sin^2\phi}\right)\e_\theta,
\end{split}
\ees
\bes
\besp
\curl\bv&=\left(\frac{1}{r\sin\phi}((v_{\theta}\sin\phi)_{,\phi}-v_{\phi,\theta})\right)\e_r\\
&+\left(\frac{1}{r\sin\phi}v_{r,\theta}-\frac{1}{r}(rv_{\theta})_{,r}\right)\e_\phi\\
&+\left(\frac{1}{r}((rv_{\phi})_{,r}-v_{r,\phi})\right)\e_\theta.
\end{split}
\ees

\section{Exercises}
\begin{enumerate}
\item Prove the relations of Eq.  (\ref{eq:propgrad}).
\item Prove the properties of the divergence in Eq.  (\ref{eq:divproperties}).
\item Prove the properties of the curl in Eq.  (\ref{eq:propcurl}).
\item Prove the identities in Eq.  (\ref{eq:propthgauss}).
\item Prove that 
\bes
\frac{d\phi}{d\n}=\grad\phi\cdot\n,\ \ \ \frac{d\bv}{d\n}=\grad\bv\ \n\ \ \ \forall\n\in\S.
\ees
\item Prove the results of Eq.  (\ref{eq:cartcoorddiffop}).
\item Consider a rigid body $\mathtt{B}$, and a point $p_0\in\mathtt{B}$. From the kinematics of rigid bodies, we now that the velocity of another point $p\in\mathtt{B}$ is given by
\bes
\bv(p)=\bv(p_0)+\bom\times(p-p_0),
\ees
with $\bom$ the angular velocity. Prove that
\bes
\bom=\frac{1}{2}\curl\bv,\ \ \ \div\bv=0.
\ees
\item In the infinitesimal theory of strain, a deformation is {\it isochoric} when $\div\bu=0,\ \bu$ being the corresponding displacement vector. Determine which, among the following ones, are locally or globally isochoric deformations:
\begin{enumerate}[i.]
\item $\bu=\alpha(x_1,x_2,x_3),\ \ \alpha\in\R,\ |\alpha|\ll1;$
\item $\bu=\beta(x_2+x_3,x_1+x_3,x_1+x_2),\ \ \beta\in\R,\ |\beta|\ll1;$
\item $\bu=\gamma(x_1x_2,x_2x_3,x_3x_1),\ \ \gamma\in\R,\ |\gamma|\ll1;$
\item $\bu=\delta(\sin x_1,-\cos x_2,\sin x_3),\ \ \delta\in\R,\ |\delta|\ll1.$
\end{enumerate}
\item In fluid mechanics, the condition $\div\bv=0$, with $\bv$ being the velocity field, characterizes {\it incompressible flows}. Verify that  the following velocity fields, given in cylindrical coordinates, correspond to incompressible flows ($\alpha\in\R$):
\begin{enumerate}[i.]
\item source or sink: $\bv=\dfrac{\alpha}{\rho}\e_\rho$;
\item vortex: $\bv=\dfrac{\alpha}{\rho}\e_\theta$;
\item doublet: $\bv=\dfrac{\alpha}{\rho^2}(\cos\theta\e_\rho+\sin\theta\e_\theta)$.
\end{enumerate}
\item A flow with $\curl\bv=\bo$ is said to be {\it irrotational}; check that the flows in the previous exercise are irrotational.

\end{enumerate}

\chapter{Curvilinear coordinates}
\label{ch:6}
\section{Introduction}
All the developments  in the previous chapters are intended for the case where algebraic and differential  operators are expressed in a Cartesian frame, i.e. with {\it rectangular coordinates}. The points of $\Eu$ are thus referred to a system of coordinates taken along straight lines that are mutually orthogonal and with the same unit along each one of the directions of the frame. Though this is a very important and common case, it is not the only possibility and in many cases  non rectangular coordinate frames are used or arise in the mathematical developments (a typical example is that of the geometry of surfaces, see Chapter \ref{ch:7}). A non rectangular coordinate frame is a frame where coordinates can be taken along non-orthogonal directions, or along some lines that intersect at right angles but that are not straight lines, or even when both of these cases occur. This situation is often denoted in the literature  as that of {\it curvilinear coordinates}; the transformations to be done to algebraic and differential operators in the case of curvilinear coordinates is the topic of this chapter. 

\section{Curvilinear coordinates, metric tensor}
\begin{figure}[ht]
	\begin{center}
         \includegraphics[scale=.6]{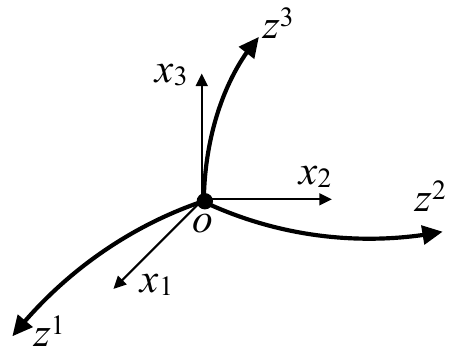}
	\caption{Cartesian and curvilinear coordinates.}
	\label{fig:21}
	\end{center}
\end{figure}

Let us consider an arbitrary origin $o$ of $\Eu$ and an orthonormal basis $e=\{\eu,\ed,\et\}$ of $\Ve$; we  indicate the coordinates of a point $p\in\Eu$ with respect to the frame $\Rep=\{o;\eu,\ed,\et\}$ by  $x_k: p=(x_1,x_2,x_3)$. Then, we also consider  another set of coordinate lines  for $\Eu$, where the position of a point $p\in\Eu$ with respect to the same arbitrary origin $o$ of $\Eu$ is now determined by a set of three numbers $z^j:p=\{z^1,z^2,z^3\}$. Nothing is {\it a-priori} required of coordinates $z^j$, namely they do not need to be a set of Cartesian coordinates, i.e. referring to an orthonormal basis of $\Ve$. In principle, the coordinates $z^j$ can be taken along non straight lines, that do not need to be mutually orthogonal at $o$ and also with different units along each line. That is why we call the $z^j$s  {\it curvilinear coordinates}, see Fig. \ref{fig:21}.
 Any point $p\in\Eu$ can be identified by either set of coordinates; mathematically, this means that there must be an isomorphism between the $x_k$s and the $z^j$s, i.e. invertible relations of the kind
 \be
 \label{eq:cc1}
\begin{array}{l}
 z^j=z^j(x_1,x_2,x_3)=z^j(x_k),\medskip\\ 
 x_k=x_k(z^1,z^2,z^3)=x_k(z^j),
 \end{array}
 \ \ \ \ \forall j,k=1,2,3
 \ee
 exist between the two sets of coordinates. The distance between two points $p,q\in\Eu$ is\footnote{The distance between two points $p$ and $q$ is still defined as the Euclidean norm of $(p-q)$, i.e., it is independent of the set of coordinates.}
 \bes
 s=\sqrt{(p-q)\cdot(p-q)}=\sqrt{(x_k^p-x_k^q)(x_k^p-x_k^q)}
 \ees
 but this is no longer true for curvilinear coordinates:
 \bes
 s\neq\sqrt{({z^j}^p-{z^j}^q)({z^j}^p-{z^j} ^q)}.
 \ees
 However, if $p\rightarrow q$, we can define
 \bes
 dx_k=x_k^p-x_k^q,\ \ \ dz^j={z^j}^p-{z^j}^q,
 \ees
 so using Eq. (\ref{eq:cc1})$_2$
 \be
 \label{eq:cc2}
 dx_k=\frac{\partial x_k}{\partial z^j}dz^j.
 \ee
 
\begin{figure}[ht]
	\begin{center}
         \includegraphics[scale=.7]{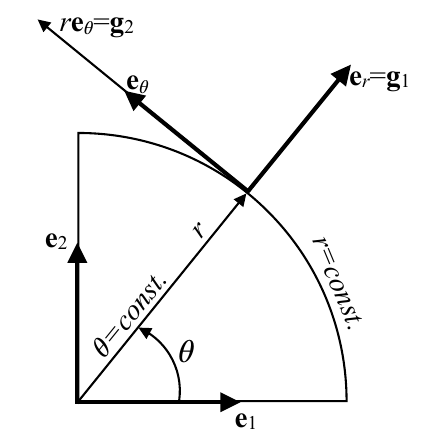}
	\caption{Tangent vectors to the polar coordinates lines.}
	\label{fig:23}
	\end{center}
\end{figure}

The (infinitesimal) distance between $p$ and $q$ will then be
\bes
ds=\sqrt{dx_kdx_k}=\sqrt{\frac{\partial x_k}{\partial z^j}\frac{\partial x_k}{\partial z^l}dz^jdz^l}=\sqrt{g_{jl}dz^jdz^l},
\ees
 where
 \be
 \label{eq:gcov1}
 g_{jl}=g_{lj}=\frac{\partial x_k}{\partial z^j}\frac{\partial x_k}{\partial z^l}
 \ee
 are the  {\it covariant\footnote{The notion of co- and contra-variant components will be detailed in the next section.} components of the metric tensor}\footnote{As usually done in the literature, we indicate the metric tensor by $\g$, i.e. using a lowercase letter, though it is a 2nd-rank tensor, not a vector.} $\g\in Sym(\Ve)$. We note that, as $\g$ defines a positive quadratic form (the square length of a vector), it is a positive definite symmetric tensor, so 
 \be
 \label{eq:signdetg}
 \det\g>0.
 \ee
 Coming back to the vector notation, from Eq. (\ref{eq:cc2}) we get\footnote{The differential $dx$ is a vector because it is the difference of two infinitely close points; that is why it is not needed to denote it in bold letters.\\}
 \bes
 dx=dx_i\e_i=\frac{\partial x_i}{\partial z^k}dz^k\e_i;
 \ees
 introducing the vector $\g_k$,
 \be
 \label{eq:vtang1}
 \g_k:=\frac{\partial x_i}{\partial z^k}\e_i,
 \ee
 we can write
 \bes
 dx=dz^k\g_k.
 \ees
 We see hence that a vector $dx$ can be expressed as a linear combination of the vectors $\g_k$; these  form therefore a basis, called the {\it local basis}. Generally speaking,  $\g_k\notin\S$ and it is clearly tangent to the lines $z^j=const$. This can be seen in Fig. \ref{fig:22} for a two-dimensional case:
 \begin{figure}[h]
	\begin{center}
         \includegraphics[scale=.6]{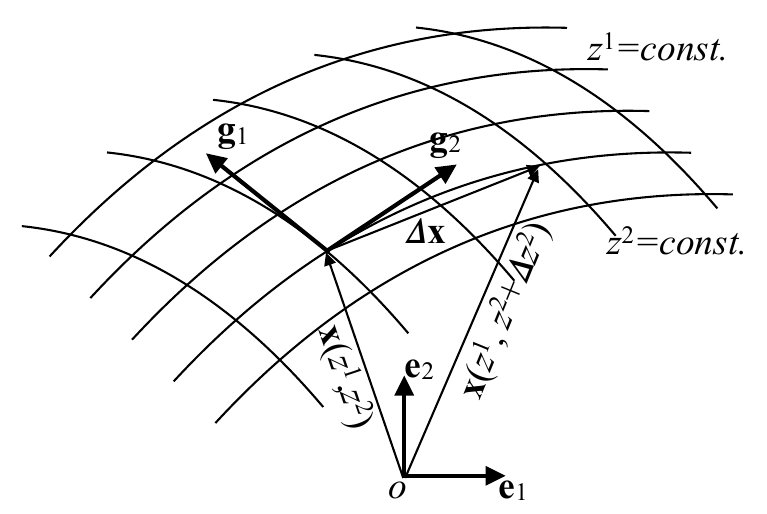}
	\caption{Tangent vectors to the curvilinear coordinates lines.}
	\label{fig:22}
	\end{center}
\end{figure}
\bes
dx=\lim_{\Delta \x\rightarrow0}\Delta \x=\lim_{\Delta \x\rightarrow0}\frac{x_i(z^1,z^2+\Delta z^2)-x_i(z^1,z^2)}{\Delta z^2}\Delta z^2\e_i=\frac{\partial x_i}{\partial z^2}\e_idz^2=\g_2dz^2.
\ees
 Then
 \be
 \label{eq:meaningghk}
 \g_k\cdot\g_l=\frac{\partial x_i}{\partial z^k}\e_i\cdot\frac{\partial x_j}{\partial z^l}\e_j=\frac{\partial x_i}{\partial z^k}\frac{\partial x_j}{\partial z^l}\delta_{ij}=g_{kl},
 \ee
 i.e. the components of the metric tensor $\g$ are the scalar products of the  tangent vectors $\g_k$s. If the curvilinear coordinates are {\it orthogonal}, i.e. if $\g_h\cdot\g_k=0\ \forall h,k=1,2,3, h\neq k$, then $\g$ is diagonal. If, in addition, $\g_k\in\S\ \forall k=1,2,3$, then $\g=\I$: It is the case of Cartesian coordinates. 
 As an example, let us consider the case of polar coordinates,
 \bes
 \left\{\begin{array}{l} x_1=r\ \cos\theta,\\ x_2=r\ \sin\theta,\end{array}\right.\ \ \ 
 \left\{\begin{array}{l} z^1=r=\sqrt{x_1^2+x_2^2},\\ z^2=\theta=\arctan\dfrac{x_2}{x_1}.\end{array}\right.
 \ees
 Hence, see Fig. \ref{fig:23},
 \bes
 \besp
 &\g_1=\frac{\partial x_1}{\partial z^1}\e_1+\frac{\partial x_2}{\partial z^1}\e_2=\cos\theta\e_1+\sin\theta\e_2=\e_r,\\
 &\g_2=\frac{\partial x_1}{\partial z^2}\e_1+\frac{\partial x_2}{\partial z^2}\e_2=-r\sin\theta\e_1+r\cos\theta\e_2=r\e_\theta.
 \end{split}
 \ees
 We remark that $|\g_1|=1$ but $|\g_2|\neq1$ and it is variable with the position.
 
 \section{Co- and contra-variant components}
\begin{figure}[ht]
	\begin{center}
         \includegraphics[scale=.8]{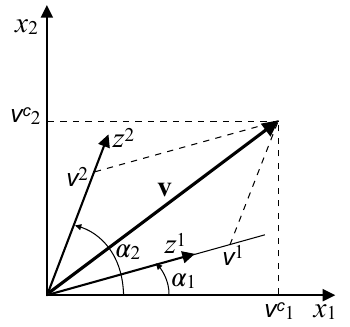}\includegraphics[scale=.8]{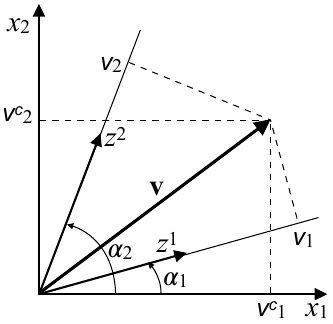}
	\caption{Contravariant (left) and covariant (right) components of a vector in a plane.}
	\label{fig:24}
	\end{center}
\end{figure}

\begin{figure}[ht]
	\begin{center}
       \includegraphics[scale=.8]{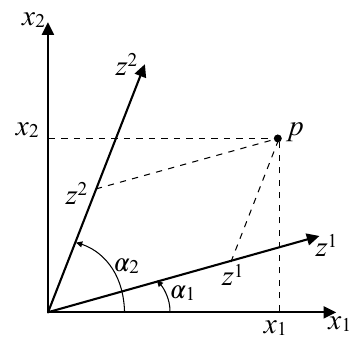}
	\caption{Relation between Cartesian and contravariant components.}
	\label{fig:26}
	\end{center}
\end{figure}

A geometrical way to introduce the concept of {\it covariant} and {\it contravariant} components is to consider how to represent a vector $\bv$ in the $z-$system. 
 There are basically two ways, cf. Fig. \ref{fig:24}, referred, for the sake of simplicity, to a planar case:
 \begin{enumerate}[i.]
 \item {\it contravariant components}: $\bv$ is projected parallel to $z^1$ and $z^2$; they are indicated by superscripts: $\bv=(v^1,v^2,v^3)$;
 \item {\it covariant components}: $\bv$ is projected perpendicularly to $z^1$ and $z^2$; they are indicated by subscripts: $\bv=(v_1,v_2,v_3)$;
 \end{enumerate}
Still referring to the planar case in Fig. \ref{fig:24},  if the Cartesian components\footnote{In the following, we use the superscript $c$ to indicate a Cartesian component: $v^c_i$ is the $i$-th Cartesian component of $\bv\in\Ve$ and $L_{ij}^c$ the $ij$-th Cartesian component of $\L\in Lin(\Ve)$.} of $\bv$ are $\bv=(v^c_1,v^c_2)$, we get
\be
\label{eq:ccont1}
\left\{\begin{array}{l}v^1=h(v_1^c\sin\alpha_2-v_2^c\cos\alpha_2),\\v^2=h(-v_1^c\sin\alpha_1+ v_2^c\cos\alpha_1),\end{array}\right.\ \ \
\left\{\begin{array}{l}v_1=v_1^c\cos\alpha_1+v_2^c\sin\alpha_1,\\v_2=v_1^c\cos\alpha_2+ v_2^c\sin\alpha_2,\end{array}\right.
\ee
and conversely
\be
\label{eq:ccov1}
\left\{\begin{array}{l}v^c_1=v^1\cos\alpha_1+v^2\cos\alpha_2,\\v^c_2=v^1\sin\alpha_1+ v^2\sin\alpha_2,\end{array}\right.\ \ \
\left\{\begin{array}{l}v_1^c=h(v_1\sin\alpha_2-v_2\sin\alpha_1),\\v^c_2=h(-v_1\cos\alpha_2+ v_2\cos\alpha_1),\end{array}\right.
\ee
 with
 \bes
 h=\frac{1}{\sin(\alpha_2-\alpha_1)}.
 \ees
 It is apparent that the Cartesian coordinates are at the same time co- and contra-variant. Still on a planar scheme, we can see how to pass from a system of coordinates to another one, cf. Fig. \ref{fig:26}.
For a point $p$ the Cartesian coordinates $(x_1,x_2)$ are related to the contravariant ones by
\bes
\besp
&x_1=z^1\cos\alpha_1+z^2\cos\alpha_2,\\
&x_2=z^1\sin\alpha_1+z^2\sin\alpha_2,
\end{split}
\ees
 and, conversely,
 \bes
 \besp
& z^1=h(x_1\sin\alpha_2-x_2\cos\alpha_2),\\
&  z^2=h(-x_1\sin\alpha_1+x_2\cos\alpha_1).
 \end{split}
 \ees
 So, differentiating, we get
 \bes
 \besp
& \frac{\partial x_1}{\partial z^1}=\cos\alpha_1,\  \frac{\partial x_1}{\partial z^2}=\cos\alpha_2,\\
& \frac{\partial x_2}{\partial z^1}=\sin\alpha_1,\  \frac{\partial x_2}{\partial z^2}=\sin\alpha_2,\\
 \end{split}
 \ees
 and
 \bes
 \besp
&  \frac{\partial z^1}{\partial x_1}=h\sin\alpha_2,\  \frac{\partial z^1}{\partial x_2}=-h\cos\alpha_2,\\
& \frac{\partial z^2}{\partial x_1}=-h\sin\alpha_1,\  \frac{\partial z^2}{\partial x_2}=h\cos\alpha_1.\\
 \end{split}
 \ees
 Injecting these expressions into Eqs. (\ref{eq:ccont1}) and (\ref{eq:ccov1}) gives
 \be
 \label{eq:cc10}
 \besp
 &v_1=v_1^c\frac{\partial x_1}{\partial z^1}+v_2^c\frac{\partial x_2}{\partial z^1},\\
 &v_2=v_1^c\frac{\partial x_1}{\partial z^2}+v_2^c\frac{\partial x_2}{\partial z^2},
 \end{split}
 \ \rightarrow\ v_i=\frac{\partial x_k}{\partial z^i}v_k^c
 \ee
 and
 \be
 \label{eq:cc11}
 \besp
 &v^1=v_1^c\frac{\partial z^1}{\partial x_1}+v_2^c\frac{\partial z^1}{\partial x_2},\\
 &v^2=v_1^c\frac{\partial z^2}{\partial x_1}+v_2^c\frac{\partial z^2}{\partial x_2},
 \end{split}
 \ \rightarrow\ v^i=\frac{\partial z^i}{\partial x_k}v_k^c.
 \ee
Now, if we calculate
\bes
g_{hi}v^i=g_{hi}\frac{\partial z^i}{\partial x_k}v_k^c,
\ees
from Eq. (\ref{eq:gcov1}) and by the chain rule\footnote{ The reader can easily see that, in practice, the chain rule allows us to handle the derivatives as fractions.} we get
\bes
\besp
g_{hi}v^i&= \frac{\partial x_j}{\partial z^h}\frac{\partial x_j}{\partial z^i}   \frac{\partial z^i}{\partial x_k}v_k^c\\
&=\frac{\partial x_j}{\partial z^h}\frac{\partial x_j}{\partial x_k}v_k^c=\frac{\partial x_j}{\partial z^h}\delta_{jk}v_k^c=\frac{\partial x_k}{\partial z^h}v_k^c=v_h,
\end{split}
\ees
 i.e. we obtain the  rule of {\it lowering of the indices} for passing from contravariant to covariant components:
 \bes
 v_h=g_{hi}v^i.
 \ees
  Introducing the {\it inverse}\footnote{To prove that the contravariant components $g^{pq}$ are the inverse of the covariant ones, $g_{pq}$, is direct: 
 \bes g^{pq}g_{pq}=\frac{\partial z^p}{\partial x_k}\frac{\partial z^q}{\partial x_k}\frac{\partial x_j}{\partial z^p}\frac{\partial x_j}{\partial z^q}=\delta_{jk}\delta_{jk}=1.
 \ees 
 } {\it to $g_{hi}$} as
 \be
  \label{eq:gcont1}
 g^{hi}=\frac{\partial z^h}{\partial x_k}\frac{\partial z^i}{\partial x_k},
 \ee
 we get, still using the chain rule,
 \bes
\besp
g^{hi}v_i&=g^{hi}\frac{\partial x_k}{\partial z^i}v_k^c=
 \frac{\partial z^h}{\partial x_j} \frac{\partial z^i}{\partial x_j}   \frac{\partial x_k}{\partial z^i}v_k^c\\
&=\frac{\partial z^h}{\partial x_j}\frac{\partial x_k}{\partial x_j}v_k^c=\frac{\partial z^h}{\partial x_j}\delta_{jk}v_k^c=\frac{\partial z^h}{\partial x_k}v_k^c=v^h,
\end{split}
\ees
which is  the rule of {\it raising of the indices} for passing from covariant to contravariant components:
 \bes
v^h= g^{hi}v_i.
 \ees 
 Again applying the chain rule, by Eq. (\ref{eq:cc10}) we get
 \bes
 \frac{\partial z^i}{\partial x_l}v_i= \frac{\partial z^i}{\partial x_l}\frac{\partial x_k}{\partial z^i}v_k^c=\frac{\partial x_k}{\partial x_l}v_k^c=\delta_{kl}v_k^c,
 \ees
 i.e.
 \be
 \label{eq:cov-cart}
 v_k^c= \frac{\partial z^i}{\partial x_k}v_i,
 \ee
 which is the converse of Eq. (\ref{eq:cc10}). In a similar way, we get the converse of Eq. (\ref{eq:cc11}):
 \be
  \label{eq:cont-cart}
v_k^c =\frac{\partial x_k}{\partial z^i}v^i.
 \ee
 Let us now calculate the norm $v$ of a vector $\bv$; starting from the Cartesian components and using the last two results, 
 \bes
 v=\sqrt{\bv\cdot\bv}=\sqrt{v_k^cv_k^c}=\sqrt{ \frac{\partial z^i}{\partial x_k}v_i  \frac{\partial z^j}{\partial x_k}v_j}=
\sqrt{ \frac{\partial z^i}{\partial x_k} \frac{\partial z^j}{\partial x_k}v_i v_j}=\sqrt{g^{ij}v_iv_j},
 \ees
 or also
 \bes
v=\sqrt{\bv\cdot\bv}=\sqrt{v_k^cv_k^c}=\sqrt{\frac{\partial x_k}{\partial z^i}v^i \frac{\partial x_k}{\partial z^j}v^j}=
 \sqrt{\frac{\partial x_k}{\partial z^i} \frac{\partial x_k}{\partial z^j}v^iv^j}=\sqrt{g_{ij}v^iv^j}
 \ees
 and even
  \bes
 v=\sqrt{\bv\cdot\bv}=\sqrt{v_k^cv_k^c}=\sqrt{\frac{\partial z^i}{\partial x_k}v_i \frac{\partial x_k}{\partial z^j}v^j}=
\sqrt{\frac{\partial z^i}{\partial x_k} \frac{\partial x_k}{\partial z^j}v_iv^j}=\sqrt{\delta^i_{\ j}v_iv^j}=\sqrt{v_iv^i}.
 \ees
 
 Through Eq.  (\ref{eq:cont-cart}) and by the definition of the tangent vectors to the lines of curvilinear coordinates, Eq. (\ref{eq:vtang1}), for a vector $\bv$  we get
 \bes
 \bv=v_i^c\e_i=v^k\frac{\partial x_i}{\partial z^k}\e_i=v^k\g_k.
 \ees
 We see hence that the contravariant components are actually the components of $\bv$ in the local basis, composed by vectors $\g_k$s,   tangent to the lines of curvilinear coordinates. In a similar manner, if we introduce the {\it dual basis} whose vectors $\g^k$ are defined as
 \be
 \label{eq:dualbasis}
 \g^k:=\frac{\partial z^k}{\partial x_i}\e_i
 \ee
 and proceeding in the same way, we obtain that 
 \bes
 \bv=v_i^c\e_i=v_k\frac{\partial z^k}{\partial x_i}\e_i=v_k\g^k,
 \ees
i.e. the covariant components are actually the components of $\bv$ in the dual basis. Finally, for a vector, we have, alternatively,
\be
\label{eq:compvectccc}
\bv=v_i^c\e_i=v^k\g_k=v_k\g^k.
\ee
 
Just as for the $\g_k$s, we have
\bes
\g^h\cdot\g^k=\left(\frac{\partial z^h}{\partial x_i}\e_i\right)\cdot\left(\frac{\partial z^k}{\partial x_j}\e_j\right)=\frac{\partial z^h}{\partial x_i}\frac{\partial z^k}{\partial x_j}\delta_{ij}=\frac{\partial z^h}{\partial x_i}\frac{\partial z^k}{\partial x_i}=g^{hk};
\ees
 moreover
\bes
\g^h\cdot\g_k=\left(\frac{\partial z^h}{\partial x_i}\e_i\right)\cdot\left(\frac{\partial x_j}{\partial z^k}\e_j\right)=\frac{\partial z^h}{\partial x_i}\frac{\partial x_j}{\partial z^k}\delta_{ij}=\frac{\partial z^h}{\partial x_i}\frac{\partial x_i}{\partial z^k}=\frac{\partial z^h}{\partial z^k}=\delta^h_{\ k},
\ees
and, by the symmetry of the scalar product,
\bes
\delta_h^{\ k}:=\g_h\cdot\g^k=\g^k\cdot\g_h=\delta^k_{\ h}.
\ees
The last equations defines the {\it orthogonality conditions} for the $\g$-vectors.
 Using these results and Eq. (\ref{eq:compvectccc}) we also have 
 \bes
 \besp
 &v^k=\delta_{\ h}^kv^h=\g^k\cdot v^h\g_h=\g^k\cdot\bv=\g^k\cdot v_h\g^h=g^{kh}v_h,\\
 &v_k=\delta_k^{\ h}v_h=\g_k\cdot v_h\g^h=\g_k\cdot\bv=\g_k\cdot v^h\g_h=g_{kh}v^h,
 \end{split}
 \ees
 thus finding again the rules of raising and lowering of the indices.

  What was done for  vectors can be transposed, using a similar approach, to tensors. In particular, for a second-rank tensor $\L$, we get
 \be
 \label{eq:transftenscartcov}
 \besp
& L^{ij}=\frac{\partial z^i}{\partial x_h}\frac{\partial z^j}{\partial x_k}L_{hk}^c,\\
& L_{ij}=\frac{\partial x_h}{\partial z^i}\frac{\partial x_k}{\partial z^j}L_{hk}^c,
 \end{split}
 \ee
 for the contravariant and covariant components, respectively, while we can also introduce the {\it mixed components}
  \be
  \label{eq:mixedtens1}
 \besp
& L^i_{\ j}=\frac{\partial z^i}{\partial x_h}\frac{\partial x_k}{\partial z^j}L_{hk}^c,\\
& L_i^{\ j}=\frac{\partial x_h}{\partial z^i}\frac{\partial z^j}{\partial x_k}L_{hk}^c.
 \end{split}
 \ee
 Conversely,
 \be
 \label{eq:touttenseurscart}
 \besp
& L_{hk}^c=\frac{\partial x_h}{\partial z^i}\frac{\partial x_k}{\partial z^j}L^{ij},\\
& L_{hk}^c=\frac{\partial z^i}{\partial x_h}\frac{\partial z^j}{\partial x_k}L_{ij},\\
& L_{hk}^c=\frac{\partial x_h}{\partial z^i}\frac{\partial z^j}{\partial x_k}L^i_{\ j},\\
& L_{hk}^c=\frac{\partial z^i}{\partial x_h}\frac{\partial x_k}{\partial z^j}L_i^{\ j}.
 \end{split}
 \ee
 Also for $\L$, the rule of lowering or raising the indices is valid:
 \be
 \label{eq:raisingloweringL}
 L^{ij}=g^{ih}g^{jk}L_{hk},\ \ \ L_{ij}=g_{ih}g_{jk}L^{hk}.
 \ee
 From Eq. (\ref{eq:touttenseurscart}) and by the same definitions of $g_{ij}$, eq.(\ref{eq:gcov1}), and $g^{ij}$, Eq. (\ref{eq:gcont1}), we get
 \bes
\L=L_{ij}^c\e_i\otimes\e_j=\frac{\partial x_i}{\partial z^h}\frac{\partial x_j}{\partial z^k}L^{hk}\e_i\otimes\e_j=L^{hk}\g_h\otimes\g_k
\ees
and
\bes
\L=L_{ij}^c\e_i\otimes\e_j=\frac{\partial z^h}{\partial x_i}\frac{\partial z^k}{\partial x_j}L_{hk}\e_i\otimes\e_j=L_{hk}\g^h\otimes\g^k.
\ees
 In a similar manner, the tensor mixed components are also found:
 \bes
\L=L_{ij}^c\e_i\otimes\e_j=\frac{\partial x_i}{\partial z^h}\frac{\partial z^k}{\partial x_j}L^h_{\ k}\e_i\otimes\e_j=L^h_{\ k}\g_h\otimes\g^k
\ees
and
 \bes
\L=L_{ij}^c\e_i\otimes\e_j=\frac{\partial z^k}{\partial x_j}\frac{\partial x_i}{\partial z^h}L_h^{\ k}\e_i\otimes\e_j=L_h^{\ k}\g^h\otimes\g_k.
\ees
 We see hence that a second-rank tensor can be given with four different combinations of  coordinates; even more complex is the case of higher order tensors, which will not be treated here.
 
Still by eqs.(\ref{eq:gcov1}) and (\ref{eq:gcont1}) and applying the chain rule to $\delta^i_j=\dfrac{\partial z^i}{\partial z^j}$, we get
 \be
 \label{eq:gedelta}
 \besp
 &g_{ij}=\frac{\partial x_k}{\partial z^i}\frac{\partial x_k}{\partial z^j}=\frac{\partial x_h}{\partial z^i}\frac{\partial x_k}{\partial z^j}\delta_{hk},\\
 &g^{ij}=\frac{\partial z^i}{\partial x_k}\frac{\partial z^j}{\partial x_k}=\frac{\partial z^i}{\partial x_h}\frac{\partial z^j}{\partial x_k}\delta_{hk},\\
 &\delta^i_{\ j}=\frac{\partial z^i}{\partial x_h}\frac{\partial x_k}{\partial z^j}\delta_{hk},\\
 &\delta_i^{\ j}=\frac{\partial x_h}{\partial z^i}\frac{\partial z^j}{\partial x_k}\delta_{hk}.
 \end{split}
 \ee
 So, applying Eq. (\ref{eq:touttenseurscart}) to the identity tensor, we get
 \bes
\I=\delta_{ij}\e_i\otimes\e_j=\frac{\partial x_i}{\partial z^h}\frac{\partial x_j}{\partial z^k}I^{hk}\e_i\otimes\e_j=I^{hk}\g_h\otimes\g_k,
\ees
but by Eqs. (\ref{eq:transftenscartcov}) and (\ref{eq:gedelta}),
\bes
I^{hk}=\frac{\partial z^h}{\partial x_i}\frac{\partial z^k}{\partial x_j}\delta_{ij}=g^{hk}
\ees
so, finally,
\bes
\I=g^{hk}\g_h\otimes\g_k.
\ees
Proceeding in a similar manner, we can also get 
\bes
\I=g_{hk}\g^h\otimes\g^k=\delta^h_{\ k}\g_h\otimes\g^k=\delta_h^{\ k}\g^h\otimes\g_k.
\ees
We see hence that the $g_{hk}$s represent $\I$ in covariant coordinates,  the $g^{hk}$s in the contravariant ones and the $\delta_h^{\ k}$s and $\delta^h_{\ k}$s in mixed coordinates.

 \section{Spatial derivatives of fields in curvilinear coordinates}
 Let $\phi$ a spatial\footnote{The term {\it spatial} here refers to differentiation with respect to spatial coordinates, which can be Cartesian or curvilinear.} scalar field, $\phi:\Eu\rightarrow\R$. Generally speaking,
 \bes
 \phi=\phi(z^j(x_i)),
 \ees
 or also
 \bes
 \phi=\phi(x_j(z^k)),
 \ees
 where the $x_j$s, $z^k$s are, respectively, Cartesian and curvilinear coordinates, related as in Eq. (\ref{eq:cc1}). By the chain rule
 \be
 \label{eq:grad1}
 \frac{\partial \phi}{\partial x_j}=\frac{\partial \phi}{\partial z^k}\frac{\partial z^k}{\partial x_j}
 \ee
 and inversely
 \bes
 \frac{\partial \phi}{\partial z^k}= \frac{\partial \phi}{\partial x_j}\frac{\partial x_j}{\partial z^k}.
 \ees
 We remark that the last quantity transforms like the components of a covariant vector, cf. Eq. (\ref{eq:cc10}). 
 
 The  gradient of $\phi$ is the vector that in the Cartesian basis, cf. Eq. (\ref{eq:cartcoorddiffop})$_1$, is given by
  \bes
 \nabla\phi= \frac{\partial \phi}{\partial x_j}\e_j; 
 \ees
so by Eqs. (\ref{eq:dualbasis}) and  (\ref{eq:grad1}) we get that, in the dual basis,
 \bes
 \nabla\phi=\frac{\partial \phi}{\partial z^k}\frac{\partial z^k}{\partial x_j}\e_j= \frac{\partial \phi}{\partial z^k}\g^k.
  \ees
 We see hence that  in curvilinear coordinates the { nabla operator}, Eq. (\ref{eq:operatorenabla}), is defined by
 \be
 \label{eq:nabla1}
 \nabla(\cdot)=\frac{\partial\ \cdot}{\partial z^k}\g^k.
 \ee
 
 The contravariant components of the gradient can be obtained by the covariant ones upon multiplication by the components of the inverse (contravariant) metric tensor, Eq. (\ref{eq:gcont1}):
 \bes
 g^{hk} \frac{\partial \phi}{\partial z^k}=
 \frac{\partial z^h}{\partial x_i}\frac{\partial z^k}{\partial x_i} \frac{\partial \phi}{\partial x_j}\frac{\partial x_j}{\partial z^k}=\delta_{ij}\frac{\partial \phi}{\partial x_j}\frac{\partial z^h}{\partial x_i}=\frac{\partial \phi}{\partial x_j}\frac{\partial z^h}{\partial x_j}\ \rightarrow  \ \nabla\phi=\frac{\partial \phi}{\partial x_j}\frac{\partial z^h}{\partial x_j}\g_h.
 \ees

 Let us now consider a vector field $\bv:\Eu\rightarrow\Ve$; we want to calculate the spatial derivative of its Cartesian components. By the chain rule and Eq. (\ref{eq:cont-cart}), we get
 \bes
 \besp
 \frac{\partial v_i^c}{\partial x_j}&= \frac{\partial v_i^c}{\partial z^k}\frac{\partial z^k}{\partial x_j}=\frac{\partial z^k}{\partial x_j}\frac{\partial}{\partial z^k}\left( \frac{\partial x_i}{\partial z^h}v^h\right)=
 \frac{\partial z^k}{\partial x_j}\left( \frac{\partial x_i}{\partial z^h}\frac{\partial v^h}{\partial z^k}+\frac{\partial^2x_i}{\partial z^k\partial z^l}v^l \right)\\
 &=\frac{\partial z^k}{\partial x_j} \frac{\partial x_i}{\partial z^h}\left(\frac{\partial v^h}{\partial z^k}+\frac{\partial z^h}{\partial x_m}\frac{\partial^2x_m}{\partial z^k\partial z^l}v^l \right),
 \end{split}
 \ees
 whence
 \be
 \label{eq:covdev1}
  \frac{\partial z^h}{\partial x_i} \frac{\partial x_j}{\partial z^k}  \frac{\partial v_i^c}{\partial x_j}= \frac{\partial v^h}{\partial z^k}+\frac{\partial z^h}{\partial x_m}\frac{\partial^2x_m}{\partial z^k\partial z^l}v^l.
 \ee
 Comparing this result with Eq. (\ref{eq:mixedtens1})$_1$ we see that the first member actually corresponds to the components of a mixed tensor field, which is the {\it gradient of the vector field} $\bv$, that we write as
 \be
 \label{eq:covdev2}
 v^h_{\ ;k}=\frac{\partial v^h}{\partial z^k}+\Gamma^h_{kl}v^l,
 \ee
 where  the functions
 \be
 \label{eq:christoffelsymb1}
 \Gamma^h_{kl}=\frac{\partial z^h}{\partial x_m}\frac{\partial^2x_m}{\partial z^k\partial z^l}
 \ee
 are the {\it Christoffel symbols}. We immediately see that $\Gamma^h_{kl}=\Gamma^h_{lk}$.
 The quantity $ v^h_{\ ;k}$ is the {\it covariant derivative of the contravariant components $v^h$}. 
 The proof  that the Christoffel symbols can also be written  as
\be
\label{eq:christoffavecg}
\Gamma^h_{kl}=\frac{1}{2}g^{hm}\left(\frac{\partial g_{mk}}{\partial z^l}+\frac{\partial g_{ml}}{\partial z^k}-\frac{\partial g_{kl}}{\partial z^m}\right)
\ee
is left to the reader as an exercise.

 Proceeding in a similar way for the covariant components of $\bv$,  but now using Eqs. (\ref{eq:cov-cart}) and (\ref{eq:mixedtens1})$_1$, we get
 \bes
 v_{h;k}=\frac{\partial v_h}{\partial z^k}-\Gamma^l_{kh}v_l, 
 \ees
 which is the {\it covariant derivative of the covariant components $v_h$}.
 
 Using Eqs. (\ref{eq:covdev1}) and (\ref{eq:covdev2}), we conclude that, cf. Eq. (\ref{eq:cartcoorddiffop})$_3$,
 \bes
 \div\bv=\frac{\partial v^c_i}{\partial x_i}=v^h_{\ ;h}.
 \ees
 
 Then, applying the operator divergence so defined to the gradient of the scalar field $\phi$ we obtain, in curvilinear coordinates $z^k$, the Laplacian $\Delta \phi$ as
 \be
 \label{eq:laplaciancurvcoord}
 \Delta \phi=\left(g^{hk}\frac{\partial \phi}{\partial z^k}\right)_{;h}=\frac{\partial}{\partial z^h}\left(g^{hk}\frac{\partial \phi}{\partial z^k}\right)+\Gamma^h_{hj}g^{jk}\frac{\partial \phi}{\partial z^k}.
 \ee
 Using the definition of the nabla operator in curvilinear coordinates, Eq. (\ref{eq:nabla1}), jointly to the fact that, cf. Section  \ref{sec:diffopercart},
 \bes
\Delta f:= \div\nabla \phi=\nabla\cdot\nabla \phi,
 \ees
we get the the following representation of the Laplace operator in curvilinear coordinates:
\bes
\besp
\Delta(\cdot)&=\nabla\cdot\nabla(\cdot)=\left(\frac{\partial}{\partial z^k}\left(\frac{\partial(\cdot)}{\partial z^h}\g^h\right)\right)\cdot\g^k
=\frac{\partial^2(\cdot)}{\partial z^k\partial z^h}\g^h\cdot\g^k+\frac{\partial \g^h}{\partial z^k}\frac{\partial(\cdot)}{\partial z^h}\cdot\g^k\\
&=\frac{\partial^2(\cdot)}{\partial z^k\partial z^h}g^{hk}+\frac{\partial \g^h}{\partial z^k}\cdot\g^k\frac{\partial(\cdot)}{\partial z^h}.
\end{split}
\ees

Let us now calculate the spatial derivatives of the components of a 2$^{nd}$-rank tensor $\L$: By Eqs. (\ref{eq:touttenseurscart})$_1$ and (\ref{eq:christoffelsymb1})  we get
\bes
\besp
\frac{\partial L_{ij}^c}{\partial x_k}&=\frac{\partial z^h}{\partial x_k}\frac{\partial}{\partial z^h}\left(\frac{\partial x_i}{\partial z^n}\frac{\partial x_j}{\partial z^p}L^{np}\right)\\
&=\frac{\partial z^h}{\partial x_k}\frac{\partial x_i}{\partial z^n}\frac{\partial x_j}{\partial z^p}\left(\frac{\partial L^{np}}{\partial z^h}+\Gamma_{hr}^nL^{rp}+\Gamma_{hr}^pL^{nr}\right),
\end{split}
\ees
which implies that
\be
\label{eq:covderivL0}
\frac{\partial L^{np}}{\partial z^h}+\Gamma_{hr}^nL^{rp}+\Gamma_{hr}^pL^{nr}=\frac{\partial x_k}{\partial z^h}\frac{\partial z^n}{\partial x_i}\frac{\partial z^p}{\partial x_j}\frac{\partial L_{ij}^c}{\partial x_k}.
\ee
So, using Eq. (\ref{eq:transftenscartcov}), we can conclude that the expression
\be
\label{eq:covderivL1}
L^{np}_{\ \ ;h}=\frac{\partial L^{np}}{\partial z^h}+\Gamma_{hr}^nL^{rp}+\Gamma_{hr}^pL^{nr}
\ee
 represents the covariant derivative of the contravariant components of the second-rank tensor $\L$. 
In a similar manner, this time by Eq. (\ref{eq:touttenseurscart})$_2$, we obtain the covariant derivatives of the covariant components of $\L$:
\bes
\frac{\partial L_{np}}{\partial z^h}-\Gamma_{nh}^rL_{rp}-\Gamma_{ph}^rL_{nr}=\frac{\partial x_k}{\partial z^h}\frac{\partial x_i}{\partial z^n}\frac{\partial x_j}{\partial z^p}\frac{\partial L_{ij}^c}{\partial x_k},
\ees
i.e.
\be
\label{eq:covderivL2}
L_{np;h}=\frac{\partial L_{np}}{\partial z^h}-\Gamma_{ph}^rL_{nr}-\Gamma_{nh}^rL_{pr}.
\ee
The same procedure with Eqs. (\ref{eq:touttenseurscart})$_{3,4}$ gives the covariant derivatives of the  mixed components\footnote{Equations (\ref{eq:covderivL1}), (\ref{eq:covderivL2}) and (\ref{eq:covderivL3}) represent the different forms of the components of an operator depending upon three indices, i.e. of a third-rank tensor: $\nabla\L$, the {\it gradient of} $\L$. } of $\L$:
\be
\label{eq:covderivL3}
\besp
L^n_{\ p;h}=\frac{\partial L^n_{\ p}}{\partial z^h}+\Gamma_{hr}^nL^r_{\ p}-\Gamma_{ph}^rL^n_{\ r},\\
L_{p\  ;h}^{\ n}=\frac{\partial L^{\ n}_{p}}{\partial z^h}-\Gamma_{ph}^rL_r^{\ n}+\Gamma_{hr}^nL_p^{\ r}.
\end{split}
\ee

If in Eqs. (\ref{eq:covderivL0}) and (\ref{eq:covderivL1}) we set $p=h$, we get
\bes
\besp
L^{nh}_{\ \ ;h}&=\frac{\partial L^{nh}}{\partial z^h}+\Gamma_{hr}^nL^{rh}+\Gamma_{hr}^hL^{nr}=\frac{\partial x_k}{\partial z^h}\frac{\partial z^n}{\partial x_i}\frac{\partial z^h}{\partial x_j}\frac{\partial L_{ij}^c}{\partial x_k}\\
&=\delta_{kj}\frac{\partial z^n}{\partial x_i}\frac{\partial L_{ij}^c}{\partial x_k}=\frac{\partial z^n}{\partial x_i}\frac{\partial L_{ij}^c}{\partial x_j},
\end{split}
\ees
which are the components of the contravariant vector field $\div\L$.

 \section{Exercises}
 \begin{enumerate}
 \item \label{ex:1ch6} Write $\g$ and $ds$ for cylindrical coordinates.
 \item \label{ex:2ch6} Write $\g$ and $ds$ for spherical coordinates.
 \item Find the length of a  helix traced on a circular cylinder of radius $R$ between the angles $\theta$ and $\theta+2\pi$.
 \item A curve is traced in a quarter  circle of radius $R$, see Fig. \ref{fig:27}, with $\rho$ proportional to $\theta$. When the quarter of circle is rolled into a cone, the curve appears as indicated in the figure. Determine the length $\ell$ of the curve, first using the polar coordinates in the plane of the quarter  circle,  then the cylindrical ones for the case of the curve on the cone (exercise given in the book by W. H. Müller, see the suggested texts). 
 \begin{figure}[ht]
	\begin{center}
         \includegraphics[scale=.7]{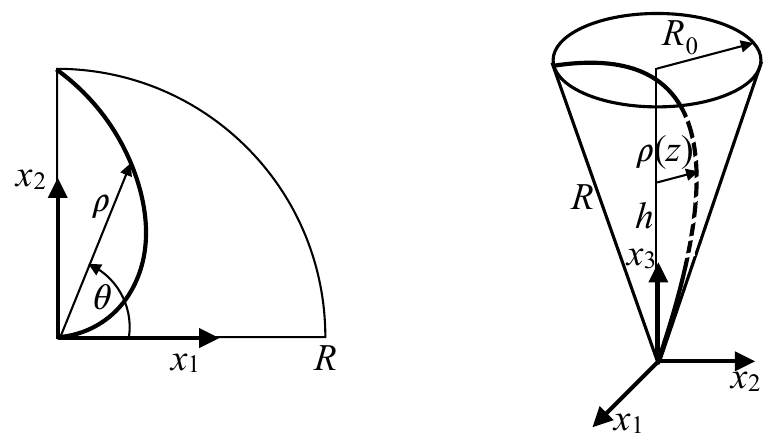}
	\caption{Curve in a plane and on a cone.}
	\label{fig:27}
	\end{center}
\end{figure}
\item \label{ex:a} Calculate $\g$   for a planar system of coordinates composed of two axes $z^1$ and $z^2$ inclined, respectively, at $\alpha_1$ and $\alpha_2$ on the axis $x_1$. Then, find the vectors $\g_k$ and $\g^k,\ k=1,2$,  check the orthogonality conditions  $\g^h\cdot\g_k=\delta^h_{\ k}$, determine the norm of these vectors and design them.
\item Calculate  the $\g_i$s for a system of spherical coordinates.
\item In the plane, elliptical coordinates are defined by the relations
\bes
x_1=c \ \cosh z^1\cos z^2, \ \ x_2=c\ \sinh z^1\sin z^2, \ \ z^1\in(0,\infty),\ z^2\in[0,2\pi);
\ees
show that the lines $z^1=const.$ and $\ z^2=const.$ are confocal ellipses and hyperbolae, determine the axes of the ellipses in terms of the parameter $c$, discuss the limit case of ellipses that degenerate into a crack and determine its length. Finally, find $\g,\g_1$ and $\g_2$.  
\item Determine the co- and contravariant components of a tensor $\L$ in cylindrical coordinates.
\item Determine the co- and contravariant components of a tensor $\L$ in spherical coordinates.
\item Show that
\bes
\tr\L=L_{ii}^c=g_{ij}L^{ij}=g^{ij}L_{ij}=L^i_{\ i}=L_j^{\ j}.
\ees
\item Prove Eq. (\ref{eq:christoffavecg}).
\item Prove the {\it Lemma of Ricci}:
\bes
\frac{\partial g_{jk}}{\partial z^h}=\Gamma^i_{jh}g_{ik}+\Gamma^i_{kh}g_{ji}.
\ees
\item Using Eq. (\ref{eq:christoffavecg}), find the Christoffel symbols for the cylindrical, spherical and elliptical (in the plane) coordinates.
\item Write the Laplacian $\Delta f$ of a spatial scalar field $f$ in cylindrical and spherical coordinates.
\item Prove that
\bes
g^{np}_{\ \ ;h}=g_{np;h}=0.
\ees

 \end{enumerate}

\chapter{Surfaces in $\Eu$}
\label{ch:7}

\section{Surfaces in $\Eu$, coordinate lines, tangent planes}
\label{sec:surf1}
\begin{figure}[ht]
	\begin{center}
         \includegraphics[scale=.7]{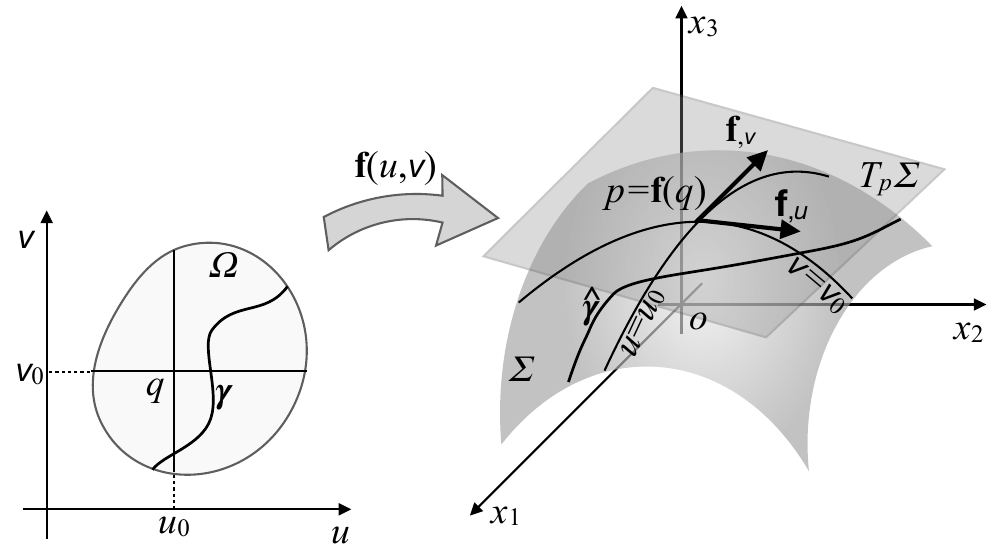}
	\caption{General scheme of a surface  and of the tangent space at a point $p$.}
	\label{fig:28}
	\end{center}
\end{figure}

A function $\f(u,v):\Omega\subset\R^2\rightarrow\Eu$ of class $\geq$C$^1$ and such that  its Jacobian 
\bes
J=\left[\begin{array}{cc}
\dfrac{\partial f_1}{\partial u}&\dfrac{\partial f_1}{\partial v}\medskip\\
\dfrac{\partial f_2}{\partial u}&\dfrac{\partial f_2}{\partial v}\medskip\\
\dfrac{\partial f_3}{\partial u}&\dfrac{\partial f_3}{\partial v}
\end{array}
\right]
\ees
has maximum rank (rank[J]=2) defines a {\it surface in $\Eu$}, see Fig. \ref{fig:28}. We say also that $\f$ is an {\it immersion} of $\Omega$ into $\Eu$ and that  the subset  $\Sigma\subset\Eu$ image of $\f$ is the {\it support} or {\it trace} of the surface $\f$. 

 As usual, we   indicate the derivatives with respect to the variables $u,v$ by, for example, $\dfrac{\partial \f}{\partial u}=\f_{,u}$ etc. The condition on the rank of $J$ is equivalent to impose that
 \be
 \label{eq:regulsurf}
 \f_{,u}(u,v)\times\f_{,v}(u,v)\neq\bo\ \ \forall (u,v)\in\Omega.
 \ee
 This allows us to introduce the {\it normal to the surface $\f$} as the vector $\N\in\S$ defined by
 \be
 \label{eq:normsurfreg}
 \N:=\frac{\f_{,u}\times\f_{,v}}{|{\f_{,u}\times\f_{,v}}|}.
 \ee
 A {\it regular point of $\Sigma$} is a point where $\N$ is defined; if $\N$ is defined $\forall p\in\Sigma$ then the surface is said to be {\it regular}.
 
 A function $\bg(t):G\subset\R\rightarrow\Omega$ whose parametric equation is $\bg(t)=(u(t),v(t))$ describes a curve in $\Omega$ whose image, through $\f$, is the curve, see Fig. \ref{fig:28},
 \bes
 \widehat{\bg}(t)=\f(u(t),v(t)):G\subset\R\rightarrow\Sigma\subset\Eu.
 \ees
 As a special case of curve in $\Omega$, let us consider the curves of the type $v=v_0$ or $u=u_0$, with $u_0,v_0$ being some constants. Then, their image through $\f$ are two curves $\f(u,v_0),\f(u_0,v)$ on $\Sigma$ called {\it coordinate lines}, see again Fig. \ref{fig:28}. The {\it tangent vectors to the coordinate lines} are, respectively, the vectors $\f_{,u}(u,v_0)$ and $\f_{,v}(u_0,v)$, while the tangent to a curve $\widehat{\bg}(t)=\f(u(t),v(t))$ is the vector 
 \be
 \label{eq:tangonsigma}
\widehat{\bg}'(t)=\f_{,u}\frac{du}{dt}+\f_{,v}\frac{dv}{dt},
 \ee
 i.e. the tangent vector to any curve on $\Sigma$ is a linear combination of the tangent vectors to the coordinate lines. We remark that the tangent vectors $\f_{,u}(u,v_0)$ and $\f_{,v}(u_0,v)$ are necessarily non-null and  linear independent as consequence of the assumption on the rank of $J$ and hence of the existence of $\N$, i.e. of the regularity of $\Sigma$. They determine a plane that contains the tangents to all the curves on $\Sigma$ passing by $p=\f(u_0,v_0)$ and form a basis on this plane, called the {\it natural basis}. Such a plane is the {\it tangent plane to $\Sigma$ in $p$} and is indicated by $T_p\Sigma$; this plane is actually the space spanned by $\f_{,u}(u,v_0)$ and $\f_{,v}(u_0,v)$ and  is also called the {\it tangent vector space}.
 
 Let us consider two open subsets $\Omega_1,\Omega_2\subset\R^2$; a {\it diffeomorphism\footnote{The definition of diffeomorphism, of course, can be given for subsets of $\R^n,n\geq1$; here, we bound the definition to the case of interest.} of class} C$^k$ between $\Omega_1$ and $\Omega_2$ is a bijective map $\vartheta:\Omega_1\rightarrow\Omega_2$ of class C$^k$ with also its inverse of class C$^k$; the diffeomorphism is {\it smooth} if $k=\infty$.
 
 Let $\Omega_1,\Omega_2$  be two open subsets of $\R^2$, $\f:\Omega_2\rightarrow\Eu$ a surface and  $\vartheta:\Omega_1\rightarrow\Omega_2$ a smooth diffeomorphism. Then the surface $\F=\f\circ\vartheta:\Omega_1\rightarrow\Eu$ is a {\it change of parameterization} for $\f$. In practice, the function defining the surface changes, but not $\Sigma$, its trace in $\Eu$. Let $(U,V)$ be the coordinates in $\Omega_1$ and $(u,v)$ those in $\Omega_2$. Then, by the chain rule,
 \bes
 \besp
& \F_{,U}=\f_{,u}\frac{\partial u}{\partial U}+\f_{,v}\frac{\partial v}{\partial U},\\
&  \F_{,V}=\f_{,u}\frac{\partial u}{\partial V}+\f_{,v}\frac{\partial v}{\partial V},
  \end{split}
 \ees
 or, denoting by $J_\vartheta$ the Jacobian of $\vartheta$, 
 \bes
 \left\{\begin{array}{c}\F_{,U}\\\F_{,V}\end{array}\right\}=\left[J_\vartheta\right]^\top\left\{\begin{array}{c}\f_{,u}\\\f_{,v}\end{array}\right\},
 \ees
 whence, making the cross product, one gets immediately
 \bes
 \F_{,U}\times\F_{,V}=\det [J_\vartheta]\ \f_{,u}\times\f_{,v}.
 \ees
 This result shows that the regularity of the surface, condition (\ref{eq:regulsurf}), the tangent plane and the tangent space vector {\it do not depend upon the parameterization of $\Sigma$}. 
 From the last equation, we get also
 \bes
 \N(U,V)=\mathrm{sgn}(\det [J_\vartheta])\ \N(u,v);
 \ees
 we say that the change of parameterization {\it preserves the orientation if $\det [J_\vartheta]>0$}, and that it {\it inverses the parameterization} in the opposite case.

 \section{Surfaces of revolution}
 A {\it surface of revolution} is a surface whose trace  is obtained by letting  a plane curve, say $\bg$, rotate around an axis, say $x_3$. To be more specific and without loss of generality, let $\bg:G\subset\R\rightarrow\R^2$ be a regular curve of the plane $x_2=0$, whose parametric equation is
 \be
 \label{eq:curvesurfrevol}
 \bg(u):\ \left\{\begin{array}{l}x_1=\phi(u),\\x_3=\psi(u),\end{array}\right.\ \ \ \phi(u)>0\ \forall u\in G.
 \ee
 Then, the subset $\Sigma_\gamma\subset\Eu$ defined by
 \bes
 \Sigma_\gamma:=\left\{(x_1,x_2,x_3)\in\Eu|x_1^2+x_2^2=\phi^2(u),x_3=\psi(u),u\in G\right\}
 \ees
 is  the trace of a surface of revolution of the curve $\bg(u)$ around the axis $x_3$. A general parameterization of such a surface is 
 \be
 \label{eq:revolsurf}
 \f(u,v):G\times(-\pi,\pi]\rightarrow\Eu|\ \ \ 
 \left\{
 \begin{array}{l}
 x_1=\phi(u)\cos v,\\
 x_2=\phi(u)\sin v,\\
 x_3=\psi(u). 
 \end{array}
 \right.
 \ee
 It is readily checked that this parameterization actually defines a regular surface:
 \bes
 \f_{,u}=\left\{\begin{array}{c}\phi'(u)\cos v\\\phi'(u)\sin v\\\psi'(u)\end{array}\right\},\ \ \f_{,v}=\left\{\begin{array}{c}-\phi(u)\sin v\\\phi(u)\cos v\\0\end{array}\right\}\ \rightarrow\ \f_{,u}\times\f_{,v}=\left\{\begin{array}{c}-\phi(u)\psi'(u)\cos v\\-\phi(u)\psi'(u)\sin v\\\phi(u)\phi'(u)\end{array}\right\}
 \ees
 so that 
 \bes
 |\f_{,u}\times\f_{,v}|=\phi^2(u)(\phi'^2(u)+\psi'^2(u))\neq0\ \forall u\in G
 \ees
 for being $\bg(u)$ a regular curve, i.e. with $\bg'(u)\neq\bo\ \forall u\in G$.
 A {\it meridian} is a curve in $\Eu$ intersection of the trace of $\f$, $\Sigma_\gamma$, with a plane containing the axis $x_3$; the equation of a meridian is obtained fixing the value of $v$, say $v=v_0$:
 \bes
 \left\{\begin{array}{l}x_1=\phi(u)\cos v_0,\\x_2=\phi(u)\sin v_0,\\x_3=\psi(u).\end{array}\right.
 \ees
 A {\it parallel} is a curve in $\Eu$  intersection  of $\Sigma_\gamma$ with a plane orthogonal to $x_3$; the equation of a parallel, which is a circle with center on the axis $x_3$, is obtained by fixing the value of $u$, say $u=u_0$:
 \bes
 \left\{\begin{array}{l}x_1=\phi(u_0)\cos v,\\x_2=\phi(u_0)\sin v,\\x_3=\psi(u_0),\end{array}\right.
 \ees
 or also
 \bes
  \left\{\begin{array}{l}x_1^2+x_2^2=\phi(u_0)^2,\\x_3=\psi(u_0);\end{array}\right.
 \ees
 the radius of the circle is $\phi(u_0)$. 
 A {\it loxodrome} or {\it rhumb line} is a curve on $\Sigma_\gamma$ crossing all the meridians  at the same angle\footnote{This concept is important for  marine and aerial navigation}. 
  \begin{figure}[ht]
	\begin{center}
         \includegraphics[height=.13\textheight]{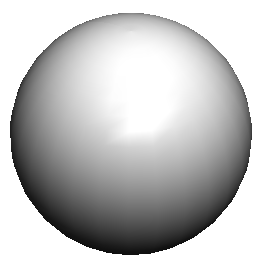}
         \includegraphics[height=.13\textheight]{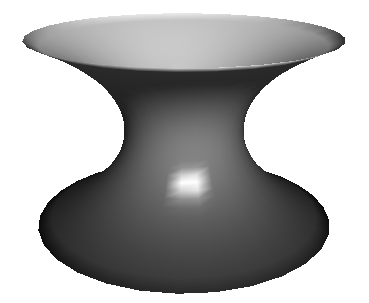}
         \includegraphics[height=.13\textheight]{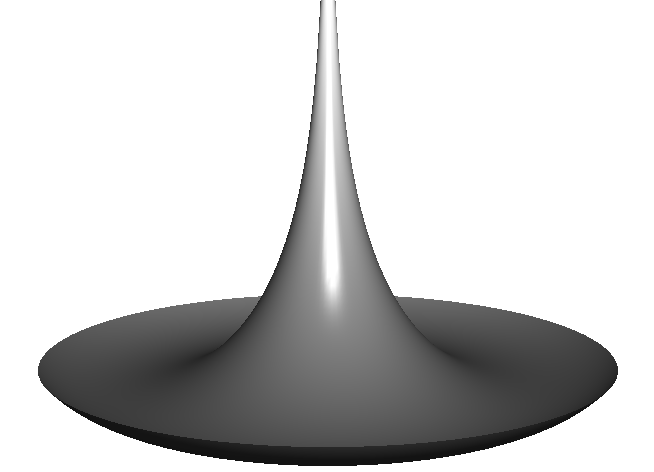}
         \includegraphics[height=.13\textheight]{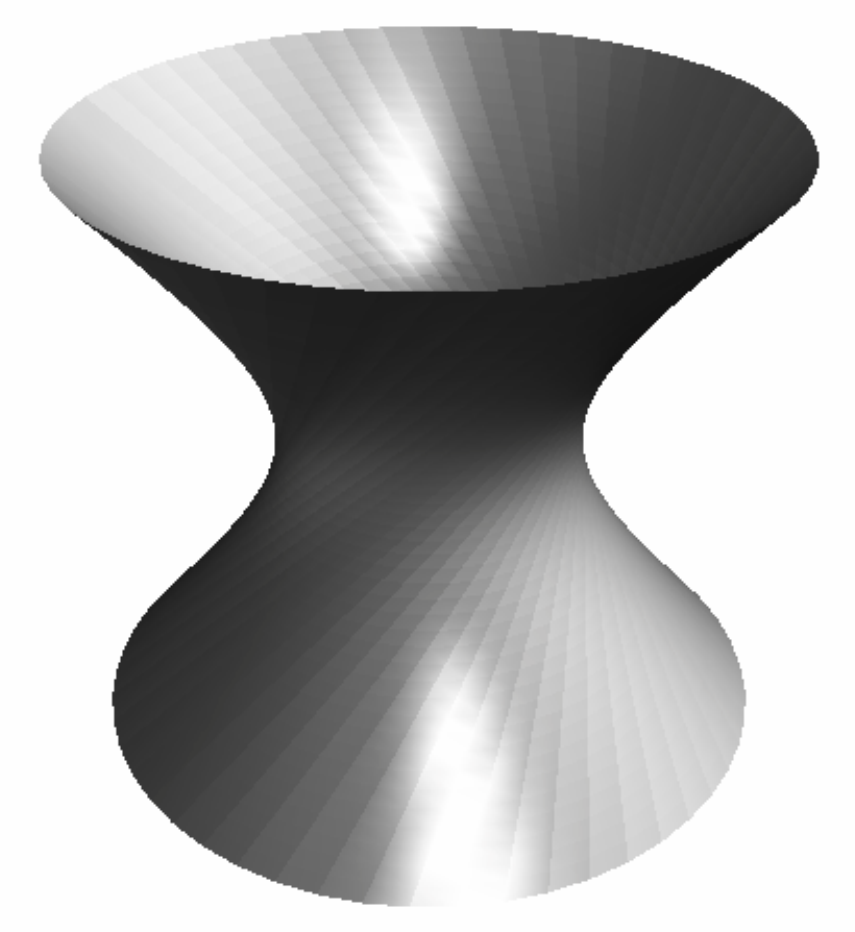}
	\caption{Surfaces of revolution. From the left: sphere, catenoid, pseudo-sphere, hyperbolic hyperboloid.}
	\label{fig:29}
	\end{center}
\end{figure}
 Some important examples of surfaces of revolution are:
 \begin{itemize}
 \item the {\it sphere}: 
 \bes
 \f(u,v):\left[-\dfrac{\pi}{2},\dfrac{\pi}{2}\right]\times(-\pi,\pi]\rightarrow\Eu|\ \ \left\{\begin{array}{l}x_1=\cos u\cos v,\\x_2=\cos u\sin v,\\x_3=\sin v;\end{array}\right.
 \ees
 \item the {\it catenoid}:
 \bes
 \f(u,v):\left[-a,a\right]\times(-\pi,\pi]\rightarrow\Eu|\ \ \left\{\begin{array}{l}x_1=\cosh u\cos v,\\x_2=\cosh u\sin v,\\x_3=u;\end{array}\right.
 \ees
 \item the {\it pseudo-sphere}:
  \be
  \label{eq:pseudosphere}
 \f(u,v):\left[0,a\right]\times(-\pi,\pi]\rightarrow\Eu|\ \ \left\{\begin{array}{l}x_1=\sin u\cos v,\\x_2=\sin u\sin v,\\x_3=\cos u+\ln\left(\tan\dfrac{u}{2}\right);\end{array}\right.
 \ee
 \item the {\it hyperbolic hyperboloid}:
  \bes
 \f(u,v):\left[-a,a\right]\times(-\pi,\pi]\rightarrow\Eu|\ \ \left\{\begin{array}{l}x_1=\cos u-v\sin u,\\x_2=\sin u+v\cos u,\\x_3= v.\end{array}\right.
 \ees
 \end{itemize}

  \section{Ruled surfaces}
 A {\it ruled surface} (also named a {\it scroll}) is a surface with the property that through every one of its points, there is a straight line that lies on the  surface. 
 A ruled surface can be seen as the set of points swept by a moving straight line.  We say that a surface is {\it doubly ruled} if through every one of its points, there are two distinct straight lines that lie on the  surface.
 
 Any ruled surface can be represented by a parameterization of the form
 \be
 \label{eq:ruledsurf}
 \f(u,v)=\bg(u)+v\bl(u),
 \ee
 where $\bg(u)$ is a regular smooth curve, the {\it directrix},  and $\bl(u)$ is a smooth curve. Fixing $u=u_0$ gives a {\it generator line} $\f(u_0,v)$ of the surface; the vectors $\bl(u)\neq\bo$ describe the directions of the generators. Some important examples of ruled surfaces are:
 \begin{itemize}
 \item {\it Cones}: For these surfaces, all the straight lines pass through a point, the {\it apex} of the cone; choosing the apex as the origin, then it must be $\bl(u)=k\bg(u),\ k\in\R\rightarrow$
 \bes
 \f(u,v)=v\bg(u);
 \ees 
\item {\it Cylinders}:  A ruled surface is a cylinder $\iff\bl(u)=const$. In this case, it is always possible to choose $\bl(u)\in\S$ and  $\bg(u)$  a planar curve lying in a plane orthogonal to $\bl(u)$; in fact, it is sufficient to choose the curve $\bg^*(u)=(\I-\bl(u)\otimes\bl(u))\bg(u)$;
 \item {\it Helicoids}:  A surface  generated by rotating and simultaneously displacing a curve, the {\it profile curve}, along an axis is a helicoid. Any point of the profile curve is the starting point of a circular helix. Generally, we get a helicoid if 
 \bes
 \bg(u)=(0,0,\phi(u)),\ \ \bl(u)=(\cos u,\sin u,0),\ \ \ \phi(u):\R\rightarrow\R.
 \ees
 \item {\it Möbius strip}: It is a ruled surface with
 \bes
 \bg(u)=(\cos 2u,\sin 2u,0),\ \ \ \bl(u)=(\cos u\cos 2u,\cos u\sin 2u,\sin u).
 \ees
  \end{itemize}
 \begin{figure}[ht]
	\begin{center}
         \includegraphics[height=.13\textheight]{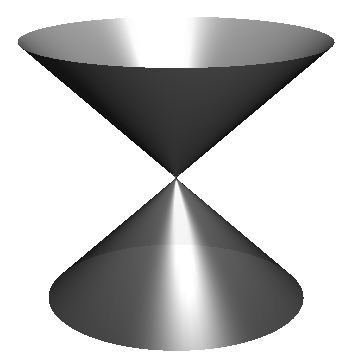}
         \includegraphics[height=.13\textheight]{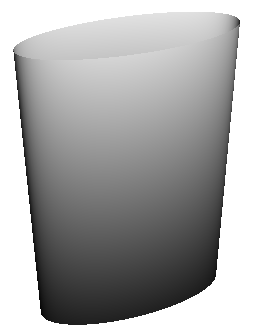}
         \includegraphics[height=.13\textheight]{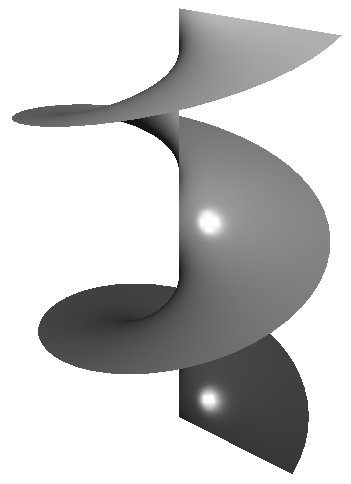}
         \includegraphics[height=.13\textheight]{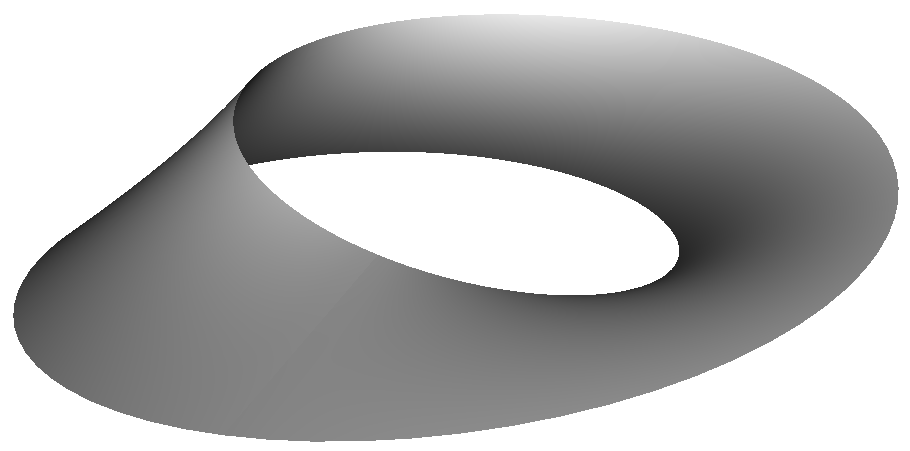}
	\caption{Ruled surfaces: (from  left) elliptical cone, elliptical cylinder, helicoid and Möbius strip.}
	\label{fig:30}
	\end{center}
\end{figure}

 \section{First fundamental form of a surface}
 Let us consider two vectors of $T_p\Sigma$, say $\bw_1,\bw_2$; we want to calculate their scalar product  in terms of their components in the natural basis $ \{\f_{,u},\f_{,v}\}$ of $T_p\Sigma$. If $\bw_1=a_1\f_{,u}+b_1\f_{,v}$ and $\bw_2=a_2\f_{,u}+b_2\f_{,v}$, then
 \bes
\bw_1\cdot\bw_2=a_1a_2\f^2_{,u}+(a_1b_2+a_2b_1)\f_{,u}\cdot\f_{,v}+b_1b_2\f^2_{,v}, 
\ees
which can be rewritten as the form
\bes
I(\bw_1,\bw_2)=\bw_1\cdot\g\bw_2,
\ees
where\footnote{\label{note:firstform}Often, in texts on differential geometry, tensor $\g$ is indicated as \bes\g=\left[\begin{array}{cc}E&F\\F&G\end{array}\right],\ees where $E:=\f_{,u}\cdot\f_{,u},F:=\f_{,u}\cdot\f_{,v},G:=\f_{,v}\cdot\f_{,v}$.}
 \bes
 \g=\left[\begin{array}{cc}\f_{,u}\cdot\f_{,u}&\f_{,u}\cdot\f_{,v}\\\f_{,v}\cdot\f_{,u}&\f_{,v}\cdot\f_{,v}\end{array}\right]
 \ees
 is precisely the metric tensor $\g$ of $\Sigma$, cf. Eq. (\ref{eq:meaningghk}). In fact, $\f_{,u}$ and $\f_{,v}$ are the tangent vectors to the coordinate lines on $\Sigma$, i.e. they coincide with the vectors $\g_k$s. 
 
 $I(\bw_1,\bw_2)$  is the {\it first fundamental form} (or simply the {\it first form}) of $\f(u,v)$. If $\bw_1=\bw_2=\bw=a\f_{,u}+b\f_{,v}$, then
  \bes
 I(\bw)=\bw^2=a^2\f^2_{,u}+2ab\f_{,u}\cdot\f_{,v}+b^2\f^2_{,v}.
 \ees
By the same definition of scalar product, $I(\bw_1,\bw_2)$ is a  positive definite, bilinear, symmetric form $\forall \bw\in T_p\Sigma$.  
 
 Through $I(\cdot,\cdot)$ we can calculate some important  quantities regarding the geometry of $\Sigma$:
 \begin{itemize}
 \item Metric on $\Sigma$: $\forall ds\in\Sigma$,
 \bes
 ds^2=ds\cdot ds=I(ds);
 \ees
 so, if 
 \bes
 ds=\f_{,u}du+\f_{,v}dv,
 \ees
 then
 \be
 \label{eq:metricsigma}
 ds^2=\f_{,u}^2du^2+2\f_{,u}\cdot\f_{,v}du\ dv+\f_{,v}^2dv^2;
 \ee
 \item Length $\ell$ of a curve $\bg:[t_1,t_2]\subset\R\rightarrow\Sigma$: We know, see Eq. (\ref{eq:esse}), that the length of a curve is the integral of the tangent vector:
 \bes
 \ell=\int_{t_1}^{t_2}|\bg'(t)|dt=\int_{t_1}^{t_2}\sqrt{\bg'(t)\cdot\bg'(t)}dt
 \ees
and hence, see  Eq. (\ref{eq:tangonsigma}), if we call $\bw=(u',v')$ the tangent vector to $\bg$, expressed by its components in the natural basis,
\be
\label{eq:lengthcurve}
\besp
\ell&=\int_{t_1}^{t_2}\sqrt{u'^2\f^2_{,u}+2u'v'\f_{,u}\cdot\f_{,v}+v'^2\f^2_{,v}}dt=\int_{t_1}^{t_2}\sqrt{(u',v')\cdot\g\ (u',v')}dt\\
&=\int_{t_1}^{t_2}\sqrt{I(\bw)}dt;
\end{split}
\ee
 \item Angle $\theta$ formed by two vectors $\bw_1,\bw_2\in T_p\Sigma$:
 \bes
 \cos\theta=\frac{\bw_1\cdot\bw_2}{|\bw_1||\bw_2|}=\frac{I(\bw_1,\bw_2)}{\sqrt{I(\bw_1)}\sqrt{I(\bw_2)}};
 \ees
 \item Area of a small surface on $\Sigma$: Let $\f_{,u}du$ and $\f_{,v}dv$ be two small vectors on $\Sigma$, forming together the angle $\theta$, that are the transformed, through\footnote{For the sake of conciseness, from now on we will indicate a surface as the function $\f:\Omega\rightarrow\Sigma$, with $\f=\f(u,v), (u,v)\in\Omega\subset\R^2$ and $\Sigma\subset\Eu$.} $\f:\Omega\rightarrow\Sigma$, of two small orthogonal vectors $du,dv\in\Omega$; then the area $d\mathcal{A}$ of the parallelogram determined by them is
 \bes
 \besp
 d\mathcal{A}&=|\f_{,u}du\times\f_{,v}dv|=|\f_{,u}\times\f_{,v}|du\ dv=\sqrt{\f^2_{,u}\f^2_{,v}\sin^2\theta}du\ dv\\
 &=\sqrt{\f^2_{,u}\f^2_{,v}(1-\cos^2\theta)}du\ dv=\sqrt{\f^2_{,u}\f^2_{,v}-\f^2_{,u}\f^2_{,v}\cos^2\theta}du\ dv\\
 &=\sqrt{\f^2_{,u}\f^2_{,v}-(\f_{,u}\cdot \f_{,v})^2}du\ dv=\sqrt{\det\g}du\ dv.
 \end{split}
 \ees
 The term $\sqrt{\det\g}$ is hence the {\it dilatation factor of the areas}; recalling Eq. (\ref{eq:signdetg}), we see that the previous expression has a sense $\forall\f(u,v)$, i.e. for any parameterization of the surface.
  \end{itemize}

 \section{Second fundamental form of a surface}
\begin{figure}[ht]
	\begin{center}
         \includegraphics[width=.8\textwidth]{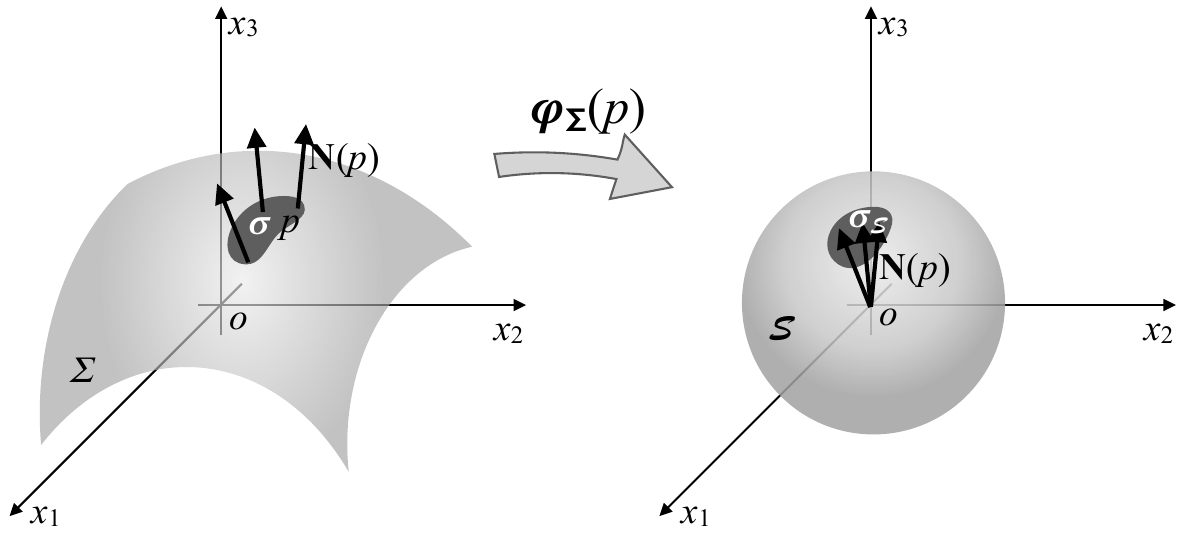}
	\caption{The map of Gauss.}
	\label{fig:37}
	\end{center}
\end{figure}

Be $\f:\Omega\rightarrow\Sigma$ a regular surface, $\{\f_{,u},\f_{,v}\}$ the natural basis for $T_p\Sigma$ and $\N\in\S$ the normal to $\Sigma$ defined as in (\ref{eq:normsurfreg}). We call {\it map of Gauss} of $\Sigma$ the map $\bphi_\Sigma:\Sigma\rightarrow\S$ that associates to each $p\in\Sigma$ its $\N:\ \bphi_\Sigma(p)=\N(p)$. To each subset $\sigma\subset\Sigma$, the map of Gauss associates hence a subset $\sigma_\S\subset\S$, Fig. \ref{fig:37} (e.g. the Gauss map of a plane is just a point on $\S$).

 We want to study how $\N(p)$ varies at the varying of $p$ on $\Sigma$. The  idea is that the change of $\N(p)$ on $\Sigma$ is related to the curvature of the surface\footnote{For curves, the curvature is linked to the change of $\btau$, but for surfaces this should not be meaningful, as $\btau$ is not unique $\forall p\in\Sigma$ while $\N$ is.}.
 For this purpose, we calculate the change in $\N$ per unit length of a curve $\bg(s)\in\Sigma$, i.e. we study how $\N$ varies along any curve of $\Sigma$ per unit of length of the curve itself; that is why we parameterize the curve with its arc-length $s$ \footnote{Actually, it is also possible to introduce the following concepts  more generally, i.e. for any parameterization of the curve but, for the sake of simplicity, in the following we  just use the parameter $s$.}. Let $\N=N_i(u,v)\e_i$; then, if $\btau\in\S$ is the tangent to the curve,
 \bes
 \besp
 \frac{d\N}{ds}&=\frac{dN_i(u(s),v(s))}{ds}\e_i=\left(\frac{\partial N_i}{\partial u}\frac{du}{ds}+\frac{\partial N_i}{\partial v}\frac{dv}{ds}\right)\e_i\\
 &= \nabla N_i\cdot\btau\e_i=(\e_i\otimes\nabla N_i)\btau=(\nabla \N)\ \btau=\frac{d\N}{d\btau}.
 \end{split}
 \ees
The change in $\N$ is hence related to the directional derivative of $\N$ along the tangent $\btau$ to $\bg(s)$, which is a linear operator on $T_p\Sigma$.  Moreover, as $\N\in\S$, then, cf. Eq. (\ref{eq:constvect}),
\bes
\N\cdot\N_{,u}=\N\cdot\N_{,v}=0\ \Rightarrow\  \N_{,u},\N_{,v}\in T_p\Sigma.
\ees
We then call {\it Weingarten operator} $\mathcal{L}_W:T_p\Sigma\rightarrow T_p\Sigma$ the  opposite of the directional derivative of $\N$:
\bes
\mathcal{L}_W(\btau):=-\frac{d\N}{d\btau}.
\ees
 Hence,
 \be
 \label{eq:Wo1}
 \mathcal{L}_W(\f_{,u})=-\N_{,u},\ \ \mathcal{L}_W(\f_{,v})=-\N_{,v}.
 \ee
 Because $\mathcal{L}_W$ is linear, then there exists a tensor $\X$ on $T_p\Sigma$ such that
 \be
 \label{eq:Wo2}
 \mathcal{L}_W(\bv)=\X\bv\ \ \ \forall\bv\in T_p\Sigma.
 \ee
 For any two vectors $\bw_1,\bw_2\in T_p\Sigma$, we define {\it second fundamental form of a surface}, denoted by $II(\bw_1,\bw_2)$ the bilinear form
 \bes
 II(\bw_1,\bw_2):=I(\mathcal{L}_W(\bw_1),\bw_2).
 \ees
 \begin{teo} {\bf (Symmetry of the second fundamental form).} $\forall \bw_1,\bw_2\in T_p\Sigma,\ II(\bw_1,\bw_2)=II(\bw_2,\bw_1)$.
 \begin{proof}
 Because $I$ and $\mathcal{L}_W$ are linear, it is sufficient to prove the thesis for the natural basis $\{\f_{,u},\f_{,v}\}$ of $T_p\Sigma$, and, by the symmetry of $I$, it is sufficient to prove that 
 \bes
 I(\mathcal{L}_W(\f_{,u}),\f_{,v})= I(\f_{,u},\mathcal{L}_W(\f_{,v})),
 \ees
 i.e. that
 \bes
 I(-\N_{,u},\f_{,v})=I(\f_{,u},-\N_{,v})
 \ees
 and in the end that
 \bes
\N_{,u}\cdot\f_{,v}=\f_{,u}\cdot\N_{,v}.
 \ees
To this purpose, we recall that
 \bes
 \N\cdot\f_{,u}=0=\N\cdot\f_{,v}.
 \ees
 So, differentiating the first equation by $v$ and the second one by $u$, we get
 \be
 \label{eq:orthsymL}
 \N_{,v}\cdot\f_{,u}=-\N\cdot\f_{,uv}=\N_{,u}\cdot\f_{,v}.
 \ee
 \end{proof}
 \end{teo}
 
The second fundamental form defines a quadratic, bilinear symmetric form:
\bes
\besp
 II(\bw_1,\bw_2)&=I(\mathcal{L}_W(\bw_1),\bw_2)= I(\bw_1,\mathcal{L}_W(\bw_2))\\ &=I(\bw_1,\X\bw_2)=\bw_1\cdot\g\X\bw_2=\bw_1\cdot\B\bw_2,
\end{split}
\ees
 where
 \be
 \label{eq:defineB}
 \B:=\g\X.
 \ee
In the natural basis $\{\f_{,u},\f_{,v}\}$ of $T_p\Sigma$, by Eq. (\ref{eq:orthsymL}), it is\footnote{\label{note:secondform}In many texts on differential geometry, the following symbols are used: 
\bes
\besp
&L=\f_{,uu}\cdot\N=-\f_{,u}\cdot\N_{,u},\\
&M=\f_{,uv}\cdot\N=-\f_{,u}\cdot\N_{,v},\\
&N=\f_{,vv}\cdot\N=-\f_{,v}\cdot\N_{,v}.
\end{split}
\ees
}
\be
\label{eq:composB}
B_{ij}=II(\f_{,i},\f_{,j})=I(\mathcal{L}_W(\f_{,i}),\f_{,j})=-\N_{,i}\cdot\f_{,j}=\N\cdot\f_{,ij};
\ee
tensor $\X$ can then be calculated by Eq. (\ref{eq:defineB}):
\be
\label{eq:XgB}
\X=\g^{-1}\B.
\ee
By Eq. (\ref{eq:composB}), because $\f_{,ij}=\f_{,ji}$ or simply because $II(\cdot,\cdot)$ is symmetric, we get that 
\bes
\B=\B^\top.
\ees

 \section{Curvatures of a surface}
 Let $\f:\Omega\rightarrow\Sigma$ be a regular surface and $\bg(s):G\subset\R\rightarrow\Sigma$ be a regular curve on $\Sigma$ parameterized with the arc length $s$. We call {\it curvature vector of $\bg(s)$} the vector $\bkappa(s)$ defined as
 \bes
 \bkappa(s):=c(s)\bnu(s)=\bg''(s),
 \ees
 where $\bnu(s)$ is the principal normal to $\bg(s)$. By Eq. (\ref{eq:frenetserret1avecs}), it is also
 \bes
  \bkappa(s)=\bg''(s).
 \ees
 Then, we call {\it normal curvature $\kappa_N(s)$ of $\bg(s)$} the projection of $\bkappa(s)$ onto $\N(s)$, the normal to $\Sigma$:
 \bes
 \kappa_N(s):=\bkappa(s)\cdot\N(s)=c(s)\ \bnu(s)\cdot\N(s)=\bg''(s)\cdot\N(s).
 \ees
 \begin{teo} The normal curvature $\kappa_N(s)$ of $\bg(s)\in\Sigma$ depends uniquely on $\btau(s)$:
 \be
 \label{eq:curvnormcalcul}
 \kappa_N(s)=\btau(s)\cdot\B\btau(s)=II(\btau(s),\btau(s)).
 \ee
 \begin{proof}
\bes
 \bg(s)=\bg(u(s),v(s))\ \rightarrow\ \btau(s)=\bg'(s)=\f_{,u}u'+\f_{,v}v',
 \ees
 therefore $\btau=(u',v')$ in the natural basis and
 \bes
 \bkappa(s)=\bg''(s)=\f_{,u}u''+\f_{,v}v''+\f_{,uu}u'^2+2\f_{,uv}u'v'+\f_{,vv}v'^2
 \ees
 and finally, by Eqs. (\ref{eq:normsurfreg}) and (\ref{eq:composB}),
 \bes
 \kappa_N(s)=\bg''(s)\cdot\N(s)=B_{11}u'^2+2B_{12}u'v'+B_{22}v'^2=\btau\cdot\B\btau=II(\btau,\btau).
 \ees
 \end{proof}
 \end{teo}
 If now $s=s(t)$ is a change of parameter for $\bg$, then
 \bes
 \bg'(t)=|\bg'(t)|\btau(t),
 \ees
 so by the linearity of $II(\cdot,\cdot)$, we get
 \bes
 II(\bg'(t),\bg'(t))=|\bg'(t)|^2II(\btau(t),\btau(t))=|\bg'(t)|^2\kappa_N(t)
 \ees
 and finally,
 \bes
 \kappa_N(t)=\frac{II(\bg'(t),\bg'(t))}{I(\bg'(t),\bg'(t))}.
 \ees
 To each point $p\in\Sigma$, it corresponds  uniquely (in the assumption of regularity of the surface $\f:\Omega\rightarrow\Sigma$) a  tangent plane and a tangent space vector $T_p\Sigma$. In $p$, there are infinite tangent vectors to $\Sigma$, all of them belonging to $T_p\Sigma$. We can associate a curvature to each direction $\bt\in T_p\Sigma$, i.e. to each tangent direction, in the following way: Let us consider the bundle  $\mathcal{H}$ of planes whose support is the straight line through $p$ and parallel to $\N$. Then, any plane $H\in\mathcal{H}$ is a {\it normal plane to $\Sigma$ in $p$}; each normal plane is uniquely determined by a tangent direction $\bt$ and the (planar) curve $\bg_{N\bt}:=H\cap\Sigma$ is called a {\it normal section of $\Sigma$}. If $\bnu$ and $\N$ are, respectively, the principal normal to $\bg_{N\bt}$ and the normal to $\Sigma$ in $p$, then
\bes
\bnu=\pm\N
\ees
 for each normal section. We have, in this way, defined a function that  to each tangent direction $\bt\in T_p\Sigma$ associates the normal curvature $\kappa_N$ of the normal section $\bg_{N\bt}$:
 \bes
 \kappa_N:\S\cap T_p\Sigma\rightarrow\R|\ \ \ \kappa_N(\bt)=\frac{II(\bt,\bt)}{I(\bt,\bt)}.
 \ees
 By the bilinearity of the second fundamental form, $\kappa_N(\bt)=\kappa_N(-\bt)$.
 
 A point $p\in\Sigma$ is said to be a {\it umbilical point} if $\kappa_{N}(\bt)=const.\ \forall\bt$, it is a {\it planar point} if $\kappa_N(\bt)=0\ \forall\bt$. In all the other points, $\kappa_N$ takes a minimum and a maximum value on distinct directions $\bt\in T_p\Sigma$. 

\begin{figure}[ht]
	\begin{center}
         \includegraphics[width=.45\textwidth]{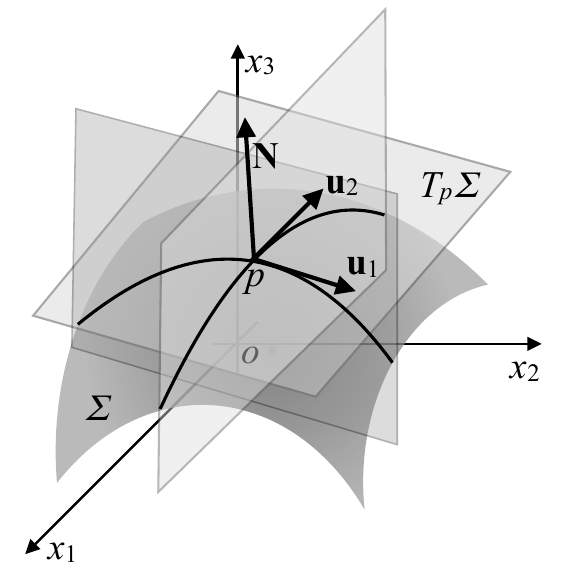}
	\caption{Principal curvatures.}
	\label{fig:38}
	\end{center}
\end{figure}

Because $\B=\B^\top$, by the spectral theorem, there exists an orthonormal basis $\{\bu_1,\bu_2\}$ of $T_p\Sigma$ such that
\bes
\B=\beta_j\bu_j\otimes\bu_j,
\ees
  with $\beta_j$ the eigenvalues of $\B$. In such a basis, by Eq. (\ref{eq:defineB}) we get 
  \bes
  \kappa_N(\bu_{i})=\frac{II(\bu_{\underline{i}},\bu_{\underline{i}})}{I(\bu_{\underline{i}},\bu_{\underline{i}})}=\frac{\bu_{\underline{i}}\cdot\B\bu_{\underline{i}}}{\bu_{\underline{i}}\cdot\g\bu_{\underline{i}}}=\frac{\bu_{\underline{i}}\cdot\g\X\bu_{\underline{i}}}{\bu_{\underline{i}}\cdot\g\bu_{\underline{i}}},\ \ \ i=1,2.
  \ees
  Then, because $\{\bu_1,\bu_2\}$ is an orthonormal basis, $\g=\I$ and
  \bes
   \kappa_N(\bu_i)=\bu_{\underline{i}}\cdot\X\bu_{\underline{i}},\ \ \ i=1,2,
  \ees
  i.e. $\X$ and $\B$ share the same eigenvectors. Moreover,  cf. Section \ref{sec:eigenvaluesvectors}, we know that the two directions $\bu_1$ and $\bu_2$ are the directions whereupon the quadratic form in the previous equation gets its maximum, $\kappa_1$, and minimum, $\kappa_2$, values, and in such a basis,
  \bes
  \X=\kappa_i\bu_i\otimes\bu_i.
  \ees
  We  call $\kappa_1$ and $\kappa_2$ the {\it principal curvatures of $\Sigma$ in $p$}  and $\bu_1,\bu_2$ the {\it principal directions of $\Sigma$ in $p$}, see Fig. \ref{fig:38}.

We call {\it Gaussian curvature $K$} the product of the principal curvatures:
 \bes
 K:=\kappa_1\kappa_2=\det\X.
 \ees
By Eq. (\ref{eq:XgB}) and the theorem of Binet, it is also
\be
\label{eq:curvgaus2}
K=\frac{\det\B}{\det\g}.
\ee

 We define {\it mean curvature $H$ of a surface}\footnote{The concept of mean curvature of a surface was introduced for the first time by Sophie Germain in her celebrated work on the elasticity of plates (1815).} $\f:\Omega\rightarrow\Sigma$ at a point $p\in\Sigma$ the mean of the principal curvatures at $p$:
 \bes
 H:=\frac{\kappa_1+\kappa_2}{2}=\frac{1}{2}\tr\X.
 \ees
 
   Of course, a change in the parameterization of a surface can change the orientation, cf. Section \ref{sec:surf1}, i.e. it can transforms $\N$ into its opposite one and, by consequence, change the sign of the second fundamental form and hence of the normal and principal curvatures. These last are hence defined  to less the sign, and  the mean curvature too, while the principal directions, umbilicality, flatness and Gaussian curvature are intrinsic to $\Sigma$, i.e. they do not depend on its parameterization. 
 
 \section{The theorem of Rodrigues}
 Then principal directions of curvature have a property which is specified by the 
 \begin{teo}{\bf (Theoreom of Rodrigues).} Let $\f(u,v)$ be a surface of class at least C$^2$ and $\bl=(\lambda_u,\lambda_v)\in T_p\Sigma$; then
 \be
 \frac{d\N(p)}{d\bl}=-\kappa_\lambda\bl 
 \ee
 if and only if $\bl$ is a principal direction; $\kappa_\lambda$ is the principal curvature relative to $\bl$.
 \begin{proof}
Let $\bl$ be a principal direction of $T_p\Sigma$. Because $\N\in\S$, then 
\be
\label{eq:rodrigues1}
\frac{d\N}{d\bl}\cdot\N=0;
\ee
 moreover, 
 \be
 \label{eq:rodrigues2}
 \frac{d\N}{d\bl}=\nabla\N\ \bl=\left[\begin{array}{ccc}0&0&0\\0&0&0\\\N_{,u}&\N_{,v}&1\end{array}\right]\left\{\begin{array}{c}\lambda_u\\\lambda_v\\0\end{array}\right\}=\N_{,u}\lambda_u+\N_{,v}\lambda_v.
 \ee
 Let $\bmu=(\mu_u,\mu_v)$ be the other principal direction of $T_p\Sigma$, then
 \bes
 \bl\cdot\bmu=0\ \rightarrow\ I(\bl,\bmu)=II(\bl,\bmu)=0.
 \ees
 Moreover,
 \bes
 \frac{d\N}{d\bl}\cdot\bmu=-II(\bl,\bmu)=0,
 \ees
 which implies, together with Eq. (\ref{eq:rodrigues1}),
 \be
 \label{eq:rodrigues3}
 \frac{d\N}{d\bl}=\alpha \bl.
\ee
 Therefore,
 \bes
 \frac{d\N}{d\bl}\cdot\bl=-II(\bl)=\alpha\bl\cdot\bl=\alpha I(\bl)
 \ees
 and finally,
 \bes
 \alpha=-\frac{II(\bl)}{I(\bl)}=-\kappa_\lambda.
 \ees
 Contrarily, if we assume Eq. (\ref{eq:rodrigues3}), as before we get $\alpha=-\kappa_\lambda$ and to end, we just need to prove that $\bl$ is a principal direction. From Eqs. (\ref{eq:rodrigues2}) and (\ref{eq:rodrigues3}), we get
 \bes
 \lambda_u\N_{,u}+\lambda_v\N_{,v}=-\kappa_\lambda(\lambda_u\f_{,u}+\lambda_v\f_{,v}).
 \ees
 Projecting this equation onto $\f_{,u}$ and $\f_{,v}$ gives the two equations
 \be
 \label{eq:rodrigues4}
 \besp
& L\lambda_u+M\lambda_v=\kappa_\lambda(E\lambda_u+F\lambda_v),\\
& M\lambda_u+N\lambda_v=\kappa_\lambda(E\lambda_u+G\lambda_v),
 \end{split}
  \ee
 with the symbols $E,F,G,L,M$ and $N$ defined in Notes \ref{note:firstform} and \ref{note:secondform} and used here for the sake of conciseness.
Let $\bw=(w_u,w_v)\in T_p\Sigma$ and consider the function 
\bes
\zeta(\bw,\kappa_\lambda)=II(\bw)-\kappa_\lambda I(\bw);
\ees
it is easy to check that $\zeta,\dfrac{\partial\zeta}{\partial w_u}$ and $\dfrac{\partial\zeta}{\partial w_v}$ take zero value for $\bw=\bl_0$, with $\bl_0$ the eigenvector of the principal direction relative to $\kappa_\lambda$, which gives the system of equations
\bes
\left\{\besp
&II(\bl_0)-\kappa_\lambda I(\bl_0)=0,\\
&\frac{\partial II(\bl_0)}{\partial w_u}-\kappa_\lambda\frac{\partial II(\bl_0)}{\partial w_u}=0,\\
&\frac{\partial II(\bl_0)}{\partial w_v}-\kappa_\lambda\frac{\partial II(\bl_0)}{\partial w_v}=0.
\end{split}
\right.
\ees
Developing the derivatives and making some standard passages, Eq.  (\ref{eq:rodrigues4}) is found again, which proves that $\bl$ is necessarily the principal direction relative to $\kappa_\lambda$.
\end{proof}
 \end{teo}
 This theorems hence states  that the derivative of $\N$ along a given direction is a vector parallel to such a direction only when this is a principal direction of curvature.

 \section{Classification of the points of a surface}
 Let $\f:\Omega\rightarrow\Sigma$ be a regular surface  and  $p\in \Sigma$ a non-planar point. Then, we say that
 \begin{itemize}
 \item $p$ is an {\it elliptic point} if $K(p)>0$;
 \item $p$ is a {\it hyperbolic point} if $K(p)<0$;
 \item $p$ is a {\it parabolic point} if $K(p)=0$. 
 \end{itemize}
 We remark that, by Eq. (\ref{eq:curvgaus2}), because $\det\g>0$, Eq. (\ref{eq:signdetg}), the value of $\det\B$ is sufficient to determine the type of a point on $\Sigma$.
 \begin{teo}
 If $p$ is an elliptical point of $\sigma$, then there exists a neighbourhood $U\in\Sigma$ of $p$ such that all the  points $q\in U$ belong to the same  half-space into which $\Eu$ is divided by the tangent plane $T_p\Sigma$. 
 \begin{proof}
For the sake of simplicity and without loss of generality, we can always chose a parameterization $\f(u,v)$  of the surface such that $p=\f(0,0)$. Expanding $\f(u,v)$ into a Taylor's series around $(0,0)$, we get the position of a point $q=\f(u,v)\in\Sigma$ in the nighbourhood of $p$ (though not indicated for the sake of brevity, all the derivatives are intended to be calculated at $(0,0)$):
\bes
\f(u,v)=\f_{,u}u+\f_{,v}v+\frac{1}{2}(\f_{,uu}u^2+2\f_{,uv}uv+\f_{,vv}v^2)+o(u^2+v^2).
\ees
The distance with sign $d(q)$ of $q\in\Sigma$ from the tangent plane $T_p\Sigma$ is the projection onto $\N$, i.e.:
\bes
\besp
d(q)&=\frac{1}{2}(\f_{,uu}u^2+2\f_{,uv}uv+\f_{,vv}v^2)\cdot\N+o(u^2+v^2)\\
&=\frac{1}{2}(B_{11}u^2+2B_{12}uv+B_{22}v^2)+o(u^2+v^2),
\end{split}
\ees
or, equivalently, once we set $\bw=u\f_{,u}+v\f_{,v}$,
\be
\label{eq:distanceq}
d(q)=\frac{1}{2}II(\bw,\bw)+o(u^2+v^2).
\ee
If $p$ is an elliptic point, the principal curvatures have the same sign because $K=\kappa_1\kappa_2>0\Rightarrow$ the sign of $II(\bw,\bw)$ does not depend upon $\bw$, i.e. upon the tangent vector. As a consequence the sign of $d(q)$ does not change with $\bw\Rightarrow\forall q\in U,\ \Sigma$ is on the same side of the tangent plane $T_p\Sigma$.
 \end{proof}
 \end{teo}
   \begin{figure}[ht]
	\begin{center}
        \includegraphics[height=.12\textheight]{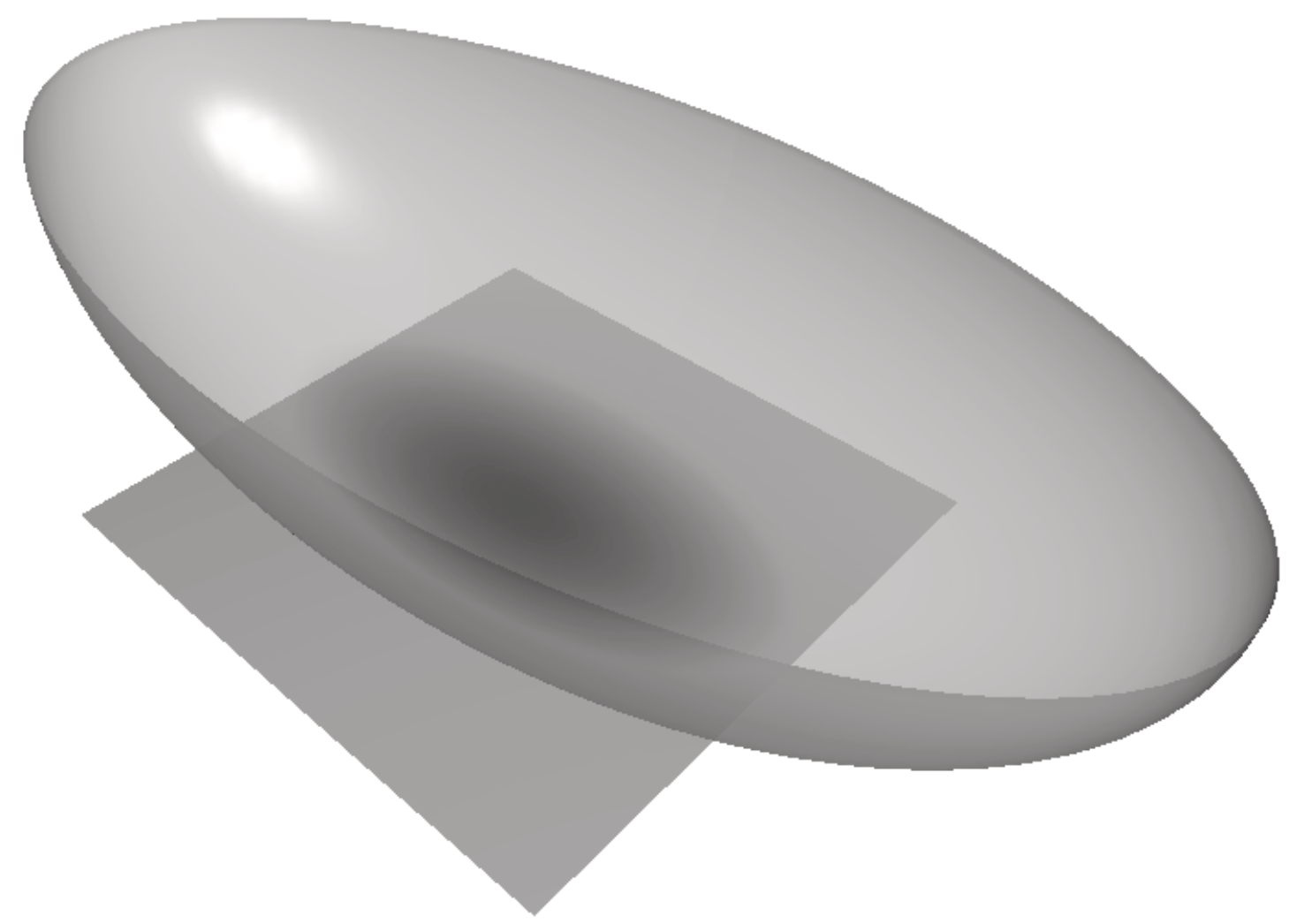}\ 
        \includegraphics[height=.12\textheight]{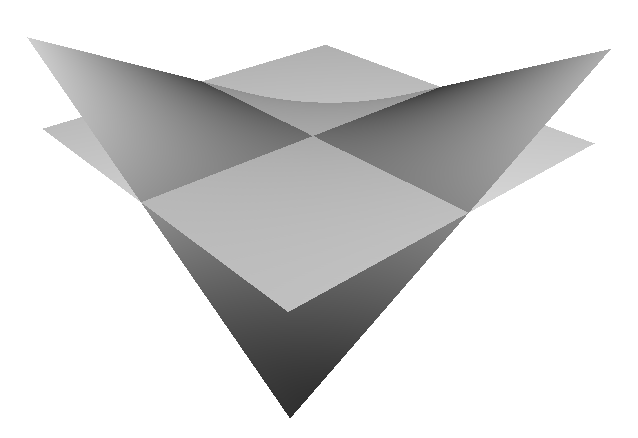}\ 
        \includegraphics[height=.12\textheight]{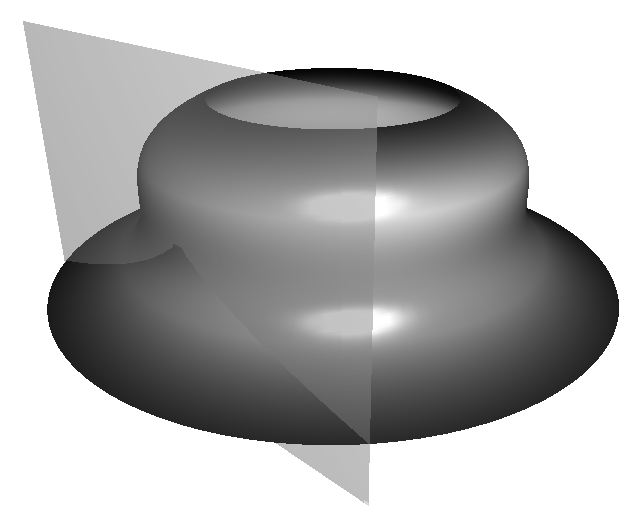}\ 
        \includegraphics[height=.12\textheight]{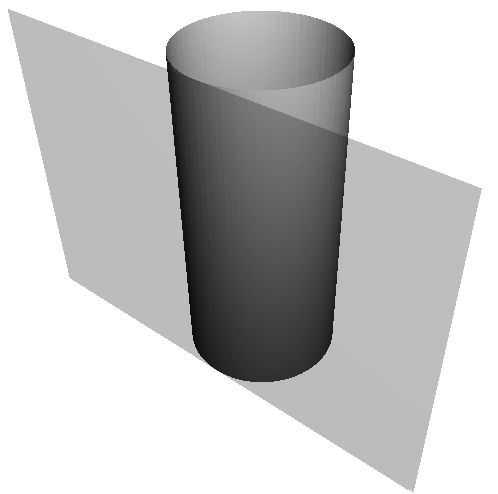}
	\caption{Elliptic (left),  hyperbolic (center), and parabolic (last two  on the right), points.}
	\label{fig:39}
	\end{center}
\end{figure}

\begin{teo}
 If $p$ is a hyperbolic point of $\Sigma$, then for each neighbourhood $U\in\Sigma$ of $p$, there are  points $q\in U$  that are in  half-spaces on the opposite sides with respect to the tangent plane $T_p\Sigma$. 
  \begin{proof}
The proof is identical to that of the previous theorem until Eq. (\ref{eq:distanceq}); now, if $p$ is a hyperbolic point, the principal curvatures have opposite  signs and by consequence $d(q)$ changes of sign at least two times in any neighbourhood $U$ of $p\Rightarrow$ there are points $q\in U$ lying in  half-spaces on the opposite sides with respect to the tangent plane $T_p\Sigma$.
 \end{proof}
 \end{teo}

 In a parabolic point, there are different possibilities: $\Sigma$ is on one side of the space with respect to $T_p\Sigma$, like for the case of a cylinder, or not, like, as an example, for the points $(0,v)$ of the surface, see Fig. \ref{fig:39},
 \bes
 \left\{\begin{array}{l}x=(u^3+2)\cos v,\\y=(u^3+2)\sin v,\\z=-u.\end{array}\right.
 \ees

\begin{figure}[ht]
	\begin{center}
        \includegraphics[height=.15\textheight]{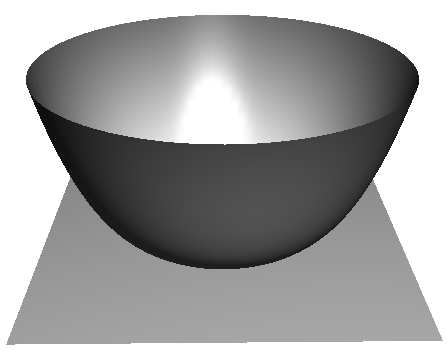}\ \ \
        \includegraphics[height=.15\textheight]{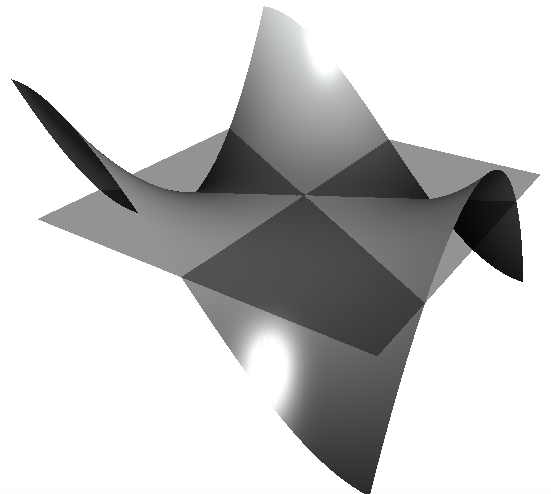}
	\caption{Two different planar points.}
	\label{fig:43}
	\end{center}
\end{figure}

This is also the case  for planar points: e.g., the point $(0,0,0)$ is a planar point for both the surfaces
\bes
z=x^4+y^4,\ \ \ z=x^3-3xy^2,
\ees
 but in the first case, all the surface is on one side of the tangent plane, while it is on both sides for the second case (the so-called {\it monkey's saddle}), see Fig.\ref{fig:43}.
 
 \section{Developable surfaces}
Let us now consider a ruled surface $\f:\Omega\rightarrow\Sigma$ as in Eq. (\ref{eq:ruledsurf});   then
\bes
\f_{,u}=\bg'+v\bl',\ \ \f_{,v}=\bl,\ \ \f_{,u}\times\f_{,v}=\bg'\times\bl+v\bl'\times\bl,\ \ \f_{,uv}=\bl',\ \ \f_{,vv}=\bo.
\ees
Consequently, $B_{22}=\N\cdot\f_{,vv}=0\ \Rightarrow\ \det\B=-B_{12}^2$: The points of $\Sigma$ are hyperbolic or parabolic. Namely, the parabolic points are those with 
 \bes
 B_{12}=\N\cdot\f_{,uv}=\frac{\f_{,u}\times\f_{,v}}{|\f_{,u}\times\f_{,v}|}\cdot\f_{,uv}=0\ \iff\ (\bg'\times\bl+v\bl'\times\bl)\cdot\bl'=\bg'\times\bl\cdot\bl'=0.
 \ees
 We remark that the ruled surfaces made of parabolic points have {\it null Gaussian curvature everywhere: $K=0$}.
 
 Let us consider ruled surfaces having only parabolic points; then, we have the following
 \begin{teo}
 \label{teo:surfruled}
 For a ruled surface $\f(u,v)=\bg(u)+v\bl(u)$, the following are equivalents:
 \begin{enumerate}[i.]
 \item $\bg',\bl,\bl'$ are linearly dependent;
 \item $\N_{,v}=\bo$.
 \end{enumerate}
\begin{proof}
Condition $ii$ implies that  $\N$ does not change along a straight line lying on the ruled surface $\Rightarrow
\f_{,u}\times\f_{,v}=\bg'\times\bl+v\bl'\times\bl$ does not depend on $v$ as well. This is possible $\iff\bg'\times\bl$ and $\bl'\times\bl$ are linearly dependent, i.e. $\iff$
\bes
(\bg'\times\bl)\times(\bl'\times\bl)=(\bl'\times\bl\cdot\bg')\bl-(\bl'\times\bl\cdot\bl)\bg'=(\bl'\times\bl\cdot\bg')\bl=\bo,
\ees
i.e. when $\bl,\bl'$ and $\bg'$ are coplanar, which proves the thesis.
\end{proof}
 \end{teo}
We say that a ruled surface is {\it developable} if one of the conditions of  Theorem \ref{teo:surfruled} is satisfied. A developable surface is a surface that can be {\it flattened without distortion onto a plane}, i.e. it can be bent without stretching or shearing or, vice-versa, it can be obtained by transforming a plane. We remark that only ruled surfaces are developable (but not all the ruled surfaces are developable).

\begin{figure}[ht]
	\begin{center}
        \includegraphics[height=.1\textheight]{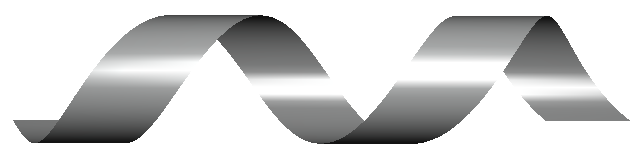}
	\caption{The ruled surface of the tangents to a cylindrical helix.}
	\label{fig:45}
	\end{center}
\end{figure}

 It is immediate to check that a cylinder or a cone are developable surfaces, while the helicoid, the hyperbolic hyperboloid or the hyperbolic paraboloid are not.  Another classical example of developable surface is the {\it ruled surface of the tangents to a curve}: Let $\bg(t):G\subset\R\rightarrow\Eu$ be a regular smooth curve; then the ruled surface of the tangents to $\bg$ is the surface $\f(u,v):G\times\R\rightarrow\Sigma$ defined by
 \bes
 \f(u,v)=\bg(u)+v\bg'(u).
 \ees
 Fig. \ref{fig:45} shows the ruled surface of the tangents to a cylindrical helix. 

\section{Points of a surface of revolution}
Let us now consider a  surface of revolution $\f:\Omega\rightarrow\Sigma_\gamma$ as in Eq. (\ref{eq:revolsurf}) and, for the sake of simplicity, let $u$ be the natural parameter of the curve $\bg(u)$ in Eq. (\ref{eq:curvesurfrevol}) generating the surface. Then
\bes
\phi'^2(u)+\psi'^2(u)=1,\ \ \ \ \psi''(u)\phi'(u)-\psi'(u)\phi''(u)=c(u).
\ees
 We can then calculate:
 \begin{itemize}
 \item the vectors of the natural basis:
 \bes
 \f_{,u}=\left\{\begin{array}{c}\phi'(u)\cos v\\\phi'(u)\sin v\\\psi'(u)\end{array}\right\},\ \ \ \ \  \f_{,v}=\left\{\begin{array}{c}-\phi(u)\sin v\\\phi(u)\cos v\\0\end{array}\right\};
 \ees
 \item the normal to the surface
 \bes
 \N=\left\{\begin{array}{c}-\psi'(u)\cos v\\-\psi'(u)\sin v\\\phi'(u)\end{array}\right\};
 \ees
 \item the metric tensor (i.e. the first fundamental form):
 \bes
 \g=\left[\begin{array}{cc}1&0\\0&\phi^2(u)\end{array}\right];
 \ees
 \item the second derivatives of $\f$:
 \bes
 \f_{,uu}=\left\{\begin{array}{c}\phi''(u)\cos v\\\phi''(u)\sin v\\\psi''(u)\end{array}\right\},\ \ \ 
  \f_{,uv}=\left\{\begin{array}{c}-\phi'(u)\sin v\\\phi'(u)\cos v\\0\end{array}\right\},\ \ \ 
   \f_{,vv}=\left\{\begin{array}{c}-\phi(u)\cos v\\-\phi(u)\sin v\\0\end{array}\right\};
 \ees
 \item tensor $\B$ (i.e. the second fundamental form):
 \bes
 \B=\left[\begin{array}{cc}c(u)&0\\0&\phi(u)\psi'(u)\end{array}\right];
 \ees
 \item the Gaussian curvature $K$:
 \bes
 K=\det\X=\frac{\det\B}{\det\g}=\frac{c(u)\psi'(u)}{\phi(u)}.
 \ees
 \end{itemize}
 Therefore, points of $\Sigma_\gamma$ where $c(u)$ and $\psi'(u)$ have the same sign are elliptic, but hyperbolic otherwise\footnote{Recall that in a revolution surface, $\phi(u)>0\ \forall u$.}. Parabolic points correspond to inflexion points of $\bg(u)$ if $c(u)=0$, or to points with horizontal tangent to $\bg(u)$ if $\psi'(u)=0$.
 
 As an example, let us consider the case of the pseudo-sphere, Eq. (\ref{eq:pseudosphere}). Then, 
 \bes
 \phi(u)=\sin u,\ \ \ \psi(u)=\cos u+\ln\tan\frac{u}{2}.
 \ees
 Some simple calculations give ($u$ is not the arc-length of $\bg(u)$, which  actually is a tractrix, see Exercise 9, Chapter 4; hence $|\bg'(u)|=\sqrt{\phi'^2(u)+\psi'^2(u)}\neq1$)
 \bes
 \psi'(u)=-\sin u+\frac{1}{\sin u},\ \ \ c(u)=-|\tan u|;
 \ees
 as a  consequence
 \bes
 K=\frac{c(u)\psi'(u)}{\phi(u)\sqrt{\phi'^2(u)+\psi'^2(u)}}=-\frac{(-\sin u+\frac{1}{\sin u})|\tan u|}{\sin u|\cot u|}=-1.
 \ees
 Finally, $K=const.=-1$, which is the reason for the name of this surface.
 
 \section{Lines of curvature, conjugated directions, asymptotic directions}
 A {\it line of curvature} is a curve on a surface with the property of being tangent, at each point, to a principal direction.  
 \begin{teo}
 The lines of curvature of a surface are the solutions to the differential equation
 \bes
 X_{21}u'^2+(X_{22}-X_{11})u'v'-X_{12}v'^2=0.
 \ees
 \begin{proof}
 A curve $\bg(t):G\subset\R\rightarrow\Sigma\subset\Eu$ is a line of curvature $\iff$
 \bes
 \bg'(t)=\f_{,u}u'+\f_{,v}v'
 \ees
 is an eigenvector of $\X(t)\ \forall t$, i.e. $\iff$ there exists a function $\mu(t)$ such that
 \bes
 \X(t)\bg'(t)=\mu(t)\bg'(t)\ \ \forall t.
 \ees
 In the natural basis of $T_p\Sigma$, this condition reads as (we omit the dependence upon $t$ for the sake of conciseness)
 \bes
 \left[\begin{array}{cc}X_{11}&X_{12}\\X_{21}&X_{22}\end{array}\right]\left\{\begin{array}{c}u'\\v'\end{array}\right\}=\mu\left\{\begin{array}{c}u'\\v'\end{array}\right\},
 \ees
 which is satisfied $\iff$ the two vectors at the left- and right-hand sides are proportional, i.e. $\iff$
 \bes
 \det\left[\begin{array}{cc}X_{11}u'+X_{12}v'&u'\\X_{21}u'+X_{22}v'&v'\end{array}\right]=0\ \rightarrow\  X_{21}u'^2+(X_{22}-X_{11})u'v'-X_{12}v'^2=0.
 \ees
 \end{proof}
 \end{teo}
 As a corollary, if $\X$ is diagonal, then the coordinate lines are at the same time,  principal directions and lines of curvature. 
 \begin{teo}
A curve $\bg(u):G\subset\R\rightarrow\Sigma$ is a line of curvature $\iff$  the surface
\be
\label{eq:rulnorm}
\f(u,v)=\bg(u)+v\N(\bg(u)),
\ee
is developable.
\begin{proof}
From Theorem \ref{teo:surfruled}, $\f(u,v)$ is developable $\iff\bg'\cdot\N\times\N'=0$. Because $\bg'$ and $\N'\in T_p\Sigma$, which is orthogonal to $\N$, the surface will be developable $\iff\ \bg'\times\N'=\bo$. Moreover, writing
\bes
\bg'=\f_{,u}u'+\f_{,v}v'
\ees
it is
\bes
\N'=\N_{,u}u'+\N_{,v}v'=-\mathcal{L}_W(\bg'),
\ees
hence $\f(u,v)$ is developable $\iff\ \mathcal{L}_W(\bg')\times\bg'=\bo$, i.e. when $\bg'$ is a principal direction.
\end{proof}
 \end{teo}
The curve in Eq. (\ref{eq:rulnorm}) is called the {\it ruled surface of the normals}.
 
 Let $p$ be a non-planar point of a surface $\f:\Omega\rightarrow\Sigma$ and $\bv_1,\bv_2$ two vectors of $T_p\Sigma$. We say that $\bv_1$ and $\bv_2$ are {\it conjugated} if $II(\bv_1,\bv_2)=0$. The directions corresponding to $\bv_1$ and $\bv_2$ are called {\it conjugated directions}. Hence, the principal directions at a point $p$ are conjugated; if $p$ is an umbilical point, any two orthogonal directions are conjugated.
 
 The direction of a vector $\bv\in T_p\Sigma$ is said to be {\it asymptotic} if it is {\it autoconjugated}, i.e. if $II(\bv,\bv)=0$. An asymptotic direction is hence a direction where the normal curvature is null. In a hyperbolic point, there are two asymptotic directions, in a parabolic point only one and in an elliptic point, there are not asymptotic directions. An {\it asymptotic line} is a curve on a surface with the property of being tangent at every point to an asymptotic direction.
 The asymptotic lines are the solution of the differential equation
 \bes
 II(\bg',\bg')=0\ \rightarrow\ B_{11}u'^2+2B_{12}u'v'+B_{22}v'^2=0;
 \ees
 in particular, if $B_{11}=B_{22}=0$ and $\B\neq\O$, then the coordinate lines are asymptotic lines. Asymptotic lines exist only in the regions where $K\leq0$.

 \section{The Dupin's conical curves}
\begin{figure}[ht]
	\begin{center}
        \includegraphics[width=.9\textwidth]{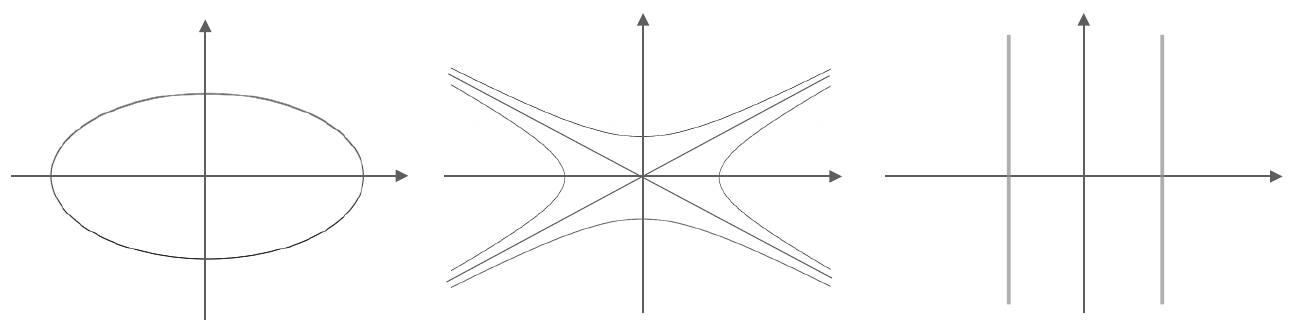}
	\caption{The conical curves of Dupin; from the left: elliptic, hyperbolic and parabolic points.}
	\label{fig:46}
	\end{center}
\end{figure}

 The {\it conical curves of Dupin} are the real curves in $T_p\Sigma$ whose equations are
 \bes
 II(\bv,\bv)=\pm1,\ \ \ \bv\in\S.
 \ees
 Let $\{\bu_1,\bu_2\}$ be the basis of the principal directions. Using polar coordinates, we can write
 \bes
 \bv=\rho\e_\rho,\ \ \ \e_\rho=\cos\theta\bu_1+\sin\theta\bu_2.
 \ees
 Therefore,
 \bes
 II(\bv,\bv)=\rho^2II(\e_\rho,\e_\rho)=\rho^2\kappa_N(\e_\rho),
 \ees
 and the conicals' equations are
 \bes
 \rho^2(\kappa_1\cos^2\theta+\kappa_2\sin^2\theta)=\pm1.
 \ees
 With the Cartesian coordinates $\xi=\rho\cos\theta,\eta=\rho\sin\theta$, we get
 \bes
 \kappa_1\xi^2+\kappa_2\eta^2=\pm1.
 \ees
 The type of conical curves depend upon the kind of point on $\Sigma$:
 \begin{itemize}
 \item {\it Elliptical points}: The principal curvatures have the same sign $\rightarrow$ one of the conical curves is an ellipse, the other one the null set (actually, it is not a real curve).
 \item {\it Hyperbolic points}: The principal curvatures have opposite signs $\rightarrow$ the conical curves are conjugated hyperbolae whose asymptotes coincide with the asymptotic directions.
 \item {\it Parabolic points}: At least one of the principal curvatures is null $\rightarrow$ one of the conical curves degenerates into a couple of parallel straight lines, corresponding to the asymptotic direction, the other one is the null set. 
 \end{itemize}
 The three possible cases are depicted in Fig. \ref{fig:46}

\section{The Gauss-Weingarten equations}
 Let $\f:\Omega\rightarrow\Sigma$ be a  surface; for any point $p\in\Sigma$, consider the basis $\{\f_{,u},\f_{,v},\N\}$, also called the {\it Gauss' basis}. It is the equivalent of the Frenet-Serret basis for the surfaces. We want to calculate the derivatives of the vectors of this basis, i.e. we want to obtain, for the surfaces, something equivalent to the Frenet-Serret equations. 
 
 $\N\in\S$ and  $\N\cdot\f_{,u}=\N\cdot\f_{,v}=0$, but, in general, $\f_{,u},\f_{,v}\notin\S$ and $\f_{,u}\cdot\f_{,v}\neq0$. In other words, we are dealing with a case of non-orthogonal (curvilinear) coordinates. So, if $w$ is the coordinate along the normal $\N$, let us call, for the sake of convenience,
 \bes
 u=z^1,\ v=z^2,
 \ees
 while, for the vectors,
 \bes
 \f_{,u}=\f_{,1}=\g_1,\ \f_{,v}=\f_{,2}=\g_2,
 \ees
 with $\g_1,\g_2$ exactly the $\g$-vectors of the coordinate lines on $\Sigma$. Then  (no summation on $i$ in the following equations),
 \bes
 \besp
 &\frac{\partial\g_i}{\partial z^j}\cdot\g_i=\frac{1}{2}\frac{\partial (\g_i\cdot\g_i)}{\partial z^j}=\frac{1}{2}\frac{\partial g_{ii}}{\partial z^j},\\
&\frac{\partial \g_i}{\partial z^i}\cdot\g_j=\frac{\partial(\g_i\cdot\g_j)}{\partial z^i}-\frac{\partial\g_i}{\partial z^j}\cdot\g_i=\frac{\partial g_{ij}}{\partial z^i}-\frac{1}{2}\frac{\partial g_{ii}}{\partial z^j};
 \end{split}
 \ \ \ \ i,j=1,2;
 \ees
for the last equation we have used the identity
 \bes
 \frac{\partial \g_j}{\partial z^i}=\f_{,ji}=\f_{,ij}=\frac{\partial\g_i}{\partial z^j}, \ \ \ \ i,j=1,2.
 \ees
Using Eq. (\ref{eq:christoffavecg}), it can be proved that it is also\footnote{The proof is rather cumbersome and it is omitted here; in many texts on differential geometry, the Christoffel symbols are just introduced in this way, as the projection of the derivatives of vectors $\g_i$s onto the same vectors, i.e. as the coefficients of the Gauss equations.}
\bes
\frac{\partial\g_i}{\partial z^j}\cdot\g_h=\Gamma^h_{ij} \ \ \ \ i,j,h=1,2.
\ees
 Moreover, by Eq. (\ref{eq:composB}),
 \bes
 \frac{\partial\g_i}{\partial z^j}\cdot\N=\f_{,ij}\cdot\N=B_{ij}\ \ \ \ i,j=1,2,
  \ees
and by Eqs. (\ref{eq:Wo1}), (\ref{eq:Wo2}),
 \bes
 \frac{\partial\N}{\partial z^i}\cdot\g_j=-\mathcal{L}_W(\g_i)\cdot\g_j=-\X\g_i\cdot\g_j=-X_{ji},\ \ \ i,j=1,2,
 \ees
 while, because $\N\in\S$, then from Eq. (\ref{eq:constvect}),
 \bes
 \frac{\partial\N}{\partial z^i}\cdot\N=0\ \ \ \forall i=1,2.
 \ees
 Finally, the decomposition of the derivatives of the vectors of the basis $\{\f_{,u},\f_{,v},\N\}$ onto these same vectors gives the  equations
 \be
 \label{eq:gaussweingarten}
 \besp
 &\frac{\partial \g_i}{\partial z^j}=\Gamma^h_{ij}\g_h+B_{ij}\N,\\
 &\frac{\partial\N}{\partial z^j}=-X_{ij}\g_i,
 \end{split}\ \ \ i,j=1,2;
 \ee
 these are the {\it Gauss-Weingarten equations} (the first one is due to Gauss and the second  to Weingarten).
 
 Moreover, if we make the scalar product of the Gauss equations by $\g_1$ and $\g_2$, i.e.
 \bes
\g_k\cdot \frac{\partial \g_i}{\partial z^j}=\g_k\cdot(\Gamma^h_{ij}\g_h+B_{ij}\N),\ \ i,j,k=1,2,
 \ees
 we get the following three systems of equations:
 \be
 \label{eq:syst1}
 \left\{
 \besp
 &\Gamma^1_{11}g_{11}+\Gamma^2_{11}g_{21}=\frac{1}{2}\frac{\partial g_{11}}{\partial z^1},\\
& \Gamma^1_{11}g_{12}+\Gamma^2_{11}g_{22}=\frac{\partial g_{12}}{\partial z^1}-\frac{1}{2}\frac{\partial g_{11}}{\partial z^2};
\end{split}
\right.
 \ee
  \smallskip
 \be
  \label{eq:syst2}
 \hspace{-13mm}
 \left\{
 \besp
 &\Gamma^1_{12}g_{11}+\Gamma^2_{12}g_{21}=\frac{1}{2}\frac{\partial g_{11}}{\partial z^2},\\
& \Gamma^1_{12}g_{12}+\Gamma^2_{12}g_{22}=\frac{1}{2}\frac{\partial g_{22}}{\partial z^1};
\end{split}
\right.
 \ee
 \medskip
 \be
  \label{eq:syst3}
\hspace{1mm}
 \left\{
 \besp
 &\Gamma^1_{22}g_{11}+\Gamma^2_{22}g_{21}=\frac{\partial g_{12}}{\partial z^2}-\frac{1}{2}\frac{\partial g_{22}}{\partial z^1},\\
& \Gamma^1_{22}g_{12}+\Gamma^2_{22}g_{22}=\frac{1}{2}\frac{\partial g_{22}}{\partial z^2}.
\end{split}
\right.
 \ee
 The determinant of each one of these systems is simply $\det\g\neq0\rightarrow$ it is possible to express the Christoffel symbols as functions of the $g_{ij}$s and of their derivatives, i.e. as functions of the first fundamental form (hence, of the metric tensor).

  \section{The \textit{Theorema Egregium}}
 
 The following theorem is a fundamental result due to Gauss:
 \begin{teo}{\bf (Theorema Egregium).} The Gaussian curvature $K$ of a surface $\f(u,v):\Omega\rightarrow\Sigma$ depends only upon the first fundamental form of $\f$.
 \begin{proof}
Let us write the identity
\bes
\frac{\partial^2 \g_1}{\partial z^1\partial z^2}=\frac{\partial^2 \g_1}{\partial z^2\partial z^1}
\ees
using the Gauss equations (\ref{eq:gaussweingarten})$_1$:
\bes
\besp
&\Gamma^1_{11}\g_{1,2}+\Gamma^2_{11}\g_{2,2}+B_{11}\N_{,2}+\Gamma^1_{11,2}\g_1+\Gamma^2_{11,2}\g_2+B_{11,2}\N=\\
&\Gamma^1_{12}\g_{1,1}+\Gamma^2_{12}\g_{2,1}+B_{12}\N_{,1}+\Gamma^1_{12,1}\g_1+\Gamma^2_{12,1}\g_2+B_{12,1}\N,
\end{split}
\ees
where, for the sake of shortness, we have abridged $\dfrac{\partial(\cdot)}{\partial z^j}$ by $(\cdot)_{,j}$.
Then, we use again Eqs. (\ref{eq:gaussweingarten}) to express $\g_{1,1},\g_{1,2},\g_{2,2},\N_{,1}$ and $\N_{,2}$; after doing that and equating to $0$ the coefficient of $\g_2$, we get
\bes
B_{11}X_{22}-B_{12}X_{21}=\Gamma^1_{11}\Gamma^2_{12}+\Gamma^2_{11}\Gamma^2_{22}+\Gamma^2_{11,2}-\Gamma^1_{12}\Gamma^2_{11}-\Gamma^2_{12}\Gamma^2_{12}-\Gamma^2_{12,1};
\ees
from Eq. (\ref{eq:defineB}), we get that
\bes
B_{11}=g_{11}X_{11}+g_{12}X_{21},\ \ \ B_{12}=g_{11}X_{12}+g_{12}X_{22},
\ees
which, injected into the previous equation, gives
\be
\label{eq:theoegreg}
g_{11}\det\X=\Gamma^1_{11}\Gamma^2_{12}+\Gamma^2_{11}\Gamma^2_{22}+\Gamma^2_{11,2}-\Gamma^1_{12}\Gamma^2_{11}-\Gamma^2_{12}\Gamma^2_{12}-\Gamma^2_{12,1}.
\ee
Equating to zero the coefficient of $\g_1$, a similar expression can also be get  for $g_{12}$. Because $\g$ is positive definite, it is not possible that $g_{11}=g_{12}=0$. So, remembering that $K=\det\X$ and the result of the previous section, we see that it is possible to express $K$ through the coefficients of the first fundamental form and of its derivatives. 
 \end{proof}
 \end{teo}

 \section{Minimal surfaces}
\begin{figure}[ht]
	\begin{center}
        \includegraphics[height=.13\textheight]{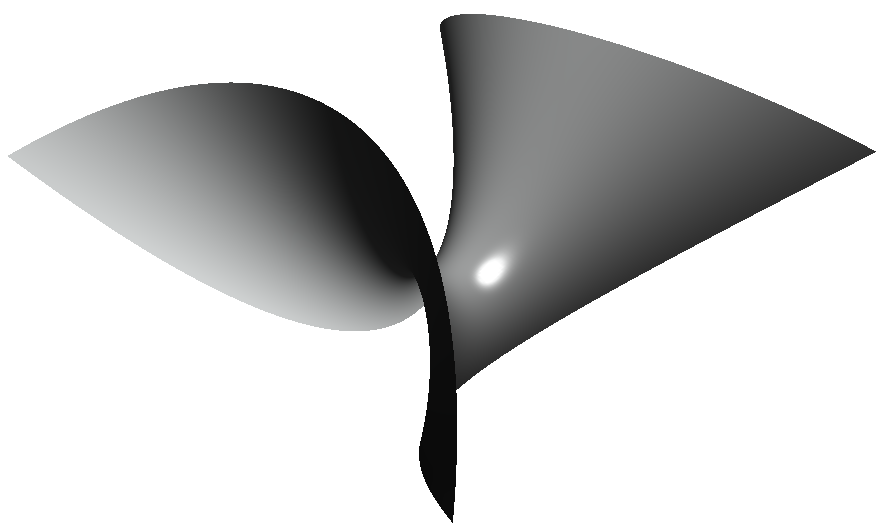}\ \
          \includegraphics[height=.13\textheight]{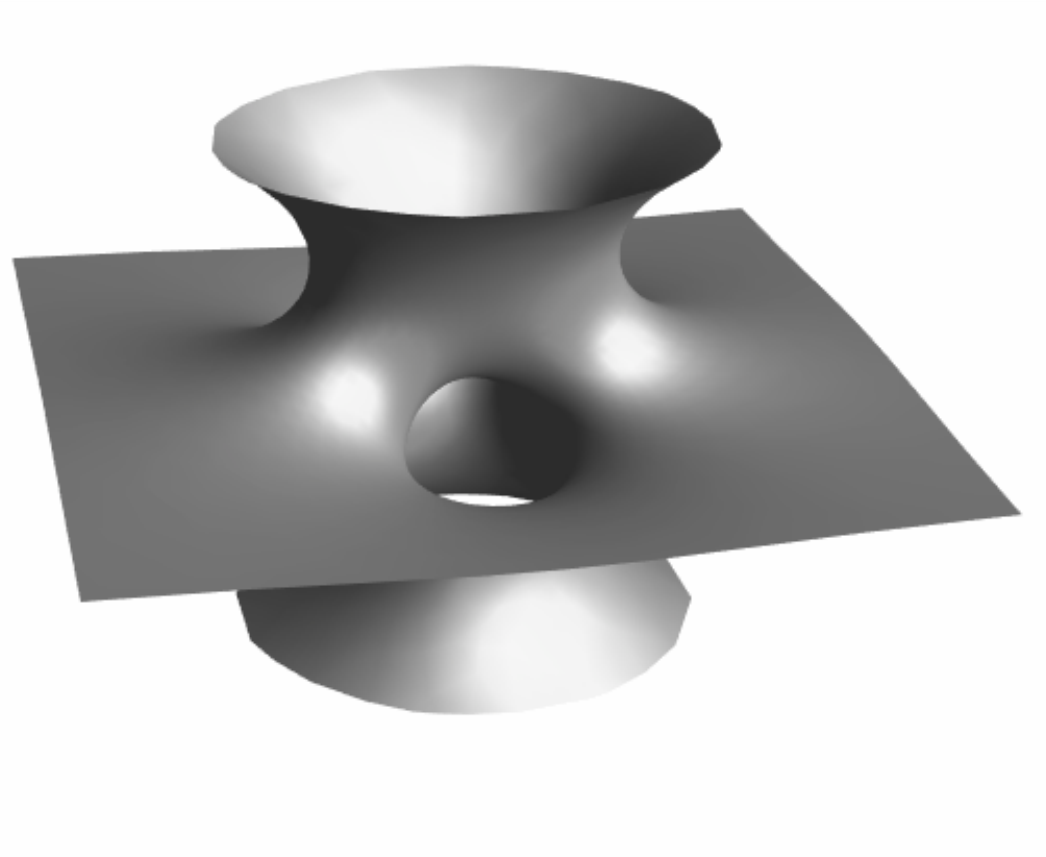}\ \
            \includegraphics[height=.13\textheight]{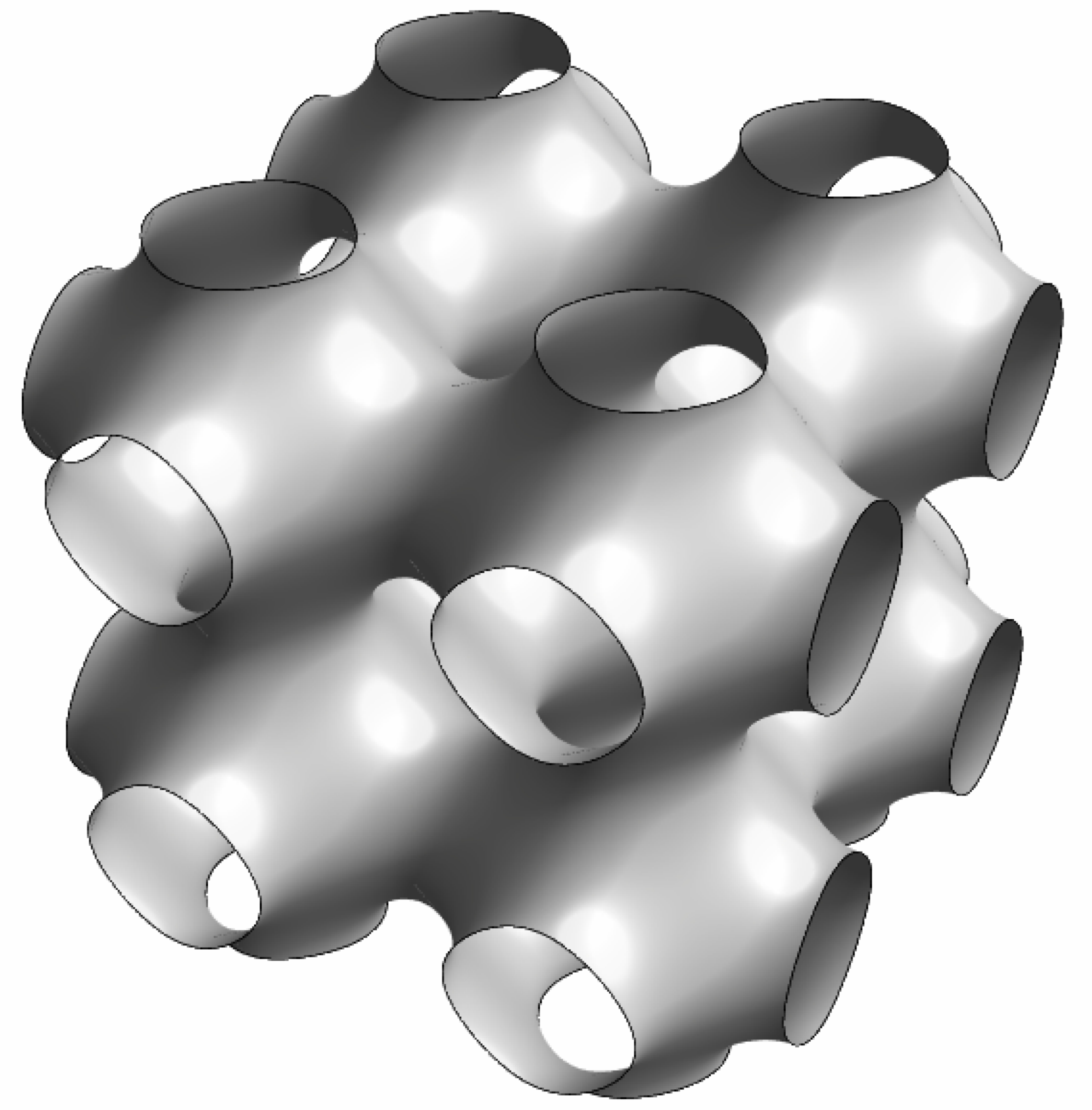}
	\caption{From the left, the minimal surfaces of Enneper, Costa and Schwarz.}
	\label{fig:47}
	\end{center}
\end{figure}

 A {\it minimal surface} is a surface $\f:\Omega\rightarrow\Sigma$ having the mean curvature $H=0\ \forall p\in\Sigma$. Typical minimal surfaces are the catenoid and the helicoid\footnote{Minimal surfaces have some interesting applications in the mechanics of tensile structures composed of prestressed membranes. Also, it can be shown that a soap film, when not bounding a closed region, takes the form of a minimal surface.}. Other minimal surfaces are the {\it Enneper's surface}
 \bes
 \left\{\begin{array}{l}x_1=u-\dfrac{u^3}{3}+uv^2,\\x_2=v-\dfrac{v^3}{3}+u^2v,\\
 x_3=u^2-v^2,\end{array}\right.
 \ees
 the {\it Costa's} and the {\it Schwarz's} surfaces, Fig. \ref{fig:47}.
\begin{teo}
The non-planar points of a minimal surface are hyperbolic.
\begin{proof}
 This is a direct consequence of the definition of mean curvature $H$ and of hyperbolic points: $H=0\iff\kappa_1\kappa_2<0$.
\end{proof}
\end{teo}

 Let $\f:\Omega\rightarrow\Sigma$ be a regular surface and $Q$ a subset of $\Omega$ with its boundary $\partial Q$ a closed regular curve in $\Omega$; then, $R=\f(Q)\subset\Sigma$ is a {\it simple region} of $\Sigma$. Let $h:Q\rightarrow\R$ be a smooth function. Then, we call {\it normal variation of R} the map $\bphi:Q\times(-\epsilon,\epsilon)\rightarrow\Eu$ defined by
 \bes
 \bphi(u,v,t)=\f(u,v)+t\ h(u,v)\N(u,v).
 \ees
 For each fixed $t$, $\bphi(u,v,t)$ is a  surface with
 \bes
 \besp
 &\bphi_{,u}(u,v,t)=\f_{,u}(u,v)+t\ h(u,v)\N_{,u}(u,v)+t\ h_{,u}(u,v)\N(u,v),\\
 &\bphi_{,v}(u,v,t)=\f_{,v}(u,v)+t\ h(u,v)\N_{,v}(u,v)+t\ h_{,v}(u,v)\N(u,v).
 \end{split}
 \ees
If the first fundamental form of $\f$ is represented by the metric tensor $\g$, we look for the metric tensor $\g^t$ representing the first fundamental form of $\bphi(u,v,t)\ \forall t$:
\bes
\besp
&g_{11}^t=\bphi_{,u}\cdot\bphi_{,u}=g_{11}+2t\ h\ \f_{,u}\cdot\N_{,u}+t^2(h^2\N_{,u}^2+h_{,u}^2),\\
&g_{12}^t=\bphi_{,u}\cdot\bphi_{,v}=g_{12}+t\ h(\f_{,u}\cdot\N_{,v}+\f_{,v}\cdot\N_{,u})+t^2(h^2\N_{,u}\cdot\N_{,v}+h_{,u}h_{,v}),\\
&g_{22}^t=\bphi_{,v}\cdot\bphi_{,v}=g_{22}+2t\ h\ \f_{,v}\cdot\N_{,v}+t^2(h^2\N_{,v}^2+h_{,v}^2),
\end{split}
\ees
and by Eq. (\ref{eq:composB}),
\bes
\besp
&g_{11}^t=g_{11}-2t\ h\ B_{11}+t^2(h^2\N_{,u}^2+h_{,u}^2),\\
&g_{12}^t=g_{12}-2t\ h\ B_{12}+t^2(h^2\N_{,u}\cdot\N_{,v}+h_{,u}h_{,v}),\\
&g_{22}^t=g_{22}-2t\ h\ B_{22}+t^2(h^2\N_{,v}^2+h_{,v}^2),
\end{split}
\ees
 whence
 \bes
 \det\g^t=\det\g-2th(g_{11}B_{22}-2g_{12}B_{12}+g_{22}B_{11})+o(t^2).
 \ees
 Then, by Eq. (\ref{eq:XgB}), we get easily that
 \bes
 g_{11}B_{22}-2g_{12}B_{12}+g_{22}B_{11}=2H\det\g,
 \ees
 so that
 \bes
  \det\g^t=\det\g(1-4thH)+o(t^2).
 \ees
 We can now calculate the area $\mathcal{A}(t)$of the simple region $R^t=\bphi(u,v,t)$ corresponding to the subset $Q$:
 \bes
 \mathcal{A}^t=\int_Q\sqrt{\det\g(1-4thH)+o(t^2)}dudv;
 \ees
 For $\epsilon\ll1\ \mathcal{A}^t$ is differentiable and its derivative for $t=0$ is
 \bes
 \left[\frac{d\mathcal{A}^t}{dt}\right]_{t=0}=-\int_Q2hH\sqrt{\det\g}dudv.
 \ees
 \begin{teo}
 A surface $\f:\Omega\rightarrow\Sigma$ is minimal $\iff\left[\dfrac{d\mathcal{A}^t}{dt}\right]_{t=0}=0\ \ \forall R\subset\Sigma$ and for each normal variation.
 \begin{proof}
 If $\f$ is minimal, the condition is clearly satisfied ($H=0$). Conversely, let us suppose that $\exists p=\f(\overline{u},\overline{v})\in\Sigma| H(p)\neq0$. Consider $r_1,r_2\in\R$ such that $|H|\neq0$ in the circle $D_2$ with center $p$ and radius $r_2$ and $|H|>\frac{1}{2}|H(p)|$ in the circle $D_1$ with center $p$ and radius $r_1$. 
 Then, we chose a smooth function $h(u,v)$ such that $i)\ h=H$ inside $D_1$, $ii)\ hH>0$ inside $D_2$ and $iii)\ h=0$ outside $D_2$. For the normal variation defined by such $h(u,v)$ we have
 \bes
 \besp
 -\left[\frac{d\mathcal{A}^t}{dt}\right]_{t=0}&=\int_{D_2}2hH\sqrt{\det\g}du\ dv\geq\int_{D_1}2H^2\sqrt{\det\g}du\ dv\\
 &\geq\int_{D_1}\frac{H(p)^2}{2}\sqrt{\det\g}du\ dv=\frac{H(p)^2}{2}\mathcal{A}(\f(D_1))\\
 &\Rightarrow\left[\dfrac{d\mathcal{A}^t}{dt}\right]_{t=0}<0
 \end{split}
 \ees 
which contradicts the hypothesis.
 \end{proof} 
 \end{teo}
 The meaning of this theorem justifies the name of minimal surfaces: These are the surfaces that have  the minimal area among all the surfaces that share the same boundary. 
 
 \section{Geodesics}
Let $\f:\Omega\rightarrow\Sigma$ be a surface and $\bg(t): G\subset\R\rightarrow\Sigma$ a curve on $\Sigma$. A vector function $\bw(t):G\rightarrow T_{\bg(t)}\Sigma$ is called a {\it vector field\footnote{More correctly, $\bw(t)$ is a curve of vectors; however, it is normally called a vector field along a curve.} along} $\bg(t)$. We call {\it covariant derivative of $\bw(t)$ along $\bg(t)$} the vector field $D_{\bg}\bw(t): G\rightarrow\Ve$ defined as\footnote{The operator that gives the projection of $\bw$ onto a vector orthogonal to  $\N\in\S$, i.e. onto $T_{\bg(t)}\Sigma$, is $\I-\N\otimes\N$, cf. Exercise 2, Chapter 2.}
\bes
D_{\bg}\bw:=(\I-\N\otimes\N)\frac{d\bw}{dt},
\ees
i.e. the  projection of the derivative of $\bw$ onto $T_{\bg(t)}\Sigma$. 
It is always possible to decompose $\bw(t)$ into its components in the natural basis $\{\f_{,u},\f_{,v}\}$:
\bes
\bw(t)=w_1(t)\f_{,u}(\bg(t))+w_2(t)\f_{,v}(\bg(t)).
\ees
Differentiating, we get (a prime here denotes the derivative with respect to $t$)
\bes
\bw'=w'_1\f_{,u}+w_1(\f_{,uu}u'+\f_{,uv}v')+w'_2\f_{,v}+w_2(\f_{,uv}u'+\f_{,vv}v')
\ees
and using the Gauss equations, Eq. (\ref{eq:gaussweingarten})$_1$, we obtain (summation on the dummy indexes and $u_1$ stands for $u$ while $u_2$ for $v$)
\bes
\bw'=\f_{,k}w'_k+(\Gamma^k_{ij}\f_{,k}+B_{ij}\N)w_iu'_j,\ \ \ i,j,k=1,2,
\ees
so that the projection onto $T_{\bg(t)}\Sigma$, i.e.  $D_{\bg}\bw(t)$, is
\be
\label{eq:derivcov1}
D_{\bg}\bw=(w'_k+\Gamma^k_{ij}w_iu'_j)\f_{,k}.
\ee
 A {\it parallel vector field $\bw$ along $\bg$} is a vector field having $\D_{\bg}\bw=\bo\ \forall t$. A regular curve $\bg$ is a {\it geodesic of} $\Sigma$ if the vector field $\bg'$ of the vectors tangent to $\bg$ is parallel along $\bg$. 
 \begin{teo}
 \label{teo: charactgeodes1}
 A curve $\bg$ is a geodesic of $\Sigma\iff\bnu\times\N=\bo$.
 \begin{proof}
 If $\bg$ is a geodesic, then the derivative of its tangent $\bg'$ has a component only along $\N$, i.e. $\bg''\times\N=\bo\Rightarrow\bg'\cdot\bg''=0$. The  principal normal to $\bg$, $\bnu$, is orthogonal to $\bg'\Rightarrow\bnu\times\N=\bo$. Vice versa, if $\bnu\times\N=\bo$, then $\bg''$ is orthogonal to $\bg'\Rightarrow\D_{\bg}\bg'=\bo$, i.e. $\bg$ is a geodesic.
  \end{proof}
 \end{teo}
\begin{teo}
If $\bg$ is a geodesic, then $|\bg'|=const.$
\begin{proof}
In a geodesic, $\bg'\cdot\bg''=0\Rightarrow\dfrac{d(\bg'\cdot\bg')}{dt}=0\Rightarrow|\bg'|=const.$
\end{proof}
\end{teo}
This result shows that, in a geodesic, the parameter is always the natural parameter $s$.

Let $\bg(s)$ be a curve on $\Sigma$ parameterized by the the arc-length $s$. We call {\it geodesic curvature of} $\bg(s)$ the function
\bes
\kappa_g:=D_{\bg}\btau\cdot(\N\times\btau),
\ees
 where $\btau=\bg'\in\S$ is the tangent vector to $\bg$. Because $\N\times\btau\in\S$ lies in $T_{\bg}\Sigma$, the component of $\btau'$ orthogonal to $T_{\bg}\Sigma$ gives a  null contribution to $\kappa_g$, so we can also write 
 \bes
 \kappa_g=\btau'\cdot(\N\times\btau).
 \ees
\begin{teo}
A regular curve $\bg(s)$ is a geodesic $\iff\kappa_g=0\ \forall s$.
\begin{proof}
If $\bg$ is a geodesic, clearly $\kappa_g=0$. Vice versa, if $\kappa_g=0$, then $\btau,\btau'$ and $\N$ are linearly dependent, i.e. coplanar. Because $\btau'\cdot\btau=\N\cdot\btau=0\Rightarrow\btau'\times\N=\bo\Rightarrow$ by Theorem \ref{teo: charactgeodes1}, $\bg$ is a geodesic.
\end{proof}
\end{teo}
 Let us now write Eq. (\ref{eq:derivcov1}) in the particular case of $\bw=\bg'$, i.e. $w_1=u',w_2=v'$:
 \bes
 D_{\bg}\bw=(u''_k+\Gamma^k_{ij}u'_iu'_j)\f_{,k};
 \ees
 therefore, the geodesics are the solutions to the system of differential equations
 \be
 \label{eq:systeqgeodesic}
 \left\{
 \besp
 &u''+\Gamma^1_{11}u'^2+2\Gamma^1_{12}u'v'+\Gamma^1_{22}v'^2=0,\\
 &v''+\Gamma^2_{11}u'^2+2\Gamma^2_{12}u'v'+\Gamma^2_{22}v'^2=0.
 \end{split}
 \right.
 \ee
 It can be shown that $\forall p\in\Sigma$ and $\forall \bw(p)\in T_p\Sigma$ the geodesic is unique. 
 
 Let $p$ be a point of a  regular surface $\f:\Omega\rightarrow\Sigma$ and $\ba(v):G\subset\R\rightarrow\Sigma$ a smooth regular curve on $\Sigma$, with $v$ being the natural parameter and such that $p=\ba(0)$. Consider the geodesic $\bg_v$ passing through $q=\ba(v)$ and such that $\bg'_v(0)=\N(\ba(v))\times\btau(v)$, with $\btau(v)$ the (unit) tangent vector to $\ba(v)$. Consider the map $\f(u,v):\Omega\rightarrow\Sigma$ defined by posing $\f(u,v)=\bg_v(u)$; this  is a surface whose  coordinates $(u,v)$ are called {\it semigeodesic coordinates}.
 
 Let us see the form that  the first fundamental form (i.e. the metric tensor $\g$), the Christoffel symbols and the Gaussian curvature take in semigeodesic coordinates. Curves $\f(u,v_0)=\bg_{v_0}(u)$ are geodesics, and $u$ is hence their natural parameter. Therefore, $\f_{,u}\in\S\Rightarrow g_{11}=1$. Then, $\f_{,uu}(u,v_0)$ is the derivative of the tangent vector to a geodesic $\f(u,v_0)=\bg_{v_0}(u)\Rightarrow\f_{,uu}(u,v_0)$ does not have a component along the tangent, hence, Eq. (\ref{eq:syst1})$_1\Rightarrow \Gamma^1_{11}=\Gamma^2_{11}=0.$ Then, by Eq. (\ref{eq:syst1})$_2$, we get $g_{12,u}=0\Rightarrow g_{12}$ does not depend upon $u\Rightarrow g_{12}(u,v)=g_{12}(0,v)\ \forall u$. Moreover, let $\theta$  be the angle between the curve $\ba$, i.e. between the coordinate line $\f(0,v)$, whose tangent vector is $\f_{,v}(0,v)$, and the geodesic $\bg_v(u)$, whose tangent vector at $(0,v)$ is $\bg'_v(u)$. Then, $\theta=\dfrac{\pi}{2}$ because $\bg'_v(0)=\N(\ba(v))\times\btau(v)$. As a consequence, $g_{12}(0,v)=0\Rightarrow g_{12}(u,v)=0\ \forall (u,v)\in\Omega$. Finally, setting $\g_{22}=g$,
 \bes
 \g=\left[\begin{array}{cc}1&0\\0&g\end{array}\right],
 \ees
 with $g>0$ because $\g$ is positive definite. Through systems (\ref{eq:syst1})$-$(\ref{eq:syst3}), we obtain 
\bes
\Gamma^1_{12}=0,\ \ \Gamma^2_{12}=\frac{g_{,u}}{2g},\ \ \Gamma^1_{22}=-\frac{g_{,u}}{2},\ \ \Gamma^2_{22}=\frac{g_{,v}}{2g},
\ees
 and using Eq. (\ref{eq:theoegreg}), we get
 \bes
 K=\det\X=-\frac{g_{,uu}}{2g}+\frac{g^2_{,u}}{4g^2}.
 \ees
 
 Given two points $p_1,p_2\in\Sigma$, we define the {\it distance $d(p_1,p_2)$} as the infimum of the lengths of the curves on $\Sigma$ joining the two points. We end with an important characterization of geodesics:
 \begin{teo}
 Geodesics are the curves of minimal distance between two points of a surface.
 \begin{proof}
 Let $\bg:G\subset\R\rightarrow\Sigma$ be  a geodesic on $\Sigma$, parameterized with the arc-length, and $\ba$ a smooth regular curve through $p$ and orthogonal to $\bg$. Through $\ba$, we set up a system of semigeodesic coordinates in a neighbourhood $U$ of $p$. With an opportune parameterization $\ba(t)$, in such coordinates we can get $p=\f(0,0)$ and  $\bg$  described by the equation $v=0$. Let $q\in U$ be a point in $\bg$, and consider a regular curve connecting $p
 $ with $q$. The length $\ell(p,q)$ of such a curve is
 \bes
 \ell(p,q)=\int_{p}^{q}\sqrt{u'^2+g\ v'^2}dt\geq\left|\int_{p}^{q}u'dt\right|=|u_q-u_p|.
 \ees
 Observing that $p=(u_p,0),q=(u_q,0)$, we remark that $|u_q-u_p|$ is exactly the length of $\bg$ between $p$ and $q$, because $\bg$ is parameterized with its arc-length.
 \end{proof}
 \end{teo}
 
 There is another, direct and beautiful way to show that geodesics are the shortest path lines: the use of the methods of the calculus of variations\footnote{The reader is addressed to texts on the calculus of variations for an insight into this  matter, cf. the suggested texts. Here, we just recall the fundamental fact to be used in the proof concerning geodesics: Let 
 \bes
 J(t)=\int_a^bF(\x,\x',t)dt
 \ees
 be a functional to be minimized by a proper choice of the function $\x(t)$ (in the case of the geodesics, $J=\ell(p,q)$); then, such a minimizing function can be found as a solution to the {\it Euler-Lagrange equations}
 \bes
 \frac{d}{dt}\frac{\partial F}{\partial \x'}-\frac{\partial F}{\partial \x}=\bo.
 \ees
 }. The length  $\ell(p,q)$ of a curve $ \bg(t)\in\Sigma$ between two points $p$ and $q$ is given by the functional (\ref{eq:lengthcurve}); it depends upon the first fundamental form, i.e. upon the metric tensor $\g$ on $\Sigma$. For the sake of conciseness, let $\bw=(w_{,t}^1,w_{,t}^2)$ be the tangent vector to the curve $\bg(w_1,w_2)\in\Sigma$. Then,
 \bes
 \ell(p,q)=\int_p^q\sqrt{I(\bw)}dt=\int_p^q\sqrt{\bw\cdot\g\bw}dt.
 \ees
The curve $\bg(t)$ that minimizes $\ell(p,q)$ is the solution to the {\it Euler-Lagrange equations}
 \bes
  \frac{d}{dt}\frac{\partial F}{\partial \bw_{,t}}-\frac{\partial F}{\partial \bw}=\bo \ \rightarrow\  \frac{d}{dt}\frac{\partial F}{\partial w^k_{,t}}-\frac{\partial F}{\partial w^k}=0, \ \ k=1,2,
 \ees
 where 
 \bes
 F(\bw,\bw_{,t},t)=\sqrt{\bw\cdot\g\bw}=\sqrt{g_{ij}w^i_{,t}w^j_{,t}}.
 \ees
 It is more direct, and equivalent, to minimize $J^2(t)$, i.e. to write the Euler-Lagrange equations for 
 \bes
 \Phi(\bw,\bw_{,t},t):=F^2(\bw,\bw_{,t},t)=g_{ij}w^i_{,t}w^j_{,t}.
 \ees 
 Therefore:
 \bes
 \besp
 &\frac{\partial\Phi}{\partial w^k_{,t}}=2g_{jk}w^j_{,t},\\
 &\frac{\partial\Phi}{\partial w^k}=\frac{\partial g_{hj}}{\partial w^k}w^h_{,t}w^j_{,t},\\
 &\frac{d}{dt}\frac{\partial\Phi}{\partial w^k_{,t}}=2\left(g_{jk}w^j_{,tt}+\frac{dg_{jk}}{dt}w^j_{,t}\right)=2\left(g_{jk}w^j_{,tt}+\frac{\partial g_{jk}}{\partial w^l}w^l_{,t}w^j_{,t}\right),
 \end{split}\ \ \ \ j,h,k=1,2.
 \ees
 The Euler-Lagrange equations are hence
 \bes
 g_{jk}w^j_{,tt}+\frac{\partial g_{jk}}{\partial w^h}w^h_{,t}w^j_{,t}-\frac{1}{2}\frac{\partial g_{hj}}{\partial w^k}w^h_{,t}w^j_{,t}=0,,\ \ \ j,h,k=1,2,
 \ees
 which can be rewritten as
\bes
 g_{jk}w^j_{,tt}+\frac{1}{2}\left(\frac{\partial g_{jk}}{\partial w^h}+\frac{\partial g_{hk}}{\partial w^j}-\frac{\partial g_{hj}}{\partial w^k}\right)w^h_{,t}w^j_{,t}=0,\ \ \ j,h,k=1,2.
 \ees 
 Multiplying by $g^{lk}$, we get
 \bes
  g^{lk}g_{jk}w^j_{,tt}+\frac{1}{2}g^{lk}\left(\frac{\partial g_{jk}}{\partial w^h}+\frac{\partial g_{hk}}{\partial w^j}-\frac{\partial g_{hj}}{\partial w^k}\right)w^h_{,t}w^j_{,t}=0,\ \ \ j,h,k,l=1,2.
 \ees
 and finally, because
 \bes
 g^{lk}g_{jk}=\delta_{lj}
 \ees
  and by Eq. (\ref{eq:christoffavecg}), we get
 \bes
 w^l_{,tt}+\Gamma^l_{jh}w^j_{,t}w^h_{,t}=0,\ \ \ j,h,l=1,2.
 \ees
 These are the differential equations whose solution is the curve of minimal length between two points of $\Sigma$; comparing these equations with those of a geodesic of $\Sigma$, Eq. (\ref{eq:systeqgeodesic}), we see that they are the same: The geodesics of a surface are hence the curves of minimal distance on the surface. 
 
 The Christoffel symbols of a plane are all null; as a consequence, the geodesic lines of a plane are straight lines. In fact, only such lines have a constant derivative.  
 
 Through systems (\ref{eq:syst1})$-$(\ref{eq:syst3}), we can calculate the Christoffel symbols for a revolution surface, Eq. (\ref{eq:revolsurf}), which are all null excepted
 \bes
 \Gamma^2_{12}=\frac{\phi'}{\phi},\ \ \ \Gamma^1_{22}=-\phi\ \phi',
 \ees
 so the system of differential equations (\ref{eq:systeqgeodesic}) becomes
 \be
 \label{eq:geodetichesurfrot}
 \left\{
 \begin{array}{l}
 u''-\phi\ \phi'v'^2=0,\\
 v''+2\dfrac{\phi'}{\phi}u'v'=0.
 \end{array}
 \right.
 \ee
 It is direct to check that the  meridians $(u=t,v=v_0)$ are geodesic lines, while the parallels $(u=u_0,v=t)$ are geodesics $\iff \phi'(u_0)=0$.
 
 \section{The Gauss-Codazzi compatibility conditions}
 Let us consider a surface $\Sigma$ whose points are determined by the vector function $\r:\Omega\subset\R^2\rightarrow\Sigma\subset\Eu,\r(\au,\ad)=x_i(\au,\ad)\bepsi_i$, with $\bepsi_i, i=1,2,3,$ the vectors of the orthonormal basis of the reference frame $\Rep=\{o;\bepsi_1,\bepsi_2,\bepsi_3\}$ and the parameters $\au,\ad$ chosen in such a way that the lines $\au=const.,\ad=const.$ are lines of  curvature, i.e. tangent at each point to the principal directions of curvature and hence mutually orthogonal\footnote{Here, the symbol $\r$ is  preferred to $\f$, like $\au$ to $u$ and $\ad$ to $v$, to recall that we have made the particular choice of coordinate lines that are lines of curvature. All the developments could be done in a more general case, but this choice is made to obtain simpler relations, which preserves anyway the generality.}. 
 With such a choice, cf. Eq. (\ref{eq:metricsigma}),
 \bes
 ds^2=A_1^2d\au^2+A_2^2d\ad^2,
 \ees
 with
 \bes
 \besp
& A_1=\sqrt{\r_{,\au}^2}=\sqrt{\frac{dx_i}{d\au}\frac{dx_i}{d\au}},\\
& A_2=\sqrt{\r_{,\ad}^2}=\sqrt{\frac{dx_i}{d\ad}\frac{dx_i}{d\ad}}
 \end{split}
 \ees
  the so-called {\it Lamé's parameters}. We remark that along the lines of curvatures, i.e. the lines $\alpha_i=const., i=1,2$, which in short, from now, on we  call  {\it the lines $\alpha_i$}, it is
  \bes
  \besp
  &ds_1=A_1d\au,\\
 & ds_2=A_2d\ad,
  \end{split}
  \ees
  and hence, 
  \be
  \label{eq:vecttanglc}
  \besp
  \bl_1=\frac{ds_1}{d\au}=A_1\eu,\\
  \bl_2=\frac{ds_2}{d\ad}=A_2\ed
  \end{split}
  \ee
  are the vectors tangent to the lines of curvature.
  Let 
  \be
  \label{eq:triedre}
  \eu=\frac{1}{A_1}\r_{,\au},\ \ \ed=\frac{1}{A_2}\r_{,\ad},\ \ \et=\eu\times\ed(=\N);
  \ee
  these three vectors form the  orthonormal (local) natural basis $\texttt{e}=\{\eu,\ed,\et\}$. We always make the choice of $\au,\ad$ such that $\e_3$ is always directed toward the convex side of $\Sigma$ if the point is elliptic or parabolic or toward the side of the centres of negative curvature, if the point is hyperbolic. 
  
  We consider a vector $\bv=\bv(p),\ p\in\Sigma$,
  \bes
  \bv=v_1\eu+v_2\ed+v_3\et,
  \ees
  and we want to calculate how it transforms  when $p$ changes. To this end, we need to calculate how $\eu,\ed,\et$ change with $\au,\ad$. Let $q\in\Sigma$ be a point in the neighborhood of $p$ on the line $\alpha_i$ and let us first consider the change of $\et$ in passing from $p$ to $q$. Because $p$ and $q$ belong to the same line $\alpha_i$, by the theorem of Rodrigues, we get (no summation on $i$ in the following equations)
  \bes
  \frac{\partial\et}{\partial\bl_i}=-\kappa_i\bl_i,\ i=1,2,
  \ees
i.e., by  Eq. (\ref{eq:vecttanglc}),
\bes
\frac{\partial\et}{\partial\ai}=\frac{A_i}{R_i}\e_i,
\ees
 with 
 \bes
 R_i=-\frac{1}{\kappa_i}
 \ees
 the (principal) radius of curvature along the line $\ai$. The minus sign  in the previous equation is due to the choice made above for orienting $\et=\N$, which gives always $\N=-\bnu$, with $\bnu$ the principal normal to the line $\ai$. This result can also be obtained  directly, see Fig. \ref{fig:52}:
      \begin{figure}[ht]
	\begin{center}
        \includegraphics[height=.2\textheight]{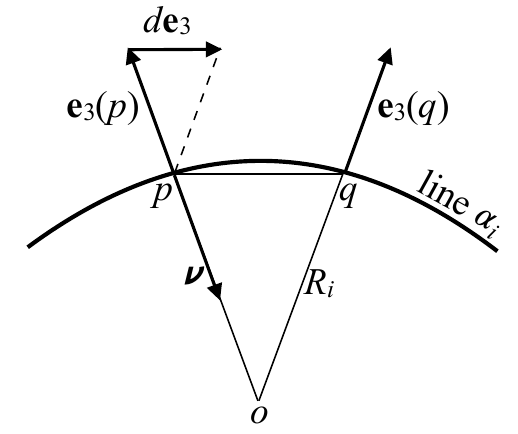}
	\caption{Variation of $\N=\et$ along a line of curvature.}
	\label{fig:52}
	\end{center}
\end{figure}
\bes
\et(q)=\et(p)+d\et
\ees
 and  in the limit of $q\rightarrow p$, $d\et$ tends to be parallel to $q-p$ and
  \bes
  \lim_{q\rightarrow p}(q-p)=\bl_i=A_i\e_i.
  \ees
 By the similitude of the triangles, it is evident that
 \bes
 \frac{|d\et|}{|\et|}=\frac{|q-p|}{R_i};
 \ees
moreover,
\bes
d\et=\frac{\partial\et}{\partial\ai}d\ai\e_i.
\ees
Finally, as $|\et|=1$, we get again
\be
\label{eq:codazzi0}
\frac{\partial\et}{\partial\ai}=\frac{A_i}{R_i}\e_i.
\ee
  Implicitly, in this last proof, we have used the theorem of Rodrigues because we have assumed that $d\et$ is parallel to $\bl_i$, as it is, because line $\ai$ is a line of curvature.  
  
  We  now move on to determine the changes in $\eu$ and $\ed$; for this purpose, we remark that
  \bes
  \frac{\partial\r_{,\au}}{\partial\ad}=\frac{\partial^2\r}{\partial\ad\partial\au}=\frac{\partial^2\r}{\partial\au\partial\ad}= \frac{\partial\r_{,\ad}}{\partial\au},
  \ees
  so by Eq. (\ref{eq:triedre}) we get
  \be
  \label{eq:codazzi1}
  \frac{\partial (A_1\eu)}{\partial\ad}= \frac{\partial (A_2\ed)}{\partial\au}.
  \ee
Let us study now $\dfrac{\partial\e_j}{\partial\ai}$; as $|\e_j|=1, j=1,2$,
  \be
  \label{eq:codazzi2}
  \frac{\partial\e_j}{\partial\ai}\cdot\e_j=0\ \ \forall i,j=1,2.
  \ee
  Because $\eu\cdot\ed=0$,
  \bes
   \frac{\partial\eu}{\partial\au}\cdot\ed=\frac{\partial(\eu\cdot\ed)}{\partial\au}-\eu\cdot\frac{\partial\ed}{\partial\au}=-\eu\cdot\frac{\partial\ed}{\partial\au}.
  \ees
 By Eq. (\ref{eq:codazzi1}), we get 
 \bes
 \frac{\partial\ed}{\partial\au}=\frac{1}{A_2}\frac{\partial( A_1\eu)}{\partial\ad}-\frac{1}{A_2}\frac{\partial A_2}{\partial\au}\ed,
 \ees 
  which when inserted into the previous equation gives, by Eq. (\ref{eq:codazzi2}),
  \bes
   \frac{\partial\eu}{\partial\au}\cdot\ed=-\frac{1}{A_2}\frac{\partial (A_1\eu)}{\partial\ad}\cdot\eu+\frac{1}{A_2}\frac{\partial A_2}{\partial\au}\ed\cdot\eu=-\frac{A_1}{A_2}\frac{\partial\eu}{\partial\ad}\cdot\eu-\frac{1}{A_2}\frac{\partial A_1}{\partial\ad}\eu\cdot\eu=-\frac{1}{A_2}\frac{\partial A_1}{\partial\ad}.
  \ees
  Then, because $\eu\cdot\et=0$,
  \bes
  \frac{\partial\eu}{\partial\au}\cdot\et=
  \frac{\partial(\eu\cdot\et)}{\partial\au}-\eu\cdot\frac{\partial\et}{\partial\au}=-\eu\cdot\frac{\partial\et}{\partial\au},
  \ees
  and by Eq. (\ref{eq:codazzi0}),
  \bes
  \frac{\partial\et}{\partial\au}=\frac{A_1}{R_1}\eu,
  \ees
  so finally,
  \bes
   \frac{\partial\eu}{\partial\au}\cdot\et=-\frac{A_1}{R_1}.
  \ees  
 Again, through Eqs. (\ref{eq:codazzi1}) and (\ref{eq:codazzi2}), we get 
 \bes
 \frac{\partial\eu}{\partial\ad}\cdot\ed=\frac{1}{A_1}\frac{\partial (A_2\ed)}{\partial\au}\cdot\ed-\frac{1}{A_1}\frac{\partial A_1}{\partial\ad}\eu\cdot\ed=\frac{A_2}{A_1}\frac{\partial \ed}{\partial\au}\cdot\ed+\frac{1}{A_1}\frac{\partial A_2}{\partial\au}\ed\cdot\ed=\frac{1}{A_1}\frac{\partial A_2}{\partial\au}
  \ees
  and also, by Eq. (\ref{eq:codazzi0}),
  \bes
  \frac{\partial\eu}{\partial\ad}\cdot\et=
  \frac{\partial(\eu\cdot\et)}{\partial\ad}-\eu\cdot\frac{\partial\et}{\partial\ad}=-\eu\cdot\frac{\partial\et}{\partial\ad}=-\frac{A_2}{R_2}\eu\cdot\ed=0.
  \ees
The derivatives of $\ed$ can be found in a similar way, and  resuming, we have
\be
\label{eq:codazzi3}
\begin{array}{ccl}
\dfrac{\partial\eu}{\partial\au}&=&-\dfrac{1}{A_2}\dfrac{\partial A_1}{\partial\ad}\ed-\dfrac{A_1}{R_1}\et,\medskip\\
\dfrac{\partial\eu}{\partial\ad}&=&\dfrac{1}{A_1}\dfrac{\partial A_2}{\partial\au}\ed,\medskip\\
\dfrac{\partial\ed}{\partial\au}&=&\dfrac{1}{A_2}\dfrac{\partial A_1}{\partial\ad}\eu,\medskip\\
\dfrac{\partial\ed}{\partial\ad}&=&-\dfrac{1}{A_1}\dfrac{\partial A_2}{\partial\au}\eu-\dfrac{A_2}{R_2}\et,\medskip\\
\dfrac{\partial\et}{\partial\au}&=&\dfrac{A_1}{R_1}\eu,\medskip\\
\dfrac{\partial\et}{\partial\ad}&=&\dfrac{A_2}{R_2}\ed.
\end{array}
\ee

  Passing now to the second-order derivatives, imposing the equality of mixed derivatives, gives some important differential relations between the Lamé's parameters $A_i$ and the radii of curvature $R_i$. In fact, from the identity
  \bes
  \frac{\partial^2\et}{\partial\au\partial\ad}=\frac{\partial^2\et}{\partial\ad\partial\au},
  \ees
  and Eqs. (\ref{eq:codazzi3})$_{5,6}$, we get
  \bes
  \frac{\partial}{\partial\ad}\left(\dfrac{A_1}{R_1}\eu\right)=\frac{\partial}{\partial\au}\left(\dfrac{A_2}{R_2}\ed\right),
  \ees
  whence
  \bes
  \frac{\partial}{\partial\ad}\left(\dfrac{A_1}{R_1}\right)\eu+\dfrac{A_1}{R_1}\frac{\partial\eu}{\partial\ad}=
  \frac{\partial}{\partial\au}\left(\dfrac{A_2}{R_2}\right)\ed+\dfrac{A_2}{R_2}\frac{\partial\ed}{\partial\au}.
  \ees
  Inserting now Eqs. (\ref{eq:codazzi3})$_{2,3}$ into the last result and rearranging the terms gives
  \bes
  \left[ \frac{\partial}{\partial\ad}\left(\dfrac{A_1}{R_1}\right)-\frac{1}{R_2}\frac{\partial A_1}{\partial\ad}\right]\eu-
  \left[ \frac{\partial}{\partial\au}\left(\dfrac{A_2}{R_2}\right)-\frac{1}{R_1}\frac{\partial A_2}{\partial\au}\right]\ed=0,
  \ees
that to be true needs that the two following conditions be identically satisfied:
\be
\label{eq:codazzi4}
\besp
&\frac{\partial}{\partial\ad}\left(\dfrac{A_1}{R_1}\right)-\frac{1}{R_2}\frac{\partial A_1}{\partial\ad}=0,\\
& \frac{\partial}{\partial\au}\left(\dfrac{A_2}{R_2}\right)-\frac{1}{R_1}\frac{\partial A_2}{\partial\au}=0.
\end{split}
\ee
The above equations are known as the {\it Codazzi conditions}. Let us now consider the other identity
\bes
\frac{\partial^2\eu}{\partial\au\partial\ad}=\frac{\partial^2\eu}{\partial\ad\partial\au};
\ees
 again using Eq. (\ref{eq:codazzi3}), with some standard operations, this identity can be transformed to
 \bes
 \left[\frac{\partial}{\partial\au}\left(\frac{1}{A_1}\frac{\partial A_2}{\partial\au}\right)+\frac{\partial}{\partial\ad}\left(\frac{1}{A_2}\frac{\partial A_1}{\partial\ad}\right)+\frac{A_1}{R_1}\frac{A_2}{R_2}\right]\ed+\left[\frac{\partial}{\partial\ad}\left(\frac{A_1}{R_1}\right)-\frac{1}{R_2}\frac{\partial A_1}{\partial\ad}\right]\et=0.
 \ees
  Also in this case, for this equation to be identically satisfied, each of the expressions in square brackets must vanish, which gives two more differential conditions, but of which only the first one is new, as the second one corresponds to Eq. (\ref{eq:codazzi4})$_1$. The new condition is hence
  \be
  \label{eq:codazzi5}
  \frac{\partial}{\partial\au}\left(\frac{1}{A_1}\frac{\partial A_2}{\partial\au}\right)+\frac{\partial}{\partial\ad}\left(\frac{1}{A_2}\frac{\partial A_1}{\partial\ad}\right)+\frac{A_1}{R_1}\frac{A_2}{R_2}=0,
  \ee
  which is known as the {\it Gauss condition}. The last identity
 \bes
\frac{\partial^2\ed}{\partial\au\partial\ad}=\frac{\partial^2\ed}{\partial\ad\partial\au}
\ees 
  does not add any  independent condition, which can be easily checked. The meaning of the {\it Gauss-Codazzi conditions}, Eqs. (\ref{eq:codazzi4}) and (\ref{eq:codazzi5}), is that of {\it compatibility conditions}: Only when these conditions are satisfied by functions $A_1,A_2,R_1$ and $R_2$, then such functions represent the Lamé's parameters and the principal radii of curvature of a surface, i.e. only in this case they define a surface, except for its position in space. The Gauss-Codazzi conditions are important in establishing the equations of the classical theory of shells.

 \section{Exercises}
 \begin{enumerate}
 \item Prove that a function of the type $x_3=f(x_1,x_2)$, with $f:\Omega\subset\R^2\rightarrow\R$ smooth, defines a surface.
 \item \label{ex:2ch7} Show that the catenoid is the rotation surface of a catenary, then find its Gaussian curvature.
 \item Show that the pseudo-sphere is the rotation surface of a tractrice and explain why the surface has this name (hint: look for its Gaussian curvature).
 \item Prove that the regularity of a cone $\f(u,v)=v\bg(u)$ is satisfied at each point except  at the apex and at the points on  the straight lines tangent to $
\bg(u)$. 
\begin{figure}[h!]
\begin{center}
\includegraphics[height=.25\textwidth]{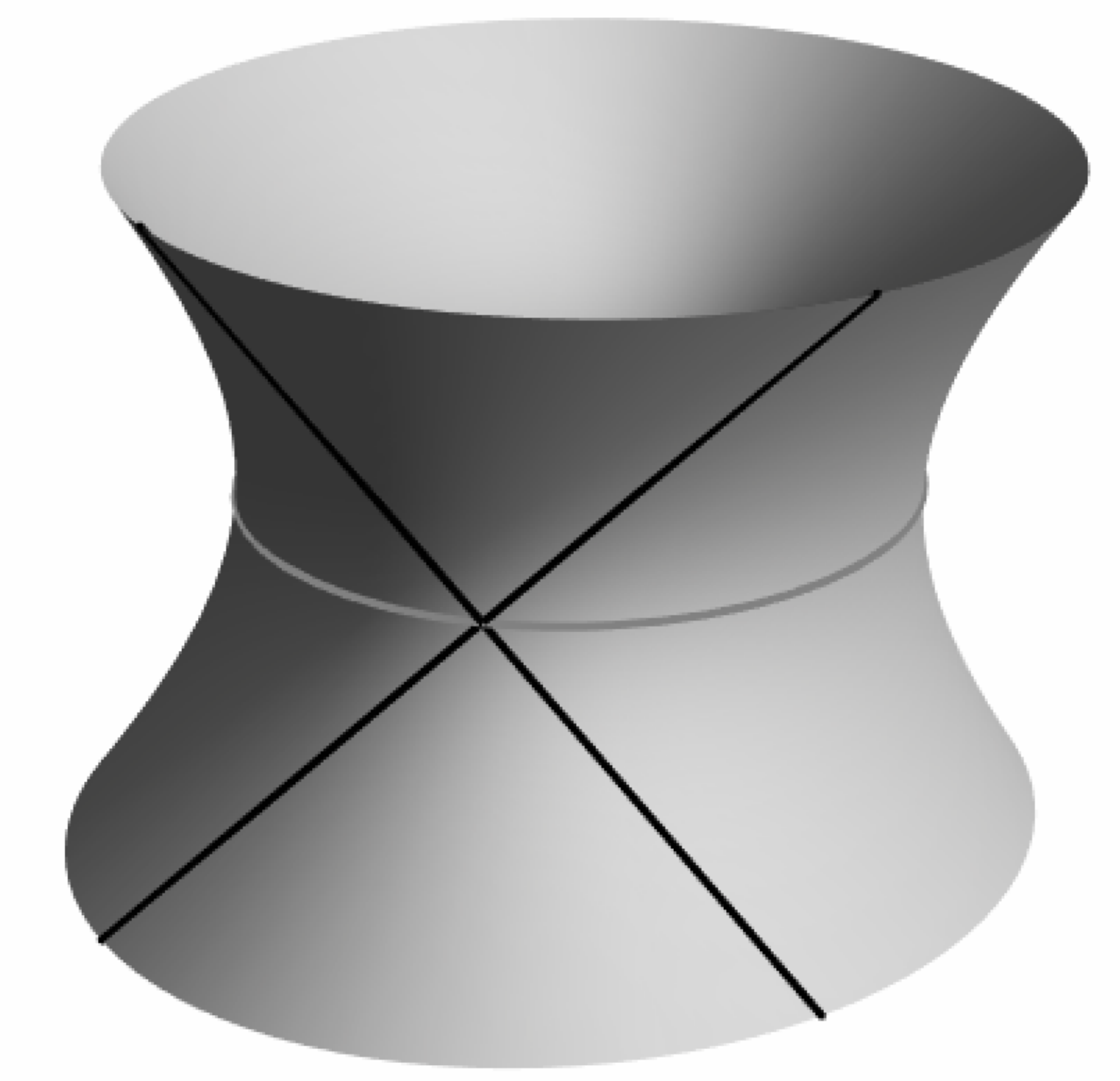}\ \ \ \ \ \
\includegraphics[height=.3\textwidth]{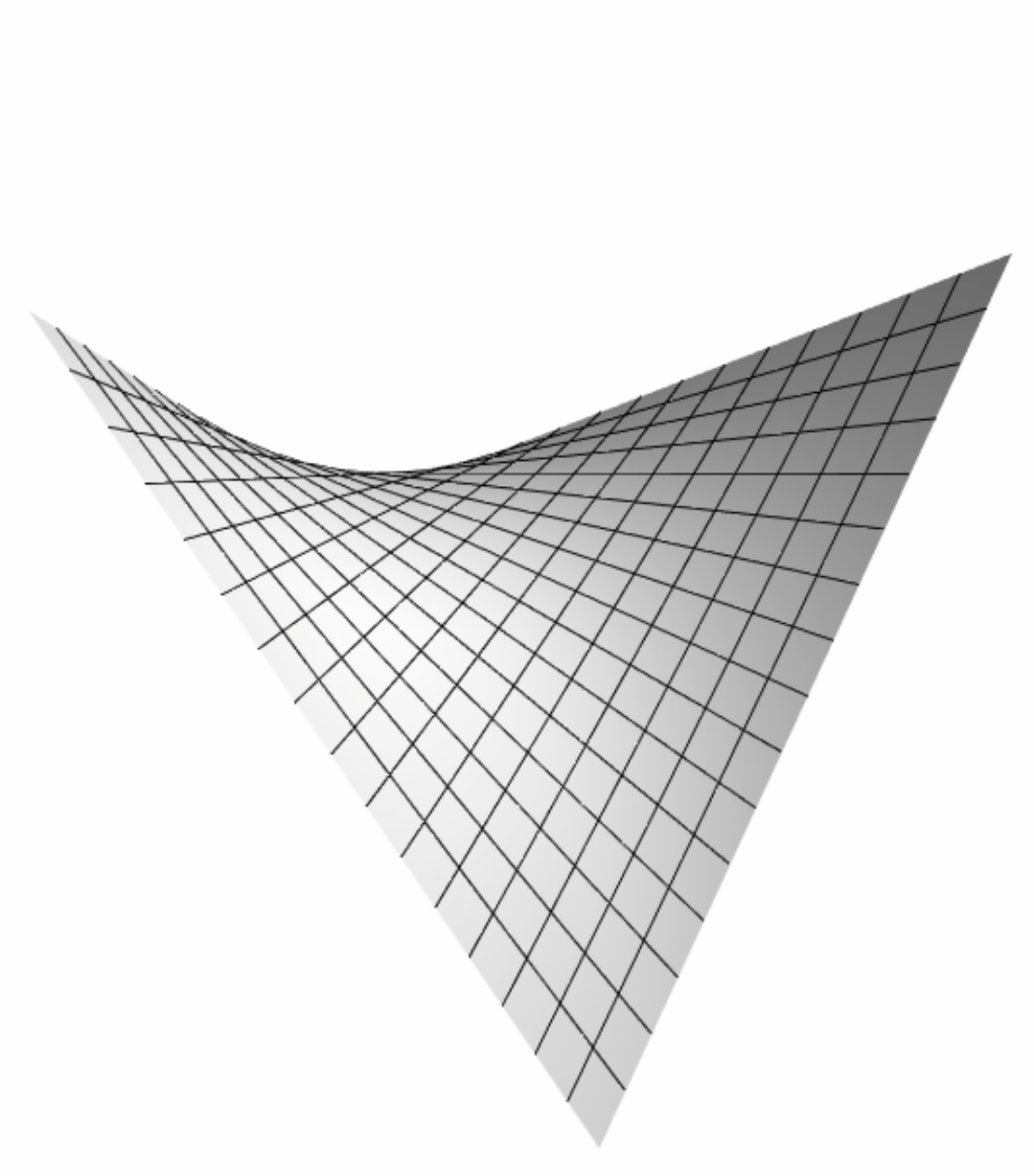}
\caption{A hyperbolic hyperboloid, left, and a hyperbolic paraboloid, right.}
\label{fig:64}
\end{center}
\end{figure}
\item Prove that the hyperbolic hyperboloid is a doubly ruled surface and determine the angle $\theta$ formed by two straight lines belonging to the two sets of lines on the surface, see the left panel of Fig. \ref{fig:64}.
\item Prove that the hyperbolic paraboloid whose Cartesian equation is $x_3=x_1x_2$, right panel of Fig. \ref{fig:64}, is a doubly ruled surface and determine the angle $\theta$ formed by two straight lines belonging to the two sets of lines on the surface. Where does $\theta=\dfrac{\pi}{2}$?
\item Consider the parameterization
\bes
 \f(u,v)=(1-v)\bg(u)+v\bl(u),
\ees
 with
 \bes
\bg(u)=(\cos(u-\alpha),\sin(u-\alpha),-1),\ \ \bl(u)=(\cos(u+\alpha),\sin(u+\alpha),1).
 \ees
 Show that:
\begin{itemize}
\item for $\alpha=0$, one gets a cylinder with equation $x_1^2+x_2^2=1$;
\item for $\alpha=\dfrac{\pi}{2}$, one gets a cone with equation $x_1^2+x_2^2=x_3^2$;
\item for $0<\alpha<\dfrac{\pi}{2}$, one gets a hyperbolic hyperboloid with equation 
\bes\dfrac{x_1^2+x_2^2}{\cos^2\alpha}-\dfrac{x_3^2}{\cot^2\alpha}=1.\ees
\end{itemize}
\item Calculate the metric tensor of a sphere of radius $R$, write its first fundamental form, determine the area of a sector of surface between the longitudes $\theta_1$ and $\theta_2$ and the length of the parallel at the latitude $\pi/4$ between these two longitudes.
\item  Prove that the surface defined by
\bes
\f(u,v):\Omega=\R\times(-\pi,\pi]\rightarrow\Eu|\ \ \f(u,v)=\left(\frac{\cos v}{\cosh u},\frac{\sin v}{\cosh u},\frac{\sinh u}{\cosh u}\right)
\ees
is a sphere. Then, show that the image of any straight line on $\Omega$ is a loxodromic line on the sphere. 
\item Calculate the vectors of the natural basis,  the tensors $\g,\B,\X$ and the first and second fundamental forms  for the catenoid.
\item \label{ex:11ch7} Calculate the same for the helicoid of parametric equation
\bes
 \f(u,v)=\bg(u)+v\bl(u),
 \ees
 with
 \bes
  \bg(u)=(0,0,u),\ \ \bl(u)=(\cos u,\sin u,0).
\ees
\item Show that the catenoid and the helicoid are made of hyperbolic points.
\item Determine the geodesic lines of a circular cylinder.
 \end{enumerate}

\cleardoublepage
\phantomsection
\addcontentsline{toc}{chapter}{Suggested texts}
\chapter*{Suggested texts}
There are many textbooks on tensors, differential geometry and calculus of variations. The style, content, language of such books greatly depend upon the scientific community the authors belong to: pure or applied mathematicians, physicists, theoretical mechanicians or engineers. It is hence difficult to suggest some readings in the domain, and ultimately, it is mostly  a matter of personal taste.  

This textbook is greatly inspired by some classical methods, style and language that are typical in the community of theoretical mechanics; the following few suggested readings, among several possible others, belong to such a kind of scientific literature. They are classical textbooks and though the list is far from being exhaustive, they constitute a solid basis for the topics briefly developed in this manuscript, where the objective is to present mathematics for mechanics.

A good introduction to tensor algebra and analysis, which greatly inspired the content of this manuscript, are the two introductory chapters of the classical textbook
\begin{itemize}
\item M. E. Gurtin: \textit{An introduction to continuum mechanics}. Academic Press, 1981,
\end{itemize}
or also, in a similar style, the long article
\begin{itemize}
\item P. Podio-Guidugli: \textit{A primer in elasticity}. Journal of Elasticity, v. 58: 1-104, 2000.
\end{itemize}
A short, effective introduction to tensor algebra and differential geometry of curves can be found in the following  text of exercises on analytical mechanics:
\begin{itemize}
\item P. Biscari, C. Poggi, E. G. Virga: \textit{Mechanics notebook}. Liguori Editore, 1999.
\end{itemize}
A classical textbook on linear algebra that is  recommended is
\begin{itemize}
\item P. R. Halmos: {\it Finite-dimensional vector spaces}. Van Nostrand Reynold, 1958.
\end{itemize}
In the previous textbooks, tensor algebra in curvilinear coordinates is not developed; an introduction to this topic, especially intended for physicists and engineers, can be found in 
\begin{itemize}
\item W. H. Müller: \textit{An expedition to continuum theory}. Springer, 2014,
\end{itemize}
which has largely influenced Chapter \ref{ch:6}.

Two modern and application-oriented textbooks on differential geometry of curves and surfaces are 
\begin{itemize}
\item V. A. Toponogov: \textit{Differential geometry of curves and surfaces - A concise guide}. Birkhäuser, 2006,
\item A. Pressley: \textit{Elementary differential geometry}. Springer, 2010.
\end{itemize}
A short introduction to the differential geometry of surfaces, oriented toward the mechanics of shells, can be found in the classical book
\begin{itemize}
\item V. V. Novozhilov: \textit{Thin shell theory}. Noordhoff LTD., 1964.
\end{itemize}

For what concerns the calculus of variations, a still valid textbook in the matter (but not only) is
\begin{itemize}
\item R. Courant, D. Hilbert: \textit{Methods of mathematical physics}. Interscience Publishers, 1953.
\end{itemize}
Two very good and classical textbooks with an introduction to the calculus of variations for engineers are
\begin{itemize}
\item C. Lanczos: \textit{The variational principles of mechanics}. University of Toronto Press, 1949,
\item H. L. Langhaar: \textit{Energy methods in applied mechanics}. Wiley, 1962.

\end{itemize}

\cleardoublepage
\phantomsection
\addcontentsline{toc}{chapter}{Solutions of the exercises}

\chapter*{Solutions to the exercises}

\section*{Chapter 1}
\begin{enumerate}
\item Suppose $\bo_1\neq\bo_2$; then, apply to a point $p$ the definition of vector null for both of them.

\item Use $\bv+\bo=\bv$ and make the scalar product with a vector $\bw$.

\item Make the norm of $\bv+\bo=\bv$ and use the above result.

\item $|\bu-\bv|=|\bu+\bv|\iff|\bu-\bv|^2=|\bu+\bv|^2\Rightarrow(\bu-\bv)\cdot(\bu-\bv)=(\bu+\bv)\cdot(\bu+\bv)\Rightarrow\bu\cdot\bv=0.$

\item By linearity, $\forall\bv=v_i\e_i,\psi(\bv)=v_i\psi(\e_i)$; moreover, $\bu\cdot\bv=u_iv_i$, so setting\\ $\bu=\psi(\e_i)\e_i, \psi(\bv)=\bu\cdot\bv$. Uniqueness: Suppose $\exists\bu_1\neq\bu_2|\psi(\bv)=\bu_1\cdot\bv=\bu_2\cdot\bv\Rightarrow(\bu_1-\bu_2)\cdot\bv=0\ \forall\bv\iff\bu_1-\bu_2=\bo\Rightarrow\bu_1=\bu_2$.

\item Let $\theta_u,\theta_v$ be the angles formed by $\bw$ with $\bu$ and $\bv$, respectively; then, $\bw\cdot\bu=\bw\cdot\bv\Rightarrow uw\cos\theta_u=vw\cos\theta_v\Rightarrow\cos\theta_u=\cos\theta_v$, as $\bu,\bv\in\S$;\\
$\theta_u=\theta_v\Rightarrow\cos\theta_u=\cos\theta_v\Rightarrow\ uw\cos\theta_u=vw\cos\theta_v\Rightarrow\bu\cdot\bw=\bv\cdot\bw$.

\item i) Coplanar vectors: Let $p$ be a point of the plane of the vectors $\Rightarrow\M^r_p\cdot\Ro=0$ because $\Ro\in$ to the plane, while $\M^r_p=(p_i-p)\times\bv^{p_i}$ is of course orthogonal to it. Then, let $q$ be a point $\notin$ to the plane of the vectors $\rightarrow\M^r_q=\M^r_p+(p-q)\times\Ro\Rightarrow\M^r_q\cdot\Ro=\M^r_p\cdot\Ro+(p-q)\times\Ro\cdot\Ro=\Ro\times\Ro\cdot(p-q)=0.$

ii) Parallel vectors: Let $\e\in\S|\ \bv^{p_i}=\alpha_i\e\ \forall i=1,...,n\Rightarrow\Ro=\sum_{i=1}^n\alpha_i\e\Rightarrow\\\forall o\in\Eu,\ \M^r_o=\sum_{i=1}^n(p_i-o)\times\alpha_i\e\Rightarrow\M^r_o\cdot\Ro=(\sum_{i=1}^n\alpha_i(p_i-o))\times\e\cdot(\sum_{i=1}^n\alpha_i)\e=0.$

\item It is a direct consequence of the theorem of reduction of the systems of applied vectors for the case $\Ro=\bo$, with $o$ being any point.

\item If $\Ro$ is applied to $p$, then the system can be reduced to $\Ro^p$ plus $\\\M^r_p=\sum_{i=1}^n(p-p)\times\bv^{p}_i=\bo.$

\item i) By the theorem of reduction, if $\Ro$ is applied to $o$, the system is reduced to $\Ro^o$ and $\M^r_o\Rightarrow$ to only $\Ro$ if $\M^r_o=\bo.$

ii) Because for coplanar or parallel vectors $\forall o\in\Ax,\ \M^r_o=\bo$, then the system is equivalent to $\Ro$ applied to any point of $\Ax$.

\item If a system is equilibrated, then, by definition, any equivalent system is equilibrated. Conversely, if it exists another system equilibrated and equivalent, then by the relation of equivalence also the system in object is equilibrated,  and this is true for any other equivalent system, which hence must be equilibrated.

\item Let $\bv^p,\bv^q=-(\bv^p)^q$ be two opposite vectors applied to $p$ and $q$, respectively $\Rightarrow\Ro=\bv^p+\bv^q=\bo$; moreover, $\M^r_p=(q-p)\otimes\bv^q\Rightarrow\forall o\neq p,\ \M^r_o=\M^r_p+(p-o)\times\Ro=\M^r_p\Rightarrow\M^r_o=\bo\ \forall o\in\Eu\iff q-p=\bo.$

\item If all the vectors pass through a point $p$, then the system is equivalent to $\Ro^p$ (Exercise \ref{ex:9ch1}) and if $\Ro=\bo$, then the system is equilibrated. 

\end{enumerate}

\section*{Chapter 2}
\begin{enumerate}
\item $\forall\bu,\L,\L\bu=\L(\bu+\bo)=\L\bu+\L\bo\iff\L\bo=\bo$

\item $\bu\in\S\rightarrow(\bu\otimes\bu)\bv=v\cos\theta\bu$; $(\I-\bu\otimes\bu)\bv=\bv-v\cos\theta\bu$, which is orthogonal to $\bu: (\I-\bu\otimes\bu)\bv\cdot\bu=\bv\cdot\bu-\bv\cdot\bu\ \bu\cdot\bu=\bo$, as $\bu\in\S$.

\item i) $\forall\a,\b\in\Ve, \a\cdot(\alpha\A)\b=(\alpha\A)^\top\a\cdot\b$ and by the linearity of the scalar product, $\a\cdot(\alpha\A)\b=\alpha\a\cdot\A\b=\alpha\A^\top\a\cdot\b\Rightarrow(\alpha\A)^\top=\alpha\A^\top$.

ii) $\forall\a,\b\in\Ve, \a\cdot(\A+\B)\b=(\A+\B)^\top\a\cdot\b$ and by the linearity of the scalar product and of tensors, $\a\cdot(\A+\B)\b=\a\cdot\A\b+\a\cdot\B\b=\A^\top\a\cdot\b+\B^\top\a\cdot\b=(\A^\top+\B^\top)\a\cdot\b\Rightarrow(\A+\B)^\top=\A^\top+\B^\top$.

iii) $\forall\bu,\bv\in\Ve\ \ \bu\cdot(\a\otimes\b)\A\bv=\bu\cdot(\a\otimes\b)(\A\bv)=(\a\otimes\b)^\top\bu\cdot(\A\bv)=(\b\otimes\a)\bu\cdot(\A\bv)=\A^\top(\b\otimes\a)\bu\cdot\bv=\A^\top(\a\cdot\bu\ \b)\cdot\bv=\a\cdot\bu\ \A^\top\b\cdot\bv=((\A^\top\b)\otimes\a)\bu\cdot\bv=\bu\cdot((\A^\top\b)\otimes\a)^\top\bv=\bu\cdot(\a\otimes\A^\top\b)\bv.$

\item By linearity and the definition of $\O:\ \forall\bu\in\Ve,\ (\L+\O)\bu=\L\bu+\O\bu=\L\bu\iff\L+\O=\L$.

\item i) $\tr\I=\tr(\delta_{ij}\ei\otimes\ej)=\delta_{ij}\tr(\ei\otimes\ej)=\delta_{ij}\ei\cdot\ej=\delta_{ij}\delta_{ij}=\delta_{ii}=3.$

ii) $\tr\L=\tr(\L+\O)=\tr\L+\tr\O\iff\tr\O=0.$

\item $\tr(\A\B)=\tr((A_{ij}\ei\otimes\ej)(B_{hk}\eh\otimes\ek))=A_{ij}B_{hk}\tr((\ei\otimes\ej)(\eh\otimes\ek))\\=A_{ij}B_{hk}\ej\cdot\eh\tr(\ei\otimes\ek)=A_{ij}B_{hk}\ej\cdot\eh\ \ei\cdot\ek=A_{ij}B_{hk}\delta_{jh}\delta_{ik}=A_{ij}B_{ji}$; in a similar way, we prove that $\tr(\B\A)=B_{ij}A_{ji}$; because $i,j$ are dummy indexes, $A_{ij}B_{ji}=B_{ij}A_{ji}\Rightarrow\tr(\A\B)=\tr(\B\A)$.

\item i) $\L^\top\cdot\M^\top=\tr((\L^\top)^\top\M^\top)=\tr(\L\M^\top)=\tr(\M^\top\L)=\M\cdot\L=\L\cdot\M$.

$\begin{array}{llll}\hspace{-2mm}\mathrm{ii)}&\hspace{-2mm}\L\M\cdot\N=\tr((\L\M)^\top\N)&=&\tr(\M^\top\L^\top\N)=\M\cdot\L^\top\N;\\
&&|\\
&&=&\tr(\N(\L\M)^\top)=\tr((\N\M^\top)\L^\top)=\tr(\L^\top(\N\M^\top))\\
&&=&\L\cdot\N\M^\top.
\end{array}$

\item i) $(\a\otimes\b)(\c\otimes\d)=((\a\otimes\b)(\c\otimes\d))_{ij}\ei\otimes\ej=(\a\otimes\b)_{ik}(\c\otimes\d)_{kj}\ei\otimes\ej\\=a_ib_kc_kd_j\ei\otimes\ej=\b\cdot\c\ a_id_j\ei\otimes\ej=\b\cdot\c\ \a\otimes\d.$

ii) $\A(\a\otimes\b)=\A(\a\otimes\b)_{ij}\ei\otimes\ej=A_{ik}(\a\otimes\b)_{kj}\ei\otimes\ej=A_{ik}a_kb_j\ei\otimes\ej\\=(\A\a)_ib_j\ei\otimes\ej=(\A\a)\otimes\b.$

\item $\L\cdot\bv\otimes\bw=\tr(\L^\top(\bv\otimes\bw))=\tr((\L^\top\bv)\otimes\bw)=\L^\top\bv\cdot\bw=\bv\cdot\L\bw.$

\item  $\A=\A^\top,\B=-\B^\top\Rightarrow \A\cdot\B=\A^\top\cdot\B^\top=\A\cdot(-\B)=-\A\cdot\B\iff\A\cdot\B=0$.

\item i) $\A=\A^\top\Rightarrow\A\cdot\L=\A\cdot(\L^s+\L^a)=\A\cdot\L^s+\A\cdot\L^a=\A\cdot\L^s$.

ii) $\B=-\B^\top\Rightarrow\B\cdot\L=\B\cdot(\L^s+\L^a)=\B\cdot\L^s+\B\cdot\L^a=\B\cdot\L^a$.

\item i) $\A\cdot(\B\C\D)=\tr(\A^\top\B\C\D)=\tr((\B^\top\A)^\top\C\D)=(\B^\top\A)\cdot(\C\D).$

ii) $\A\cdot(\B\C\D)=(\B\C\D)\cdot\A=\tr((\B\C\D)^\top\A)=\tr(\A(\D^\top\C^\top\B^\top))\\=\tr((\A\D^\top)(\C^\top\B^\top))=\tr((\C^\top\B^\top)(\A\D^\top))=\tr((\B\C)^\top(\A\D^\top))=\B\C\cdot\A\D^\top\\=\A\D^\top\cdot\B\C.$

\item $\L\in Sym(\Ve)\Rightarrow\L\cdot\W=0$, as already proved. Now, if $\L\cdot\W=0\ \forall\W\in Skw(\Ve)$, suppose $\L\notin Sym(\Ve)\Rightarrow\L=\L^s+\L^a\Rightarrow\L\cdot\W=\L^s\cdot\W+\L^a\cdot\W=\L^a\cdot\W=0$; if in $Skw(\Ve)$, we chose $\W=\L^a$, we get $\W\cdot\L^a=\L^a\cdot\L^a=0\iff\L^a=\O\Rightarrow\L\in Sym(\Ve)$.

\item $I_2=\dfrac{1}{2}(\tr^2\L-\tr\L^2)=\dfrac{1}{2}(L_{ii}L_{jj}-L_{ij}L_{ji})\\=L_{11}L_{22}+L_{11}L_{33}+L_{22}L_{33}-L_{12}L_{21}-L_{13}L_{31}-L_{23}L_{32}.$

\item $\a\times\b\cdot\c=\left[\begin{array}{ccc}0 & -a_3 & a_2 \\a_3 & 0 &-a_1 \\-a_2& a_1 & 0\end{array}\right]\left\{\begin{array}{c}b_1 \\b_2 \\b_3\end{array}\right\}\cdot\left\{\begin{array}{c}c_1 \\c_2 \\c_3\end{array}\right\}=a_2b_3c_1-a_3b_2c_1+a_3b_1c_2-a_1b_3c_2+a_1b_2c_3-a_2b_1c_3=\det\left[\begin{array}{ccc}0 & -a_3 & a_2 \\a_3 & 0 &-a_1 \\-a_2& a_1 & 0\end{array}\right].$

\item Let $\L_1\neq\L_2$ be two distinct inverse tensors of $\L$; then $\L_1\L=\I=\L_2\L\Rightarrow\\ \L_1\L-\L_2\L=\O\Rightarrow(\L_1-\L_2)\L=\O\ \forall\L\iff\L_1-\L_2=\O\Rightarrow\L_1=\L_2.$

\item $(\a\otimes\b)_{ij}=a_ib_j$; it is then sufficient to write the matrix representing $(\a\otimes\b)$ and to compute its determinant.

\item $(\alpha\L)^{-1}(\alpha\L)=\alpha(\alpha\L)^{-1}\L=\I\Rightarrow(\alpha\L)^{-1}\L=\dfrac{1}{\alpha}\I\Rightarrow(\alpha\L)^{-1}\L\L^{-1}=\dfrac{1}{\alpha}\I\L^{-1}\Rightarrow(\alpha\L)^{-1}=\dfrac{1}{\alpha}\L^{-1}.$

\item $|\W^2|=\W\cdot\W=\tr(\W^\top\W)=-\tr(\W\W)=\tr(\I-\bw\otimes\bw)=3-1=2\Rightarrow\\\W\W=-\dfrac{1}{2}|\W^2|(\I-\bw\otimes\bw).$

\item Let $\bw_1=(a_1,b_1,c_1),\bw_2=(a_2,b_2c_2)$; then, form $\W_1,\W_2$ and compute the two scalar products.

\item $\bu\times\bv=\bo\iff\bv=k\bu,k\in\R$. So, $\bu\times\bv=\bo\Rightarrow\bu\otimes\bv=k\bu\otimes\bu\in Sym(\Ve)$. Conversely, if $\bu\otimes\bv\in Sym(\Ve)$ then $\forall\bw,\ \bw\cdot\bv\ \bu=(\bu\otimes\bv)\bw=(\bv\otimes\bu)\bw=\bw\cdot\bu\ \bv\Rightarrow\bv=\dfrac{\bw\cdot\bv}{\bw\cdot\bu}\bu\Rightarrow\bu\times\bv=\bo.$

\item $\L=\L^\top\Rightarrow(\Rb\L\Rb^\top)^\top=\Rb\L^\top\Rb^\top=\Rb\L\Rb^\top$; moreover, $\bu\cdot\L\bu>0\ \forall\bu\Rightarrow\\\bu\cdot(\Rb\L\Rb^\top)\bu=(\Rb^\top\bu)\cdot\L(\Rb^\top\bu)>0$.

\item i) $\det(\L^{sph}-\lambda\I)=\det\left(\dfrac{1}{3}\tr\L\I-\lambda^{sph}\I\right)=\left(\dfrac{1}{3}\tr\L\I-\lambda^{sph}\I\right)^3\det\I\\=\left(\dfrac{1}{3}\tr\L\I-\lambda^{sph}\I\right)^3=0\Rightarrow\lambda^{sph}_i=\dfrac{1}{3}\tr\L,\ i=1,2,3$.

ii) $(\L^{sph}-\lambda^{sph}_i\I)\bv=\bo\ \forall i=1,2,3\Rightarrow\left(\L^{sph}-\dfrac{1}{3}\tr\L\I\right)\bv=\bo\Rightarrow\\(\L^{sph}-\L^{sph})\bv=\bo\Rightarrow\O\bv=\bo$, which is true $\forall\bv$.

\item $\det(\L^{dev}-\lambda^{dev}\I)=\det(\L-\L^{dev}-\lambda^{dev}\I)=\det\left(\L-\left(\dfrac{1}{3}\tr\L+\lambda^{dev}\right)\I\right)\\=\det\left(\L-\left(\lambda^{sph}+\lambda^{dev}\right)\I\right)=0\Rightarrow\lambda=\lambda^{sph}+\lambda^{dev}$ is an eigenvalue of $\L\Rightarrow\\\lambda^{dev}=\lambda-\lambda^{sph}.$

\end{enumerate}

\section*{Chapter 3}
\begin{enumerate}
\item $\forall\L\in Lin(\Ve),\ (\ei\otimes\ej)\boxtimes(\ek\otimes\el)\L=(\ei\otimes\ej)\L(\ek\otimes\el)^\top=(\ei\otimes\ej)\L(\el\otimes\ek)=(\ei\otimes\ej)((\L\el)\otimes\ek)=\ej\cdot(\L\el)\ei\otimes\ek=L_{jl}\ei\otimes\ek$; moreover $(\ei\otimes\ek\otimes\ej\otimes\el)\L=(\ei\otimes\ek)\otimes(\ej\otimes\el)\L=((\ej\otimes\el)\cdot\L)(\ei\otimes\ek)=((\ej\otimes\el)\cdot(L_{pq}\e_p\otimes\e_q))(\ei\otimes\ek)=L_{pq}\delta_{jp}\delta_{lq}(\ei\otimes\ek)=L_{jl}\ei\otimes\ek\Rightarrow(\ei\otimes\ej)\boxtimes(\ek\otimes\el)=\ei\otimes\ek\otimes\ej\otimes\el$.

\item $\forall\L,\M\in Lin(\Ve),\ \L\cdot(\Aq\Bq)\M=\Aq^\top\L\cdot\Bq\M=\Bq^\top\Aq^\top\L\cdot\M\Rightarrow(\Aq\Bq)^\top=\Bq^\top\Aq^\top.$

\item $\forall \C\in Lin(\Ve),\ (\A\otimes\B\Lq)\C=(\A\otimes\B)\Lq\C=\B\cdot\Lq\C\A=\Lq^\top\B\cdot\C\A=(\A\otimes\Lq^\top\B)\C.$

\item $\forall\L\in Lin(\Ve),\ ((\A\boxtimes\B)(\C\boxtimes\D))\L=\A\boxtimes\B\C\L\D^\top=\A\C\L\D^\top\B^\top\\=(\A\C)\boxtimes(\D^\top\B^\top)^\top\L=(\A\C)\boxtimes(\B\D)\L.$

\item Let $\Aq=A_{ijkl}\ei\otimes\ej\otimes\ek\otimes\el=A_{ijkl}(\ei\otimes\ek)\boxtimes(\ej\otimes\el)$ and $\Bq=B_{pqrs}\ep\otimes\eq\otimes\er\otimes\es=B_{pqrs}(\ep\otimes\er)\boxtimes(\eq\otimes\es)$. \\ Then, $ \Aq\Bq=A_{ijkl}B_{pqrs}((\ei\otimes\ek)\boxtimes(\ej\otimes\el))((\ep\otimes\er)\boxtimes(\eq\otimes\es))\\=A_{ijkl}B_{pqrs}((\ei\otimes\ek)(\ep\otimes\er))\boxtimes((\ej\otimes\el)(\eq\otimes\es))\\=A_{ijkl}B_{pqrs}\delta_{kp}\delta{lq}(\ei\otimes\er)\boxtimes(\ej\otimes\es)=A_{ijkl}B_{klrs}(\ei\otimes\ej\otimes\er\otimes\es).$

\item $\forall\L\in Lin(\Ve),\ (\A\otimes\B)(\C\boxtimes\D)\L=(\A\otimes\B)\C\L\D^\top=\B\cdot(\C\L\D^\top)\A=\B\D\cdot(\C\L)\A=\C^\top\B\D\cdot\L\ \A=(\C^\top\boxtimes\D^\top)\B\cdot\L\ \A=\A\otimes((\C^\top\boxtimes\D^\top)\B)\L.$

\item $\forall\L\in Lin(\Ve),\ (\A\boxtimes\B)(\C\otimes\D)\L=\D\cdot\L(\A\boxtimes\B)\C=(\D\cdot\L)\A\C\B^\top=\A\C\B^\top(\D\cdot\L)=((\A\boxtimes\B)\C)(\D\cdot\L)=(((\A\boxtimes\B)\C)\otimes\D)\L.$

\item $(\P\otimes\P)_{ijhk}=P_{ij}P_{hk}=(\p\otimes\p)_{ij}(\p\otimes\p)_{hk}=p_ip_jp_hp_k\\=p_ip_hp_jp_k=(\p\otimes\p)_{ih}(\p\otimes\p)_{jk}=P_{ih}P_{jk}=(\P\boxtimes\P)_{ijhk}.$

\item i) $\Iq\Aq=(\I\boxtimes\I)\Aq=(\I\boxtimes\I)A_{ijkl}(\ei\otimes\ej)\otimes(\ek\otimes\el)=A_{ijkl}((\I\boxtimes\I)(\ei\otimes\ej))\otimes(\ek\otimes\el)=A_{ijkl}(\I(\ei\otimes\ej)\I^\top)\otimes(\ek\otimes\el)=A_{ijkl}(\ei\otimes\ej\otimes\ek\otimes\el)=\Aq.$

ii) $\Aq\Iq=A_{ijkl}(\ei\otimes\ej)\otimes(\ek\otimes\el)(\I\boxtimes\I)=A_{ijkl}(\ei\otimes\ej)\otimes((\I\boxtimes\I)\ek\otimes\el)=A_{ijkl}(\ei\otimes\ej)\otimes(\I(\ek\otimes\el))\I^\top)=A_{ijkl}(\ei\otimes\ej\otimes\ek\otimes\el)=\Aq.$

\item $(\A\otimes\B)\cdot(\C\otimes\D)=\tr_4((\A\otimes\B)^\top(\C\otimes\D))=\tr_4((\B\otimes\A)(\C\otimes\D))\\=\tr_4((\B\otimes\A)\C)
\otimes\D=\tr_4(\A\cdot\C\ \B\otimes\D)=\A\cdot\C\ \tr_4(\B\otimes\D)=\A\cdot\C\ \B\cdot\D.$

\item $\dfrac{\I}{|\I|}\otimes\dfrac{\I}{|\I|}=\dfrac{\I}{\sqrt{3}}\otimes\dfrac{\I}{\sqrt{3}}=\dfrac{1}{3}\I\otimes\I=\Sq^{sph}.$

\item $\forall\L\in Lin(\Ve),\ \L=\L^{sph}+\L^{dev}$ and $\L^{sph}=\dfrac{1}{3}\tr\L\ \I\rightarrow$ just one number is sufficient to determine $\L^{sph}\Rightarrow\dim(Sph(\Ve))=1$. Then, $\L^{dev}=\L-\L^{sph}$ is determined by 5 numbers : $\dim(Dev(\Ve))=\dim(Lin(\Ve)-Sph(\Ve))=6-1=5$.

\item i) $\Sq^{sph}\Sq^{sph}=\left(\dfrac{1}{3}\I\otimes\I\right)\left(\dfrac{1}{3}\I\otimes\I\right)=\dfrac{1}{9}(\I\otimes\I)=\dfrac{1}{9}\I\cdot\I\ \I\otimes\I=\dfrac{1}{3}\I\otimes\I=\Sq^{sph}.$

ii) $\Dq^{dev}\Dq^{dev}=(\Iq^s-\Sq^{sph})(\Iq^s-\Sq^{sph})=\Iq^s-2\Sq^{sph}+\Sq^{sph}\Sq^{sph}=\Iq^s-\Sq^{sph}=\Dq^{dev}.$

iii) $\Sq^{sph}\Dq^{dev}=\Sq^{sph}(\Iq^s-\Sq^{sph})=\Sq^{sph}-\Sq^{sph}\Sq^{sph}=\Sq^{sph}-\Sq^{sph}=\Oq.$

\item i) $S^{sym}_{ijkl}=(\ei\otimes\ej)\cdot\Sq^{sym}(\ek\otimes\el)=(\ei\otimes\ej)\cdot\dfrac{\ek\otimes\el+\el\otimes\ek}{2}\\=\dfrac{1}{2}(\delta_{ik}\delta_{jl}+\delta_{il}\delta_{jk}),\ 
W^{skw}_{ijkl}=(\ei\otimes\ej)\cdot\Wq^{skw}(\ek\otimes\el)=(\ei\otimes\ej)\cdot\dfrac{\ek\otimes\el-\el\otimes\ek}{2}\\=\dfrac{1}{2}(\delta_{ik}\delta_{jl}-\delta_{il}\delta_{jk}),\rightarrow S^{sym}_{ijkl}+W^{skw}_{ijkl}=\delta_{ik}\delta_{jl}=I_{ijkl}\Rightarrow\Sq^{sym}+\Wq^{dev}=\Iq.$

ii) $S^{sym}_{ijkl}-W^{skw}_{ijkl}=\delta_{il}\delta_{jk}=T^{trp}_{ijkl}\Rightarrow\Sq^{sym}-\Wq^{dev}=\Tq^{trp}.$

\item i) $\Sq^{sph}\cdot\Sq^{sph}=\left(\dfrac{1}{3}\I\otimes\I\right)\cdot\left(\dfrac{1}{3}\I\otimes\I\right)=\dfrac{1}{9}\tr_4((\I\otimes\I)^\top(\I\otimes\I))=\dfrac{1}{9}\tr_4((\I\otimes\I)(\I\otimes\I))\\=\dfrac{1}{9}\tr_4((\I\otimes\I)\I)\otimes\I=\dfrac{1}{9}\tr_4(\I\cdot\I\ \I\otimes\I)=\dfrac{1}{3}\tr_4(\I^\top\otimes\I)=\dfrac{1}{3}\I\cdot\I=1.$

ii) $\Dq^{dev}\cdot\Dq^{dev}=(\Iq^s-\Sq^{sph})\cdot(\Iq^s-\Sq^{sph})=\Iq^s\cdot\Iq^s-2\Iq^s\cdot\Sq^{sph}+\Sq^{sph}\cdot\Sq^{sph}=\dfrac{1}{4}(\delta_{ih}\delta_{jk}+\delta_{ik}\delta_{jh})(\delta_{ih}\delta_{jk}+\delta_{ik}\delta_{jh})-2\dfrac{\delta_{ih}\delta_{jk}+\delta_{ik}\delta_{jh}}{2}\dfrac{1}{3}\delta_{ij}\delta_{hk}+1=\dfrac{1}{4}(\delta_{ih}\delta_{ih}\ \delta_{jk}\delta_{jk}+2\delta_{ik}\delta_{jh}\delta_{ih}\delta_{jk}+\delta_{ik}\delta_{ik}\ \delta_{jh}\delta_{jh})-\dfrac{1}{3}(\delta_{ij}\delta_{hk}\delta_{ih}\delta_{jk}+\delta_{ij}\delta_{hk}\delta_{ik}\delta_{jh})+1.$ Now, one should check that $\delta_{ih}\delta_{ih}= \delta_{jk}\delta_{jk}=\delta_{ik}\delta_{ik}=\delta_{jh}\delta_{jh}=\delta_{ik}\delta_{jh}\delta_{ih}\delta_{jk}=\delta_{ij}\delta_{hk}\delta_{ih}\delta_{jk}=\delta_{ij}\delta_{hk}\delta_{ik}\delta_{jh}=3\Rightarrow\\\Dq^{dev}\cdot\Dq^{dev}=5.$ 

iii) $\Sq^{sph}\cdot\Dq^{dev}=\Sq^{sph}\cdot(\Iq^s-\Sq^{sph})=\Iq^s\cdot\Sq^{sph}-\Sq^{sph}\cdot\Sq^{sph}=\dfrac{\delta_{ih}\delta_{jk}+\delta_{ik}\delta_{jh}}{2}\dfrac{1}{3}\delta_{ij}\delta_{hk}-1=\dfrac{1}{6}(\delta_{ih}\delta_{jk}\delta_{ij}\delta_{hk}+\delta_{ik}\delta_{jh}\delta_{ij}\delta_{hk})-1=\dfrac{1}{6}(3+3)-1=0.$

\item $\Sy_R=\I-2\n\otimes\n\rightarrow\Sq_R=(\I-2\n\otimes\n)\boxtimes(\I-2\n\otimes\n)\\=\I\boxtimes\I-2(\n\otimes\n\boxtimes\I+\I\boxtimes\n\otimes\n+4\n\otimes\n\boxtimes\n\otimes\n)\\
=\Iq-2(\n\otimes\n\boxtimes\ei\otimes\ei+\ei\otimes\ei\boxtimes\n\otimes\n+4\n\otimes\n\otimes\n\otimes\n)\\
=\Iq-2(\n\otimes\ei\otimes\n\otimes\ei+\ei\otimes\n\otimes\ei\otimes\n+4\n\otimes\n\otimes\n\otimes\n)$. \\By components, $(\Sq_R)_{ijkl}=(\I-2\n\otimes\n)_{ik}(\I-2\n\otimes\n)_{jl}=(\delta_{ik}-2n_in_k)(\delta_{jl}-2n_jn_l)\\=\delta_{ik}\delta_{jl}-2(\delta_{ik}n_jn_l+\delta_{jl}n_in_k)+4n_in_jn_kn_l.$

\item i) $R_0=0$: This is the case of the so-called $R_0$-orthotropy. 

ii) $R_1=0$: This is the case of the square symmetry (all the components depend upon $4\theta$).

iii) $\Phi_0-\Phi_1=k\dfrac{\pi}{4},\ k\in\{0,1\}$: This is the case of the two ordinary orthotropies.

iv) $R_0=R_1=0$: This is the condition for  isotropy. Nothing depends upon $\theta\Rightarrow$ all the directions are equivalent and thus,  at the same time, axes of elastic symmetry. 

\end{enumerate}

\section*{Chapter 4}

\begin{enumerate}
\item i) $(\bu(t)+\bv(t))-(\bu(t_0)+\bv(t_0))=(t-t_0)(\bu+\bv)'+o(t-t_0)$ and also $\bu(t)-\bu(t_0)+\bv(t)-\bv(t_0)=(t-t_0)\bu'+(t-t_0)\bv'+o(t-t_0)\Rightarrow(\bu+\bv)'=\bu'+\bv'.$

The proof that $(\L+\M)'=\L'+\M'$ and that $(\Lq+\Mq)'=\Lq'+\Mq'$ can be done in a similar way.

ii) We indicate, in short, $\alpha(t)=\alpha,\bv(t)=\bv,\alpha(t_0)=\alpha_0,\bv(t_0)=\bv_0,\alpha'(t_0)=\alpha'_0,\bv'(t_0)=\bv'_0\rightarrow\alpha\bv-\alpha_0\bv_0=(\alpha\bv)'_0(t-t_0)+o(t-t_0)$; \\moreover $\alpha=\alpha_0+\alpha'_0(t-t_0)+o(t-t_0),\bv=\bv_0+\bv'_0(t-t_0)+o(t-t_0)\Rightarrow\alpha\bv-\alpha_0\bv_0=(\alpha_0+\alpha'_0(t-t_0)+o(t-t_0))(\bv_0+\bv'_0(t-t_0)+o(t-t_0))-\alpha_0\bv_0=\alpha'_0\bv_0(t-t_0)+\bv_0 o(t-t_0)+\alpha_0\bv'_0(t-t_0)+\alpha'_0\bv'_0(t-t_0)^2+\bv'_0(t-t_0)o(t-t_0)+\alpha_0o(t-t_0)+\alpha'_0(t-t_0)o(t-t_0)+o(t-t_0)^2=(\alpha'_0\bv_0+\alpha_0\bv'_0)(t-t_0)+o(t-t_0)$, so by comparison, $(\alpha\bv)'=\alpha'\bv+\alpha\bv'$.

By the same technique, one can easily prove the differentiation rule for all the product-like quantities: 
$(\bu\cdot\bv)',(\bu\times\bv)',(\bu\otimes\bv)',(\alpha\L)',(\L\bv)',(\L\M)',(\L\cdot\M)',\\(\L\otimes\M)',
(\L\boxtimes\M)',(\alpha\Lq)',(\Lq\L)',(\Lq\Mq)',(\Lq\cdot\Mq)'.$

\item The proof is the same for all the cases, we just write that for $\bv(t)$. Using the two properties shown in the previous exercise, we get: $\bv(t)=v_i(t)\ei\Rightarrow\\\bv'(t)=(v_i(t)\ei)'=v'_i(t)\ei+v_i(t)\ei'=v'_i(t)\ei$, because $\ei$ does not depend on $t$.

\item i) $p(\theta)=(a\theta\cos\theta,a\theta\sin\theta)\rightarrow c=\dfrac{2+\theta^2}{a(1+\theta^2)^\frac{3}{2}}.$

ii) $\ell=\dfrac{a}{2}(\theta\sqrt{1+\theta^2}+\ln(\theta+\sqrt{1+\theta^2}))$. 

iii) Let $p_i,p_{i+1}$ be two consecutive intersection points of the spiral ($i$ denotes the order of the intersection point) with a straight line passing through the origin and inclined at $\theta$; their distances from the origin are $r_i=a(\theta+2\pi i),r_{i+1}= a(\theta+2\pi(i+1))\Rightarrow |p_{i+1}-p_i|=r_{i+1}-r_i=2\pi a$ that does not depend upon $\theta$.

\item i) $r=a\ e^{b\theta}\Rightarrow r=0\iff\theta\rightarrow+\infty$, for $b<0,-\infty$ for $b>0$.

ii) $p(\theta)=(a\ e^{b\theta}\cos\theta,a\ e^{b\theta}\sin\theta)\rightarrow c=\dfrac{1}{a\ e^{b\theta}\sqrt{1+b^2}}.$

iii) $\ell=\dfrac{a}{b}\sqrt{1+b^2}e^{b\theta}$.

iv) With the same meaning as in the previous exercise, $|p_{i+1}-p_i|=r_{i+1}-r_i=a (e^{2\pi b}-1)e^{b(\theta+2\pi i)}$, i.e. the distance depend upon the order of the intersection: $\dfrac{r_{i+1}}{r_i}=e^{2\pi b}$, which is a geometric progression. 

v) $\btau=\dfrac{1}{\sqrt{1+b^2}}(b\cos\theta-\sin\theta,b\sin\theta+\cos\theta)\Rightarrow\\(p-o)\cdot\btau=\dfrac{ab\ e^{b\theta}}{\sqrt{1+b^2}}\Rightarrow\cos\phi=\dfrac{(p-o)\cdot\btau}{|p-o||\btau|}=\dfrac{b}{\sqrt{1+b^2}}\Rightarrow\phi$, the angle between $\tau$ and $p-o$, is constant.

vi) $\bnu=\dfrac{1}{\sqrt{1+b^2}}(-b\sin\theta-\cos\theta,b\cos\theta-\sin\theta)\Rightarrow q=p+\dfrac{1}{c}\bnu=ab\ e^{b\theta}(-sin\theta,\cos\theta)$ is a point of the evolute, whose polar equation is hence $r=ab\ e^{b\theta}$, which is still a logarithmic spiral.

\item i) $c=\dfrac{1}{a\theta}.$

ii) $ \ell=\dfrac{1}{2}a\theta^2$.

iii) $\bnu=(-\sin\theta,\cos\theta)\Rightarrow p+\rho\bnu=p+\dfrac{1}{c}\bnu=a(\cos\theta,\sin\theta)$, which is the parametric equation of a circle of center $o$ and radius $a$.

\item i) $\btau=\dfrac{1}{\sqrt{a^2+b^2}}(-a\sin\omega\theta,a\cos\omega\theta,b)\Rightarrow\cos\phi=\btau\cdot\et=\dfrac{b}{\sqrt{a^2+b^2}}$, which is independent of $\theta$.

ii) $\ell=\omega\theta\sqrt{a^2+b^2}$.

iii) $c=\dfrac{a}{a^2+b^2}.$

iv) $\vartheta=-\dfrac{b}{a^2+b^2}.$

v) Let $p_i, p_{i+1}$ be two points, intersection of a same generatrix of the cylinder with the helix, i.e. for, say, $\theta+i\dfrac{2\pi}{\omega}$ and $\theta+(i+1)\dfrac{2\pi}{\omega}\Rightarrow d=(p_{i+1}-p_i)\cdot\et=2\pi b$.

vi) By definition, a curve is a helix $\iff \btau\cdot\et=const$; differentiating gives $\btau'\cdot\et+\btau\cdot\et=\btau'\cdot\et=0\Rightarrow$ by the first equation of Frenet-Serret $c\bnu\cdot\et=0\Rightarrow\bnu\cdot\et=0\Rightarrow\bb$ is tangent to the cylinder and $\bb\cdot\et=const.$ Moreover, differentiating again, $\bnu\cdot\et+\bnu\cdot\et'=\bnu'\cdot\et=0\Rightarrow$ by the third equation of Frenet-Serret $(-c\btau-\vartheta\bb)\cdot\et=0\Rightarrow-c\btau\cdot\et=\vartheta\bb\cdot\et\Rightarrow\dfrac{c}{\vartheta}=-\dfrac{\bb\cdot\et}{\btau\cdot\et}=const.$\\
Conversely, if $p(s)$ is a curve with $\dfrac{c}{\vartheta}=\alpha=const.$ through the first and second equations  of Frenet-Serret, we get $\bnu=\dfrac{1}{c}\btau'=\dfrac{1}{\vartheta}\bb'\Rightarrow\btau'=\dfrac{c}{\vartheta}\bb'=\alpha\bb'\Rightarrow(\btau-\alpha\bb)'=0\Rightarrow \bv=\btau-\alpha\bb=const.\Rightarrow\btau\cdot\bv=1\Rightarrow\btau$ forms with $\bv$ a constant angle, and because $\bv$ is a constant vector, $\btau\cdot\et=const.\Rightarrow$ the curve is a helix.

vii) $p'\times p''=\omega^3(ab\sin\omega\theta,-ab\cos\omega\theta,a^2)\Rightarrow A=ab\omega^3,B=a^2\omega^3$.

\item i) $p(\theta)=R(\theta-\sin\theta,1-\cos\theta)$.

ii) $\ell(\theta)=4R\left(1-\cos\dfrac{\theta}{2}\right)\Rightarrow\ell(2\pi)=8R.$

iii) $c=\dfrac{1}{2R\sqrt{2(1-\cos\theta)}}.$

iv) $\bnu=\dfrac{1}{\sqrt{2(1-\cos\theta)}}(\sin\theta,\cos\theta-1)\Rightarrow q(\theta)=p(\theta)+\dfrac{1}{c}\bnu=R(\theta+\sin\theta,\cos\theta-1)$; this curve is the evolute of the cycloid and it can be obtained also as $\\q(\theta)=p(\theta+\pi)-(\pi,2)R$, i.e., it is the same cycloid $p(\theta)$ translated by $-(\pi,2)R$.

\item i) $c=\dfrac{1}{\cosh^2t}.$

ii) $\bnu=\left(-\dfrac{\sinh t}{\cosh t},\dfrac{1}{\cosh t}\right)\Rightarrow q(t)=p(t)+\dfrac{1}{c}\bnu=(t-\sinh t\cosh t,2\cosh t)$ is the evolute.

iii) $s=\int|p'(t)|dt=\sinh t, \ \btau=\left(\dfrac{1}{\cosh t},\dfrac{\sinh t}{\cosh t}\right)\Rightarrow b(t)=p(t)+(a-s)\btau=\left(t+\dfrac{a-\sinh t}{\cosh t},\cosh t+\sinh t\dfrac{a-\sinh t}{\cosh t}\right)$ is the equation of the involutes.

\item i) $\btau=(\cos t,\sin t)$; tangent to the tractrix at $p$: $z=p+w\btau=\\\left((1+w)\cos t+\ln\left(\tan\dfrac{t}{2}\right),(1+w)\sin t\right)$. Intersection point $g$ with $x_1$ axis for $\\w=-1\rightarrow g=\left(\ln\left(\tan\dfrac{t}{2}\right),0\right)\Rightarrow|p-g|=1\ \forall t$.

ii) $\ell=\ln\dfrac{\sin t_2}{\sin t_1}$.

iii) $c=\tan t$.

iv) $\bnu=(-\sin t,\cos t)\Rightarrow q=p+\dfrac{1}{c}\bnu=\left(\ln\left(\tan\dfrac{t}{2}\right),\dfrac{1}{\sin t}\right)$. Setting $\sigma=\ln\left(\tan\dfrac{t}{2}\right)\Rightarrow\tan\dfrac{t}{2}=e^\sigma, \dfrac{1}{\sin t}=\cosh\sigma\Rightarrow q=(\sigma,\cosh\sigma)$, which is the equation of a catenary.

\item i) $p(\theta)=(\cos\theta,\sin\theta,\sin\theta)\Rightarrow p'(=-\sin\theta,\cos\theta,\cos\theta), p''=(-\cos\theta,-\sin\theta,-\sin\theta)\\ \Rightarrow c=\dfrac{\sqrt{2}}{(1+\cos^2\theta)^\frac{3}{2}}\Rightarrow c_{max}=\sqrt{2}.$

ii) $p'''=(\sin\theta,-\cos\theta,-\cos\theta)\Rightarrow p'\times p''\cdot p'''=\cos\theta-\cos\theta=0\Rightarrow\vartheta=0\Rightarrow$ the curve is planar.

\item i) $\bv=\dot{p},\btau=\dfrac{\dot{p}}{|\dot{p}|}\Rightarrow\dot{p}=|\dot{p}|\btau\Rightarrow\bv=v\btau,v=|\dot{p}|$ is the scalar velocity. 

ii) $\a=\ddot{p}=\dot{\bv}=(v\btau)^\cdot=\dot{v}\btau+v\dot{\btau}; \dot{\btau}=|\dot{\btau}|\bnu;c=\dfrac{|\dot{\btau}|}{|\dot{p}|}\Rightarrow|\dot{\btau}|=|\dot{p}|c=\dfrac{v}{\rho}\Rightarrow\\ \a=\dot{v}\btau+\dfrac{v^2}{\rho}\bnu$; $\dot{v}$ is the {\it tangential acceleration} and $\dfrac{v^2}{\rho}$ the {\it centripetal acceleration}.

iii) $f=m\a=m\dot{v}\btau+\dfrac{mv^2}{\rho}\bnu=f_\tau\btau+f_\nu\bnu; f_\tau=m\dot{v}$ is the {\it tangential force}, which is responsible for the change of the scalar velocity; $f_\nu=\dfrac{mv^2}{\rho}$ is the {\it centripetal force}, which is responsible for the path change.

\end{enumerate}

\section*{Chapter 5}
\begin{enumerate}
\item i) By Eq. (\ref{eq:gradproducts})$_3$: $ \mathrm{grad}(\bv\cdot\bw)=(\mathrm{grad}\bw)^\top\bv+(\mathrm{grad}\bv)^\top\bw= (\mathrm{grad}\bw)^\top\bv+(\grad\bw)\bv-(\grad\bw)\bv+(\mathrm{grad}\bv)^\top\bw+(\grad\bv)\bw-(\grad\bv)\bw=(\grad\bw)\bv+(\grad\bv)\bw+((\mathrm{grad}\bw)^\top-\grad\bw))\bv+((\mathrm{grad}\bv)^\top-\grad\bv))\bw=(\grad\bw)\bv+(\grad\bv)\bw-(\curl\bw)\times\bv-\\(\curl\bv)\times\bw=(\grad\bw)\bv+(\grad\bv)\bw+\bv\times(\curl\bw)+\bw\times(\curl\bv).$

ii) By  Eqs. (\ref{eq:gradproducts})$_{2,3}$ and Exercise 3 iii), Chapter 2,
 $\grad(\bu\cdot\bv\ \bw)=\bu\cdot\bv\ \grad\bw+\\\bw\otimes\grad(\bu\cdot\bv)=\bu\cdot\bv\ \grad\bw+\bw\otimes((\grad\bu)^\top\bv+(\grad\bv)^\top\bu)=\bu\cdot\bv\ \grad\bw+\\\bw\otimes(\grad\bu)^\top\bv+\bw\otimes(\grad\bv)^\top\bu=\bu\cdot\bv\ \grad\bw+(\bw\otimes\bv)\grad\bu+(\bw\otimes\bu)\grad\bv.$

iii) By Theorem \ref{teo:proddiv} i), $\div((\grad\bv)\bv-(\div\bv)\bv)+(\div\bv)^2=\div((\grad\bv)\bv)-\\\div((\div\bv)\bv)+(\div\bv)^2=(\grad\bv)^\top\cdot\grad\bv+\bv\cdot\div(\grad\bv)^\top-(\div\bv)^2-\\\bv\cdot\grad(\div\bv)+(\div\bv)^2=(\grad\bv)^\top\cdot\grad\bv+\bv\cdot\div(\grad\bv)^\top-\bv\cdot\div(\grad\bv)^\top=(\grad\bv)^\top\cdot\grad\bv.$

\item i) By Theorem \ref{teo:proddiv} i), $\div(\grad\bv)^\top=\div(v_{j,i}\ei\otimes\ej)=v_{j,i}\div(\ei\otimes\ej)+\\\ei\otimes\ej\grad v_{j,i}=\ei\otimes\ej v_{j,ik}\ek=v_{j,ik}\delta_{jk}\ei=v_{j,ij}\ei=v_{j,ji}\ei=(\div\bv)_{,i}\ei=\grad(\div\bv).$

ii) By iii) of the previous exercise and Theorem \ref{teo:proddiv} i), $\div((\grad\bv)\bv-(\div\bv)\bv)=\grad\bv\cdot(\grad\bv)^\top-(\div\bv)^2$ and also $\div((\grad\bv)\bv-(\div\bv)\bv)=\div(\grad\bv)\bv)-\div((\div\bv)\bv)=\div((\grad\bv)\bv)-(\div\bv)^2-\bv\cdot\grad(\div\bv)$, so comparing the two results $\div((\grad\bv)\bv)=\grad\bv\cdot(\grad\bv)^\top+\bv\cdot\grad(\div\bv)$.

iii) By Theorem \ref{teo:proddiv} i) and iv), $\div(\phi\L\bv)=\phi\div(\L\bv)+\L\bv\cdot\grad\phi=\L\bv\cdot\grad\phi+\phi\L^\top\cdot\grad\bv+\phi\bv\cdot\div\L^\top.$

\item i) By Theorem \ref{teo:propgrad} ii) and Eq. (\ref{eq:doublecrossprod}), $\forall\a=const.\in\Ve,\ (\curl(\phi\bv))\times\a=(\grad(\phi\bv)-(\grad(\phi\bv)^\top)\a=(\phi\grad\bv+\bv\otimes\grad\phi-(\phi\grad\bv+\bv\otimes\grad\phi)^\top)\a=
(\phi\grad\bv+\\\bv\otimes\grad\phi-\phi(\grad\bv)^\top-\grad\phi\otimes\bv)\a=\phi(\curl\bv)\times\a+\a\cdot\grad\phi\ \bv-\a\cdot\bv\ \grad\phi=\phi(\curl\bv)\times\a+\a\times(\bv\times\grad\phi)=\phi(\curl\bv)\times\a-(\bv\times\grad\phi)\times\a=(\phi\curl\bv-\\\bv\times\grad\phi)\times\a\Rightarrow\curl(\phi\bv)=\phi\curl\bv+\grad\phi\times\bv.$

ii) Using the Ricci's alternator for the cross product, $\bv\times\bw=\epsilon_{pqr}v_qw_r\e_p$ and $\curl(\bv\times\bw)=\epsilon_{ijk}(\bv\times\bw)_{k,j}\ei=\epsilon_{ijk}\epsilon_{kqr}(v_qw_r)_{,j}\ei=\epsilon_{ijk}\epsilon_{kqr}(v_{q,j}w_r+v_qw_{r,j})\ei=\epsilon_{kij}\epsilon_{kqr}(v_{q,j}w_r+v_qw_{r,j})\ei$. Then, because $\epsilon_{kij}\epsilon_{kqr}=\delta_{iq}\delta_{jr}-\delta_{ir}\delta_{jq}$ we get $\curl(\bv\times\bw)=(\delta_{iq}\delta_{jr}-\delta_{ir}\delta_{jq})(v_{q,j}w_r+v_qw_{r,j})\ei=\delta_{iq}\delta_{jr}(v_{q,j}w_r+v_qw_{r,j})\ei-\delta_{ir}\delta_{jq}(v_{q,j}w_r+v_qw_{r,j})\ei=(v_{i,j}w_j+v_iw_{j,j})\ei-(v_{j,j}w_i+v_jw_{i,j})\ei=\grad\bv\ \bw-\grad\bw\ \bv+\bv\div\bw-\bw\div\bv.$

\item i) Through Theorem \ref{teo:proddiv} iv) we get $\int_{\partial\Omega}\bv\cdot\L\n\ dA=\int_{\partial\Omega}\L^\top\bv\cdot\n\ dA=\int_\Omega\div(\L^\top\bv)dV=\int_\Omega\L\cdot\grad\bv+\bv\cdot\div\L\ dV.$

ii) By Theorems \ref{teo:rotoremecaflu} and \ref{teo:proddiv} i),  $\forall\a=const.\in\Ve, \int_{\partial\Omega}(\L\n)\otimes\bv\ \a\ dA=\int_{\partial\Omega}\a\cdot\bv\ \L\n\ dA=\int_\Omega\div(\a\cdot\bv\ \L)dV=\int_\Omega\a\cdot\bv\ \div\L+\L\grad(\a\cdot\bv)dV=\int_{\Omega}\a\cdot\bv\ \div\L+\L(\grad\a)^\top\bv+\L(\grad\bv)^\top\a\ dV=\int_\Omega((\div\L)\otimes\bv+\L(\grad\bv)^\top)\a\ dV.$

iii) By Theorem \ref{teo:proddiv} ii), $\int_{\partial\Omega}(\bw\cdot\n)\bv\ dA=\int_{\partial\Omega}(\bv\otimes\bw)\n\ dA=\int_\Omega\div(\bv\otimes\bw)dV=\int_\Omega\bv\ \div\bw+(\grad\bv)\bw\ dV.$

\item i) Take $\bu=\alpha\n,\ \n\in\S\rightarrow\phi(p+\alpha\n)=\phi(p)+\alpha\ \grad\phi\cdot\n+o(\alpha)\Rightarrow\\\dfrac{d\phi}{d\n}:=\lim_{\alpha\rightarrow0}\dfrac{\phi(p+\alpha\n)-\phi(\p)}{\alpha}=\grad\phi\cdot\n$.

ii) In a similar way, $\bv(p+\alpha\n)=\bv(p)+\alpha\ \grad\bv\ \n+o(\alpha)\Rightarrow\\\dfrac{d\bv}{d\n}:=\lim_{\alpha\rightarrow0}\dfrac{\bv(p+\alpha\n)-\bv(\p)}{\alpha}=\grad\bv\ \n$.

\item i) Applying the first proof of the previous exercise to $\n=\ei,\ i=1,2,3$, we get immediately $\dfrac{df}{d\ei}:=f_{,i}=\grad f\cdot\ei:=(\grad f)_{,i}\Rightarrow\grad f=f_{,i}\ei.$

ii) Applying the second proof of the previous exercise to $\n=\ej,\ j=1,2,3$, we obtain $\dfrac{d\bv}{d\ej}:=\bv_{,j}=v_{i,j}\ei=\grad\bv\ \ej=(\grad\bv)_{ik}\ei\otimes\ek \ \ej=\delta_{kj}(\grad\bv)_{ik}\ei=(\grad\bv)_{ij}\ei\\\Rightarrow(\grad\bv)_{ij}=v_{i,j}\Rightarrow\grad\bv=v_{i,j}\ei\otimes\ej.$

iii) $\div\bv:=\tr(\grad\bv)=\tr(v_{i,j}\ei\otimes\ej)=v_{i,j}\tr(\ei\otimes\ej)=\delta_{ij}v_{i,j}=v_{i,i}.$

iv) $\forall\bu=const.\in\Ve,\ (\div\L)\cdot\bu:=\div(\L^\top\bu)\Rightarrow(\div\L)_iu_i=(\L^\top\bu)_{j,j}=(L_{ij}u_i)_{,j}=L_{ij,j}u_i\Rightarrow(\div\L)_i=L_{ij,j}\Rightarrow\div\L=L_{ij,j}\ei.$

v) $\Delta f:=\div(\grad f)=\div(f_{,i}\ei)=f_{,ii}.$

vi) $\Delta\bv:=\div(\grad \bv)=\div(v_{i,j}\ei\otimes\ej)=v_{i,jj}.$

vii) For the sake of brevity, let $\bw=\curl\bv;\ \forall\bu\in\Ve,\\ (\curl\bv)\times\bu:=(\grad\bv-\grad\bv^\top)\bu\Rightarrow
\left[\begin{array}{ccc}0 & -w_3 & w_2 \\w_3 & 0 & -w_1 \\-w_2& w_1 & 0\end{array}\right]\left\{\begin{array}{c}u_1 \\u_2 \\u_3\end{array}\right\}=\\
\left(\left[\begin{array}{ccc}v_{1,1} & v_{1,2} & v_{1,3} \\v_{2,1} & v_{2,2} & v_{2,3} \\v_{3,1} & v_{3,2} & v_{3,3}\end{array}\right]-\left[\begin{array}{ccc}v_{1,1} & v_{2,1} & v_{3,1}  \\v_{1,2} &  v_{2,2} & v_{3,2} \\v_{1,3}  & v_{2,3} & v_{3,3}\end{array}\right]\right)\left\{\begin{array}{c}u_1 \\u_2 \\u_3\end{array}\right\}=\\
\left[\begin{array}{ccc}0& v_{1,2}-v_{2,1}& v_{1,3}-v_{3,1}  \\v_{2,1}-v_{1,2} &  0&v_{2,3}- v_{3,2} \\v_{3,1}-v_{1,3}  & v_{3,2}-v_{2,3} & 0\end{array}\right]\left\{\begin{array}{c}u_1 \\u_2 \\u_3\end{array}\right\}\Rightarrow \curl\bv=\left\{\begin{array}{c}v_{3,2}-v_{2,3}\\v_{1,3}-v_{3,1} \\v_{2,1}-v_{1,2} \end{array}\right\}.$

\item i) $\bv(p)=\bv(p_0)+\bom\times(p-p_0)\Rightarrow\ \exists\W_{\omega}\in Skw(\Ve)|\ \bv(p)=\bv(p_0)+\W_\omega(p-p_0)$, with $\W_\omega$ the axial tensor of $\bom$. Moreover, by the definition of gradient, $\bv(p)=\bv(p_0)+(\grad\bv)(p-p_0)\Rightarrow\W_\omega=\grad\bv\Rightarrow\grad\bv=\dfrac{\grad\bv+\grad\bv^\top}{2}+\dfrac{\grad\bv-\grad\bv^\top}{2}=\W_\omega=-\W_\omega^\top\iff\dfrac{\grad\bv+\grad\bv^\top}{2}=\O\Rightarrow\grad\bv=\dfrac{\grad\bv-\grad\bv^\top}{2}\Rightarrow\\\bv(p)=\bv(p_0)+\dfrac{\grad\bv-\grad\bv^\top}{2}(p-p_0)$ so, by the definition of curl, $\\\bv(p)=\bv(p_0)+\dfrac{1}{2}\curl\bv\times(p-p_0)$, and comparing the two results we get $\bom=\dfrac{1}{2}\curl$.

ii) By the definition of divergence and Eq. (\ref{eq:tracciaskew}), $\div\bv=\tr(\grad\bv)=\tr\W_\omega=0.$

\item i) $\div\bu=3\alpha\rightarrow$ nowhere isochoric.

ii) $\div\bu=0\rightarrow$ globally isochoric.

iii) $\div\bu=\gamma(x_1+x_2+x_3)=0$: isochoric on the points of the plane $x_1+x_2+x_3=0$.

iv) $\div\bu=\delta(\cos x_1+\sin x_2+\cos x_3)=0$: isochoric on the points of the surface $\cos x_1+\sin x_2+\cos x_3=0$.

\item Using Eq. (\ref{eq:divvectcyl}), we get: 

i) $v_\theta=v_z=0, \ v_\rho=\dfrac{\alpha}{\rho}\Rightarrow\div\bv=-\dfrac{\alpha}{\rho^2}+\dfrac{\alpha}{\rho^2}=0.$

ii) $v_\rho=v_z=0, \ v_\theta=\dfrac{\alpha}{\rho}\Rightarrow\div\bv=0.$

iii) $v_\rho=\dfrac{\alpha\cos\theta}{\rho^2},v_\theta=\dfrac{\alpha\sin\theta}{\rho^2},v_z=0\Rightarrow\div\bu=-\dfrac{2\alpha\cos\theta}{\rho^3}+\dfrac{2\alpha\cos\theta}{\rho^3}=0.$

\item Using Eq. (\ref{eq:rotorecoordcylind}), we get:

i) $v_{\rho,\theta}=v_{\rho,z}=v_{\theta,\rho}=v_{\theta,z}=v_{z,\rho}=v_{z,\theta}=0\Rightarrow\curl\bv=\bo.$

ii) $v_{\rho,\theta}=v_{\rho,z}=v_{\theta,z}=v_{z,\rho}=v_{z,\theta}=0, v_{\theta,\rho}=-\dfrac{\alpha}{\rho^2}\Rightarrow\curl\bv=\left(0,0,\dfrac{\alpha}{\rho^2}-\dfrac{\alpha}{\rho^2}\right)=\bo.$

iii) $v_{\rho,\theta}=-\dfrac{\alpha\cos\theta}{\rho^2},v_{\rho,z}=0,v_{\theta,\rho}=-\dfrac{2\alpha\sin\theta}{\rho^3},v_{\theta,z}=0,v_{z,\rho}=v_{z_,\theta}=0\Rightarrow\\\curl\bv=\left(0,0,\dfrac{\alpha\cos\theta}{\rho^3}-2\dfrac{\alpha\cos\theta}{\rho^3}+\dfrac{\alpha\cos\theta}{\rho^3}\right)=\bo.$

\end{enumerate}

\section*{Chapter 6}

\begin{enumerate}
\item Setting $\rho=z^1,\theta=z^2,z=z^3$, by Eq. (\ref{eq:gcov1}) we get $\g=\left[\begin{array}{ccc}1 & 0 & 0 \\0 & \rho^2 & 0 \\0 & 0 & 1\end{array}\right]\Rightarrow\\ ds=\sqrt{g_{hk}dz^hdz^k}=\sqrt{d\rho^2+\rho^2d\theta^2+dz^2}.$

\item Setting $r=z^1,\theta=z^2,\phi=z^3$, proceeding in a similar way, we get $\\\g=\left[\begin{array}{ccc}1 & 0 & 0 \\0 & r^2\sin^2\phi & 0 \\0 & 0 & r^2\end{array}\right]\Rightarrow ds=\sqrt{g_{hk}dz^hdz^k}=\sqrt{dr^2+r^2\sin^2\phi d\theta^2+r^2d\phi^2}.$

\item For cylindrical coordinates, cf. Exercise  \ref{ex:1ch6}, $ds=\sqrt{d\rho^2+\rho^2d\theta^2+dz^2}$, and for a curve on a circular cylinder, $\rho=R\Rightarrow d\rho=0\Rightarrow ds=\sqrt{R^2d\theta^2+dz^2}$; if the equation of the helix is $p(\theta)=R\cos\theta\eu+R\sin\theta\ed+b\theta\et$, then $\dfrac{dz}{d\theta}=b\Rightarrow dz=b\ d\theta\Rightarrow\\ ds=\sqrt{R^2+b^2}d\theta\Rightarrow\ell=\int_\theta^{\theta+2\pi}\sqrt{R^2+b^2}d\theta=2\pi\sqrt{R^2+b^2}.$

\item i) $r=\dfrac{2R}{\pi}\theta\Rightarrow dr=\dfrac{2R}{\pi}d\theta;\ ds=\sqrt{dr^2+r^2d\theta^2}=\dfrac{2R}{\pi}\sqrt{1+\theta^2}d\theta\Rightarrow\ell=\int_0^{\frac{\pi}{2}}ds=\dfrac{2R}{\pi}\int_0^{\frac{\pi}{2}}\sqrt{1+\theta^2}d\theta=\dfrac{R}{\pi}[\theta\sqrt{1+\theta^2}+\textrm{arcsinh}\theta]_0^{\frac{\pi}{2}}=\left(\dfrac{1}{4}\sqrt{4+\pi^2}+\dfrac{1}{\pi}\mathrm{arcsinh}\dfrac{\pi}{2}\right)R\sim1.324R.$

ii) $\dfrac{\pi}{2}R=2\pi R_0\Rightarrow R_0=\dfrac{R}{4}\Rightarrow h=\sqrt{R^2-R_0^2}=\dfrac{15}{4}R;\ \rho(z)=\dfrac{R_0}{h}z,\ z=\dfrac{h}{2\pi}\theta\Rightarrow\rho(\theta)=\dfrac{R_0}{2\pi}\theta\Rightarrow\rho(\theta)=\dfrac{R}{8\pi}\theta,\ z(\theta)=\dfrac{\sqrt{15}}{8\pi}R\theta\Rightarrow d\rho=\dfrac{R}{8\pi}d\theta,\ dz=\dfrac{\sqrt{15}}{8\pi}R\ d\theta$; equation of the conical helix: $p(\theta)=\rho(z)\cos\theta\eu+\rho(z)\sin\theta\ed+\sqrt{15}\theta\et=\dfrac{R}{8\pi}(\theta\cos\theta\eu+\theta\sin\theta\ed+\sqrt{15}\theta\et)\Rightarrow ds=\sqrt{d\rho^2+\rho^2d\theta^2+dz^2}=\dfrac{R}{8\pi}\sqrt{d\theta^2+\theta^2d\theta^2+15d\theta^2}=\dfrac{R}{8\pi}\sqrt{16+\theta^2}d\theta\Rightarrow\ell=\int_0^{2\pi}ds=\dfrac{R}{8\pi}\int_0^{2\pi}\sqrt{16+\theta^2}d\theta=\\\dfrac{R}{8\pi}\left[\dfrac{1}{2}\theta\sqrt{16+\theta^2}+8\textrm{arcsinh}\dfrac{\theta}{4}\right]_0^{2\pi}\sim1.324R.$

\item i) Referring to Fig. \ref{fig:26} and by Eq. (\ref{eq:gcov1}), $x_1=z^1\cos\alpha_1+z^2\cos\alpha_2,x_2=z^1\sin\alpha_1+z^2\sin\alpha_2\Rightarrow\g=\left[\begin{array}{cc}1 & \cos(\alpha_2-\alpha_1)  \\\cos(\alpha_2-\alpha_1)  & 1\end{array}\right].$

ii) By Eq. (\ref{eq:vtang1}), $\g_1=\cos\alpha_1\eu+\sin\alpha_1\ed,\ \g_2=\cos\alpha_2\eu+\sin\alpha_2\ed.$

iii) $z^1=h(x_1\sin\alpha_2-x_2\cos\alpha_2),\ z^2=h(-x_1\sin\alpha_1+x_2\cos\alpha_1),\ h=\dfrac{1}{\sin(\alpha_2-\alpha_1)}\Rightarrow$ by Eq. (\ref{eq:dualbasis}), $\g^1=\dfrac{\sin\alpha_2}{\sin(\alpha_2-\alpha_1)}\eu-\dfrac{\cos\alpha_2}{\sin(\alpha_2-\alpha_1)}\ed,\ \g^2=-\dfrac{\sin\alpha_1}{\sin(\alpha_2-\alpha_1)}\eu+\dfrac{\cos\alpha_1}{\sin(\alpha_2-\alpha_1)}\ed.$

iv) $\g_1\cdot\g^1=\dfrac{\cos\alpha_1\sin\alpha_2-\sin\alpha_1\cos\alpha_2}{\sin(\alpha_2-\alpha_1)}=1,\\
\g_2\cdot\g^2=\dfrac{-\sin\alpha_1\cos\alpha_2+\sin\alpha_2\cos\alpha_1}{\sin(\alpha_2-\alpha_1)}=1,\\
\g_1\cdot\g^2=\dfrac{-\cos\alpha_1\sin\alpha_1+\sin\alpha_1\cos\alpha_1}{\sin(\alpha_2-\alpha_1)}=0,\\
\g_2\cdot\g^1=\dfrac{\cos\alpha_2\sin\alpha_2-\sin\alpha_2\cos\alpha_2}{\sin(\alpha_2-\alpha_1)}=0.$

v) $|\g_1|=|\g_2|=1,\ |\g^1|=\sqrt{\sin^2\alpha_1+\dfrac{\sin^2\alpha_2}{\sin^2(\alpha_2-\alpha_1)}},÷ |\g^2|=\dfrac{1}{|\sin(\alpha_2-\alpha_1)|}.$

vi) \begin{center}\includegraphics[scale=.8]{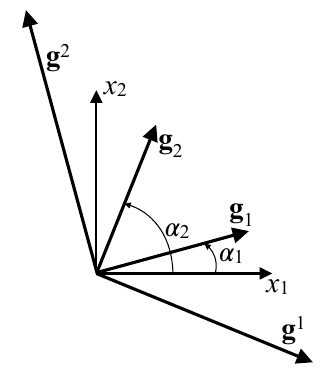}\end{center}

\item Referring to Exercise \ref{ex:2ch6} and by Eq. (\ref{eq:vtang1}), we get $\g_1=\cos\theta\sin\phi\eu+\sin\theta\sin\phi\e2+\cos\phi\et,\ \g_2=-r\sin\theta\sin\phi\eu+r\cos\theta\sin\phi\ed,\ \g_3=r\cos\theta\cos\phi\eu+r\sin\theta\cos\phi\ed-r\sin\phi\et.$

\item i) If $z^1=const.\Rightarrow x_1=a \cos z^2,\ x_2=b \sin z^2,$ with $a=c\ \cosh z^1=const.,\ b=c\ \sinh z^1=const.\Rightarrow\dfrac{x_1^2}{a^2}+\dfrac{x_2^2}{b^2}=1$: family of ellipses all with the same focuses $x_e=\pm\sqrt{a^2-b^2}=\pm c.$

ii) If $z^2=const.\Rightarrow x_1=A\cosh z^1,\ x_2=B\sinh z^1,$ with $A=c\ \cos z^2=const.,\ B=c\ \sin z^2=const.\Rightarrow\dfrac{x_1^2}{A^2}-\dfrac{x_2^2}{B^2}=1$: family of hyperbolae all with the same focuses $x_h=\pm\sqrt{A^2+B^2}=\pm c\Rightarrow x_e=x_h$.

iii) The axes of the ellipses are $2a=2c\ \cosh z^1$ and $2b=2c\ \sinh z^1$

iv) A crack along the horizontal axis corresponds to $b\rightarrow0$, which happens $\iff z^1\rightarrow0\Rightarrow\cosh z^1\rightarrow1$ and $a\rightarrow c\Rightarrow$ length of the crack: $2c$.

v) Applying  Eq. (\ref{eq:gcov1}), we get $\g=\dfrac{c^2}{2}(\cosh 2z^1-\cos2z^2)\I.$

vi) By Eq. (\ref{eq:vtang1}), $\g_1=c\ \sinh z^1\cos z^2\eu+c\ \cosh z^1\sin z^2\ed,\ \g_2=-c\ \cosh z^1\sin z^2\eu+c\ \sinh z^1\cos z^2\ed$. We note that $\g_1\cdot\g_2=0$.

\item i) By Eq. (\ref{eq:transftenscartcov})$_1$, setting $z^1=\rho,z^2=\theta,z^3=z$, we get:
$\medskip\\\begin{array}{l}
L^{11}=L_{11}^x\cos^2\theta+(L_{12}^x+L_{21}^x)\sin\theta\cos\theta+L_{22}^x\sin^2\theta,\\
L^{12}=\dfrac{1}{\rho}((L_{22}^x-L_{11}^x)\sin\theta\cos\theta+L_{12}^x\cos^2\theta-L_{21}^x\sin^2\theta),\\
L^{13}=L_{13}^x\cos\theta+L_{23}^x\sin\theta,\\
L^{21}=\dfrac{1}{\rho}((L_{22}^x-L_{11}^x)\sin\theta\cos\theta-L_{12}^x\sin^2\theta+L_{21}^x\cos^2\theta),\\
L^{22}=\dfrac{1}{\rho^2}(L_{11}^x\sin^2\theta-(L_{12}^x+L_{21}^x)\sin\theta\cos\theta+L_{22}^x\cos^2\theta),\\
L^{23}=-L_{13}^x\sin\theta+L_{23}^x\cos\theta,\smallskip\\
L^{31}=L_{31}^x\cos\theta+L_{32}^x\sin\theta,\smallskip\\
L^{32}=-L_{31}^x\sin\theta+L_{32}^x\cos\theta,\smallskip\\
L^{33}=L_{33}^x.
\end{array}$

ii) The covariant components can alternatively be found by  Eq. (\ref{eq:transftenscartcov})$_2$ or, using the results of the previous point, by Eq. (\ref{eq:raisingloweringL})$_2$; by this latter way, using the result of Exercise \ref{ex:1ch6}, we get easily
$L_{11}=L^{11},\  
L_{12}=\rho^2L^{12},\ 
L_{13}=L^{13},\ 
L_{21}=\rho^2L^{21},\\ 
L_{22}=\rho^4L^{22},\ 
L_{23}=\rho^2L^{23},\ 
L_{31}=L^{31},\ 
L_{32}=\rho^2L^{32},\ 
L_{33}=L^{33}.$

\item i) In this case, we first calculate the covariant components: By Eq. (\ref{eq:transftenscartcov})$_2$, setting $z^1=r,z^2=\theta,z^3=\phi$, we get:
$\medskip\\\begin{array}{l}
L_{11}=L_{11}^x\cos^2\theta\sin^2\phi+(L_{12}^x+L_{21}^x)\sin\theta\cos\theta\sin^2\phi+(L_{13}^x+L_{31}^x)\cos\theta\sin\phi\cos\phi\\\hspace{7mm}+L_{22}^x\sin^2\theta\sin^2\phi+(L_{23}^x+L_{32}^x)\sin\theta\sin\phi\cos\phi+L_{33}^x\cos^2\phi,\smallskip\\
L_{12}=-r\cos\theta\sin\theta\sin^2\phi L_{11}^x+r\sin^2\phi(L_{12}^x\cos^2\theta-L_{21}^x\sin^2\theta)\\\hspace{7mm}+r \sin\theta\cos\theta\sin^2\phi L_{22}^x+r\sin\phi\cos\phi(L_{32}^x\cos\theta-L_{31}^x\sin\theta),\smallskip\\
L_{13}=r\cos^2\theta\sin\phi\cos\phi L_{11}^x+r\sin\theta\cos\theta\sin\phi\cos\phi(L_{12}^x+L_{21}^x)\\\hspace{7mm}+r \sin^2\theta\cos\phi\sin\phi L_{22}^x-r\sin^2\phi\sin\theta L_{23}^x+r\cos\theta (L_{31}^x\cos^2\phi-L_{13}^x\sin^2\phi)\\\hspace{7mm}-r\cos\phi\sin\phi L_{33}^x,\smallskip\\
L_{21}=-r\sin\theta\cos\theta\sin^2\phi L_{11}^x+r\sin^2\phi(L_{21}^x\cos^2\theta-L_{12}^x\sin^2\theta)\\\hspace{7mm}+r\sin\theta\cos\theta\sin^2\phi L_{22}^x+r\sin\phi\cos\phi(L_{23}^x\cos\theta-L_{13}^x\sin\theta),\smallskip\\
L_{22}=r^2 \sin^2\theta\sin^2\phi L_{11}^x-r^2\sin\theta\cos\theta\sin^2\phi(L_{12}^x+L_{21}^x)+r^2\cos^2\theta\sin^2\phi L_{22}^x,\smallskip\\
L_{23}=-r^2\sin\theta\cos\theta \sin\phi\cos\phi (L_{22}^x-L_{11}^x)+r^2\sin\phi\cos\phi(L_{21}^x\cos^2\theta-L_{12}^x\sin^2\theta)\\\hspace{7mm}+r^2\sin^2\phi(L_{13}^x\sin\theta-L_{23}^x\cos\theta),\smallskip\\
L_{31}=r\sin\phi\cos\phi(L_{11}^x\cos^2\theta+L_{22}^x\sin^2\theta)+r\sin\theta\cos\theta\sin\phi\cos\phi(L_{12}^x+L_{21}^x)\\\hspace{7mm}+r\cos\theta(L_{13}^x\cos^2\phi-L_{31}^x\sin^2\phi)+r\sin\theta(L_{23}^x\cos^2\phi-L_{32}^x\sin^2\phi)\\\hspace{7mm}-r\sin\phi\cos\phi L_{33}^x,\smallskip\\
L_{32}=r^2\sin\theta\cos\theta\sin\phi\cos\phi(L_{22}^x-L_{11}^x)+r^2\sin\phi\cos\phi(L_{12}^x\cos^2\theta-L_{21}^x\sin^2\theta)\\\hspace{7mm}+r^2\sin^2\phi(L_{31}^x\sin\theta-L_{32}^x\cos\theta),\smallskip\\
L_{33}=r^2\cos^2\phi(L_{11}^x\cos^2\theta+L_{22}^x\sin^2\theta)+r^2\sin\theta\cos\theta\cos^2\phi(L_{12}^x+L_{21}^x)\\\hspace{7mm}-r^2\cos\theta\sin\phi\cos\phi(L_{13}^x+L_{31}^x)-r^2\sin\theta\sin\phi\cos\phi(L_{23}^x+L_{32}^x)+r^2\sin^2\phi L_{33}^x.
\end{array}$

ii) For the contravariant components, we use Eq. (\ref{eq:raisingloweringL})$_1$, after having calculated $\g^{cont}$; this can be done either using Eq. (\ref{eq:gcont1}) or simply observing that $\g^{cov}$ is diagonal (see Exercise \ref{ex:2ch6}) and that $g^{pq}=\dfrac{1}{g_{pq}}\Rightarrow\g^{cont}=\left[\begin{array}{ccc}1 & 0 & 0 \\0 & \dfrac{1}{r^2\sin^2\phi} & 0 \\0 & 0 & \dfrac{1}{r^2}\end{array}\right]\Rightarrow\\
L^{11}=L_{11},\ L^{12}=\dfrac{L_{12}}{r^2\sin^2\phi},\ L^{13}=\dfrac{L_{13}}{r^2},\ L^{21}=\dfrac{L_{21}}{r^2\sin^2\phi},\ L^{22}=\dfrac{L_{22}}{r^4\sin^4\phi},\\ L^{23}=\dfrac{L_{23}}{r^4\sin^2\phi},\ L^{31}=\dfrac{L_{31}}{r^2},\ L^{32}=\dfrac{L_{32}}{r^4\sin^2\phi},\ L^{33}=\dfrac{L_{33}}{r^4}.$

\item i) $\tr\L=L_{hh}^x$,  Eq. (\ref{eq:tracciaL}).

ii) By Eq. (\ref{eq:tracciaL})$_1\Rightarrow L_{hh}^x=\dfrac{\partial x_h}{\partial z^i}\dfrac{\partial x_k}{\partial z^j}L^{ij}\delta_{hk}=g_{ij}L^{ij}$.

iii) By Eq. (\ref{eq:tracciaL})$_2\Rightarrow L_{hh}^x=\dfrac{\partial z^i}{\partial x_h}\dfrac{\partial z^j}{\partial x_k}L_{ij}\delta_{hk}=g^{ij}L_{ij}$.

iv) By Eq. (\ref{eq:tracciaL})$_3\Rightarrow L_{hh}^x=\dfrac{\partial x_h}{\partial z^i}\dfrac{\partial z^j}{\partial x_k}L^i_{\ j}\delta_{hk}=\delta_i^{\ j}L^i_{\ j}=L^i_{\ i} $.

v) By Eq. (\ref{eq:tracciaL})$_4\Rightarrow L_{hh}^x=\dfrac{\partial z^i}{\partial x_h}\dfrac{\partial x_k}{\partial z^j}L_i^{\ j}\delta_{hk}=\delta^i_{\ j}L_i^{\ j}=L_j^{\ j} $.

\item $\dfrac{1}{2}g^{hm}\left(\dfrac{\partial g_{mk}}{\partial z^l}+\dfrac{\partial g_{ml}}{\partial z^k}-\dfrac{\partial g_{kl}}{\partial z^m}\right)=\smallskip\\\dfrac{1}{2}\dfrac{\partial z^h}{\partial x_p}\dfrac{\partial z^m}{\partial x_p}
\left(\dfrac{\partial }{\partial z^l}\dfrac{\partial x_p}{\partial z^m}\dfrac{\partial x_p}{\partial z^k}
+\dfrac{\partial }{\partial z^k}\dfrac{\partial x_p}{\partial z^m}\dfrac{\partial x_p}{\partial z^l}
-\dfrac{\partial }{\partial z^m}\dfrac{\partial x_p}{\partial z^k}\dfrac{\partial x_p}{\partial z^l}\right)=\smallskip\\
\dfrac{1}{2}\dfrac{\partial z^h}{\partial x_p}\dfrac{\partial z^m}{\partial x_p}
\left(\dfrac{\partial^2x_p }{\partial z^l\partial z^m}\dfrac{\partial x_p}{\partial z^k}
+\dfrac{\partial x_p}{\partial z^m}\dfrac{\partial^2 x_p}{\partial z^l\partial z^k}
+\dfrac{\partial^2x_p }{\partial z^k\partial z^m}\dfrac{\partial x_p}{\partial z^l}
+\dfrac{\partial x_p}{\partial z^m}\dfrac{\partial^2 x_p}{\partial z^k\partial z^l}\right.\smallskip\\
\left.-\dfrac{\partial^2x_p }{\partial z^m\partial z^k}\dfrac{\partial x_p}{\partial z^l}
-\dfrac{\partial^2x_p }{\partial z^m\partial z^l}\dfrac{\partial x_p}{\partial z^k}\right)
=\dfrac{\partial z^h}{\partial x_p}\dfrac{\partial z^m}{\partial x_p}\dfrac{\partial x_p}{\partial z^m}\dfrac{\partial^2x_p }{\partial z^k\partial z^l}
=\dfrac{\partial z^h}{\partial x_p}\dfrac{\partial^2x_p }{\partial z^k\partial z^l}=\Gamma^h_{kl}.$

\item First, we remark that $g^{im}g_{ik}=\dfrac{\partial z^i}{\partial x_p}\dfrac{\partial z^m}{\partial x_p}\dfrac{\partial x_q}{\partial z^i}\dfrac{\partial x_q}{\partial z^k}=\delta_{pq}\dfrac{\partial z^m}{\partial x_p}\dfrac{\partial x_q}{\partial z^k}=\dfrac{\partial z^m}{\partial z^k}=\delta_{mk}$, and similarly $g^{im}g_{ij}=\delta_{mj}$. Then,
$\Gamma^i_{jh}g_{ik}+\Gamma^i_{kh}g_{ji}=\dfrac{1}{2}\left(g^{im}\left(\dfrac{\partial g_{mj}}{\partial z^h}+\dfrac{\partial g_{mh}}{\partial z^j}-\dfrac{\partial g_{jh}}{\partial z^m}\right)g_{ik}+\right.\\
\left.g^{im}\left(\dfrac{\partial g_{mk}}{\partial z^h}+\dfrac{\partial g_{mh}}{\partial z^k}-\dfrac{\partial g_{kh}}{\partial z^m}\right)g_{ji}\right)=
\dfrac{1}{2}\left(g^{im}g_{ik}\left(\dfrac{\partial g_{mj}}{\partial z^h}+\dfrac{\partial g_{mh}}{\partial z^j}-\dfrac{\partial g_{jh}}{\partial z^m}\right)+\right.\\
\left.g^{im}g_{ij}\left(\dfrac{\partial g_{mk}}{\partial z^h}+\dfrac{\partial g_{mh}}{\partial z^k}-\dfrac{\partial g_{kh}}{\partial z^m}\right)\right)=
\dfrac{1}{2}\left(\delta_{mk}\left(\dfrac{\partial g_{mj}}{\partial z^h}+\dfrac{\partial g_{mh}}{\partial z^j}-\dfrac{\partial g_{jh}}{\partial z^m}\right)+\right.\\
\left.\delta_{mj}\left(\dfrac{\partial g_{mk}}{\partial z^h}+\dfrac{\partial g_{mh}}{\partial z^k}-\dfrac{\partial g_{kh}}{\partial z^m}\right)\right)=
\dfrac{1}{2}\left(\dfrac{\partial g_{kj}}{\partial z^h}+\dfrac{\partial g_{kh}}{\partial z^j}-\dfrac{\partial g_{jh}}{\partial z^k}+\dfrac{\partial g_{jk}}{\partial z^h}+\dfrac{\partial g_{jh}}{\partial z^k}-\dfrac{\partial g_{kh}}{\partial z^j}\right)\\=\dfrac{\partial g_{jk}}{\partial z^h}.$

\item i) $g_{\rho\rho}=g_{zz}=1,g_{\theta\theta}=\rho^2$ and the other components are null $\Rightarrow g^{\rho\rho}=g^{zz}=1,g^{\theta\theta}=\dfrac{1}{\rho^2}\Rightarrow\\
\Gamma_{\theta\theta}^\rho=\dfrac{1}{2}g^{\rho m}\left(\dfrac{\partial g_{m\theta}}{\partial\theta}+\dfrac{\partial g_{m\theta}}{\partial\theta}-\dfrac{\partial g_{\theta\theta}}{\partial z^m}\right)=-\dfrac{1}{2}g^{\rho\rho}\dfrac{\partial g_{\theta\theta}}{\partial \rho}=-\rho,\\
\Gamma_{\rho\theta}^\theta=\dfrac{1}{2}g^{\theta m}\left(\dfrac{\partial g_{m\rho}}{\partial\theta}+\dfrac{\partial g_{m\theta}}{\partial\rho}-\dfrac{\partial g_{\rho\theta}}{\partial z^m}\right)=\dfrac{1}{2}g^{\theta\theta}\dfrac{\partial g_{\theta\theta}}{\partial \rho}=\dfrac{1}{\rho},$
the other $\Gamma_{ij}^k$ are null.

ii) $g_{rr}=1,g_{\phi\phi}=r^2\sin^2\phi,g_{\theta\theta}=r^2\Rightarrow g^{rr}=1,g^{\phi\phi}=\dfrac{1}{r^2\sin^2\phi},g^{\theta\theta}=\dfrac{1}{r^2}$ and the other components are null $\Rightarrow\\
\Gamma_{\phi r}^\phi=\dfrac{1}{2}g^{\phi m}\left(\dfrac{\partial g_{m\phi}}{\partial r}+\dfrac{\partial g_{m r}}{\partial\phi}-\dfrac{\partial g_{\phi r}}{\partial z^m}\right)=\dfrac{1}{2}g^{\phi\phi}\dfrac{\partial g_{\phi\phi}}{\partial r}=\dfrac{1}{r},\\
\Gamma_{\theta r}^\theta=\dfrac{1}{2}g^{\theta m}\left(\dfrac{\partial g_{\theta m}}{\partial r}+\dfrac{\partial g_{r m}}{\partial\theta}-\dfrac{\partial g_{\theta r}}{\partial z^m}\right)=\dfrac{1}{2}g^{\theta\theta}\dfrac{\partial g_{\theta\theta}}{\partial r}=\dfrac{1}{r},\\
\Gamma_{\phi \phi}^r=\dfrac{1}{2}g^{r m}\left(\dfrac{\partial g_{\phi m}}{\partial \phi}+\dfrac{\partial g_{\phi m}}{\partial \phi}-\dfrac{\partial g_{\phi \phi}}{\partial z^m}\right)=-\dfrac{1}{2}g^{rr}\dfrac{\partial g_{\phi\phi}}{\partial r}=-r,\\
\Gamma_{\theta \theta}^r=\dfrac{1}{2}g^{r m}\left(\dfrac{\partial g_{\theta m}}{\partial \theta}+\dfrac{\partial g_{\theta m}}{\partial \theta}-\dfrac{\partial g_{\theta \theta}}{\partial z^m}\right)=-\dfrac{1}{2}g^{rr}\dfrac{\partial g_{\theta\theta}}{\partial r}=-r\sin^2\phi,\\
\Gamma_{\theta \phi}^\theta=\dfrac{1}{2}g^{\theta m}\left(\dfrac{\partial g_{\theta m}}{\partial \phi}+\dfrac{\partial g_{\phi m}}{\partial \theta}-\dfrac{\partial g_{\theta \phi}}{\partial z^m}\right)=\dfrac{1}{2}g^{\theta\theta}\dfrac{\partial g_{\theta\theta}}{\partial \phi}=\cot\phi,\\
\Gamma_{\theta \theta}^\phi=\dfrac{1}{2}g^{\phi m}\left(\dfrac{\partial g_{\theta m}}{\partial \theta}+\dfrac{\partial g_{\theta m}}{\partial \theta}-\dfrac{\partial g_{\theta \theta}}{\partial z^m}\right)=-\dfrac{1}{2}g^{\phi\phi}\dfrac{\partial g_{\theta\theta}}{\partial \phi}=-\sin\phi\cos\phi,\\
\Gamma_{r\phi}^\phi=\Gamma_{\phi r}^\phi, \Gamma_{r\theta}^\theta=\Gamma_{\theta r}^\theta,\Gamma_{\phi\theta}^\theta=\Gamma_{\theta\phi}^\theta,
$
 the other $\Gamma_{ij}^k$ are null.

iii) $g_{11}=g_{22}=\dfrac{c^2}{2}(\cosh2z^1-\cos2z^2)\Rightarrow g^{11}=g^{22}=\dfrac{2}{c^2(\cosh2z^1-\cos2z^2)}$ and the other components are null $\Rightarrow\\
\Gamma_{11}^1=\dfrac{1}{2}g^{1m}\left(\dfrac{\partial g_{1 m}}{\partial z^1}+\dfrac{\partial g_{1 m}}{\partial z^1}-\dfrac{\partial g_{11}}{\partial z^m}\right)=\dfrac{1}{2}g^{11}\dfrac{\partial g_{1 1}}{\partial z^1}=\dfrac{\sinh2z^1}{\cosh2z^1-\cos2z^2},\\
\Gamma_{12}^1=\dfrac{1}{2}g^{1m}\left(\dfrac{\partial g_{1 m}}{\partial z^2}+\dfrac{\partial g_{2 m}}{\partial z^1}-\dfrac{\partial g_{12}}{\partial z^m}\right)=\dfrac{1}{2}g^{11}\dfrac{\partial g_{1 1}}{\partial z^2}=\dfrac{\sin2z^2}{\cosh2z^1-\cos2z^2},
$

and $\Gamma_{12}^2=\Gamma_{21}^2=\Gamma_{11}^1=-\Gamma_{22}^1,
\Gamma_{22}^2= \Gamma_{12}^1=\Gamma_{21}^1=-\Gamma_{11}^2.$

\item Applying Eq. (\ref{eq:laplaciancurvcoord}), we get:

 i) $\Delta f=\dfrac{\partial}{\partial \rho}\left(g^{\rho k}\dfrac{\partial f}{\partial z^k}\right)+\Gamma_{\rho j}^\rho g^{jk}\dfrac{\partial f}{\partial z^k}+
 \dfrac{\partial}{\partial \theta}\left(g^{\theta k}\dfrac{\partial f}{\partial z^k}\right)+\Gamma_{\theta j}^\theta g^{jk}\dfrac{\partial f}{\partial z^k}+
 \dfrac{\partial}{\partial z}\left(g^{zk}\dfrac{\partial f}{\partial z^k}\right)+\Gamma_{zj}^zg^{jk}\dfrac{\partial f}{\partial z^k}=
 \dfrac{\partial}{\partial \rho}g^{\rho\rho}\dfrac{\partial f}{\partial \rho}+\Gamma_{\rho\theta}^\rho g^{\theta\theta}\dfrac{\partial f}{\partial \theta}+\dfrac{\partial}{\partial \theta}g^{\theta\theta}\dfrac{\partial f}{\partial \theta}+\Gamma_{\theta\theta}^\theta g^{\theta\theta}\dfrac{\partial f}{\partial \theta}+ \Gamma_{\theta \rho}^\theta g^{\rho\rho}\dfrac{\partial f}{\partial \rho}+\dfrac{\partial}{\partial z}g^{zz}\dfrac{\partial f}{\partial z}=
 \dfrac{\partial^2f}{\partial\rho^2}+\dfrac{1}{\rho^2}\dfrac{\partial^2f}{\partial\theta^2}+\dfrac{1}{\rho}\dfrac{\partial f}{\partial \rho}+\dfrac{\partial^2f}{\partial z^2}=\dfrac{1}{\rho}(\rho\ f_{,\rho})_{,\rho}+\dfrac{1}{\rho^2}f_{,\theta\theta}+f_{,zz},
 $ which is the same already found in Section \ref{sec:cylcoord}.

ii) $\Delta f=\dfrac{\partial}{\partial r}\left(g^{rk}\dfrac{\partial f}{\partial z^k}\right)+\Gamma_{rj}^rg^{jk}\dfrac{\partial f}{\partial z^k}+
\dfrac{\partial}{\partial \theta}\left(g^{\theta k}\dfrac{\partial f}{\partial z^k}\right)+\Gamma_{\theta j}^\theta g^{jk}\dfrac{\partial f}{\partial z^k}+
\dfrac{\partial}{\partial \phi}\left(g^{\phi k}\dfrac{\partial f}{\partial z^k}\right)+\Gamma_{\phi j}^\phi g^{jk}\dfrac{\partial f}{\partial z^k}=\hspace{-1mm}
\dfrac{\partial}{\partial r}g^{rr}\dfrac{\partial f}{\partial r}+
\Gamma_{\theta r}^\theta g^{rk}\dfrac{\partial f}{\partial z^k}+\dfrac{\partial}{\partial \theta}g^{\theta\theta}\dfrac{\partial f}{\partial \theta}+\Gamma_{\theta \phi}^\theta g^{\phi k}\dfrac{\partial f}{\partial z^k}+
\dfrac{\partial}{\partial \phi}g^{\phi\phi}\dfrac{\partial f}{\partial \phi}+\Gamma_{ \phi r}^\phi g^{r k}\dfrac{\partial f}{\partial z^k}\hspace{-1mm}=
\dfrac{\partial^2f}{\partial r^2}+\dfrac{1}{r^2\sin^2\phi}\dfrac{\partial^2f}{\partial \theta^2}+\dfrac{2}{r}\dfrac{\partial f}{\partial r}+\cot\phi\dfrac{1}{r^2}
\dfrac{\partial f}{\partial \phi}+\dfrac{1}{r^2}\dfrac{\partial^2f}{\partial \phi^2}=\\\dfrac{1}{r^2}(r^2f_{,r})_{,r}+\dfrac{1}{r^2\sin\phi}\left(\dfrac{f_{,\theta\theta}}{\sin\phi}+(f_{,\phi}\sin\phi)_{,\phi}\right)$, which is to be compared to the one given in Section \ref{sec:sphercoord}.

\item i) By Eqs. (\ref{eq:gcont1}), (\ref{eq:christoffelsymb1}) and (\ref{eq:covderivL1}), $g^{np}_{;h}=\dfrac{\partial g^{np}}{\partial z^h}+\Gamma_{hr}^ng^{rp}+\Gamma_{hr}^pg^{nr},\ g^{np}=\dfrac{\partial z^n}{\partial x_k}\dfrac{\partial z^p}{\partial x_k}, \\ \Gamma_{hr}^n=\dfrac{\partial z^n}{\partial x_m}\dfrac{\partial^2 x_m}{\partial z^h\partial z^r},\ \Gamma_{hr}^p=\dfrac{\partial z^p}{\partial x_t}\dfrac{\partial^2 x_t}{\partial z^h\partial z^r}\Rightarrow g^{np}_{;h}=\dfrac{\partial}{\partial x_k}\dfrac{\partial z^n}{\partial x_h}\dfrac{\partial z^p}{\partial x_k}+
\dfrac{\partial z^n}{\partial x_k}\dfrac{\partial}{\partial x_k}\dfrac{\partial z^p}{\partial z^h}+\dfrac{\partial z^n}{\partial x_m}\dfrac{\partial}{\partial z^r}\dfrac{\partial x_m}{\partial z^h}\dfrac{\partial z^r}{\partial x_q}\dfrac{\partial z^p}{\partial x_q}+
\dfrac{\partial z^p}{\partial x_t}\dfrac{\partial}{\partial z^r}\dfrac{\partial x_t}{\partial z^h}\dfrac{\partial z^n}{\partial x_s}\dfrac{\partial z^r}{\partial x_s}\hspace{-1mm}=\hspace{-1mm}
\dfrac{\partial z^n}{\partial x_m}\dfrac{\partial}{\partial z^h}\dfrac{\partial x_m}{\partial x_q}\dfrac{\partial z^p}{\partial x_q}+
\dfrac{\partial z^p}{\partial x_t}\dfrac{\partial}{\partial z^h}\dfrac{\partial x_t}{\partial x_s}\dfrac{\partial z^n}{\partial x_s}\hspace{-1mm}=\hspace{-1mm}0$ because, e.g., $\dfrac{\partial z^n}{\partial z^h}=\delta_{nh}, \dfrac{\partial x_m}{\partial x_q}=\delta_{mq}$ etc., so their derivatives are null.

ii) By Eqs. (\ref{eq:gcov1}), (\ref{eq:christoffelsymb1}) and (\ref{eq:covderivL2}), $g_{np;h}=\dfrac{\partial^2x_k}{\partial z^h\partial z^n}\dfrac{\partial x_k}{\partial z^p}+\dfrac{\partial x_k}{\partial z^n}\dfrac{\partial^2x_k}{\partial z^h\partial z^p}-\\\dfrac{\partial z^r}{\partial x_m}
\dfrac{\partial^2x_m}{\partial z^p\partial z^h}\dfrac{\partial x_q}{\partial z^n}\dfrac{\partial x_q}{\partial z^r}-
\dfrac{\partial z^r}{\partial x_t}
\dfrac{\partial^2x_t}{\partial z^n\partial z^h}\dfrac{\partial x_s}{\partial z^p}\dfrac{\partial x_s}{\partial z^r}=
\dfrac{\partial^2x_k}{\partial z^h\partial z^n}\dfrac{\partial x_k}{\partial z^p}+\dfrac{\partial x_k}{\partial z^n}\dfrac{\partial^2x_k}{\partial z^h\partial z^p}-\\\delta_{qm}\dfrac{\partial^2x_m}{\partial z^p\partial z^h}\dfrac{\partial x_q}{\partial z^n}- 
\delta_{st}\dfrac{\partial^2x_t}{\partial z^n\partial z^h}\dfrac{\partial x_s}{\partial z^p}=
\dfrac{\partial^2x_k}{\partial z^h\partial z^n}\dfrac{\partial x_k}{\partial z^p}+\dfrac{\partial x_k}{\partial z^n}\dfrac{\partial^2x_k}{\partial z^h\partial z^p}-\dfrac{\partial^2x_q}{\partial z^p\partial z^h}\dfrac{\partial x_q}{\partial z^n}-\\\dfrac{\partial^2x_s}{\partial z^n\partial z^h}\dfrac{\partial x_s}{\partial z^p}=0.$

\end{enumerate}

\section*{Chapter 7}

\begin{enumerate}
\item It is sufficient to pose $x_1=u,x_2=v,x_3=f(u,v)\Rightarrow p(u,v)$ defines a surface because as $f(u,v)$ is smooth,  $p(u,v)$ is also smooth, and because the Jacobian is $[J]=\left[\begin{array}{cc}1 & 0 \\0 & 1 \\f_{,u} & f_{,v}\end{array}\right]$, then rank$[J]=2$.

\item Catenoid: $\f(u,v):\left\{\begin{array}{l}x_1=\cosh u\cos v,\\x_2=\cosh u\sin v,\\x_3=u;\end{array}\right.$ Meridians: $v=const.$; if, for example, $v=0\Rightarrow\left\{\begin{array}{l}x_1=\cosh u,\\x_2=0,\\x_3=u,\end{array}\right.$ is a catenary in the plane $(x_1,x_3). \\
\f_{,u}=\left\{\begin{array}{c}\sinh u\cos v\\\sinh u\sin v\\1\end{array}\right\},\f_{,v}=\left\{\begin{array}{c}-\cosh u\sin v\\\cosh u\cos v\\0\end{array}\right\}\Rightarrow\g=\left[\begin{array}{cc}\cosh^2u & 0 \\0 & \cosh^2u\end{array}\right].\\
\f_{,u}\times \f_{,v}=\left\{\begin{array}{c}-\cosh u\cos v\\-\cosh u\sin v\\\sinh u\cosh u\end{array}\right\},|\f_{,u}\times \f_{,v}|=\cosh^2u\Rightarrow
\N=\dfrac{1}{\cosh u}\left\{\begin{array}{c}-\cos v\\-\sin v\\\sinh u\end{array}\right\};\\\f_{,uv}=\f_{,vu}=\left\{\begin{array}{c}-\sinh u\sin v\\\sinh u\cos v\\0\end{array}\right\},\f_{,uu}=\left\{\begin{array}{c}\cosh u\cos v\\\cosh u\sin v\\0\end{array}\right\},\\\f_{,vv}=\left\{\begin{array}{c}-\cosh u\cos v\\-\cosh u\sin v\\0\end{array}\right\}\Rightarrow\B=\left[\begin{array}{cc}-1 & 0 \\0 & 1\end{array}\right]\Rightarrow K=-\dfrac{1}{\cosh^4 u}.$

\item Pseudo-sphere: $\f(u,v): \left\{\begin{array}{l}x_1=\sin u\cos v,\\x_2=\sin u\sin v,\\x_3=\cos u+\ln\left(\tan\dfrac{u}{2}\right).\end{array}\right.$ Meridians: $v=const.$; if, for example, $v=0\Rightarrow\left\{\begin{array}{l}x_1=\sin u,\\x_2=0,\\x_3=\cos u+\ln\left(\tan\dfrac{u}{2}\right),\end{array}\right.$ is a tractrix in the plane $(x_1,x_3).\\
\f_{,u}=\left\{\begin{array}{c}\cos u\cos v\\\cos u\sin v\\-\sin u+\dfrac{1}{\sin u}\end{array}\right\},\f_{,v}=\left\{\begin{array}{c}-\sin u\sin v\\\sin u\cos v\\0\end{array}\right\}\Rightarrow\g=\left[\begin{array}{cc}\dfrac{\cos^2u}{\sin^2u} & 0 \\0 & \sin^2u\end{array}\right].\\
\f_{,u}\times\f_{,v}=\left\{\begin{array}{c}-\cos^2 u\cos v\\-\cos^2 u\sin v\\\sin u\cos u\end{array}\right\},|\f_{,u}\times\f_{,v}|=|\cos u|\Rightarrow\N=\dfrac{1}{|\cos u|}\left\{\begin{array}{c}-\cos^2 u\cos v\\-\cos^2 u\sin v\\\sin u\cos u\end{array}\right\};\\
\f_{,uu}=\left\{\begin{array}{c}-\sin u\cos v\\-\sin u\sin v\\-\cos u-\dfrac{\cos u}{\sin^2 u}\end{array}\right\},
\f_{,uv}=\f_{,vu}=\left\{\begin{array}{c}-\cos u\sin v\\\cos u\cos v\\0\end{array}\right\},\\
\f_{,vv}=\left\{\begin{array}{c}-\sin u\cos v\\-\sin u\sin v\\0\end{array}\right\}\Rightarrow \B=\left[\begin{array}{cc}-\dfrac{\cos^2u}{\sin u|\cos u|} & 0 \\0 & \dfrac{\cos^2u\sin u}{|\cos u|}\end{array}\right]\Rightarrow K=-1.$

\item Cone: $\f(u,v)=v\bg(u)\Rightarrow\f_{,u}=v\bg',\f_{,v}=\bg\rightarrow \N\neq\bo\iff v\neq0$ and $\bg\neq\alpha\bg'$, i.e. everywhere except at the apex of the cone and on straight lines tangent to $\bg$.

\item The most general equation of the hyperbolic hyperboloid is $\\ \f(u,v): \left\{\begin{array}{l}x_1=a(\cos u-v\sin u),\\x_2=b(\sin u+v\cos u),\\x_3= c\ v,\end{array}\right. a,b,c\in\R\Rightarrow\f(u,v)=\bg(u)+v\bl(u)$, with $\bg(u)=a\cos u\eu+b\sin u\ed,\bl(u)=-a\sin u\eu+b\cos u\ed+c\et\Rightarrow\f(u,v)$ is a ruled surface.
Fixing $u=u_0\Rightarrow\left\{\begin{array}{l}x_1=a(\cos u_0-v\sin u_0),\\x_2=b(\sin u_0+v\cos u_0),\\x_3=c\ v,\end{array}\right.$ equation of a bundle of straight lines belonging to $\f(u,v)$, as well as $\left\{\begin{array}{l}x_1=a(\cos u_0-v\sin u_0),\\x_2=b(\sin u_0+v\cos u_0),\\x_3= -c\ v.\end{array}\right.$ The angle formed by two straight lines of the two sets is $\theta=\arccos\dfrac{a^2\sin^2u_0+b^2\cos^2u_0-c^2}{a^2\sin^2u_0+b^2\cos^2u_0+c^2}.$

\item $x_3=x_1x_2;$ setting $u=x_1,v=x_2\Rightarrow\f(u,v):\left\{\begin{array}{l}x_1=u,\\x_2=v,\\x_3=u\ v,\end{array}\right.$ is of the type $\\ \f(u,v)=\bg(u)+v\bl(v)$, with $\bg(u)=u\eu,\bl(u)=\ed+u\et.\\$ The straight lines $\left\{\begin{array}{l}x_1=u_0,\\x_2=v,\\x_3=u_0\ v,\end{array}\right.$ and $\left\{\begin{array}{l}x_1=u,\\x_2=v_0,\\x_3=u\ v_0,\end{array}\right.$ belong of course to $\f(u,v)$; they form the angle $\theta=\arccos\dfrac{u_0v_0}{\sqrt{(1+u_0^2)(1+v_0^2)}}; \theta=\dfrac{\pi}{2}\iff 0=v_0=0$, i.e. at $(0,0,0)$.

\item i) $\bg(u)=(\cos u,\sin u,-1),\bl(u)=(\cos u,\sin u,1)\Rightarrow\f(u,v)=(\cos u,\sin u,2v-1)$, of the form $\f(u,v)=\bg_1(u)+v\bl_1(u),$ with $\bg_1(u)=(\cos u,\sin u,-1),\bl_1(u)=(0,0,2)=const.\Rightarrow\f(u,v)$ is a cylinder whose Cartesian equation is $x_1^2+x_2^2=1$.

ii) $\bg(u)=(\sin u,-\cos u,-1),\bl(u)=(-\sin u,\cos u,1)\Rightarrow\\\f(u,v)=(2v-1)(-\sin u,\cos u,1)$, of the form $ \f(u,v)=\bg_2(u)+v\bl_2(u),$ with $\bg_2(u)=(\sin u,-\cos u,-1),\bl_1(u)=(-2\sin u,2\cos u,2)=-2\bg_2(u)\Rightarrow\f(u,v)$  is a cone whose Cartesian equation is $x_1^2+x_2^2=x_3^2.$

iii) $\left\{\begin{array}{l}x_1=(1-v)\cos(u-\alpha)+v\cos(u+\alpha),\\x_2=(1-v)\sin(u-\alpha)+v\sin(u+\alpha),x_3=-(1-v)+v,\end{array}\right.\Rightarrow\\
\left\{\begin{array}{l}x_1=\cos u\cos\alpha-\sin u\sin\alpha(2v-1),\\x_2=\sin u\cos\alpha+\cos u\sin\alpha(2v-1),x_3=2v-1.\end{array}\right.\\$ Change in parameter $w=\dfrac{\sin\alpha}{\cos\alpha}(2v-1)\Rightarrow\left\{\begin{array}{l}x_1=\cos\alpha(\cos u-w\sin u),\\x_2=\cos\alpha(\sin u+w\cos u),\\x_3=w\dfrac{\cos\alpha}{\sin\alpha},\end{array}\right.\Rightarrow\left\{\begin{array}{l}x_1=a(\cos u-w\sin u),\\x_2=a(\sin u+w\cos u),\\x_3=c w,\end{array}\right. a=\cos\alpha,\ c=\dfrac{\cos\alpha}{\sin\alpha}$, which is the parametric equation  of a hyperbolic hyperboloid with Cartesian equation $\dfrac{x_1^2}{a^2}+\dfrac{x_2^2}{a^2}-\dfrac{x_3^2}{c^2}=1\rightarrow \dfrac{x_1^2+x_2^2}{\cos^2\alpha}-\dfrac{x_3^2}{\cot^2\alpha}=1.$

\item i) Sphere of radius $R:x_1^2+x_2^2+x_3^2=R^2\Rightarrow$ using the spherical coordinates $\theta=u,\phi=v$  for expressing the $x_i$s, we get the parametric equation $\f(u,v):\left\{\begin{array}{l}x_1=R\cos\theta\sin\phi,\\x_2=R\sin\theta\sin\phi,\\x_3=R\cos\phi,\end{array}\right.\Rightarrow
\f_{,u}=R(-\sin u\sin v,\cos u\sin v,0),\\\f_{,v}=R(\cos u\cos v,\sin u\cos v,-\sin v)\Rightarrow\g=\left[\begin{array}{cc}R^2\sin^2v&0\\0&R^2\end{array}\right].$

ii) If $\bw=a\f_{,u}+b\f_{,v}\in T_p\Sigma,I(\bw)=\bw\cdot\g\bw=R^2(a^2\sin^2v+b^2)$.

iii) $A=\int_{\theta_1}^{\theta_2}\int_0^\pi\sqrt{\det\g}du\ dv=\int_{\theta_1}^{\theta_2}\int_0^\pi\sqrt{R^4\sin^2v}du\ dv=2R^2(\theta_2-\theta_1).$

iv) Parallel: putting $u=t,v=\dfrac{\pi}{4}, \bg(t):\left\{\begin{array}{l}x_1=\dfrac{R}{\sqrt{2}}\cos t,\\x_2=\dfrac{R}{\sqrt{2}}\sin t,\\x_3=\dfrac{R}{\sqrt{2}},\end{array}\right.\Rightarrow\\\bg'(t):\left\{\begin{array}{l}x_1=-\dfrac{R}{\sqrt{2}}\sin t,\\x_2=\dfrac{R}{\sqrt{2}}\cos t,\\x_3=0,\end{array}\right.\Rightarrow
\bg'(t)=\dfrac{du}{dt}\f_{,u}+\dfrac{dv}{dt}\f_{,v}=\f_{,u}\Rightarrow$ in the natural basis of $T_p\Sigma,\bw=(1,0)$ is the tangent vector to the parallel $\bg(t)\Rightarrow\\ I(\bw)=\bw\cdot\g\bw=R^2\sin^2v=\dfrac{R^2}{2}\Rightarrow\ell=\int_{\theta_1}^{\theta_2}\sqrt{I(\bw)}dt=\dfrac{R}{\sqrt{2}}(\theta_2-\theta_1).$

\item i) $x_1^2+x_2^2+x_3^2=\dfrac{\cos^2v}{\cosh^2u}+\dfrac{\sin^2v}{\cosh^2u}+\dfrac{\sinh^2u}{\cosh^2u}=\dfrac{1}{\cosh^2u}+\dfrac{\sinh^2u}{\cosh^2u}=\dfrac{\cosh^2u}{\cosh^2u}=1\rightarrow$ Cartesian equation of a sphere of centre $(0,0,0)$ and radius $R=1$.

ii) Straight line in $\Omega:\left\{\begin{array}{l}u=u_0+a\ t,\\v=v_0+b\ t;\end{array}\right.\Rightarrow$ curve on $\Sigma: \bg(t):\left\{\begin{array}{l}x_1=\dfrac{\cos(v_0+bt)}{\cosh(u_0+at)},\\x_2=\dfrac{\sin(v_0+bt)}{\cosh(u_0+at)},\\x_3=\dfrac{\sinh(u_0+at)}{\cosh(u_0+at)},\end{array}\right.$ or also $\bg(t)=\f(u(t),v(t))\Rightarrow\bg'(t)=\dfrac{du}{dt}\f_{,u}+\dfrac{dv}{dt}\f_{,v}=a\f_{,u}+b\f_{,v}.$

 Meridians: setting $v=const.=\hat{v}\Rightarrow\bmu(u)=\f(u,\hat{v})$; in fact $\dfrac{x_2}{x_1}=\tan\hat{v}=const.\rightarrow$ equation of a vertical plane. Tangent to the meridian $\bmu(u): \bmu'(u)=\f_{,u}\Rightarrow$ in the natural basis $\{\f_{,u},\f_{,v}\},\bg'(t)=(a,b),\bmu'(u)=(1,0);\cos\theta=\dfrac{I(\bg',\bmu')}{\sqrt{I(\bg')I(\bmu')}}.\\
 \f_{,u}=\left(-\cosh v\dfrac{\sinh u}{\cosh^2u},-\sin v\dfrac{\sinh u}{\cosh^2 u},\dfrac{1}{\cosh^2u}\right),\f_{,v}=\left(-\dfrac{\sin v}{\cosh u},\dfrac{\cos v}{\cosh u},0\right)\Rightarrow\\\g=\dfrac{1}{\cosh^2u}\I\Rightarrow I(\bg',\bmu')=\bg'\cdot\g\bmu'=\dfrac{a}{\cosh^2u},I(\bg')=\bg'\cdot\g\bg'=\dfrac{a^2+b^2}{\cosh^2u},\\I(\bmu')=\bmu'\cdot\g\bmu'=\dfrac{1}{\cosh^2u}\Rightarrow\cos\theta=\dfrac{a}{\sqrt{a^2+b^2}}=const.\Rightarrow\bg(t)$ is a loxodromic line on the sphere.

\item i) Catenoid $\f(u,v):\left\{
 \begin{array}{l}
 x_1=\phi(u)\cos v,\\
 x_2=\phi(u)\sin v,\\
 x_3=\psi(u). 
 \end{array}
 \right.$ with $\phi(u)=\cosh u,\psi(u)=u\Rightarrow\\\phi'(u)=\sinh u,\phi''(u)=\cosh u,\psi'(u)=1,\psi''(u)=0\Rightarrow\\
 \f_{,u}=(\sinh u\cos v,\sinh u\sin v,1),\f_{,v}=(-\cosh u\sin v,\cosh u\cos v,0).$
 
 ii) $\g=\cosh^2u\I$.
 
 iii) $\f_{,u}\times\f_{,v}=(-\cos v,-\sin v,-\sinh u),|\f_{,u}\times\f_{,v}|=\cosh u\Rightarrow\\\N=\left(-\dfrac{\cos v}{\cosh u},-\dfrac{\sin v}{\cosh u},-\dfrac{\sinh u}{\cosh u}\right).\\
 \f_{,uu}=(\cosh u\cos v,\cosh u\sin v,0),\f_{,uv}=\f_{,vu}=(-\sinh u\sin v,\sinh u \cos v,0),\\\f_{,vv}=-\f_{,uu}\Rightarrow
 \B=\left[\begin{array}{cc}-1&0\\0&1\end{array}\right]$.

iv) $\g^{-1}=\dfrac{1}{\cosh^2u}\I\Rightarrow\X=\g^{-1}\B=\dfrac{1}{\cosh^2u}\B.$

v) Let $\bw=(a,b)\in T_p\Sigma\Rightarrow I(\bw)=\bw\cdot\g\bw=\cosh^2 u(a^2+b^2).$

vi) $II(\bw)=\bw\cdot\B\bw=b^2-a^2.$

\item i) Helicoid $\f(u,v):\left\{\begin{array}{l}x_1=v\cos u,\\x_2=v\sin u,\\x_3=u,\end{array}\right.\Rightarrow\f_{,u}=(-v\sin u,v \cos u,1),\f_{,v}=(\cos u,\sin u,0).$

ii) $\g=\left[\begin{array}{cc}1+v^2&0\\0&1\end{array}\right].$

iii) $\f_{,u}\times\f_{,v}=(-\sin u,\cos u,-v),|\f_{,u}\times\f_{,v}|=\sqrt{1+v^2}\Rightarrow\N=\dfrac{1}{\sqrt{1+v^2}}(-\sin u,\cos u,-v).\\
\f_{,uu}=(-v\cos u,-v\sin u,0),\f_{,uv}=\f_{,vu}=(-\sin u,\cos u,0),\f_{,vv}=(0,0,0)\Rightarrow\\\B=\dfrac{1}{\sqrt{1+v^2}}\left[\begin{array}{cc}0&1\\1&0\end{array}\right]$.

iv) $\g^{-1}=\left[\begin{array}{cc}\dfrac{1}{1+v^2}&0\\0&1\end{array}\right].\Rightarrow\X=\g^{-1}\B=\left[\begin{array}{cc}0&\dfrac{1}{(1+v^2)^{\frac{3}{2}}}\\\dfrac{1}{(1+v^2)^{\frac{1}{2}}}&0\end{array}\right].$

v) Let $\bw=(a,b)\in T_p\Sigma\Rightarrow I(\bw)=\bw\cdot\g\bw=(1+v^2)a^2+b^2.$

vi) $II(\bw)=\bw\cdot\B\bw=\dfrac{2ab}{\sqrt{1+v^2}}.$

\item i) Catenoid (see Exercise \ref{ex:2ch7}): $K=-\dfrac{1}{\cosh^4u}<0\ \forall u\Rightarrow$ hyperbolic points.

ii) Helicoid (see Exercise \ref{ex:11ch7}): $K=\dfrac{\det\B}{\det\g}=-\dfrac{1}{(1+v^2)^2}<0\ \forall v\Rightarrow$ hyperbolic points.

\item Parametric equation of a circular cylinder of radius $R\rightarrow\f(u,v):\left\{\begin{array}{l}x_1=R\cos v,\\x_2=R\sin v,\\x_3=u,\end{array}\right.$ with $u=z,v=\theta$ of a system of cylindrical coordinates. Referring to  Eq. (\ref{eq:ruledsurf}),  $ \phi(u)=R,\psi(u)=u\Rightarrow$ Eq. (\ref{eq:geodetichesurfrot}) is $\left\{\begin{array}{l}u''=0,\\v''=0,\end{array}\right.\Rightarrow\left\{\begin{array}{l}u(t)=\alpha t+\alpha_1,\\v(t)=\beta t+\beta_1,\end{array}\right.$ with $\alpha,\alpha_1,\beta,\beta_1=const.$ If $\alpha_1=\beta_1=0$, we get the geodesic $\bg(t)$ passing through $(R,0,0)$ for $t=0\Rightarrow\bg(t):\left\{\begin{array}{l}x_1=R\cos v(t),\\x_2=R\sin v(t),\\x_3=u(t),\end{array}\right.\Rightarrow\left\{\begin{array}{l}x_1=R\cos (\beta t),\\x_2=R\sin (\beta t),\\x_3=\alpha t,\end{array}\right.\Rightarrow$ equation of a helix if $\alpha,\beta\neq0$, of a circle (cross section) if $\alpha=0,\beta\neq0$ and of a straight line on the cylinder (generatrix) if $\alpha\neq0,\beta=0$.

\end{enumerate}

\label{ch:sol}

\end{document}